\documentclass[11pt, a4paper,reqno]{amsart}
\pdfoutput=1
\usepackage{tikz,genyoungtabtikz,enumitem,rotating,latexsym,bm,stmaryrd,caption,textgreek}
\usetikzlibrary{positioning,intersections,decorations}
\usepackage[all]{xy}
\usepackage{bbm}
\usetikzlibrary{shapes}

\usepackage{array}

\definecolor{ao(english)}{rgb}{0.0, 0.5, 0.0}

\definecolor{darkgreen}{rgb}{0.0, 0.5, 0.0}

\usetikzlibrary{decorations.pathmorphing}

	\definecolor{eng}{rgb}{0.0, 0.5, 0.0}
\definecolor{apple}{rgb}{0.55, 0.71, 0.0}
\definecolor{cadmium}{rgb}{0.0, 0.42, 0.24}
\definecolor{darkspringgreen}{rgb}{0.09, 0.45, 0.27}
\definecolor{amethyst}{rgb}{0.6, 0.4, 0.8}
\definecolor{ao}{rgb}{0.0, 0.0, 1.0}
\definecolor{atomictangerine}{rgb}{1.0, 0.6, 0.4}
\definecolor{carmine}{rgb}{0.59, 0.0, 0.09}

\definecolor{toggle}{rgb}{1.0, 0.94, 0.96}

\usepackage[normalem]{ulem}

 \newcommand{\dom}{{\sf dom}}
  \newcommand{\cod}{{\sf cod}}

 \newcommand{\al}{{{{\color{magenta}\boldsymbol\alpha}}}}
  \newcommand{\gam}{{{{\color{ao(english)}\boldsymbol\gamma}}}}
 \newcommand{\bet}{{{\color{cyan}\boldsymbol\beta}}}
 \newcommand{\emp}{{{\color{magenta}\mathbf{\boldsymbol\clock}}}}
  \newcommand{\empb}{{{\color{cyan}\mathbf{\boldsymbol\clock}}}}
    \newcommand{\empg}{{{\color{ao(english)}\mathbf{\boldsymbol\clock}}}}

\newcommand{\exx}{{b_\al }}

   \renewcommand{\tau}{h}

\newcommand{\reflectpath}{\SSTP_\al^\flat}

\newcommand{\eps}{ \varepsilon}

\newcommand{\isit}{{i}}
\newcommand{\isitone}{\al(i+1)}
\newcommand{\Shl}{\widehat{\mathfrak{S}}_{\aatch}}
\newcommand{\w}{{\underline{w}}}
\newcommand{\x}{{\underline{x}}}

\newcommand{\Alc}{\text{\bf Alc}}
 \newcommand{\Pdiptwo}{{\sf M}_{i,i+2}}
 \newcommand{\Pdipj}{{{\sf M}_{i,j}}}

\usepackage{standalone}

\usepackage{etex}

\usepackage{todonotes}

\usepackage{tikz}
\usetikzlibrary{matrix,  intersections, calc, decorations.pathreplacing} 

\usepackage{tikz-cd}

      \newcommand{\aatch}{h}
      \newcommand{\aatchpair}{{\underline{h}}}
      \newcommand{\enn}{{h}}

\newlength{\superthick}
\newlength{\cornerradius}
\tikzstyle{corner}=[rounded corners=\cornerradius]
\tikzstyle{dot}=[circle, inner sep=0pt, minimum size=4.8pt]
\tikzstyle{string}=[line width=\superthick]
\tikzstyle{std}=[string,dash pattern=on 0.9pt off 0.9pt]
\definecolor{realcyan}{rgb}{0,1,1}

 \usepackage{amsmath,amsthm,amssymb,mathrsfs,pb-diagram}
\usepackage[
bookmarks=true,colorlinks=true,linktoc=page,citecolor=darkgreen,linkcolor=darkgreen,urlcolor=darkgreen]{hyperref}
\usepackage{caption}
\usepackage{lipsum,wasysym}
\usepackage{mathtools}
 \usepackage[a4paper,margin=1in]{geometry}
\usepackage{cleveref}
 \usepackage{amsmath}
\mathchardef\mhyphen="2D
\usepackage{color}
\usepackage{xcolor}
\usepackage{ifthen}
\usepackage{sidecap}   
\definecolor{mediumblue}{rgb}{0.0, 0.0, 0.8}

\newcommand\mptn[2]{\mathscr{P}_{#1}({#2})}
\renewcommand{\geq}{\geqslant}
\renewcommand{\leq}{\leqslant}
\renewcommand{\trianglerighteq}{\trianglerighteqslant}
\renewcommand{\trianglelefteq}{\trianglelefteqslant}

\tikzset{wei/.style= 
{red,double=red,double
distance=0.5pt}}

\tikzset{wei2/.style={red,double=red,double
distance=0.5pt}}

\allowdisplaybreaks
\numberwithin{equation}{section}
\usepackage{scalefnt}

\newtheorem{thm}{Theorem}[section]
\newtheorem{cor}[thm]{Corollary}

\newtheorem{lem}[thm]{Lemma}
\newtheorem{prop}[thm]{Proposition}

\newtheorem*{prop*}{Proposition}
\newtheorem*{thmA*}{Theorem A}
\newtheorem*{corA}{Corollary A}
\newtheorem*{corB}{Corollary B}
\newtheorem*{thmB*}{Theorem B}
\newtheorem*{thmC*}{Theorem C}\newtheorem*{thm*}{Theorem D}
\newtheorem*{cor*}{Corollary}

\newtheorem*{conj*}{Conjecture A}

\newtheorem*{conj1*}{Conjecture B}
\newtheorem*{Acknowledgements*}{Acknowledgements}

\theoremstyle{rmk}

\theoremstyle{defn}
\newtheorem{rmk}[thm]{Remark}
\newtheorem{defn}[thm]{Definition}
\newtheorem{eg}[thm]{Example}

\newtheorem{notn}[thm]{Notation}
 \newtheorem{conventioning}[thm]{Convention}
\newcommand{\great}{>}
\newcommand{\less}{<}

\newcommand{\rad}{\mathrm{rad}}
\newcommand{\res}{\mathrm{res}}

\newcommand{\Std}{{\rm Std}}

\newcommand{\Shape}{\operatorname{Shape}} 
\newcommand{\Path}{{\rm Path}}

\newcommand{\la}{\lambda}
\newcommand{\I}{i}
\newcommand{\J}{j}

\newcommand{\M}{m}

\newcommand{\SSTS}{\mathsf{S}}

\newcommand{\SSTT}{\mathsf{T}}  
\newcommand{\SSTP}{\mathsf{P}}  
\newcommand{\SSTU}{\mathsf{U}}  
\newcommand{\SSTV}{\mathsf{V}}  
\newcommand{\SSTQ}{\mathsf{Q}}  
\newcommand{\sts}{\mathsf{s}}  
\newcommand{\stt}{\mathsf{t}}  
\newcommand{\ZZ}{{\mathbb Z}}

\newcommand{\CC}{{\mathbb{C}}}
\newcommand{\RR}{{\mathbb R}}

\DeclareMathOperator{\Hom}{Hom}

\let\<=\langle
\let\>=\rangle

\newcommand{\north}{top }
\newcommand{\northT}{{\sf T}}
\newcommand{\south}{bottom } 
\newcommand{\southT}{{\sf B}}

\newcommand\mydots{\makebox[1em][c]{.\hfil.\hfil.}}
\tikzset{
ultra thin/.style= {line width=0.05pt},
very thin/.style=  {line width=0.2pt},
thin/.style=       {line width=0.1pt},
semithick/.style=  {line width=0.6pt},
thick/.style=      {line width=0.8pt},
very thick/.style= {line width=1.2pt},
ultra thick/.style={line width=1.6pt}
}

\crefname{defn}{Definition}{Definitions}
\crefname{thm}{Theorem}{Theorems}
\crefname{prop}{Proposition}{Propositions}
\crefname{lem}{Lemma}{Lemmas}
\crefname{cor}{Corollary}{Corollaries}
\crefname{conj}{Conjecture}{Conjectures}
\crefname{section}{Section}{Sections}
\crefname{subsection}{Subsection}{Subsections}
\crefname{eg}{Example}{Examples}
\crefname{figure}{Figure}{Figures}
\crefname{rem}{Remark}{Remarks}
\crefname{rmk}{Remark}{Remarks}
\crefname{equation}{equation}{equation}

\Crefname{defn}{Definition}{Definitions}
\Crefname{thm}{Theorem}{Theorems}
\Crefname{prop}{Proposition}{Propositions}
\Crefname{lem}{Lemma}{Lemmas}
\Crefname{cor}{Corollary}{Corollaries}
\Crefname{conj}{Conjecture}{Conjectures}
\Crefname{section}{Section}{Sections}
\Crefname{subsection}{Subsection}{Subsections}
\Crefname{eg}{Example}{Examples}
\Crefname{figure}{Figure}{Figures}
\Crefname{rem}{Remark}{Remarks}
\Crefname{rmk}{Remark}{Remarks}

\usepackage[hang,flushmargin]{footmisc}

\newcommand\REMOVETHESE[2]{{{{\mathsf{M}}_{#1}^{#2}	}}}
\newcommand\ADDTHIS[2]{{{{\mathsf{P}}_{#1}^{#2}}}}

 \newcommand{\Spotzero}{\SSTS_{0,\al}}
  \newcommand{\Spotone}{\SSTS_{1,\al}}
   \newcommand{\Spottwo}{\SSTS_{2,\al}}
    \newcommand{\Spotthree}{\SSTS_{3,\al}}
 \newcommand{\Spotq}{\SSTS_{q,\al}}
  \newcommand{\Spotb}{\SSTS_{\exx,\al}}
  \newcommand{\Spotqplus}{\SSTS_{q+1,\al}}

 \newcommand{\Forkq}{{\sf F}_{q,\emp\al}}
  \newcommand{\Forkqplus}{{\sf F}_{q+1,\emp\al}}
 \newcommand{\Forkone}{{\sf F}_{1,\emp\al}}
 \newcommand{\Forkzero}{{\sf F}_{0,\emp\al}} 
  \newcommand{\Forktwo}{{\sf F}_{2,\emp\al}}
    \newcommand{\Forkthree}{{\sf F}_{3,\emp\al}}
  \newcommand{\Forkb}{{\sf F}_{\exx,\emp\al}}

\begin{document}
  
 \title[Diagrammatic Hecke   and Bott--Samelson endomorphism algebras]{Path isomorphisms between quiver Hecke  and  diagrammatic Bott--Samelson  endomorphism algebras}
 \author{Chris Bowman}
       \address{Department of Mathematics, 
University of York, Heslington, York, YO10 5DD, UK}
\email{Chris.Bowman-Scargill@york.ac.uk}

  \author{Anton Cox}

 	\address{Department of Mathematics, City, University of London,   London, UK}
\email{A.G.Cox@city.ac.uk}

   \author{Amit Hazi}
        \address{Department of Mathematics, 
University of York, Heslington, York, YO10 5DD, UK}

 \email{Amit.Hazi@york.ac.uk}


   \maketitle
   

\vspace{-0.75cm}
  \begin{abstract}
  We construct an explicit isomorphism between (truncations of) quiver Hecke algebras and 
Elias--Williamson's  diagrammatic endomorphism algebras of  Bott--Samelson  
bimodules. As a corollary, we deduce that the  decomposition numbers of  these algebras (including as examples the symmetric groups and generalised blob algebras) are tautologically equal to 
 the associated  $p$-Kazhdan--Lusztig polynomials,  provided that the  characteristic is greater than the Coxeter number.  
 We hence give an elementary and more explicit proof of the main theorem of Riche--Williamson's recent monograph and  extend their  categorical equivalence to cyclotomic quiver Hecke algebras, thus solving Libedinsky--Plaza's categorical blob conjecture.     
 
  \end{abstract}

\section{Introduction}
 The  
  symmetric group lies at the intersection of two great categorical theories.  
 The first   is Khovanov--Lauda and Rouquier's  categorification of quantum groups and their knot invariants \cite{MR2525917,ROUQ}; 
  this setting has provided powerful new graded presentations of the symmetric group and its affine Hecke algebra \cite{MR2551762}.  
  The second  is  Elias--Williamson's diagrammatic categorification in terms of  endomorphisms of Bott--Samelson bimodules; it was in this setting that the counterexamples to Lusztig's  conjecture  were first found \cite{MR3671935} and  
 that the first general character formulas for  decomposition numbers of symmetric groups were discovered \cite{MR3805034} (in characteristic $p>h$, the Coxeter number).

 The purpose of this paper is to     construct an explicit isomorphism between   these two  diagrammatic worlds.  
This allows us to provide an elementary  algebraic proof 
of \cite[Theorem 1.9]{MR3805034}  and to vastly generalise this theorem 
   to the quiver  Hecke (or KLR) algebras $\mathcal{H}_n$; we hence 
   settle Libedinsky--Plaza's categorical blob conjecture \cite{blobby}.  
    Understanding its  simple modules is equivalent to understanding those of its  cyclotomic quotients     $\mathcal{H}_n^\sigma $ for $\sigma=(\sigma_0,\sigma_1,\dots, \sigma_{\ell-1}) \in \ZZ^\ell$.  
We prove that   
    $\mathcal{H}_n^\sigma $  
      has graded 
    decomposition  numbers   $d_{\la,\mu}(t)$  equal to the 
    $p$-Kazhdan--Lusztig polynomials  of type 
$$      A_{h_0}\times \mydots \times A_{h_{\ell-1}} \backslash \widehat{A}_{h_0+\dots+h_{\ell-1}} $$
provided that $\la$ and $\mu$ have at most $h_m$
 columns in the $m$th component 
(where $h_m\leq   \sigma_{m+1}-\sigma_{m} $ for $0\leq m< \ell-1$ and  $h_{\ell-1}<  e+\sigma_0-\sigma_{\ell-1}  $).  
We denote the set of such $\ell$-multipartitions by $\mathscr{P}_{\underline{h}}(n)$ 
 for $\underline{h}=(h_0,\dots, h_{\ell-1})\in \ZZ_{\geq0}^\ell$ and refer to such an $\underline{h}\in \ZZ^\ell $  
 as being {\em $(\sigma,e)$-admissible}.  
 This is the broadest possible generalisation, in the context of  the quiver Hecke algebra, of studying the category of tilting modules of the principal block of the general linear group, ${\rm GL}_h(\Bbbk)$, in characteristic $p>h$.

   \begin{thmA*}
\color{black}  Let     $\sigma\in \ZZ^\ell $  and $e\in \ZZ_{>1}$ and suppose that 
$\underline{h}\in \ZZ^\ell_{\geq0}$ is $(\sigma,e)$-admissible.  
We have a canonical isomorphism of  graded $\ZZ$-algebras  between  certain  subquotients of the quiver Hecke algebra 
  $\mathcal{H}_n^\sigma $  and 
 Elias--Williamson's diagrammatic   category under which the simple and standard modules labelled by $ \mathscr{P}_{\underline{h}}(n)$ are preserved.  
 The isomorphism  is  defined in  \cref{thisistheisomorphismfor the intro}.

  \end{thmA*}

   Perhaps most importantly, our isomorphism 
     allows one to pass information back and forth between these two diagrammatic categorifications for the first time.   
Combining our result with \cite{MR2551762}  allows one to import Soergel calculus to calculate decomposition numbers
directly within the setting of the symmetric group (and more generally, within the cyclotomic quiver Hecke algebras). 
 For instance, the key to the counterexamples of \cite{MR3671935} 
 are the mysterious 
``intersection forms" controlling decompositions of Bott--Samelson 
bimodules; in light of our isomorphism, these intersection forms can 
be seen simply as an   efficient version  of James' classical 
bilinear form on the 
Specht modules of $\Bbbk\mathfrak{S}_n$, and the efficiency  arises by way of idempotent 
truncation (in particular, the Gram matrices of these forms are equal).
   In other words, by virtue of our isomorphism, one can view the current state-of-the-art regarding  $p$-Kazhdan--Lusztig theory  (in type $A$) 
  entirely within  the context of the group algebra of the symmetric group, 
without the need for calculating intersection cohomology groups, or working with
  parity sheaves, or appealing to the deepest results of 2-categorical Lie theory.  In \cref{whyyyy} we will explain 
   that the regular decomposition numbers of cyclotomic quiver Hecke algebras are {\em tautologically equal to}  $p$-Kazhdan--Lusztig polynomials, simply by the categorical definition of these polynomials.  

\begin{thmB*} \color{black} 
The isomorphism of Theorem A  maps each choice of light leaves cellular basis  to a cellular basis element of $\mathcal{H}^\sigma_n$. 
Thus the Gram matrix of  
the intersection form associated to the fibre  of a Bott--Samelson resolution  of a Schubert variety coincides  with the Gram matrix of 
  James' bilinear  form on  the idempotent truncated Specht module   for $\la \in \mathscr{P}_{\underline{h}}(n)$.  
 \end{thmB*}

In the other  direction:  Soergel diagrammatics is, at present, confined to regular blocks ---  
whereas  quiver Hecke diagrammatics is not so restricted --- 
we expect our isomorphism to offer insight toward  constructing Soergel diagrammatics for singular blocks.  
  In particular, our isomorphism  interpolates  between the (well-understood) {LLT-style combinatorics} of KLR algebras and the (more mysterious) {Kazhdan--Lusztig-style combinatorics} of diagrammatic Bott--Samelson endomorphism algebras.  
       
\smallskip\noindent
{\bf Symmetric groups. }
For $\ell=1$ our Theorem A has the  immediate  corollary of reproving  the  famous result of   Riche--Williamson (and later Elias--Losev) which states that regular decomposition numbers of symmetric groups are equal to  $p$-Kazhdan--Lusztig polynomials 
\cite{MR3805034,ELpaper}.  Our proof is 
conceptually simpler than both existing proofs, as it does not require any higher  categorical Lie theory.  
 Once one has developed the appropriate combinatorial framework, our proof simply verifies  that the two diagrammatically defined algebras are isomorphic by checking the  relations.  
 In this regard, our proof is   akin to the work of Brundan--Kleshchev  
 \cite{MR2551762} and extends their ideas to the world of Soergel diagrammatics.
 We state the simplified version of Theorem A now, for ease of reference.

   \begin{corA}\color{black} 
   For $\Bbbk$ a field of characteristic $p> h$, we have an isomorphism of  graded $\Bbbk$-algebras between   certain  subquotients of 
$ 
 \Bbbk\mathfrak{S}_n $ and   Elias--Williamson's diagrammatic   category of type $ A_{h-1}\backslash \widehat{A}_{h -1}$. 
The decomposition numbers of symmetric groups labelled by partitions with at most $h< p $ columns are
{tautologically} equal to the $p$-Kazhdan--Lusztig polynomials of type $ A_{h-1}\backslash \widehat{A}_{h -1}  $.  
   
   \end{corA}

\smallskip\noindent
{\bf Blob algebras and statistical mechanics. }
The (generalised) blob algebras first arose as the transfer matrix algebras for the one-boundary Potts models in statistical mechanics.  
In a series of beautiful and prophetic papers \cite{MS94,MW00,MW03}, Paul Martin and his collaborators 
conjectured that these algebras would be controlled by non-parabolic affine  Kazhdan--Lusztig polynomials and verified this conjecture for level $\ell=2$.  It was the  advent of quiver Hecke and Cherednik algebras   that provided the necessary perspective  for solving this conjecture \cite{manycell}.  
This perspective    allowed  Libedinsky--Plaza to  push these ideas still further (into the modular setting) in the form of a beautiful conjecture
which brings together ideas from statistical mechanics, diagrammatic algebra, and $p$-Kazhdan--Lusztig theory for the first time \cite{blobby}.  
For $h=(1^\ell)$ our Theorem A verifies their conjecture, as follows:   
  
     \begin{corB}[Libedinsky--Plaza's categorical blob conjecture]
     \color{black} 
 For $ \Bbbk$ a field, we have  an isomorphism of  graded $\Bbbk$-algebras, 
between  certain  subquotients of the generalised blob algebra of level $\ell$  and   Elias--Williamson's diagrammatic   category of type $ \widehat{A}_{\ell -1}$. 
 In particular the  decomposition numbers of generalised blob algebras  are {tautologically} equal to the $p$-Kazhdan--Lusztig polynomials of type $ \widehat{A}_{\ell-1}    $.  
   
   \end{corB}

\smallskip\noindent
{\bf Weightings and gradings on cyclotomic quiver Hecke algebras. }  
Recently, Elias--Losev generalised  
\cite[Theorem 1.9]{MR3805034}   to calculate    decomposition numbers of cyclotomic quiver Hecke algebras.  
However, we emphasise that 
   our Theorem A and 
  Elias--Losev's work intersect only in the case of the symmetric group (providing two independent proofs of  \cite[Theorem 1.9]{MR3805034}).  
In particular, Elias--Losev's work {\em does   not} imply 
Libedinsky--Plaza's conjecture (as explained in detail in Libedinsky--Plaza's paper \cite{blobby}).    
This lack of overlap arises from different choices of weightings on the cyclotomic  quiver Hecke algebra,  
we refer the reader to \cite{blobby,manycell,STEEN} for more details.

  \color{black} 
\smallskip\noindent
{\bf The structure and ideas  of the paper. }
The isomorphism of this paper was a surprise to many of the experts in this field.  This is because of the fundamental differences in the ways 
 we think of Bott--Samelson endomorphism algebras versus quiver Hecke algebras.  
The elements of the former algebras are thought of as morphisms between words (in the Coxeter generators of $\widehat{\mathfrak{S}}_{\aatch}$), their complex representation theory is controlled by Soergel's algorithm, which can be thought of in terms of paths in the Bruhat graph of $ {\mathfrak{S}}_{\aatch}\leq \widehat{\mathfrak{S}}_{\aatch}$.  
The elements of the latter algebras algebras are thought of as ``graded versions" of permutations, the complex representation theory of these algebras is  controlled by the LLT algorithm, which can be thought of in terms of graded standard tableaux \cite{KN10}.  
Of course the LLT algorithm and Soergel's algorithm produce the same results, even though the steps involved appear quite different.  
One can think of this as being because the LLT algorithm has many  more ``degree zero steps" which simply ``pad out" the tableaux.  
This is a good heuristic for this paper, which we now expound section  by section.
  
Sections 2 and 3 introduce the combinatorics and basic definitions of  quiver Hecke and  diagrammatic Bott--Samelson endomorphism algebras in tandem.
  We provide a dictionary for passing between standard tableaux (of the former world) 
and   expressions in  cosets of affine Weyl group (of the latter world) 
by means of coloured  paths in our   alcove geometries.  
We subtly tweak the classical perspective for  quiver Hecke algebras 
by recasting each element of the algebra as a morphism   between a  pair  of paths in the alcove geometry.   
Heuristically, we  ``equate the combinatorics" of the LLT and Soergel algorithms
 by  writing 
  tableaux/paths as the  concatenation of  component paths (each of which corresponds to a single reflection hyperplane).


One of the  core principles  of this paper is that 
  diagrammatic Bott--Samelson endomorphisms   are 
 simply  a  ``condensed shorthand" for  KLR path-morphisms.
 Section 4 details the reverse process by which we ``dilate" simple elements of the KLR algebra and hence   construct these path-morphisms. 
 Section 4 also  provides a translation principle by which we can see that   a path-morphism  depends only on the series of hyperplanes in the path's trajectory, not the individual steps taken within the path.  Heuristically, this translation principle says that ``the degree zero steps in the LLT algorithm are unimportant".  

In Section 5, we recast the generators of the diagrammatic Bott--Samelson endomorphism algebra within the setting of the quiver Hecke algebra; this allows us to explicitly state the isomorphism, $\Psi$,  of  Theorem A.    
%
In Section 6 we verify that $\Psi$ is a graded $\ZZ$-algebra homomorphism by recasting the relations of the diagrammatic Bott--Samelson endomorphism within the setting of the quiver Hecke algebra.  
This involves rewriting  products of the path-morphisms in the KLR algebra  one step at a time --- for the products involving  forks and spots there is a single ``important step" in this procedure with the others corresponding to ``LLT padding".

 Finally, in Section 7 we match-up the light leaves bases of these algebras under the map $\Psi$ and hence prove that $\Psi$ is bijective and thus complete the proofs of Theorems A and B.

In Appendix A we provide a coherence theorem for weakly graded monoidal categories which allows us to relate the classical Bott-Samelson endomorphism algebras to certain breadth-enhanced  versions which are more convenient for the purposes of this paper.  
   The reader can think of this as inserting   ``extra monoidal  identity padding"      into the diagrammatic Bott--Samelson endomorphisms algebras  which corresponds (on the KLR side of the isomorphism) to the  {\em steps of degree zero} in paths/tableaux. 

  Finally we emphasise that the LLT/Soergel analogy above is motivated by the situation over $\CC$. This is merely a heuristic and our results work over a field of arbitrary characteristic (indeed, the isomorphism is actually proven to hold over the integers).

  For the convenience of the reader we provide three tables summarising the notation used throughout the paper in Appendix B. 
 
    \color{black}

\section{Parabolic and non-parabolic alcove geometries and path combinatorics}  \label{sec3}
\renewcommand{\sts}{\mathsf{S}}   
\renewcommand{\stt}{\mathsf{T}} 

Without loss of generality,  we  assume that  $\sigma \in \ZZ^\ell$ is 
weakly increasing and $e>h\in \ZZ_{\geq1}$.  We say that     $\underline{h}= (h_0,\dots,h_{\ell-1})  \in \ZZ_{\geq0}^\ell$
with $h_0+h_1+\dots +h_{\ell-1}=h$  is {\sf$(\sigma,e)$-admissible} if 
  $h_m\leq   \sigma_{m+1}-\sigma_{m} $ for $0\leq m< \ell-1$ and  $h_{\ell-1}<  e+\sigma_0-\sigma_{\ell-1}  $.  
  (This condition on $\underline{h}= (h_0,\dots,h_{\ell-1})  \in \ZZ_{\geq0}^\ell$ is equivalent to   the empty partition   not lying on any hyperplane of our alcove geometry, so that the resulting Kazhdan--Lusztig theory is ``non-singular".)

\subsection{Multipartitions, residues and tableaux} \label{necessaryevil}\color{black}
  We define a {\sf composition}, $\lambda$,  of $n$ to be a   finite   sequence  of non-negative integers $ (\lambda_1,\lambda_2, \ldots)$ whose sum, $|\lambda| = \lambda_1+\lambda_2 + \mydots$, equals $n$.   We say that $\la$ is a partition if, in addition, this sequence is weakly decreasing.  
An    {\sf $\ell $-multicomposition} (respectively {\sf $\ell$-multipartition})  $\lambda=(\lambda^{(0)},\mydots,\lambda^{(\ell-1)})$ of $n$ is an $\ell $-tuple of   compositions (respectively of partitions)   such that $|\lambda^{(0)}|+\mydots+ |\lambda^{(\ell-1)}|=n$. 
We will denote the set of $\ell $-multicompositions (respectively $\ell$-multipartitions) of $n$ by
$\mathscr{C}_{\ell}(n)$ (respectively by $\mathscr{P}_{\ell}(n)$).   
Given  $\lambda=(\lambda^{(0)},\lambda^{(1)},\ldots ,\lambda^{(\ell-1)}) \in \mathscr{P}_{\ell}(n)$, the {\sf (dual) Young diagram} of $\lambda$    is defined to be the set of nodes, 
\[\color{black} 
[\la]=\{(r,c,m) \mid  1\leq  r\leq 
(\lambda^{(m)})_c, 0
\leq m <\ell \}.
\]
Notice that we have taken the transpose-dual of the
 usual conventions so that the multipartitions are the 
 sequences  whose columns are weakly decreasing (this is a
  trivial, if unfortunate, 
   relabelling inherited from our earlier work \cite{CoxBowman,cell4us}).  
We do not distinguish between the multipartition and its (dual) Young diagram.  
We refer to a node $(r,c,m)$ as being in the $r$th row and $c$th column of the $m$th component of $\lambda$.  
Given a node, $(r,c,m)$,  
we define the {\sf content} of this node  to be  ${\rm ct}(r,c,m)  = \sigma_m+  c - r $ and we define its  {\sf residue} to be ${\rm res}(i,j,m)= {\rm ct}(i,j,m) \pmod e$.    
We refer to a node of residue  $i\in \ZZ/e\ZZ$ as an $i$-node.      
   Given $\la \in \mathscr{C}_ \ell( n)$  or $\mptn \ell n$, we let   ${\rm Rem} (\la)$ (respectively ${\rm Add}  (\la)$) 
  denote the set of all    removable  (respectively addable) nodes   of the Young diagram of $\la$ so that the resulting diagram is the Young diagram of an $\ell$-composition or an $\ell$-partition. 
 
Given $\lambda\in  \mathscr{C}_\ell(n)$, we define a {\sf tableau} of shape $\lambda$ to be a filling of the nodes  of $\la $ with the numbers 
$\{1,\mydots , n\}$. 
We define a {\sf standard tableau} of shape $\lambda$ to be a tableau of shape $\lambda$ such that  entries increase along the rows and down the columns of each component.  
We let $\Std (\lambda)$ denote the set of all standard tableaux of shape $\lambda\in\mathscr{P}_\ell(n)$.  We let $\varnothing$ denote the empty multipartition.  
 
\color{black}

\begin{defn} 
  Given a pair of $i$-nodes $(r,c,m), (r',c',m')$, we write 
   $(r,c,m) \lhd  (r',c',m')$ if either 
   ${\rm ct}(r,c,m) < {\rm ct}(r',c',m')$ or 
   ${\rm ct}(r,c,m) = {\rm ct}(r',c',m')$  and $m>m'.$  
 For $\la,\mu\in \mptn \ell n$, we write   $ \mu \trianglelefteq \la $  if there is a 
  bijective   map ${\sf A}:[\lambda] \to [\mu]$ 
 such that either $  {\sf A}(r,c,m)\vartriangleleft (r,c,m) $  or $  {\sf A}(r,c,m)=(r,c,m) $ for all $(r,c,m)\in \lambda$.  

   \end{defn}
  
   Given $\sts\in \Std(\la)$ a , we write $\sts{\downarrow}_{\leq k} $ or $\sts{\downarrow}_{\{1,\dots, k\}} $ (respectively $\sts{\downarrow}_{\geq k} $) for the subtableau of $\sts$ consisting solely of the entries $1$ through $k$ (respectively of the entries $k$ through $n$).   Given $\la \in  \mptn {\ell}n$, we let $\stt_\la$ denote the     $\la$-tableau in which we place the  entry $n$ in the minimal (under the $\rhd$-ordering) removable  node of $\la$, then continue in this fashion inductively.  
 Given $1\leq k\leq n$, we let $(r_k,c_k,m_k)\in \la$ be the node such that 
 $\SSTT(r_k,c_k,m_k)=k$.  
We let  ${\mathcal A} _\stt(k)$ (respectively ${\mathcal R} _\stt(k)$) 
 denote the set of   all addable (respectively removable) $\res (r_k,c_k,m_k)$-nodes   
   of the multipartition $\Shape(\stt{\downarrow}_{\{1,\dots ,k\}})$ which
   are less than  $(r_k,c_k,m_k)$ in the $\rhd $-order.     We define the $(\rhd)$-degree of $\stt\in \Std(\la)$ for $\la \in  \mathscr{P}_\ell  (n)$  as follows,
$$
\deg (\stt) = \sum_{k=1}^n \left(	|{\mathcal A} _\stt(k)|-|{\mathcal R} _\stt(k)|	\right).
$$

\begin{defn}\color{black}Given $\underline{h}\in \ZZ_{\geq0}^\ell$,   we let    $\mathscr{P}_{\underline{h}}(n)\subseteq
\mathscr{C}_{\underline{h}}(n) $ denote the subsets of  $\ell$-multipartitions and $\ell$-multicompositions with at most $h_m$ columns in the $m$th component for $0\leq m <\ell$.     
\end{defn}
   
  If  $\underline{h}\in \ZZ_{\geq0}^\ell$ is $(\sigma,e$)-admissible, then   ${\rm deg}(\SSTT_
  \la)=0$ for $\la \in \mptn {\underline{h}}n$.  
   
\begin{eg}\label{continuation}\color{black}
Let  $\sigma=(0,3,8)\in \ZZ^3$ and $e=13$.  
We note that $\aatchpair = (3,5,4)$ is $(\sigma,e)$-admissible.  We depict 
$\la= ((5,4,2 ),(5,4,3,2^2 ),(5,3^2, 2 ))\in    \mathscr{P}_{\underline{h}}(n)$ along with the residues of this multipartition as follows: 
  $$   \left( \ 
 \begin{minipage}{2.1cm}
\scalefont{0.9} 
  \end{minipage}  \right).
   $$
   Notice that any given residue  $i\in \ZZ/e\ZZ$ appears  at most once in a fixed  row of the multipartition.    \end{eg}

\subsection{Alcove geometry}   
 For ease of notation, we  set $H_m=h_0+\dots+h_m$ for $0\leq m <\ell$, 
and   $\aatch=\tau_0+\dots+\tau_{\ell-1}$.  
 For each $ 1\leq \I \leq h_m$ and $ 0\leq \M < \ell$ we let
  $\eps_{i,m}:=\eps_{(\tau_0+\dots+\tau_{m-1}) + i}$    denote a
formal symbol, and define an   $\aatch  $-dimensional real vector space 
\[
{\mathbb E}_{\aatchpair   }
=\bigoplus_{
	\begin{subarray}c 
	0 \leq m <  \ell   \\ 
		1\leq i \leq  \tau_m    
	\end{subarray}
} \mathbb{R}\varepsilon_{i,m}
\]
and $\overline{\mathbb E}_{\aatchpair  }$ to be the quotient of this space by the one-dimensional subspace spanned by 
\[\sum_{
	\begin{subarray}c 
	0 \leq m <  \ell   \\ 
		1\leq i \leq  \tau_m    
	\end{subarray}
} \varepsilon_{i,m}.\]
We have an inner product $\langle \; , \; \rangle$ on ${\mathbb E}_{\aatchpair  }$ given by extending
linearly the relations 
\[
\langle \varepsilon_{i,p} , \varepsilon_{j,q} \rangle= 
\delta_{\I,\J}\delta_{p,q}
\]
for all $1\leq \I, \J \leq n$ and $0 \leq p,q< \ell $, where
$\delta_{i,j}$ is the Kronecker delta.  
We identify $\lambda \in  \mptn {\aatchpair  }n $ with an element of the integer lattice inside $\mathbb{E}_{\aatchpair  }  $ via the map
\begin{equation}\label{pprrpr}
\lambda\longmapsto
 \sum_{\begin{subarray}c   0\leq {m}<\ell \\ 1\leq \I\leq \tau_m  \end{subarray}}
 \lambda^{(m)} 
 _\I \varepsilon_{i,m}. 
 \end{equation}
%
We let $\Phi$ denote the root system of type $A_{\aatch  -1}$ consisting of the roots 
$$\{\varepsilon_{i,p}-\varepsilon_{j,q}:  \ 0\leq p,q<\ell, \
1\leq i \leq \tau_p ,1\leq j \leq \tau_q, 
 \text{with}\ (i,p)\neq (j,q)\}$$
and $\Phi_0$ denote the root system of type $A_{\tau_0-1}\times\cdots\times A_{\tau_{\ell-1}-1}$ consisting of the roots 
 $\{\varepsilon_{i,m}-\varepsilon_{j,m}: 
  0\leq m<\ell, 1\leq i \neq j\leq \tau_m \}.$ 
 We choose $\Delta$ (respectively $\Delta_0$) to be the set of simple roots inside $\Phi$ (respectively  $\Phi_0$) of the form $\varepsilon_t-\varepsilon_{t+1}$ for some $1\leq t\leq h$, 
 and write $\Phi^+$ (respectively $\Phi_0^+$) for the set of positive roots with respect to this choice of simple roots.
Given $r\in\ZZ$ and $\alpha\in\Phi$ we define $s_{\alpha,re}$ to be the reflection which acts on ${\mathbb E}_{\aatchpair  }$ by
$$s_{\alpha,re}x=x-(\langle x,\alpha\rangle -re)\alpha.$$
The group generated by the $s_{\alpha,0}$ with $\alpha\in\Phi$ (respectively $\alpha\in\Phi_0$) is isomorphic to the symmetric group $\mathfrak{S}_{\enn   }$ (respectively to $\mathfrak{S}_f:=\mathfrak{S}_{\tau_0}\times\cdots\times\mathfrak{S}_{\tau_{\ell-1}}$), while the group generated by the $s_{\alpha,re}$ with $\alpha \in\Phi$ and $r\in\ZZ$ is isomorphic to $\widehat{\mathfrak{S}}_{\enn   }$, the affine Weyl group of type $A_{\enn   -1}$. 
    We set $\alpha_0=\varepsilon_{\enn}-\varepsilon_1$ and $\Pi=\Delta\cup\{\alpha_0\}$.
The elements
 $S=
\{s_{\alpha,0}:\alpha\in\Delta\}\cup\{s_{\alpha_0,-e}\}
$ 
generate $\widehat{\mathfrak S}_{\enn}$. 
We have chosen $\alpha_0=\varepsilon_{\enn}-\varepsilon_1$ (rather 
than $\alpha_0=\varepsilon_{1}-\varepsilon_\enn$) as this 
is compatible with out path combinatorics.

\begin{notn}
	We shall frequently find it convenient to refer to the generators in $S$ in terms of the elements of $\Pi$, and will abuse notation in two different ways. First, we will write $s_{\alpha}$ for $s_{\alpha,0}$ when $\alpha\in\Delta$ and $s_{\alpha_0}$ for $s_{\alpha_0,-e}$. This is unambiguous except in the case of the affine reflection $s_{\alpha_0,-e}$, where this notation has previously been used for the element $s_{\alpha,0}$. As the element $s_{\alpha_0,0}$ will not be referred to hereafter this should not cause confusion.
	Second, we will write $\alpha=\eps_i-\eps _{i+1}$ in all cases; if $i=\enn$ then all occurrences of $i+1$ should be interpreted modulo $\enn$ to refer to the index $1$.
\end{notn}

We shall consider a shifted action of the affine Weyl group $\widehat{\mathfrak{S}}_{\aatch  }$  on ${\mathbb E}_{\aatchpair}$ 
by the element
$ 
 \rho:= (\rho_{0}, \rho_{2}, \ldots, \rho_{\ell-1}) \in \ZZ^{\aatch  } $  where $
\rho_m := (  \sigma_m+\aatch_m-1 ,  \sigma_m+\aatch_m-2,\dots,   \sigma_m ) \in \ZZ^{\tau_m},$    
that is, given an element $w\in \widehat{\mathfrak{S}}_{\aatch} $,  we set 
$
w\cdot x=w(x+\rho)-\rho.
$
  This shifted action induces a well-defined action on $\overline{\mathbb E}_{\underline{h}}$; we will define various geometric objects in ${\mathbb E}_{\underline{h}}$ in terms of this action, and denote the corresponding objects in the quotient with a bar without further comment.
We let ${\mathbb E} ({\alpha, re})$ denote the affine hyperplane
consisting of the points  
$${\mathbb E} ({\alpha, re}) = 
\{ x\in{\mathbb E}_{\underline{h}} \mid  s_{\alpha,re} \cdot x = x\} .$$
Note that our assumption that $\underline{h}\in  \ZZ_{\geq0} ^\ell$ is $(\sigma,e)$-admissible implies that the origin does not lie on any hyperplane.   Given a
hyperplane ${\mathbb E} ( \alpha,re)$ we 
remove the hyperplane from ${\mathbb E}_{\underline{h}}$ to obtain two
distinct subsets ${\mathbb 
	E}^{\great}(\alpha,re)$ and ${\mathbb E}^{\less}(\alpha,re)$
where the origin  lies in $  {\mathbb E}^{\less }(\alpha,re)$. The connected components of 
$$\overline{\mathbb E}_{\underline{h}}  \setminus (\cup_{\alpha \in \Phi_0}\overline{\mathbb E}(\alpha,0))$$ are called chambers.
 The dominant chamber, denoted
$\overline{\mathbb E}_{\aatchpair}^+ $, is defined to be 
$$\overline{\mathbb E}_{\aatchpair}^+=\bigcap_{ \begin{subarray}c
	\alpha \in \Phi_0
	\end{subarray}
} \overline{\mathbb E}^{\less} (\alpha,0).$$
The connected components of $$\overline{\mathbb E}_{\underline{h}}  \setminus (\cup_{\alpha \in \Phi,r\in \ZZ}\overline{\mathbb E}(\alpha,re))$$
are called alcoves, and any such alcove is a fundamental domain for the action of the group $  \widehat{\mathfrak{S}}_{\aatch}$ on
the set $\Alc$ of all such alcoves. We define the {\sf fundamental alcove} $A_0$ to be the alcove containing the origin (which is inside the dominant chamber). 
 We note that the   map  
 $\mptn \aatchpair n \to  {\mathbb E}_{\underline{h}} \cap\ZZ_{\geq0}\{\eps_1,\dots,\eps_h\}$   
  restricts to be surjective when 
  we restrict  the codomain  to the dominant Weyl chamber.

%
%

We have a bijection from $\widehat{\mathfrak{S}}_{\aatch}$ to $\Alc$ given by $w\longmapsto wA_0$. Under this identification $\Alc$ inherits a right action from the right action of $\widehat{\mathfrak{S}}_{\aatch}$ on itself.
Consider the subgroup 
$$
\mathfrak{S}_f:=
\mathfrak{S}_{\tau_0}\times\cdots\times\mathfrak{S}_{\tau_{\ell-1}} \leq \widehat{\mathfrak{S}}_{\aatch  }.
$$The dominant chamber is a fundamental domain for the action of $\mathfrak{S}_f$ on the set of chambers in $\overline{\mathbb E}_{\underline{h}}$.  
We let $ \mathfrak{S}^f$ denote the set of minimal length representatives for right cosets $\mathfrak{S}_f \backslash \widehat{\mathfrak{S}}_{\aatch}$.  So multiplication gives a bijection $\mathfrak{S}_f\times \mathfrak{S}^f \to  \widehat{\mathfrak{S}}_{\aatch}$. 
 This induces a bijection between right cosets and the alcoves in our dominant chamber. 
  Under this identification,   the alcoves are partially ordered  by the Bruhat-ordering on $ \mathfrak{S}^f$.  (This is the opposite of the    ordering, $\trianglelefteq $, on multipartitions belonging to these alcoves.)

If the intersection of a hyperplane $\overline{\mathbb E}(\alpha,re)$ with the closure of an alcove $A$ is generically  of codimension one in $\overline{\mathbb E}_{\underline{h}}$ then we call this intersection a {\sf wall} of $A$. The fundamental alcove $A_0$ has walls corresponding to $\overline{\mathbb E}(\alpha,0)$ with $\alpha\in\Delta$ together with an affine wall $\overline{\mathbb E}(\alpha_0,e)$. We will usually just write $\overline{\mathbb E}(\alpha)$ for the walls $\overline{\mathbb E}(\alpha,0)$ (when $\alpha\in\Delta$) and $\overline{\mathbb E}(\alpha,e)$ (when $\alpha=\alpha_0$). We regard each of these walls as being labelled by a distinct colour (and assign the same colour to the corresponding element of $S$). Under the action of $\widehat{\mathfrak{S}}_{\aatch}$ each wall of a given alcove $A$ is in the orbit of a unique wall of $A_0$, and thus inherits a colour from that wall. We will sometimes use the right action of $\widehat{\mathfrak{S}}_{\aatch}$ on $\Alc$.  
\color{black} 
Given an alcove $A$ and an element $s\in S$  we have that $A=wA_0$ for some $w$ under the identification above (that is, $\widehat{\mathfrak{S}}_{\aatch}$ to $\Alc$ given by $w\longmapsto wA_0$). Thus the right action of 
$s$ on $A$ gives the element $wsA_0$ in $\Alc$, and this can easily be seen to be obtained   by reflecting $A$ in the wall of $A$ with colour corresponding to the colour of $s$. 
With this observation it is now easy to see that if $w=s_{1}\ldots s_{t}$ where the $s_i$ are in $S$ then $wA_0$ is the alcove obtained from $A_0$ by successively reflecting through the walls corresponding to $s_1$ up to $s_t$.
\color{black}

 If a wall of an alcove $A$ corresponds to $\overline{\mathbb E}(\alpha,re)$ and $A\subset \overline{\mathbb E}^>(\alpha,re)$ then we call this a {\sf lower alcove wall} of $A$; otherwise we call it an {\sf upper alcove wall} of $A$. We will call a multipartition {\sf $\sigma$-regular} (or just {\sf regular}) if its image in $\overline{\mathbb E}_{\aatchpair}$ lies in some alcove; those multipartitions whose images lies on one or more walls will be called {\sf $\sigma$-singular}.

      Let $\lambda \in \overline{\mathbb E }_{\underline{h}}  $.     There are only finitely many hyperplanes   
  ${\mathbb E}  (\alpha,re)  $ for $\alpha 
  \in \Pi$ and $r\in \ZZ$  lying   between  the   points   $\lambda  \in {\mathbb E} _{\underline{h}}$
 and the point  $ \varnothing\in \overline{\mathbb E}  _{\underline{h}}$.  
We let $\ell_\alpha (\la )$ denote the total number of these  hyperplanes for a given $\alpha \in \Pi$ 
 (including any hyperplane upon which $\lambda$ 
lies).   

\newcommand{\Po}{\operatorname{Po}}

\subsection{Paths in the geometry}
We now  
develop the  combinatorics of paths  inside our geometry.  
Given   a map  $p: 
\{1,\mydots
, n\}\to \{1,\mydots ,   {\aatch}    \}$ we define points $\SSTP(k)\in
{\mathbb E}_{\underline{h}}$ by
\[
\SSTP(k)=\sum_{1\leq i \leq k}\varepsilon_{p(i)} 
\]
for $1\leq i \leq n$. 
We define the associated path of length $n$    by 
$$
\SSTP=\left(
\varnothing=\SSTP(0),\SSTP(1),\SSTP(2), \ldots, \SSTP(n) \right) 
$$
 and we say that 
the path has shape $\pi= \SSTP(n) \in   {\mathbb E}_{\underline{h}}$.   
 We also denote this path by 
 $\SSTP=(\varepsilon_{p(1)},\ldots,\varepsilon_{p(n)}) $ 
\color{black} and call
$\varepsilon_{p(i)}$ the $i$th step in this path. \color{black}
Given $\lambda\in  {\mathbb E}_{\underline{h}} \cap\ZZ_{\geq0}\{\eps_1,\dots,\eps_h\}$ we let $\Path(\lambda)$ denote the set of paths of length $n$ with shape $\lambda$. We define $\Path_{\underline{h}}(\lambda)$ to be the subset of $\Path(\lambda)$ consisting of those paths lying entirely inside the dominant chamber, ${\mathbb E}  _{\underline{h}}^+$; in other words, those $\SSTP$ such that $\SSTP(i)$ is dominant for all $0\leq i\leq n$. 

 Given a path $\mathsf{P}$ defined by such a map $p$ of length $n$ and shape $\lambda$ we can write each $p(j)$ uniquely in the form $\eps_{p(j)}=\eps_{c_j,m_j}$ where $0\leq m_j<\ell$ and $1\leq c_j\leq \tau_j$. We record these elements in a tableau of shape $\lambda^T$ by induction on $j$, where we place the positive integer $j$ in the first empty node in the $c_j$th column of component $m_j$.
By definition, such a tableau will have entries increasing down columns; if $\lambda$ is a multipartition then the entries also increase along rows if and only if the given path is in $\Path_{\underline{h}}(\lambda)$, and hence there is a bijection between $\Path_{\underline{h}}(\lambda)$ and $\Std(\lambda)$. For this reason we will sometimes refer to paths as tableaux, to emphasise that what we are doing is generalising the classical tableaux combinatorics for the symmetric group. 

\begin{notn}
	Given a path $\SSTP$ we will let 
	$ 
	\SSTP^{-1}(r,\eps_{c,m})
	$ 	with $0\leq m<\ell$ and $1\leq c\leq h_{m}$ denote the $(r,c)$-entry of the $m$th component of the tableau corresponding to $\SSTP$.  In terms of our path this is the  point at which the $r$th  step of the form $+\eps_{c,m}$   occurs in $\SSTP$.  
Given a path $\SSTP$ we define $$\res(\SSTP)=(\res_\SSTP(1),\ldots,\res_\SSTP(n))$$ where $\res_\SSTP(i)$ denotes the residue of the node labelled by $i$ in the tableau corresponding to $\SSTP$.
\end{notn}

\begin{figure}[ht] 
 $$ 
  \begin{minipage}{10.5cm}  
\end{minipage}
  $$
 \caption{An alcove path  in $\Path_{(3)}(20,5^2)$ and the corresponding tableau in $\Std(20,5^2)$. The black vertices denote vertices on the path in the orbit of the origin.}
\label{diag2}
\end{figure}

\begin{eg}\label{eg1}
	We will illustrate our various definitions with an example in $\overline{\mathbb E}_{3,1}^+$ with $e=5$. This space is the projection of $\RR^3$ in two dimensions, which we shall represent as shown in Figure \ref{diag2}. Notice in particular that $\varepsilon_1+\varepsilon_2+\varepsilon_3=0$ in this projection, as required. Only the dominant chamber is illustrated, 
 	with the origin marked in the fundamental alcove $A_0$. 
	
	The affine Weyl group $\widehat{\mathfrak{S}}_3$ has generating set $S$ corresponding to the green and blue (non-affine) reflections $s_{\color{darkgreen}\varepsilon_2-\varepsilon_3,0}$ and $s_{\color{cyan}\varepsilon_1-\varepsilon_2,0}$
	about the lower walls of the fundamental alcove, together with the (affine) reflection $s_{\color{magenta}\varepsilon_3-\varepsilon_1,-5}$ about the red wall of that alcove. Recall that we will abuse notation, and refer to these simply as $s_{\color{darkgreen}\varepsilon_2-\varepsilon_3,}$, $s_{\color{cyan}\varepsilon_1-\varepsilon_2}$, and
	$s_{\color{magenta}\varepsilon_3-\varepsilon_1}$.
	The associated colours for the remaining alcove walls are as shown.
	
	Given $\lambda=(3^5,1^{15})$ we have illustrated a path $\SSTP$ from the origin to $\lambda$ with a black line. Recall that we embed partitions via the transpose map (as in \cref{pprrpr})  and so the final point in the path corresponds to the point $(20,5,5)\in E_{3,1}$. The corresponding steps in the path are recorded in the standard tableau at the bottom of the figure, where an entry $i$ in column $j$ of the tableau (again, note the transpose) corresponds to the $i$th step of the path being in the direction $\varepsilon_j$. This is an element of $\Path_{\underline{h}}(\lambda)$ as it never leaves the dominant region. 
	
	The path passes through the sequence of alcoves obtained from the fundamental alcove by reflecting through the walls labelled ${\color{magenta}R}$ then ${\color{darkgreen}G}$ then ${\color{cyan}B}$ then ${\color{magenta}R}$ then ${\color{darkgreen}G}$ then ${\color{cyan}B}$, and so the final alcove corresponds to the element $s_{\color{magenta}\varepsilon_3-\varepsilon_1}s_{\color{darkgreen}\varepsilon_2-\varepsilon_1}
	s_{\color{cyan}\varepsilon_3-\varepsilon_2}s_{\color{magenta}\varepsilon_3-\varepsilon_1}
	s_{\color{darkgreen}\varepsilon_2-\varepsilon_1}s_{\color{cyan}\varepsilon_3-\varepsilon_1}A_0$.
		If $\sigma=(0)$ then we have 
	$$\res(\SSTP)=(0,1,4,0,3,4,2,1,0,2,4,\ldots,1).$$	
\end{eg}

   \begin{eg}\color{black} 
Further  examples of paths and tableaux are given   in   \cref{paths-first,tableau-second,figure1}. 
\end{eg}

Given paths $\SSTP=(\varepsilon_{p(1)},\ldots,\varepsilon_{p(n)})$ and $\SSTQ=(\varepsilon_{q(1)},\ldots,\varepsilon_{q(n)})$ we say that $\SSTP\sim\SSTQ$ if there exists an $\alpha
 \in\Phi$ and $r\in\ZZ$ and $s\leq n$ such that
$$\SSTP(s)\in{\mathbb E}(\alpha,re)\qquad \text{ 
and  }
\qquad \varepsilon_{q(t)}=\left\{\begin{array}{ll}
\varepsilon_{p(t)}& \ \text{for}\ 1\leq t\leq s\\
s_{\alpha}\varepsilon_{p(t)}& \ \text{for}\ s+1\leq t\leq n.\end{array}\right.$$ 
In other words the paths $\SSTP$ and $\SSTQ$ agree up to some point $\SSTP(s)=\SSTQ(s)$ which lies on ${\mathbb E}(\alpha,re)$, after which each $\SSTQ(t)$ is obtained from $\SSTP(t)$ by reflection in ${\mathbb E}(\alpha,re)$. We extend $\sim$ by transitivity to give an equivalence relation on paths, and say that two paths in the same equivalence class are related by a series of {\sf wall reflections of paths}.
    We say that $\SSTP=(\varepsilon_{p(1)},\ldots,\varepsilon_{p(n)})$  is a {\sf reduced path} if 
$\ell_\alpha(\SSTP(s+1))\geq\ell_\alpha\SSTP(s )) $ for  all $1\leq s <n$ and  $\alpha \in \Pi$ .  
There exist a unique reduced path in each $\sim$-equivalence class.

\begin{lem}\label{inlightof}
\color{black}	We have $\SSTP\sim\SSTQ$ if and only $\res(\SSTP)=\res(\SSTQ)$. 
\end{lem}
\begin{proof}
	Let $\alpha=\varepsilon_{i,a}-\varepsilon_{j,b}$. We first note that a path of shape $\lambda$ lies on ${\mathbb E}(\alpha,re)$ if and only if the addable nodes in $\lambda $ in the $i$th column of the $a$th component and in the $j$th column of the $b$th component have the same residue. 
	 (This is straightforward from the definition of the inner product, see for example \cite[Lemma 6.19]{CoxBowman}.) Also $s_{\alpha}\varepsilon_t=\varepsilon_t$ for all $t\notin\{H_{a-1}+i,H_{b-1}+j\}$ and $s_\alpha$ permutes the elements of this set. So if two paths coincide up to some point and then diverge, but have the same sequence of residues, then the point where they diverge must lie on some ${\mathbb E}(\alpha,re)$ and the divergence must initially be by a reflection in this hyperplane. From this the result easily follows by induction on the number of hyperplanes which the two paths cross.
\end{proof}


We recast the degree of a tableau in the path-theoretic setting as follows.

\begin{defn}\label{Soergeldegreee}
Given a path $\sts=(\sts(0),\sts(1),\sts(2), \ldots, \sts(n))$  we set
$\deg(\sts(0))=0$ and define  
 \[
 \deg(\sts ) = \sum_{1\leq k \leq n} d (\sts(k),\sts(k-1)), 
 \]
for $d(\sts(k),\sts(k-1))$  defined as follows. 
For $\alpha\in	\Phi^+$ we set $d_{\alpha}(\sts(k),\sts(k-1))$ to be
\begin{itemize}
\item $+1$ if $\sts(k-1) \in 
   {\mathbb E}(\alpha,re)$ and 
   $\sts(k) \in 
   {\mathbb E}^{\less}(\alpha,re)$ for some $r \in \ZZ$;
   
\item $-1$ if $\sts(k-1) \in 
   {\mathbb E}^{\great}(\alpha,re)$ and 
   $\sts(k) \in 
   {\mathbb E}(\alpha,re)$ for some $r \in \ZZ$;
\item $0$ otherwise.  
   \end{itemize}
We let 
$$ 
 {\rm deg}(\SSTS)= \sum _{1\leq k \leq n }\sum_{\alpha \in \Phi^+}d_\alpha(\sts(k-1),\sts(k)).$$  
  \end{defn}

\begin{figure}[ht!]
 $$  \scalefont{0.8}   
 
 $$
\caption{\color{black} Two paths $\SSTS$ and $\SSTT$ in an alcove geometry.  These paths are used in \cref{egfortheref}. }
\label{paths-first}

\end{figure}

\vspace{-0.1cm}
\begin{figure}[ht!]
$$\scalefont{0.8}    
 \begin{minipage}{1cm}
  %
\end{minipage}   $$
  \caption{\color{black} The two tableaux $\SSTS$ and $\SSTT$ corresponding to the paths in \cref{paths-first}. These paths are used in \cref{egfortheref}.}
  \label{tableau-second}
  \end{figure}

\vspace{-0.15cm}

\subsection{Alcove paths}
When passing from multicompositions to our geometry $\overline{\mathbb E}_{\aatchpair}$, many non-trivial elements map to the origin. One such element is $\delta=((h_0),\mydots,(h_{\ell-1}))\in{\mathscr P}_{\underline{h}} ({\aatch})$.
 (Recall our transpose convention for embedding multipartitions into our geometry,   as in \cref{pprrpr}.) We will sometimes refer to this as the {\sf determinant} as (for $\ell=1$) it corresponds to the determinant representation of the associated general linear group. We will also need to consider elements corresponding to powers of the determinant, namely $\delta_n=((h_0^n),\mydots,(h_{\ell-1}^n))\in{\mathscr P}_\ell(n{\aatch})$.


We now restrict our attention to paths between points in the principal linkage class, in other words to paths between points in $\Shl\cdot 0$. Such points can be represented by the   $\mu$ in the orbit $\Shl\cdot\delta_n$ for some choice of $n$.


\begin{defn}\label{alphabet} We will associate alcove paths to certain words in the {\sf alphabet}
	$$
	S\cup\{1\}=\{ s_\alpha \mid \alpha \in \Pi \cup \{\emptyset\}\}
	$$
	where $s_\emptyset =1$.  That is, we will consider words in the generators of the affine Weyl group, but enriched with explicit occurrences of the identity 
	in these expressions.  
	When we wish to consider a particular expression for an element $w\in\widehat{\mathfrak{S}}_{\aatch}$ in terms of our alphabet we will denote this by $\underline{w}$.  
\end{defn}  

Our aim is to define certain distinguished paths from the origin to multipartitions in the principal linkage class; for this we will need to proceed in stages. In order to construct our path we want to proceed inductively. There are two ways in which we shall do this.

\begin{defn} Given two paths 
	$$\SSTP=(\eps_{i_1},\eps_{i_2},\dots, \eps_{i_p}) \in \Path(\mu)
	\quad\text{and}\quad
	\SSTQ=(\eps_{j_1},\eps_{j_2},\dots, \eps_{j_q}) 
	\in \Path(\nu)$$  we define the {\sf naive concatenated path} 
	$$\SSTP\boxtimes \SSTQ =
	(\eps_{i_1},\eps_{i_2},\dots, \eps_{i_p}, \eps_{j_1},\eps_{j_2},\dots, \eps_{j_q}) 
	\in \Path (\mu+\nu).$$ 
\end{defn}

There are several problems with naive concatenation. Most seriously, the naive concatenation   of two paths between points in the principal linkage class will not in general itself connect points in that class. Also, if we want to associate to our path the coloured sequence of walls through which it passes, then this is not compatible with naive concatentation. To remedy these failings, we will also need to define a {\sf contextualised concatenation}.

Given a path $\SSTP$ between points in the principal linkage class, the end point lies in the interior of an alcove of the form $wA_0$ for some $w\in\widehat{\mathfrak S}_{\aatch}$. If we write $w$ as a word in our alphabet, and then replace each element $s_{\alpha}$ by the corresponding non-affine reflection $s_{\alpha}$ in ${\mathfrak S}_{\aatch}$ to form the element $\overline{w}\in{\mathfrak S}_{\aatch}$ then the basis vectors $\varepsilon_i$ are permuted by the corresponding action of $\overline{w}$ to give $\varepsilon_{\overline{w}(i)}$, and there is an isomorphism from $\overline{\mathbb E}_{\aatchpair}$ to itself which maps $A_0$ to $wA_0$ such that $0$ maps to $w\cdot 0$, coloured walls map to walls of the same colour, and each basis element $\varepsilon_i$ map to $\varepsilon_{\overline{w}(i)}$. Under this map we can transform a path $\SSTQ$ starting at the origin to a path starting at $w\cdot 0$ which passes through the same sequence of coloured walls as $\SSTQ$ does. 

More generally, the end point of a path $\SSTP$ may lie on one or more walls. In this case, we can choose a distinct transformation as above for each alcove $wA_0$ whose closure contains the endpoint. We can now use this to define our contextualised concatenation.

\begin{defn}\label{contextuals}
	Given two paths 
	 $\SSTP=(\eps_{i_1},\eps_{i_2},\dots, \eps_{i_p}) \in \Path(\mu)$
and $	\SSTQ=(\eps_{j_1},\eps_{j_2},\dots, \eps_{j_q}) 
	\in \Path(\nu)$   with the endpoint of $\SSTP$ lying in the closure of some alcove $wA_0$
	we define the {\sf contextualised concatenated path} 
	$$\SSTP\otimes_w \SSTQ =
	(\eps_{i_1},\eps_{i_2},\dots, \eps_{i_p})\boxtimes  
	(\eps_{\overline{w}(j_1)},\eps_{\overline{w}(j_2) },\dots, \eps_{\overline{w}(j_q) }) 
	\in \Path (\mu+(w \cdot \nu)).$$ 
	If there is a unique such $w$ then we may simply write $\SSTP\otimes \SSTQ$.
	If $w=s_\al $ we will simply write $\SSTP\otimes_\al\SSTQ$.    \end{defn}

{\color{black} 
It is not difficult to understand contextualised concatenation in terms of tableaux.  
 Each symbol $\eps_i$ for $1\leq i \leq h$ labels a column of a partition. 
 Contextualised concatenation is then given by permuting the columns 
 (according to the rule in \cref{contextuals})
  and then vertically stacking the tableaux (and shifting the entries), see \cref{herewegoooooo}.  }

Our next aim is to define the building blocks from which all of our distinguished paths will be constructed.  We begin by defining certain integers that describe the position of the origin in our fundamental alcove.
  
\begin{defn} Given ${\al}\in \Pi$ we define 
	$\exx$ to be the distance from the origin to the wall corresponding to 
	${\al} $, and let $b_\emptyset =1$.    
	Given our earlier conventions this corresponds to setting
 $$
	b_{\varepsilon_{i,p}-\varepsilon_{ j,q }} = \sigma_q-\sigma_p + j-i
	 	\qquad
	b_{\varepsilon_{\aatch }-\varepsilon_{1 }} = e - \sigma_{0} + \sigma_{\ell-1} + h_{\ell-1} -1 
	$$
 	for $0\leq p\leq q<\ell$ and 	  $1\leq i \leq h_{p}$,  $1\leq j \leq h_{q}$.
 	We sometimes write $\delta_{\al}$ for the element $\delta_{\exx}$. 
	Given ${\al},\bet\in \Pi$  we set $b_{\al\bet}=b_\al + b_\bet$.   
\end{defn}

\begin{eg}\label{diag22point4}
	Let $e=5$, $h=3$ and $\ell=1$ as in \cref{diag2}.  
	Then $b_{\color{darkgreen}\varepsilon_2-\varepsilon_3}$ and 
	$b_{\color{cyan}\varepsilon_1-\varepsilon_2}$ both equal $1$, while $b_{\color{magenta}\varepsilon_3-\varepsilon_1}=3$ and $b_{\emptyset}=1$.  
\end{eg}

\begin{eg}
	Let $e=7$, $h=2$ and $\ell=2$ and $\sigma=(0,3)\in\ZZ^2$.  
	Then  
	$b_{\color{cyan}\varepsilon_1-\varepsilon_2}$ and $b_{\color{magenta}\varepsilon_3-\varepsilon_4}$ both equal $1$, while 
	$b_{\color{darkgreen}\varepsilon_4-\varepsilon_1}=3$, $b_{\color{orange}\varepsilon_2-\varepsilon_3}=2$, and $b_{\emptyset}=1$.  
\end{eg}

We can now define our basic building blocks for paths.

\begin{defn} \label{base} Given ${\al}  =\varepsilon_i -\varepsilon_{i+1}\in \Pi$, we consider the multicomposition 
	$ 
	s_\al\cdot \delta_{\al} 
	$ 
	with all columns of length $\exx$, with the exception of the $i$th and $(i+1)$st columns, which are of length  $0$ and $2\exx$, respectively.  
	We set 
	$$
	\REMOVETHESE {i} {} =  (\varepsilon_{1 }, \mydots, \varepsilon_{i-1 }
	, \widehat{\varepsilon_{ i  }},   \varepsilon_{ i+1  },\mydots, \varepsilon_{{\aatch} })
	\quad \text{and}\quad 
	\ADDTHIS {i} {} =(+\eps_i)$$
	where $\widehat{.}$ denotes omission of a coordinate.
	Then our distinguished path corresponding to $s_\al$ is given by
	$$ \SSTP_\al=
 	\REMOVETHESE{i} {\exx} \boxtimes \ADDTHIS {i+1}{\exx} \in \Path(s_\al\cdot\delta _{\al}).
	$$
	The distinguished path corresponding to $\emptyset$ is given by
	$$
	\SSTP_{\emptyset}=  (\varepsilon_{1 }, \mydots, \varepsilon_{i-1 }
	,\varepsilon_{ i  } ,   \varepsilon_{ i+1  },\mydots, \varepsilon_{{\aatch} })
	\in \Path(\delta)=\Path(s_\emptyset\cdot\delta) 
	$$ 
and set $\SSTP_\emp=(\SSTP_\emptyset)^{\exx}$.  	We will also find it useful to have the following variant of $\REMOVETHESE{i}{ }$. We set  $$
 	\Pdipj =
	(\eps_1,\dots, {\eps_{i-1}},\widehat{\eps_{i}}, {\eps_{i+1}},\dots
	, {\eps_{j-1}},\widehat{\eps_{j}}, {\eps_{j+1}},\dots,
	\eps_{\aatch}).
	$$
\end{defn}

\begin{figure}[ht!]
  $$\color{black} 
  \begin{minipage}{2cm}
\end{minipage} 
 $$
\caption{
Three paths and their corresponding tableaux.  
The leftmost two paths are the path $\SSTP_\al$ which  walks through an ${\al}$-hyperplane in $\overline{\mathbb E}_{1,3}^+$, and the path $\reflectpath $  which reflects the former path  through the same ${\al}$-hyperplane.
The rightmost   path is $\SSTP_\emptyset$ (which we repeat thrice).  
  } \label{figure1}
\end{figure}

\begin{figure}[ht!]  
  $$
  \begin{minipage}{1.2cm}  %
  \end{minipage}
  $$
\caption{\color{black} The tableau $\SSTP_\al\otimes_\al \SSTP_\al$ obtained by contextualised concatenation from the path/tableau $\SSTP_\al$ in \cref{figure1}. 
The reflection $s_\al $ for $\al=\eps_1-\eps_3$ permutes the first and third columns of $\SSTP_\al$.  
The entries of tableaux are coloured to facilitate comparison. 
The reader is invited to draw the corresponding path.}
\label{herewegoooooo}
\end{figure}

%
   \color{black}

\begin{eg}\color{black} 
The paths/tableaux $\SSTS$ and $\SSTT$ from \cref{paths-first,tableau-second} are equal to 
$\SSTP_\al\otimes_\al \SSTP_\bet \otimes_\bet \SSTP_\gam \otimes _\gam\SSTP_\bet $ 
and 
$\SSTP_\al\otimes_\al     \SSTP_\gam \otimes _\gam\SSTP_\bet 
\otimes_\bet
  \SSTP_\gam
$
 respectively for $\al=\eps_1-\eps_3$, 
 $\bet=\eps_1-\eps_2$,  $\gam=\eps_2-\eps_3$.   
\end{eg}

Given all of the above, we can finally define our distinguished paths for general words in our alphabet. There will be one such path for each word in our alphabet, and they will be defined by induction on the length of the word, as follows.

\begin{defn}\label{thepathwewatn}
	We now define a {\sf distinguished path} $\SSTP_{\w}$ for each word $\w$ in our alphabet $S\cup\{1\}$ by induction on the length of $\w$. If $\w$ is $s_\emptyset$ or a simple reflection $s_\alpha$ we have already defined the distinguished path in Definition \ref{base}. Otherwise if $\w=s_\alpha\w'$ then we define  
	$$
	\SSTP_{\underline{w}}:= \SSTP_{\alpha}
	\;\otimes_{\color{magenta}\al}  \; \SSTP_{\underline{w}'}	.$$ 
  \end{defn}

	If $\w$ is a reduced word in $\widehat{\mathfrak S}_{\aatch }$,  then   the  
	 path $\SSTP_{\w}$ is a  reduced path.

\begin{rmk}
	Contextualised concatenation is not associative (if we wish to decorate the tensor products with the corresponding elements $w$). As we will typically be constructing paths as in Definition \ref{thepathwewatn} we will adopt the convention that an unbracketed concatenation of $n$ terms corresponds to bracketing from the right: 
	$$\SSTQ_1\otimes\SSTQ_2\otimes \SSTQ_3\otimes\cdots \SSTQ_n=\SSTQ_1\otimes(\SSTQ_2\otimes(\SSTQ_3\otimes(\cdots \otimes \SSTQ_n)\cdots)).$$
\end{rmk}

\noindent We will also need certain reflections of our distinguished paths corresponding to elements of $\Pi$. 

\begin{defn}
	{ \renewcommand{\SSTU}{{\sf P}}\renewcommand{\mu}{{\pi}} 
		Given ${\al} \in \Pi$ we set 
		$$
		\reflectpath =
		\REMOVETHESE {i} {b_\al} \boxtimes \ADDTHIS {i} {\exx} =
		\REMOVETHESE {i} {b_\al} \otimes_\al \ADDTHIS {i+1} {\exx} =
		(+\varepsilon_{1 }, \mydots, +\varepsilon_{i-1 }
		, \widehat{+\varepsilon_{ i  }},   +\varepsilon_{ i+1  },\mydots, +\varepsilon_{{\aatch} })^\exx 
		\boxtimes  (\eps_i)^{\exx}
		$$
		the path obtained by reflecting the second part of $\SSTP_\al$ in the wall through which it passes.  
	} 
	
\end{defn}

\begin{eg}
We illustrate these various constructions in a series of examples. In the first two diagrams of Figure \ref{figure1}, we illustrate the basic path $\SSTP_{\al}$ and the path $\reflectpath$ 
 and in the rightmost diagram of Figure \ref{figure1},  we illustrate the path $\SSTP_\emptyset$. A more complicated example is illustrated in Figure \ref{diag2}, where we show the distinguished path $\SSTP_{\underline{w}}$ for $\underline{w}=
 s_{\color{magenta}\varepsilon_3-\varepsilon_1}
 s_{\color{darkgreen}\varepsilon_1-\varepsilon_2}
	s_{\color{cyan}\varepsilon_2-\varepsilon_3}
	s_{\color{magenta}\varepsilon_3-\varepsilon_1}
	s_{\color{darkgreen}\varepsilon_1-\varepsilon_2}
	s_{\color{cyan}\varepsilon_2-\varepsilon_3}$ as in Example \ref{eg1}. The components of the path between consecutive black nodes correspond to individual $\SSTP_{\al}$s.   \end{eg}

\begin{rmk}
\color{black} 
 There are plenty of other paths we could have chosen.  For example, we could replace the leftmost path in \cref{figure1} with the path 
 $$(\eps_1,\eps_1,\eps_1,\eps_2,\eps_2,\eps_2,\eps_1,\eps_1,\eps_1) \in \Path(6,3).$$
 In \cref{adjust1} we will see that it does not matter which path we pick, providing it ``does not hit any extra hyperplanes".  Our ``zig zagging" paths are merely the  easiest to define  such general paths.
\end{rmk}

\begin{rmk}\label{idemp-remark-for-the-refsssss}\color{black}
By \cref{inlightof} we have $\res(\SSTP_{\al})=\res(\reflectpath)$. 
This fact is key to our construction of the KLR versions of the diagrammatic Bott--Samelson generators using step-preserving permutations.
\end{rmk}

\begin{defn}
We say that a   word $\w=s_{\alpha^{(1)}}\mydots s_{\alpha^{(p)}} $
in either of
the alphabets 
$S$ or $S \cup\{1\}$ 
has breadth 
$${\sf breadth}_\sigma(\w)=\sum_{1\leq i \leq p}  b_{\alpha^{(i)}} $$ which we denote simply by $b_\w$ when the context is clear.    
  We let  
 $  \Lambda ({n,\sigma})$ (respectively   $\Lambda ^+({n,\sigma})$)
  denote the set of words $\w$ in the alphabet $S\cup\{1\}$ (respectively the alphabet  $S$)
such that 
$  {\sf breadth}_{\sigma}(\w)= n$.  
We define $$ \mathscr{P}_{\underline{h}}(n,\sigma)=\{ \la \in\mathscr{P}_{\underline{h}}(n)\mid \text{there exists } \SSTP_\w \in \Std(\la) , \w \in \Lambda(n,\sigma)\} .
$$  \end{defn}

\begin{eg} \label{exampleforanton}
We can  insert  the path   $\SSTP_\emptyset=(+\varepsilon_1,+\varepsilon_2,+\varepsilon_3) $ 
into the path in  \cref{diag2}  at seven distinct points to obtain a new alcove path. 
For example, we can insert two copies of this path (in two distinct ways) to obtain 
$\SSTP_\w$ and $\SSTP_{\w'}$ for  
$ 
\w=  s_\emptyset   s_\emptyset  s_{\color{magenta}\varepsilon_3-\varepsilon_1}
s_{\color{darkgreen}\varepsilon_2-\varepsilon_3}
s_{\color{cyan}\varepsilon_1-\varepsilon_2}
s_{\color{magenta}\varepsilon_3-\varepsilon_1} 
s_{\color{darkgreen}\varepsilon_2-\varepsilon_3}
s_{\color{cyan}\varepsilon_1-\varepsilon_2}
 $ and $
\w'=    s_{\color{magenta}\varepsilon_3-\varepsilon_1}s_\emptyset  
s_{\color{darkgreen}\varepsilon_2-\varepsilon_3}
s_{\color{cyan}\varepsilon_1-\varepsilon_2}s_\emptyset 
s_{\color{magenta}\varepsilon_3-\varepsilon_1}
s_{\color{darkgreen}\varepsilon_2-\varepsilon_3}
s_{\color{cyan}\varepsilon_1-\varepsilon_2}
$ 
respectively. Then ${\sf res}(\SSTP_\w)$ and ${\sf res}(\SSTP_{\w'})$ are equal to 
\begin{align*}
(
0,1,2,4,0,1, 
\color{magenta} 
3,4,
2,3,
1,2, 
0,
4,
3,
\color{darkgreen}	
0,
2,
1,
\color{cyan}
0,
1,
4,
\color{magenta}
3,  	4,
2,     3,
1,	2,
0,
4,
3,
\color{darkgreen}
0, 
2,		
1,\color{cyan}
0,    1,
4
\color{black}),
\\
(
\color{magenta} 
0,1,
4,0,
3,4, 
2,
1,
0,
\color{black} 
2,3,4,
\color{darkgreen}	
4,
1,
0,
\color{cyan}
4,
0,
3,
\color{black} 
2,3,4,
\color{magenta}
3,  	4,
2,     3,
1,	2,
0,
4,
3,
\color{darkgreen}
0, 
2,		
1,\color{cyan}
0,    1,
4
\color{black}).
\end{align*}
\end{eg}

For any $ \la  \in {\mathscr P}_{\underline{h}}(n)$, we define the  set of 
{\sf alcove-tableaux}, $\Std_{n,\sigma}(\la)$,  to consist  of all standard  tableaux which can be obtained by contextualised concatenation of paths from the set 
$$\{ \SSTP_\al \mid \al \in \Pi  \}\cup\{ \SSTP_\al^\flat  \mid \al \in \Pi  \}\cup \{ \SSTP_\emptyset   \} .$$
We let $\Std_{n,\sigma}^+ (\la) \subseteq \Std_{n,\sigma} (\la)$ denote the subset of {\sf strict alcove-tableaux} 
of the form $(\SSTP_\emptyset)^{\otimes p} \otimes \SSTQ$ for $\SSTQ$ obtained by contextualised concatenation of paths from the set $\{ \SSTP_\al \mid \al \in \Pi  \}\cup\{ \SSTP_\al^\flat  \mid \al \in \Pi  \}$ and some $p\geq0$.

\begin{eg}
The   tableau 
of shape $ (20,5^2)$ in \cref{diag2} is the strict alcove tableau given by $ \SSTP_\al \otimes_\al \SSTP_\gam \otimes_\gam \SSTP_\bet \otimes_\bet 
\SSTP_\al \otimes_\al \SSTP_\gam \otimes_\gam \SSTP_\bet.$ 
\end{eg}

Clearly any such (strict) alcove tableau terminates at a regular partition in the principal linkage class of the algebra.  
By definition, we have that there is   precisely one alcove-tableau
$\SSTP_ {\underline{w}}$ 
for each expression $\underline{w}$ in the simple reflections (and the emptyset). 
Similarly,  we have that there is   precisely one strict alcove-tableau
$\SSTP_ {\underline{w}}$ 
for each expression $\underline{w}$ in the simple reflections.

 \begin{eg}
Let $h=3$ and $\ell=1$ and $e=5$ and  $\al=\varepsilon_3-\varepsilon_1$.
We have that $\exx=3$.  
  We have that 
\begin{align*}\SSTP_{\al \emp} &=(  \varepsilon_1,  \varepsilon_2, 
  \varepsilon_1,  \varepsilon_2, 
  \varepsilon_1,  \varepsilon_2, 
   \varepsilon_1,  \varepsilon_1,  \varepsilon_1)
\otimes 
 (  \varepsilon_1,  \varepsilon_2, 
  \varepsilon_3,  \varepsilon_1, 
  \varepsilon_2,  \varepsilon_3, 
   \varepsilon_1,  \varepsilon_2 ,  \varepsilon_ 3)
\\ &=
  (  \varepsilon_1,  \varepsilon_2, 
  \varepsilon_1,  \varepsilon_2, 
  \varepsilon_1,  \varepsilon_2, 
   \varepsilon_1,  \varepsilon_1,  \varepsilon_1,
    \varepsilon_3,  \varepsilon_2, 
  \varepsilon_1,  \varepsilon_3, 
  \varepsilon_2,  \varepsilon_1, 
   \varepsilon_3,  \varepsilon_2 ,  \varepsilon_1)
 \\
 \SSTP_{\emp \al} &=(  \varepsilon_1,  \varepsilon_2,   \varepsilon_3, 
  \varepsilon_1,  \varepsilon_2,   \varepsilon_3, 
  \varepsilon_1,  \varepsilon_2  , \varepsilon_3, 
   \varepsilon_1,  \varepsilon_2,   
  \varepsilon_1,  \varepsilon_2,   
    \varepsilon_1,  \varepsilon_2,   
  \varepsilon_1,  \varepsilon_1,   \varepsilon_1 )
 \end{align*} 
are both dominant paths of shape   $  (3^3,2^3,1^3) $. 
\end{eg}

 \subsubsection{Permutations as morphisms between paths }\label{perms as morphisms}
We now discuss how one can  think of a permutation as a  morphism between pairs of paths in the alcove geometries of \cref{sec3}. 
\color{black}  This shift in perspective, from permutations  acting on tableaux (the usual combinatorics of $\mathfrak{S}_n$) to ``morphisms between paths" is a  central idea of this paper.

\color{black}

  \begin{defn}
Let $\la \in   \ZZ_{\geq0}\{\eps_1,\dots,\eps_h\}$.  
 Given a pair of paths $\SSTS,\SSTT\in \Path(\la)$ we write the steps $\varepsilon_i$ in 
 $\SSTS$ and $\SSTT$ in sequence along the   \north    and   \south       edges of a frame, respectively.
 We define   $w^{\SSTS}_{\SSTT}\in \mathfrak{S}_n$ to be the 
 unique step-preserving permutation   with the minimal number of crossings.  
  \end{defn}
  
 \color{black} 
Recall that a step  $\eps_i$ in a path corresponds to adding a node in the $i$th column (indexed from left to right) of the multi-partition tableau.  
Thus one can rewrite the above for  pairs of 
  {\em column standard} tableaux  as follows: $w^\SSTS_\SSTT$ is the unique element   such that $w^\SSTS_\SSTT (\SSTS)= \SSTT$ (under the usual action of the symmetric group on tableaux).  
 An example is given in \cref{egfortheref}.  
 
 \color{black}

\begin{eg}\label{concanetaion}
We consider $\Bbbk\mathfrak{S}_9$ in the case of $p=5$.  
We set $\al=\eps_3-\eps_1\in \Pi $.   
Here we have  $$\SSTP_\emp=(\eps_1,\eps_2,\eps_3,\eps_1,\eps_2,\eps_3,\eps_1,\eps_2,\eps_3)
\quad  \text{  
and }
\quad \reflectpath  =(\eps_1,\eps_2,\eps_1,\eps_2,\eps_1,\eps_2,\eps_3,\eps_3,\eps_3) $$
(the corresponding tableaux are given in \cref{figure1}).  
The unique step-preserving permutation of minimal length  is given by 
\begin{equation} 
w^{\SSTP_\emp}_
{\reflectpath }
= \begin{minipage}{6cm}
  
  \end{minipage} \end{equation}
  Notice that if two strands have the same step-label, then they do not cross.  This is, of course, exactly what it means for a step-preserving permutation to be  of minimal length.  \end{eg}
  
\begin{eg}\label{egfortheref} 
\color{black} 
We depict two paths $\SSTS,\SSTT\in \Path(11,6,1)$ in   \cref{paths-first} and the corresponding tableaux in \cref{tableau-second}.  The path-morphism $w^\SSTS_\SSTT$ is as follows
$$  w^\SSTS_\SSTT 
 = 
 \begin{minipage}{11cm}
 %
 \end{minipage}\;.$$
Notice that the sequence of $\eps_i$ along the top (bottom) of the word simply record the columns of the entries of the tableaux $\SSTS, \SSTT$ read in order according to the entries $1\leq i \leq 18$.  We always use 
$\eps_i$ as our labels of strand 
(dropping the epsilons would cause confusion later on, when we further attach KLR residues to these strands).

\end{eg}

   When we wish to explicitly  write down a specific reduced expression for  $w^{\SSTS}_{\SSTT}$ for concreteness, we will find the following notation incredibly useful.

\begin{defn}\label{whatabouti}
	Given $t$ an integer, we let $r_{\aatch}(t)$  denote the remainder of $t$ 
	modulo ${{\aatch}}$. 
Given $p,q \geq 1$ such that   $r_{\aatch}(p)\neq i$  and $\al=\eps_i-\eps_{i+1}\in\Pi$, we  set 
	$$
	\al({p}) =\SSTP_\al ^{-1}(1,r_{\aatch}(p))
 \qquad\text{and}\qquad	\emptyset({q}) =  \SSTP_\emptyset ^{-1}(1,r_{\aatch}(q))  	$$
 	This notation allows us to implicitly  use the cyclic ordering on the  labels of  roots without further ado.  
\end{defn}


\begin{conventioning}\label{conventioning}
  Throughout the paper, we let  $\al=\eps_i-\eps_{i+1}$, $\bet =\eps_j-\eps_{j+1}$, $\gam = \eps_k-\eps_{k+1}$, ${\color{orange}\delta}=\eps_m-\eps_{m+1}$.   
  We will assume that  $\bet,\gam,{\color{orange}\delta}$ label distinct
   commuting reflections.   We will assume  throughout that $\bet$ and $\al$ label non-commuting reflections.  
Here we read these subscripts in the obvious cyclotomic ordering, without further ado (in other words, we read occurrences of ${\aatch}+1$ simply as 1).  
\end{conventioning}

\section{The diagrammatic algebras}
\label{section-trois}
We now introduce the two  protagonists of this paper:
   the   {diagrammatic Bott--Samelson} endomorphism algebras and the quiver Hecke algebras ---  these   can be defined either
  as monoidal (tensor) categories or as finite-dimensional diagrammatic algebras.  
We favour the latter perspective for aesthetic reasons,   but we  borrow the notation from the former world  by letting $\otimes$ denote
  horizontal concatenation of diagrams --- in the quiver Hecke case, we must  first ``contextualise" before concatenating as we shall explain in \cref{truncation}. 
  (We refer to   \cite{BSBS} for a more detailed discussion of the interchangeability of
   these two languages.)  
 The relations for both algebras are entirely local (here a local relation means one that 
 can be specified by its effect on a  sufficiently small region of the wider diagram).  
 We then consider the cyclotomic quotients of these algebras: these can be viewed as quotients by right-tensor-ideals, or equivalently (as we do in this paper) as quotients by a {\em non-local} diagrammatic relation concerning the leftmost strand in the ambient concatenated diagram.  (We remark that the cyclotomic relations break the monoidal structure of both   categories.)   
We continue with the notation of \Cref{conventioning}.

  \begin{rmk} \color{black} 
The cyclotomic quotients of (anti-spherical) Hecke categories are small categories  
with finite-dimensional morphism spaces given by the light leaves basis of \cite{MR3555156,antiLW}.     
Working with such a category is equivalent working to with a locally unital algebra, as defined in \cite[Section 2.2]{BSBS}, see  \cite[Remark 2.3]{BSBS}.  
   Throughout this paper we will work in the latter setting.  
   The reader who prefers to think of categories can equivalently phrase everything in the this paper  in terms of categories and representations of categories.

   \end{rmk}

   \subsection{The   {diagrammatic Bott--Samelson} endomorphism algebras}  \label{soergel}
These algebras were  defined  by Elias--Williamson in \cite{MR3555156}.  
 In this section, all our words will be in the alphabet $S$.

\begin{defn}
Given  $\al=\eps_i-\eps_{i+1}$  we define the corresponding Soergel idempotent, 
${\sf 1}_{ {\al}}$   to be a frame of width 1 unit, containing a single vertical strand   coloured with $\al  \in \Pi  $.   
For $\w=s_{\alpha^{(1)}}\mydots s_{\alpha^{(p)}}$ an expression with 
  $\alpha^{(i)}\in \Pi  $ simple roots, we set 
  $${\sf 1}_\w= {\sf 1}_{\alpha^{(1)}}\otimes 
  {\sf 1}_{\alpha^{(2)}}\otimes \dots 
\otimes    {\sf 1}_{\alpha^{(p)}} $$
to be the diagram obtained by horizontal concatenation.  \end{defn}

\begin{eg}
Consider the colour-word from the path in \cref{diag2}.  Namely,  $$\w= 
 s_{\color{magenta}\varepsilon_3-\varepsilon_1}
 s_{\color{darkgreen}\varepsilon_2-\varepsilon_3}
 s_{\color{cyan}\varepsilon_1-\varepsilon_2}
 s_{\color{magenta}\varepsilon_3-\varepsilon_1}
 s_{\color{darkgreen}\varepsilon_2-\varepsilon_3}
 s_{\color{cyan}\varepsilon_1-\varepsilon_2}
 \in\widehat{\mathfrak {S}}_3 .
$$ The  corresponding Soergel idempotent is  as follows 
$$ 
 {\sf 1}_{\w}=\begin{minipage}{4.5cm}
\end{minipage} 
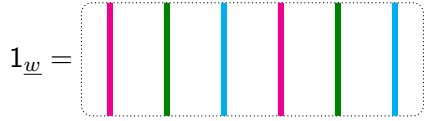
 $$
\end{eg}

\begin{defn} Given $\w=s_{\alpha^{(1)}}\mydots s_{\alpha^{(p)}},\w'=s_{\beta^{(1)}}\mydots s_{\beta^{(q)}}\in \mathfrak{S}_{\aatch}$, a $(\w,\w')$-Soergel diagram  $D$ 
is defined to be any diagram obtained by horizontal and vertical concatenation of the 
following diagrams $$
 %
$$
their   flips through the horizontal axis and their isotypic deformations  such that the    \north    and    \south edges  of the graph are given by the idempotents 
${\sf 1}_{\underline{w}}$ and  
${\sf 1}_{\underline{w}'}$ respectively.  
 Here the vertical concatenation of a $(\w,\w')$-Soergel diagram on top of a 
$(\underline{v},\underline{v}')$-Soergel diagram is zero if $\underline{v}\neq \w'$.  
We define the degree of these generators (and their   flips) to be $0, 1, -1, 0$, and $0$ respectively.  
    \end{defn}

\begin{eg}
Examples  of   $(\w,\w')$-Soergel diagrams, for $$\w= 
 s_{\color{magenta}\varepsilon_3-\varepsilon_1}
 s_{\color{darkgreen}\varepsilon_2-\varepsilon_3}
 s_{\color{cyan}\varepsilon_1-\varepsilon_2}
  s_{\color{magenta}\varepsilon_3-\varepsilon_1}
 s_{\color{darkgreen}\varepsilon_2-\varepsilon_3}
 s_{\color{cyan}\varepsilon_1-\varepsilon_2}
  s_{\color{cyan}\varepsilon_1-\varepsilon_2},  
$$
$$
\w' =    s_{\color{magenta}\varepsilon_3-\varepsilon_1} 
 s_{\color{darkgreen}\varepsilon_2-\varepsilon_3}
 s_{\color{cyan}\varepsilon_1-\varepsilon_2} 
 s_{\color{magenta}\varepsilon_3-\varepsilon_1}
 s_{\color{darkgreen}\varepsilon_2-\varepsilon_3}
 s_{\color{cyan}\varepsilon_1-\varepsilon_2}    $$
are as follows
$$ 
 \begin{tikzpicture}[scale=1.2]
\draw[densely dotted, rounded corners] (-0.25,0) rectangle (3.25,1);
\draw[magenta,line width=0.08cm](0,1)to [out=-90,in=90] (0.5,0);
\draw[darkgreen,line width=0.08cm](0.5,1) to [out=-90,in=90] (1,0);
\draw[cyan,line width=0.08cm,fill=cyan](1,1)  --++(-90:0.35) circle (2pt);

\draw[cyan,line width=0.08cm](1.5,0)--++(90:1);
\draw[magenta,line width=0.08cm](2,0)--++(90:1);
\draw[darkgreen,line width=0.08cm](2.5,0)--++(90:1);
\draw[cyan,line width=0.08cm](3,0)--++(90:1);
\end{tikzpicture}
\qquad\qquad
 \begin{tikzpicture}[scale=1.2]
\draw[densely dotted, rounded corners] (-0.25,0) rectangle (3.25,1);
\draw[magenta,line width=0.08cm](0,1)to [out=-90,in=90] (0.5,0);
\draw[darkgreen,line width=0.08cm](0.5,1) to [out=-90,in=90] (1,0);
\draw[cyan,line width=0.08cm ](1,1) to  [out=-90,in=150] (1.5,0.5);

\draw[cyan,line width=0.08cm](1.5,0)--++(90:1);
\draw[magenta,line width=0.08cm](2,0)--++(90:1);
\draw[darkgreen,line width=0.08cm](2.5,0)--++(90:1);
\draw[cyan,line width=0.08cm](3,0)--++(90:1);
\end{tikzpicture}
 $$ 
 \end{eg}

 We let $\ast$ denote the  map  which   flips a diagram   through its  horizontal axis.

\begin{defn}\label{ranknsoergl}
Let $\Bbbk$ be an arbitrary commutative ring. 
We define the  
  {diagrammatic Bott--Samelson} endomorphism algebra, $\mathscr{S} ({n,\sigma})$ to be the span of all $(\w,\w')$-Soergel diagrams for $\w,\w' \in \Lambda(n,\sigma)$, with $\Bbbk$-associative multiplication given by vertical concatenation and 
subject to   isotypic deformation  and the following local relations:  
 For each colour (i.e. each generator   $s_\al$ for 
 ${\al}\in \Pi $) we have
 \begin{equation} \tag{S1} \label{rel1}
 \begin{minipage}{1.3cm}
\end{minipage}=\; 0
\end{equation}along with their horizontal and vertical flips and the Demazure relation 
 \begin{equation}\tag{S2}  \label{rel2}
\qquad \begin{minipage}{1.3cm}%
\end{minipage} 
 \end{equation} 
 We now picture the two-colour relations for  
non-commuting reflections $s_\al, s_{\bet}\in \Shl$.  We have 
\begin{align}\tag{S3} \label{rel3}
 & \begin{minipage}{1.3cm}%
\end{minipage} \end{align}
along with their flips through the horizontal  and vertical axes. 
 We also have the cyclicity relation \begin{align}
 \tag{S4}  \label{rel4}
 \begin{minipage}{2.8cm}%
\end{minipage}
\end{align}
 and the two-colour barbell relations
\begin{align}\tag{S5}\label{rel5} \begin{minipage}{1.3cm}%
\end{minipage}
\end{align} 
for $\Phi $ of rank greater than 1 (or double the righthand-side if $\Phi $ has rank 1).
For commuting reflections 
 $s_{\bet}, s_{\color{darkgreen}\gamma}\in \Shl$  we have the following relations 
\begin{equation}\tag{S6} \label{rel6}
   \begin{minipage}{1.3cm}%
 \end{minipage}
 \end{equation} along with  their flips through the horizontal  and vertical axes.  
 In order to picture the three-colour commuting relations we require 
 a fourth root  
 $s_{\color{orange}\boldsymbol\delta}\in \Shl$ which commutes with all other roots (such that
  $s_{\color{orange}\boldsymbol\delta}s_\al 
  =
  s_\al  s_{\color{orange}\boldsymbol\delta}$,
   $s_{\color{orange}\boldsymbol\delta}s_\bet
  =
  s_\bet  s_{\color{orange}\boldsymbol\delta}$,
   $s_{\color{orange}\boldsymbol\delta}s_\gam 
  =
  s_\gam  s_{\color{orange}\boldsymbol\delta}$)
  and we have the following,
\begin{align}\tag{S7}\label{S7}
 \begin{minipage}{1.3cm}
\end{minipage}\end{align}Finally, 
 we require     the tetrahedron relation for which we make the additional assumption on $\gam$ that it does not commute with $\al$.  This relation is as follows, 
\begin{equation}\tag{S8} \label{dragrace}
 \begin{minipage}{4.1cm}%
 \end{minipage}
\end{equation}
\end{defn}

\begin{rmk} 
\color{black}
%
%
%
The diagrammatic Bott--Samelson category of $\widehat{\mathfrak{S}}_h$ is normally defined using an underlying reflection representation $\mathfrak{h}=(V, \{\alpha_{\al }^\vee : \al  \in S\}, \{\alpha_{\al } : \al  \in S\})$ of $\widehat{\mathfrak{S}}_h$ called a {\sf realisation}. 
Our construction of the diagrammatic Bott--Samelson endomorphism algebra implicitly assumes that the roots $\{\alpha_{\al} : \al \in S\} \subset V^\ast$ form a basis, and that the pairing between roots and coroots is given by the usual Cartan matrix of type $\widehat{A}_{h-1}$. 
These two conditions uniquely determine the realisation, which we call the {\sf universal realisation} of $\widehat{\mathfrak{S}}_{h}$ with respect to this Cartan matrix \cite{withemily}.
It coincides with the modular reduction of the ``dual geometric'' realisation of $\widehat{\mathfrak{S}}_h$ (which can be defined over $\ZZ$ as $\widehat{\mathfrak{S}}_h$ is simply laced) \cite{antiLW}.
%
%
%
%
\color{black}
\end{rmk}

\begin{rmk}\color{black} 
We do not include ``isotopy" as an explicit relation here (unlike in \cite{MR3555156}) 
as it follows from the one-colour relations and   cyclicity of the braid generator 
(see \cite[Proposition 8.6]{MR4220642}).  This is the more modern definition, see for example \cite[Section 2.3]{modern}
 \end{rmk}

\begin{defn}
We define the 
 cyclotomic {diagrammatic Bott--Samelson} endomorphism algebra, $$ \mathscr{S}_{\underline{h}}({n,\sigma}) := {\rm End}_{\mathcal{D}^{\rm asph,\oplus}_{\rm BS} ( A_{h-1}\times \mydots \times A_{h-1}\backslash \widehat{A}_{{\aatch}-1} )}\left(\oplus_{\underline{w} \in \Lambda({n,\sigma})} B_{\underline{w}} \right) 
 $$ to be the quotient of $\mathscr{S}(n,\sigma)$ by the relations 
 \begin{equation}\tag{S9} \label{cyclotomic}{\sf1}_{\al}\otimes {\sf1}_{\w}=0\qquad\quad \text{and} \qquad \quad 
  \begin{minipage}{0.7cm}
 \begin{tikzpicture}[scale=0.8]
\draw[densely dotted,rounded corners](-0.4cm,-0.62cm)  rectangle (0.4cm,0.62cm);

 \draw[line width=0.08cm, darkgreen](0,-0.4)--(0,0.4);
\fill[darkgreen] (0,0.35)  circle (5pt);
\fill[darkgreen] (0,-0.35)  circle (5pt);
\end{tikzpicture} \end{minipage} 
\otimes {\sf1}_{\w}=0
\end{equation}
for $\gam\in\Pi$ arbitrary,    $\al\in \Pi$ corresponding to a wall of the dominant chamber,   and $\w$ any word in the alphabet $S$.  
  \end{defn}

 \subsection{The breadth-enhanced {diagrammatic Bott--Samelson} endomorphism algebra}
We now  use the notion of a weakly graded monoidal category (see \Cref{Appendix}) to 
    introduce the breadth-enhanced 
 {diagrammatic Bott--Samelson} endomorphism algebra. 
  On one level this definition and construction is utterly superficial.  
  It merely allows us to keep track of occurrences of the identity of  $\Shl$ in a given expression.  
   The  occurrences  of $s_\emptyset =1$ 
are usually   ignored in the world of Soergel diagrammatics and so   this will seem very foreign to some.  
 We ask these readers to be patient as    this extra ``blank space" will be very important in this paper: each occurrence of $s_\emptyset$ corresponds to adding ${\aatch}$ additional strands in the quiver Hecke algebra or, if you prefer, corresponds to ``tensoring with the determinant".    
For this reason, in this section  all our words will be in the alphabet $S\cup\{1\}$.  
\begin{defn}
Given  $\al=\eps_i-\eps_{i+1}$  we   define the breadth-enhanced  Soergel idempotent, 
${\sf 1}_{ {\al}},$   to be a frame of width $2b_{\al }$ with a single  vertical strand   coloured with $\al  \in \Pi  $ placed in the centre.     
We define the breadth-enhanced  Soergel idempotent 
${\sf 1}_{ {\emptyset}}$   to be an empty frame of width $2 $.     
For $\w=s_{\alpha^{(1)}}\mydots s_{\alpha^{(p)}}$ an expression with 
  $\alpha^{(i)}\in \Pi   \cup \{\emptyset\}$, we set 
  $${\sf 1}_\w= {\sf 1}_{\alpha^{(1)}}\otimes 
  {\sf 1}_{\alpha^{(2)}}\otimes \dots 
\otimes    {\sf 1}_{\alpha^{(p)}} $$
to be the diagram obtained by horizontal concatenation.   
In order that we better illustrate this idea, we colour the   \north    and   \south       edges of a frame with the corresponding element of $\Pi \cup\{\emptyset\}$.    
 \end{defn}

\begin{eg}\label{continuity}
 Continuing with \cref{diag2,diag22point4}, we    
let   $$\w= 
s_\emptyset s_\emptyset s_{\color{magenta}\varepsilon_3-\varepsilon_1}
 s_{\color{darkgreen}\varepsilon_2-\varepsilon_3}
 s_{\color{cyan}\varepsilon_1-\varepsilon_2}
 s_{\color{magenta}\varepsilon_3-\varepsilon_1}
 s_{\color{darkgreen}\varepsilon_2-\varepsilon_3}
 s_{\color{cyan}\varepsilon_1-\varepsilon_2}
$$  $$ \w'= 
  s_{\color{magenta}\varepsilon_3-\varepsilon_1}s_\emptyset
 s_{\color{darkgreen}\varepsilon_2-\varepsilon_3}
 s_{\color{cyan}\varepsilon_1-\varepsilon_2} s_\emptyset
 s_{\color{magenta}\varepsilon_3-\varepsilon_1}
 s_{\color{darkgreen}\varepsilon_2-\varepsilon_3}
 s_{\color{cyan}\varepsilon_1-\varepsilon_2}
 .$$   
  The breadth-enhanced Soergel idempotents  are as follows
\begin{equation}\label{twoidemps}
 {\sf 1}_{\w}=\begin{minipage}{4.5cm}
\end{minipage}
\end{equation}
\end{eg}

\begin{defn}
Let   $w\in  \mathfrak{S}_{\aatch}$ and  suppose  $\underline{w}=s_{\alpha^{(1)}}\mydots s_{\alpha^{(p)}} $ and  $\underline{w}'=s_{\beta^{(1)}}\mydots s_{\beta^{(p)}} $ for $\alpha^{(i)},\beta^{(j)}\in \Pi \cup\{\emptyset\}$ are two expressions which differ only by occurrences of $s_\emptyset$ within the word.   We define the corresponding Soergel adjustment
 ${\sf 1}_{\underline{w}'}^{\underline{w}}$, to be the diagram with 
 ${\sf 1}_{\underline{w}}$ along the top  
 and  ${\sf 1}_{\underline{w}'}$ along the bottom and no crossing strands.  
\end{defn}

\begin{eg}
 Continuing with \cref{continuity}, we have that 
$$ {\sf 1}_{\w'}^{\w} =\begin{minipage}{5cm}
\end{minipage}
   $$
\end{eg}

\begin{defn} Given $\w=s_{\alpha^{(1)}}\mydots s_{\alpha^{(p)}},\w'=s_{\beta^{(1)}}\mydots s_{\alpha^{(q)}} $ for  
$\alpha^{(i)},\beta^{(j)}\in \Pi \cup\{\emptyset\}$,
 a breadth-enhanced $(\w,\w')$-Soergel diagram  $D$ 
is defined to be any diagram obtained by  horizontal and vertical concatenation of 
the following generators 
\begin{equation}\label{Soergelgens} \begin{minipage}{1.1cm} %
\end{minipage}  \end{equation}
and     their flips through the horizontal  axes   such that the    \north    edge of the graph is given by the breadth-enhanced idempotent 
${\sf 1}_{\underline{w}}$ and the   \south       edge given by the breadth-enhanced idempotent 
${\sf 1}_{\underline{w}'}$.  
Here the vertical concatenation of a $(\w,\w')$-diagram on top of a 
$(\underline{v},\underline{v}')$-diagram is zero if $\underline{v}\neq \w'$.  
The degree of these generators (and their   flips)     are  $0,0,0,  1, -1, 0$, and $0$ respectively.  
 When we wish to avoid drawing diagrams, we will denote the above diagrams by
$${\sf 1}_\al \quad   {\sf 1}_\emptyset \quad  {\sf 1}_{\emptyset\al}^{\al\emptyset}
\quad  {\sf SPOT}_\al^\emp 
\quad  {\sf FORK}_{\al\al}^{\emp\al} 
\quad {\sf HEX}_{\al\bet\al}^{\bet\al\bet} 
\quad \text{ and }\quad  {\sf COMM}_{\bet\gam}^{\gam\bet}.$$   
These diagrams are known as  ``single strand", ``blank space", ``single adjustment", ``spot", ``fork", 
``hexagon" (in order to distinguish from the symmetric group braid)
and   ``commutator".    

\end{defn}

\begin{defn} 
We define the breadth-enhanced diagrammatic Bott--Samelson endomorphism algebra, 
$\mathscr{ S}^{\rm br}({n,\sigma})$
 (respectively, its cyclotomic quotient  
 $\mathscr{S}^{\rm br}_{\aatchpair}({n,\sigma})$)  to be the span of all $(\w,\w')$-breadth enhanced Soergel diagrams for $\w,\w' \in \Lambda(n,\sigma)$, with multiplication given by vertical concatenation, subject to the breadth-enhanced analogues of the relations \ref{rel1} to \ref{dragrace}
(plus the additional cyclotomic relation \ref{cyclotomic}, respectively) which are explicitly pictured in \cref{relations}, 
the adjustment inversion and naturality   relations pictured in \cref{nature,nature2} and their flips through the horizontal axis.  

\end{defn}

 \begin{figure}[ht!]
 $$  
 \begin{minipage}{1.1cm} 

\end{minipage}
$$
\caption{The adjustment-inversion relations and the naturality relation for the spot diagram (we also require their flips through horizontal axis).}
\label{nature2}
\end{figure}

 \begin{figure}[ht!]
$$
   \begin{minipage}{1.8cm}\; %
\end{minipage}
$$
\caption{The remaining naturality relations (we also require their flips through horizontal axis).}
\label{nature}
\end{figure}

 We are free to use the breadth-enhanced form of the diagrammatic Bott--Samelson endomorphism algebra instead of the usual one because of the following result.
 We  let  $\phi: \cup_{0\leq m \leq n}\Lambda^+(m,\sigma)\hookrightarrow \Lambda(n,\sigma) $
denote the map which takes $\w\in \Lambda^+ ({m,\sigma})$ to $(s_\emptyset)^{n-m} \w\in \Lambda ({n,\sigma})$.  
We refer to the image,  ${\rm im}(\phi)=\Lambda^+( \leq n,\sigma)$, 
as the subset of {\sf left-adjusted} words in $\Lambda(n,\sigma)$ and we define an associated idempotent 
\begin{equation*}
{\sf 1}_{n, \sigma}^+=\sum
 _{\w \in \Lambda^+(\leq n,\sigma)}{\sf 1}_{\w}.  
\end{equation*}
 \begin{prop}\label{amit}
We have the following isomorphisms of graded $\Bbbk$-algebras, 
\begin{align*}
\mathscr{S}(n,\sigma)   \cong  {\sf 1}_{n,\sigma}^+\mathscr{S}^{\rm br}(n,\sigma){\sf 1}_{n,\sigma}^+   \qquad \mathscr{S}_{\underline{h}}(n,\sigma)   \cong  {\sf 1}_{n,\sigma}^+\mathscr{S}_{\underline{h}}^{\rm br}(n,\sigma){\sf 1}_{n,\sigma}^+.
\end{align*}
 
\end{prop}

\begin{proof}
This is the one point in the paper in which we require the notions from \Cref{Appendix} and all references within this proof are to the appendix.  
Thus for this proof only, we   briefly switch perspectives and think of the algebras above as {\em categories} $\mathscr{S}$ and $\mathscr{S}^{\rm br}$ and use the notation in \Cref{Appendix}. 
The category $\mathscr{S}$ (resp.~$\mathscr{S}^{\rm br}$) has objects given by expression in the alphabet $S$ (resp.~$S \cup \{1\}$) and $\Hom$-spaces given by ${\sf 1}_{\w}\mathscr{S}(n,\sigma){\sf 1}_{\w'}$ (resp.~${\sf 1}_{\w}\mathscr{S}(n,\sigma){\sf 1}_{\w'}$) for some sufficiently large $n$ (resp.~for some $n$).

We will establish the first isomorphism; the second isomorphism is similar.   
Let $b : \mathrm{Ob}(\mathscr{S}) \rightarrow \ZZ_{\geq 0}$ be a monoidal homomorphism given by $b(s_{\al})=b_{\al}$ for all $\al \in \Pi$. 
We now apply Theorem~\ref{thm:wkgradingcoherence} to show that $\mathscr{S}^{\rm br}(n,\sigma)$ is isomorphic to the weak grading of $\mathscr{S}(n,\sigma)$ concentrated in breadth $b$.
Most of the hypotheses of this result follow by design. 
For example, since $\mathscr{S}$ is already defined by generators and relations, it's enough to add breadth-enhanced versions of the relations to ensure the composition and tensor product axioms in the theorem.
Also, adjustments on objects are defined monoidally, so the weak grading axioms \ref{item:mult} and \ref{item:identity} automatically hold.
Finally \ref{item:naturality} follows from the adjustment inversion and naturality relations above.
\end{proof}

 \subsection{The quiver Hecke algebra}\label{quiver}

We now introduce the second star of the paper, the  quiver Hecke or KLR algebras.  
Given $\underline{i}=(i_1,\dots, i_n)\in (\ZZ/e\ZZ)^n$ and $s_r=(r,r+1)\in \mathfrak{S}_n$ we set 
 $s_r(\underline{i})= (i_1,\dots, i_{r-1}, i_{r+1}, i_r ,i_{r+2}, \dots ,i_n)$.  
 
\begin{defn}[\cite{MR2551762,MR2525917,ROUQ}]\label{defintino1}
Fix $e > 2 $.  
The {\sf quiver Hecke algebra} (or KLR algebra),  $\mathcal{H}_n $,   is defined to be the unital, associative $\ZZ$-algebra with generators
$$ 
\{e_{\underline{i}}  \ | \ {\underline{i}}=(i_1,\mydots,i_n)\in   (\ZZ/e\ZZ)^n\}\cup\{y_1,\mydots,y_n\}\cup\{\psi_1,\mydots,\psi_{n-1}\},
$$
subject to the  following  relations, 
for all $r,s,{\underline{i}},{\underline{j}}$ we have that  
\begin{align}
\tag{R1}  
\textstyle\sum
 e_{{\underline{i}}} =1 _{\mathcal{H} _n};  
 \;\;\;    e_{\underline{i}} e_{\underline{j}} =\delta_{{\underline{i}},{\underline{j}}} e_{\underline{i}} 
\;\;\; y_r		e_{{\underline{i}}}=e_{{\underline{i}}}y_r 
\;\;\;
\psi_r e_{{\underline{i}}} = e_{s_r{\underline{i}}} \psi_r 
\;\;\;
y_ry_s =y_sy_r
 \end{align}
 where the sum is over all ${{\underline{i}} \in   (\ZZ/e\ZZ )^n } $ 
  and 
\begin{align}
\tag{R2}\label{rel1.7}
\psi_ry	_s  = y_s\psi_r \  \text{ for } s\neq r,r+1&
\qquad  
&\psi_r\psi_s = \psi_s\psi_r \ \text{ for } |r-s|>1
 \\ \tag{R3}\label{rel1.8} 
y_r \psi_r e_{\underline{i}}
  =
(\psi_r y_{r+1} \color{black}-\color{black} 
\delta_{i_r,i_{r+1}})e_{\underline{i}}  &
 \qquad 
&y_{r+1} \psi_r e_{\underline{i}}   =(\psi_r y_r \color{black}+\color{black} 
\delta_{i_r,i_{r+1}})e_{\underline{i}} 
\end{align}
\begin{align}
\tag{R4}\label{rel1.10}
\psi_r \psi_r  e_{\underline{i}} &=\begin{cases}
0 & \text{if }i_r=i_{r+1},\\
e_{\underline{i}}  & \text{if }i_{r+1}\neq i_r, i_r\pm1,\\
(y_{r+1} - y_r) e_{\underline{i}}  & \text{if }i_{r+1}=i_r +1  ,\\
(y_r - y_{r+1}) e_{\underline{i}} \qquad\qquad& \text{if }i_{r+1}=i_r -1   
\end{cases}\\
\tag{R5}\label{rel1.11}
\psi_r \psi_{r+1} \psi_r &=\begin{cases}
(\psi_{r+1}\psi_r\psi_{r+1} - 1)e_{\underline{i}}\qquad\qquad& \text{if }i_r=i_{r+2}=i_{r+1}+1  ,\\
(\psi_{r+1}\psi_r\psi_{r+1} + 1)e_{\underline{i}}& \text{if }i_r=i_{r+2}=i_{r+1}-1    \\
 \psi_{r+1} \psi_r\psi_{r+1}e_{\underline{i}}   &\text{otherwise} 
\end{cases} 
\end{align}
 for all  permitted $r,s,i,j$.  We identify such elements with decorated permutations  and the multiplication with vertical concatenation, $\circ$,  of these diagrams  in the  standard fashion of \cite[Section 1]{MR2551762}. 
   We let $\ast$ denote the anti-involution which fixes the generators (this can be visualised as a flip through the horizontal axis of the diagram).  
 
\end{defn}

We identify an undecorated single strand with the sum over all possible residues on that strand, as in  
  $\sum_{{\underline{i}} \in   (\ZZ/e\ZZ )^n } e_{{\underline{i}}} =1 _{\mathcal{H} _1}$.  
We freely identify an element   $d\in \mathcal{H}_n $ with
an element 
  of $\mathcal{H}_{n+1}$  
  by adding such an  undecorated vertical strand to the right; we extend this to all $ \mathcal{H}_m $ with  $m>n$.   
The $y_k$ elements are visualised as dots on strands; we hence refer to them as {\sf KLR dots}.
   Given $\stt \in \Std(\la)$, we set $e_\SSTT:= e_{{\rm res} (\SSTT)}\in \mathcal{H}_n $.  Using the notation of \cref{necessaryevil}, we define 
\begin{align}  
\label{tableau not} 
y_{\SSTT}  = \prod _{k=1}^n y_k^{ |{\mathcal A} _{\SSTT }(k)|}e_{\SSTT},  
\end{align} 
  such elements should be familiar to those working in KLR algebras, see for example \cite[Section 4.3]{hm10}.
Given   $p<q$ we set 
$$
w^p_q= s_p s_{p+1}\dots s_{q-1}\quad
w^q_p=  s_{q-1} \dots s_{p+1} s_{p }
\qquad \psi^p_q= \psi_{p}\psi_{p+1}\mydots \psi_{q-1} 
\qquad \psi^q_p=   \psi_{q-1}  \mydots
 \psi_{p+1}\psi_{p}.
 $$  
and  given an expression $\w=s_{i_1}\dots s_{i_p}\in \mathfrak{S}_n$ we set 
$\psi_\w= \psi_{i_1}\dots \psi_{i_p}\in \mathcal{H}_n$.      

\begin{defn}Fix $e > 2 $ and $\sigma\in  \ZZ ^\ell$.  
The {\sf cyclotomic quiver Hecke algebra}, $\mathcal{H}_n^\sigma$,   is defined to be the
 quotient of $\mathcal{H}_n$ by the relation
\begin{align}
\label{rel1.12} y_1^{\sharp\{\sigma_m | \sigma_m= i_1 ,1\leq m \leq \ell 	\}} e_{\underline{i}} &=0 
\quad 
\text{ for    ${\underline{i}}\in (\ZZ/e\ZZ)^n$.}
\end{align}  
 \end{defn}

\begin{defn}We define the degree    on the generators  as follows, 
$$
{\rm deg}(e_{\underline{i}})=0 \quad 
{\rm deg}( y_r)=2\quad 
{\rm deg}(\psi_r e_{\underline{i}})=
\begin{cases}
-2		&\text{if }i_r=i_{r+1} \\
1		&\text{if }i_r=i_{r+1}\pm 1 \\
 0 &\text{otherwise} 
\end{cases}. $$
\end{defn}

 \begin{defn}
Given a  pair of  paths  $\SSTS,\SSTT \in \Path(\la)$,  
 and a fixed choice of reduced expression for $w^\SSTS_\SSTT=s_{i_1}s_{i_2}\dots s_{i_k}$ we define 
 $ \psi ^\SSTS_\SSTT=e_\SSTS \psi_{i_1}\psi_{i_2}\dots \psi_{i_k}e_\SSTT$.  
 \end{defn}

 \!
 \begin{rmk}
 The quiver Hecke algebra and its cyclotomic quotients are isomorphic (over a field) to the classical affine Hecke algebra 
 and its cyclotomic quotients (at a root of unity) by \cite[Main Theorem]{MR2551762}.   
 Setting $e=p$ and $\sigma=(0)\in\ZZ^1$ we have that $\Bbbk\mathfrak{S}_n$ is isomorphic to $\mathcal{H}_n^\sigma$ and we freely identify these algebras  
 without further mention.  
\end{rmk}

     \begin{figure}[ht!]$$
 \begin{minipage}{6.5cm}
 \begin{tikzpicture}
  [xscale=0.7,yscale=  
 - 0.7]

  \draw  (0.2+14.8-2  ,4.5)  node  {$ \SSTP_{{\al}} ^\flat $};

     \draw  (0.2+14.8-2 ,3) node {$  \SSTP_\emp  $};
     \draw(3,3) rectangle (12,4.5);
        \foreach \i in {3.5,4.5,...,11.5}
  {
   \fill(\i,3) circle(1.5pt) coordinate (a\i);
          \fill(\i,4.5)circle(1.5pt)  coordinate (d\i);
    }

\draw(5.5,3) --(9.5,4.5);
\draw (4+4.5,3) --(10.5,4.5);

\draw(3.5,3) --(3.5,4.5);
\draw(4.5,3) --(4.5,4.5);
\draw(6.5,3)  --(5.5,4.5);

\draw(3+4.5,3) --(6.5,4.5);

\draw(11.5,3) --(11.5,4.5);

 \draw (4+5.5,3)  --(4+3.5,4.5);
\draw(4+6.5,3) --(4+4.5,4.5);

 \scalefont{0.9} 

\draw (3.5,4.5) node[below] {$ 0$};
\draw (4.5,4.5) node[below] {$ 1$};
\draw (5.5,4.5) node[below] {$ 4$};
\draw (6.5,4.5) node[below] {$ 0$};
\draw (7.5,4.5) node[below] {$ 3$};
\draw (8.5,4.5) node[below] {$ 4$};
\draw (9.5,4.5) node[below] {$ 2$};
\draw (10.5,4.5) node[below] {$ 1$};
\draw (11.5,4.5) node[below] {$ 0$};

\draw (3.5,3) node[above] {$ 0$};
\draw (4.5,3) node[above] {$ 1$};
\draw (5.5,3) node[above] {$ 2$};
\draw (6.5,3) node[above] {$ 4$};
\draw (7.5,3) node[above] {$ 0$};
\draw (8.5,3) node[above] {$ 1$};
\draw (9.5,3) node[above] {$ 3$};
\draw (10.5,3) node[above] {$ 4$};
\draw (11.5,3) node[above] {$ 0$};

  \end{tikzpicture} 
 \end{minipage}
 $$

\!\!\!
  \caption{The element $  \psi ^{\SSTP_\emp}_
{\reflectpath }
$ for    $\Bbbk\mathfrak{S}_9$ in the case $p=5$ and  $\al=\eps_3-\eps_1\in \Pi $ (see also  \cref{concanetaion}). }
  \end{figure}
  \!\!

\!\!\!
\subsubsection{Our   quotient algebra and regular blocks} \label{quotient}
A long-standing belief in modular Lie theory  is that we should (first) restrict our attention to fields whose characteristic, $p$, is greater than the Coxeter number, $ h$, of the algebraic  group we are studying.   
This allows one to consider a ``regular" or ``principal block"     of 
the  algebraic group in question.
For example, the {diagrammatic Bott--Samelson} endomorphism algebras categorify the endomorphisms between  
 tilting modules for the {\em principal} block of the algebraic group, ${\rm GL}_h(\Bbbk)$, and this is the crux of the proof of  \cite[Theorem 1.9]{MR3805034}.    
   Extending this ``Soergel diagram calculus" 
  to singular blocks is a difficult problem.   
As such, all   results in
 \cite{MR3805034,MR3868004} and  
  this paper are restricted to {\em regular} blocks.  
 In the language of 
 \cite{MR3805034,MR3868004} this   restricts the study of the algebraic group in question to primes $p>h$.

 What does this mean on the other side of the Schur--Weyl duality  relating ${\rm GL}_h(\Bbbk)$ and $\Bbbk\mathfrak{S}_n$?  
By the second fundamental theorem of invariant theory, the kernel of the group algebra of the symmetric group acting on $n$-fold $h$-dimensional tensor space is the element 
$ 
\textstyle \sum _{g \in \mathfrak{S}_{h+1}\leq \mathfrak{S}_n} {{\rm sgn}(g)} g\in \Bbbk\mathfrak{S}_n
$.   Modulo ``more dominant terms"  this  element is equal to 
 $ {\sf y}_{ \SSTT^{(h+1)} }$ (the element introduced in \cref{tableau not}).  
       The module category of  $\Bbbk\mathfrak{S}_n/ \Bbbk\mathfrak{S}_n{\sf y}_{ \SSTT^{(h+1)} }\Bbbk\mathfrak{S}_n$ 
       is     the Serre subcategory of 
$\Bbbk\mathfrak{S}_n$-${\rm mod}$ whose 
 simple  modules   are indexed by 
partitions with at most $h$ columns. 
      For   $p>h$, the algebra $\Bbbk\mathfrak{S}_n/ \Bbbk\mathfrak{S}_n{\sf y}_{ \SSTT^{(h+1)} }\Bbbk\mathfrak{S}_n$ 
    is the largest quotient of $\Bbbk\mathfrak{S}_n$ controlled by  the  
  {diagrammatic Bott--Samelson} endomorphism algebra with $h$ distinct colours.  
    Combinatorially, the condition that $p>h$ ensures that $\varnothing$ does not lie on any hyperplane in the alcove geometry (and so the $p$-Kazhdan--Lusztig theory is ``regular" not ``singular"). 
     The  importance of this 
 Serre subcategory and the condition $p>h$
  can also be explained in the context of  calibrated/unitary modules \cite[Introduction]{BNS2}.  
  The main theorem of 
       \cite{MR3805034} calculates 
       decomposition numbers  of $\Bbbk\mathfrak{S}_n/ \Bbbk\mathfrak{S}_n{\sf y}_{ \SSTT^{(h+1)} }\Bbbk\mathfrak{S}_n$.

   There is a canonical manner in which the above situation    generalises to  cyclotomic  Hecke algebras.  
 For a given $e>h$, one can ask {\em ``what is the largest quotient of $\mathcal{H}_n^\sigma$ controlled by the {diagrammatic Bott--Samelson} endomorphism algebra with $h$ distinct colours?''}    Assuming that $\aatchpair\in \ZZ_{>0}^\ell$  is $(\sigma,e)$-admissible, we 
   define 
   $$ {\sf y}_{\underline{h} }=
 \sum_{\begin{subarray}c
\alpha= (\varnothing,\dots,\varnothing,(h_a+1),\varnothing,\dots,\varnothing)
   \\
 0\leq a <\ell    \end{subarray}} 
 y_{\SSTT^{\alpha}}
$$
and we  claim that the answer to the question is provided by the quotient algebras 
$\mathcal{H}^\sigma_n/ \mathcal{H}^\sigma_n {\sf y}_{\aatchpair }\mathcal{H}^\sigma_n $ for $(\sigma,e)$-admissible $\aatchpair\in \ZZ_{\geq0}^\ell$.  
%
 Our claim is justified as follows: for $e>h$ the  condition that $\aatchpair\in \ZZ_{>0}^\ell$  is $(\sigma,e)$-admissible is equivalent to requiring that $\varnothing$ does not lie on any hyperplane in the alcove geometry (so that our $p$-Kazhdan--Lusztig theory is  ``regular" not ``singular" as required).  
We further remark that the  
importance  of the Serre subquotient  
with regards to calibrated/unitary modules  goes through verbatim to our setting, see  \cite[Introduction]{BNS2}.    
 
\begin{eg}\color{black} 
  Let  $e=3$ and $h=3\in \ZZ$ (and let $\sigma=(0)\in \ZZ$). We have that 
  ${\sf y}_{\underline{h} }=y_{\SSTT^{(3)}}= y_4e(0,1,2,3)$.  
\end{eg}   

\begin{eg}\color{black} 
Continuing with \cref{continuation}, we let  $\sigma=(0,3,8)\in \ZZ^3$ and $e=13$.  
 We have that 
  ${\sf y}_{\underline{h} }= y_4e(0,1,2,3)+y_6e(3,4,5,6,7,8)+ e(8,9,10,11,12)$.  The reader should compare these residue sequences with the   residues appearing in the first row of the tableau in  \cref{continuation}.
\end{eg}

\begin{rmk}
\color{black}  The tableaux $\SSTT^\alpha$ for $0\leq a <\ell$ all have different residue sequences, in particular the corresponding $e_{\SSTT^\alpha}$ are pairwise orthogonal idempotents.
For $h_a <  \sigma_{a+1}-\sigma_{a} $ and  $0\leq a \leq \ell-2$, we have that 
  $y_{\SSTT^\alpha}=e_{\SSTT^\alpha}$.  
  Similarly, for $a=\ell-1$ and  $h_{a} < e+\sigma_0-\sigma_{a-1}-1$, we have that 
  $y_{\SSTT^\alpha}=e_{\SSTT^\alpha}$.  
 If we replace either of the strict inequalities above with an equality, then we obtain     $y_{\SSTT^\alpha}=y_{h_a+1}e_{\SSTT^\alpha}$.  
  Thus the element ${\sf y}_{\aatchpair }$ need not be homogenous, however 
   each element  $y_{\SSTT^\alpha}$ is homogeneous in the grading (of degree 0 or 1).  
We have that the ideal generated by ${\sf y}_{\aatchpair }$  is the same as the ideal generated by the set of  homogeneous elements $\{ y_{\SSTT^\alpha} \mid 0\leq a <\ell\}$ and  therefore the quotient 
 is a graded algebra.
\end{rmk}

\begin{rmk}\label{inwhatfollows}
In \cite[4.1 Lemma]{hm10} it is proven that relation \ref{rel1.12} is equivalent to 
 $e_{\underline{i}}=0$ for any  $\underline{i}\neq \res(\SSTS)$ for some $\SSTS \in \Std(\la)$  with $\la \in \mathscr{P}_{\ell}(n)  $.   
In $\mathcal{H}^\sigma_n/ \mathcal{H}^\sigma_n {\sf y}_{\aatchpair }\mathcal{H}^\sigma_n$ we have that  $e_{\underline{i}}=0$ for any  $\underline{i}\neq \res(\SSTS)$ for some $\SSTS \in \Std(\la)$   with $\la \in \mathscr{P}_{\underline{h}}(n)  $.  For more details, see \cite[Theorem 1.19(a)]{cell4us}.  
\end{rmk}

\subsubsection{The Bott--Samelson truncation} \label{truncation}

 In the previous section, we defined the Bott--Samelson endomorphism algebra and its breadth-enhanced counterpart.  The idempotents in the  former (respectively latter) algebra  were  
   indexed by expressions $\underline{w}$ in the simple reflections (respectively, 
   in the simple reflections {\em and} the empty set).   
   We define 
$$
   {\sf f}_{n,\sigma}^+=\sum_{\begin{subarray}c
   \sts\in \Std_{n,\sigma}^+(\la)
   \\
   \la  \in \mathscr{P}_{\underline{h}}(n)
   \end{subarray}} { e}_{\SSTS}
\qquad  {\sf f}_{n,\sigma}=\sum_{\begin{subarray}c
   \sts\in \Std_{n,\sigma}(\la)
   \\
   \la  \in \mathscr{P}_{\underline{h}}(n )
   \end{subarray}} { e}_{\SSTS}
 $$
and the bulk of this paper will be dedicated to proving that 
$$
  {\sf f}_{n,\sigma}^+(\mathcal{H}^\sigma_n/ \mathcal{H}^\sigma_n {\sf y}_{\aatchpair }\mathcal{H}^\sigma_n )   {\sf f}_{n,\sigma}^+ \quad \text { and }\quad 
    {\sf f}_{n,\sigma} (\mathcal{H}^\sigma_n/ \mathcal{H}^\sigma_n {\sf y}_{\aatchpair }\mathcal{H}^\sigma_n )      {\sf f}_{n,\sigma}
  $$
are isomorphic to the cyclotomic  Bott--Samelson endomorphism algebra and its breadth-enhanced counterpart, respectively.   
 For the most part, we work in the breadth-enhanced Bott--Samelson endomorphism algebra where the isomorphism 
 is more natural (and we then finally truncate at the end of the paper to deduce our Theorem A).  
 
\subsubsection{Concatenation } We now discuss horizontal concatenation of diagrams in (our truncation of) 
 the quiver Hecke algebra.  
 First we let $\boxtimes $ denote the ``naive concatenation" of  KLR diagrams  
  side-by-side as illustrated in \cref{concanetaion2}. 
 Now, given two quiver Hecke diagrams $  \psi ^{\SSTP }_{\SSTQ} $ and $ \psi^ {\SSTP'}_{\SSTQ'} $ we define 
 $$
 \psi ^{\SSTP }_{\SSTQ}
  \otimes 
\psi^ {\SSTP'}_{\SSTQ'}
  =
      e _{\SSTP' \otimes   \SSTQ' }
   \circ \psi^{\SSTP\otimes \SSTQ}_
  {\SSTP'\otimes \SSTQ'}
     \circ      e _{\SSTP' \otimes   \SSTQ' }
   .
 $$ We refer to this as the {\sf contextualised} concatenation of diagrams
 (as the the residue sequences appearing along the bottom 
  of the diagram are not obtained by simple concatenation, but rather from considering the residue sequence of the concatenated path).

 \!\! \begin{figure}[ht!]
$$ \begin{minipage}{12.5cm}
 \end{minipage} $$

 \!\!  \caption{Continuing    \cref{concanetaion}, we depict   $\psi^{\SSTP_\emp}_{\reflectpath }
\boxtimes 
\psi^{\SSTP_\emp}_{\reflectpath }
$ and 
$  
 \psi^{\SSTP_\emp}_{\reflectpath }
\otimes 
\psi^{\SSTP_\emp}_{\reflectpath }
$ respectively.}
 \label{concanetaion2}
 \end{figure}

 \section{Translation and dilation}\label{transdilatesection}
In this     section  we prove some  technical  results  for   KLR elements which will appear repeatedly in what follows.  The reader should feel free to skip this section on first reading. 
We continue with the  notation of \Cref{conventioning}.  

\subsection{The translation principle for paths}
Our first result of this section says that our choice of distinguished path $\SSTP_\w$ in \cref{thepathwewatn} for $\w=\alpha_1 \alpha_2\dots \alpha_p$  was entirely arbitrary (the only thing that matters is that  the path crosses the hyperplanes $\alpha_1,  \alpha_2, \dots \alpha_p$ in sequence).  
 \begin{lem}\label{resconsider}
 Let $\SSTP$ denote any path which terminates at a regular point and let $r\in \ZZ/e\ZZ$.   Then $$  {e_{\SSTP }} 
  \boxtimes e_{r,r}
   =0.  $$  
\end{lem}
\begin{proof}
 The result follows from \cref{inwhatfollows} in light of the proof of \cref{inlightof}.   
   \end{proof}

\!\!
   \begin{figure}[ht!]
   $$
    \begin{minipage}{2.3cm}
\end{minipage}
  $$

\caption{ A series of paths   $\SSTP$, $\SSTQ$, ${\sf R}$, $\SSTS$, $\SSTT$ and $\SSTU$. The paths $\SSTP, \SSTQ, \SSTU$  are   $\al$-crossing paths.      }
\label{pathbendingfork333}
\end{figure}

\!
For $\al \in \Pi$, we say that a path $\SSTP$ of length $n$ is an $\al$-crossing path if $(i)$ 
 there exists  $1<p_1 \leq p_2 <n$ such that 
 $ \SSTP(k)  \in \mathbb{E}(\al )  
 $ 
if and only if   $k\in [p_1,p_2]$ and $(ii)$    
   $\SSTP(k)  \not   \in \mathbb{E}(\bet,se)\not = \mathbb{E}(\al ) $ 
   for any $1\leq k\leq n$.     
We say that  $\SSTP$ is  an $\emptyset$-crossing path  if $\SSTP(k)$ is a regular point for all $1\leq k \leq n$.  
We say  a path is $\al$-bouncing if it is obtained from an  $\al$-crossing path by reflection through the $\al$-hyperplane. 
 
\begin{eg}\label{albounceex}
Let $e=5$, $\ell=1$,   $h=3$, and $\al=\eps_3-\eps_1$. For the paths in \cref{pathbendingfork333}, we have that 
$
\res(\SSTP)= (0,1,4,0,3,4,2,1,0,2)
$, 
$\res(\SSTQ)= (0,1,4,0,3,4,			 2,	1,2, 0)$, $
\res({\sf R})= (0,1,4,0,3,4,			2,	2,	1,0)$, 
$
\res(\SSTS)= (0,1,4,0,3,4,			2,	2,	1,0) $, 
$\res(\SSTT)= (0,1,4,0,3,2,4,			 	2,	1,0)$, 
 and $
\res(\SSTU)= (0,1,4,0,2, 3,4,			 	2,	1,0)$
and we have that   
$$
\res_\SSTP(\SSTP^{-1}(1,\eps_3))=2 \quad 
\res_\SSTP(\SSTP^{-1}(3,\eps_1))=3 \quad 
\res_\SSTP(\SSTP^{-1}(4,\eps_1))=2 \quad 
\res_\SSTP(\SSTP^{-1}(5,\eps_1))=1. 
$$
It is not difficult to see that the elements 
$\psi^\SSTP_\SSTQ$, $\psi^\SSTP_{\sf R}$,  $\psi^\SSTP_\SSTS$, $\psi^\SSTP_\SSTT$, and 
$\psi^\SSTP_\SSTU$    have 0, 1, 2, 3, 3,  crossings of non-zero degree respectively.  We will see that 
$ {\sf e}_\SSTP=\psi^\SSTP_\SSTQ\psi_\SSTP^\SSTQ =\psi^\SSTP_\SSTT \psi_\SSTP^\SSTT = \psi^\SSTP_\SSTU \psi_\SSTP^\SSTU $.  
\end{eg}

\begin{rmk}
Given $\SSTP$ and $\SSTU$ 
 two ($\al$-crossing) paths, we can 
 pass  between them inductively, this lifts  to a factorisation of  $w^\SSTP_\SSTU$ as a product of Coxeter generators.  An example is given by the sequence of paths  $\SSTP$, $\SSTQ$, ${\sf R}$, $\SSTS$, $\SSTT$ and $\SSTU$ in \cref{pathbendingfork333} (for example 
 $w^\SSTS_\SSTT=(6,7)$).  
 The degree of each of these crossings is determined by whether 
 we are stepping onto or off-of a wall.  
For example, the elements $ 
 \psi_{\sf Q}^{\sf R}= e_{\sf R}\psi_{8} e_\SSTQ$, 
  $ 
 \psi_{\sf R}^{\sf S}=e_{\sf S}  \psi_{7} e_{\sf R}$, 
 and 
   $ 
 \psi_{\sf S}^{\sf T}=e_{\sf T}\psi_{6} e_{\sf S}$  have degrees $1, -2, $ and $1$ respectively.   \end{rmk}

 \begin{prop}  
 \label{adjust1}
  Fix   $\al \in \Pi \cup\{\emptyset\} $.  Let 
  $\SSTP, \SSTQ$ be a pair of $\al$-crossing/bouncing paths of length $n$ 
 from $\varnothing \in  A_0$ to $\la \in s_\al   A_0$.  
We have that  \begin{equation}
 \label{trans}
 \psi^{\SSTP}_{\SSTQ}\psi^{\SSTQ}_{\SSTP}
=e_{\SSTP}
\quad \text{ and }\quad 
\psi^{\SSTQ}_{\SSTP}\psi^{\SSTP}_{\SSTQ}
=e_{\SSTQ}.\end{equation}

\end{prop}

\begin{proof} 
The $\al=\emptyset$ case is trivial, and so we set $\al=\eps_i-\eps_{i+1}$.  
We fix   $\SSTP=(\eps_{j_1},\dots,\eps_{j_n})$ and $
\SSTQ=(\eps_{k_1},\dots,\eps_{k_n})$.     Recall that $w^\SSTP_\SSTQ$ is minimal and step-preserving
 and  that the paths   $\SSTP$ and $\SSTQ$ only cross the hyperplane $\al\in\Pi$.  
  This implies, for any pair of strands from $1\leq x < y\leq n$ to $1\leq w^\SSTP_\SSTQ(y) < w^\SSTP_\SSTQ(x) \leq n$ whose crossing has {\em non-zero degree}, that   $ \eps_{j_x}=\eps_{i+1} $ and $\eps_{j_y} = \eps_{i}$
   and   $\SSTP( y)\in   s_\al   A_0$  and $\SSTQ(w^\SSTP_\SSTQ(y))\in    A_0$ 
    (one can swap $\SSTP$ and $\SSTQ$ and hence reorder so that $1\leq y < x \leq n$).      
We let $1\leq y\leq n$ be minimal such that $\SSTP( y)\in   s_\al   A_0$  and we suppose that  $\res_\SSTP(y)= r\in \ZZ/e\ZZ$.  We let  $Y$ denote  this $r$-strand from $y$ to $w^\SSTP_\SSTQ(y)$.

We   recall  our assumption that $\SSTP$ and $\SSTQ$   cross the $\al$-hyperplane precisely once.   This implies that there exists a unique $1\leq z \leq n$ such that 
$ \SSTP^{-1}(z,\eps_{i+1})\in [p_1,p_2]$.   
We have that 
${\sf res}_\SSTP(\SSTP^{-1}(z,\eps_{i+1}))=r+1 $, 
${\sf res}_\SSTP(\SSTP^{-1}(z+1,\eps_{i+1}))=r $, 
and ${\sf res}_\SSTP(\SSTP^{-1}(z+2,\eps_{i+1}))=r-1$.  
  The 
$Y$ strand 
crosses each of the  strands 
connecting the points 
 $\SSTP^{-1}(z ,\eps_{i+1})$, 
$\SSTP^{-1}(z+1,\eps_{i+1})$, 
and $\SSTP^{-1}(z+2,\eps_{i+1})$
 to  the points 
  $\SSTQ^{-1}(z ,\eps_{i+1})$, 
$\SSTQ^{-1}(z+1,\eps_{i+1})$, and $\SSTQ^{-1}(z+2,\eps_{i+1})$
and these are all the crossings  involving the $Y$-strand which are of non-zero degree.  
 We refer to these strands as $Z_{+1}$, $Z_{0}$, $Z_{-1}$.  
 
We are ready to  consider the product $\psi^\SSTP_\SSTQ \psi_\SSTP^\SSTQ$. 
We use case 4 of relation \ref{rel1.10} to resolve the double-crossing  of the $Y$ and $Z_{+1}$ strands, which yields two terms with KLR-dots on these strands.
The term with a KLR-dot on the $Z_{+1}$ strand vanishes after applying case 1 of \ref{rel1.10} to the like-labelled double-crossing $r$-strands $Y$ and $Z_{0}$. 
The remaining term has a KLR-dot on the $Y$ strand. 
We next use \ref{rel1.8} to pull this KLR-dot through one of the like-labelled crossings of $Y$ and $Z_{0}$. 
Again we obtain the difference of two terms, one of which vanishes by applying case 1 of \ref{rel1.10}.
This remaining term has the $r$-strands $Y$ and $Z_{0}$ crossing only once.  
We then pull the $Z_{-1}$-strand through this crossing using 
the second case of relation \ref{rel1.11},  
to obtain another sum of two terms. 
The term with more crossings is zero by \cref{resconsider}, while the remaining term has no non-trivial double-crossings involving the $Y$ strand. 
As the $Y$ strand was chosen to be minimal, we now repeat the above argument with the next such strand; we proceed until all double-crossings of non-zero degree have been undone.     \end{proof}

\begin{rmk}
More generally, given $\SSTP $ and $\SSTQ$ two   $\al$- and $\bet$-crossing/bouncing paths, 
we can apply \cref{adjust1} to any local regions
$\SSTS \otimes \SSTP \otimes \SSTT$ 
and 
$\SSTS \otimes \SSTQ \otimes \SSTT$ 
 of a wider pair of paths.     
The proof again follows simply by applying the same sequence of relations as in the proof of  \cref{adjust1}.  
Indeed,  $\SSTP$ and $\SSTQ$ can be said to be ``translation-equivalent" if  the non-zero double-crossings in $\psi^\SSTP_\SSTQ \psi^\SSTQ_\SSTP$  are precisely those detailed in the proof of \cref{adjust1} (and so are in bijection with the   crossings  of non-zero degree in \cref{transeg}).
\end{rmk}

\begin{eg}\label{transeg}
We now go through the steps of the above proof for the product $\psi^\SSTP_\SSTU
\psi_\SSTP^\SSTU=e_{(0,1,4,0,3,4,2,1,0,2)}$ from Example \ref{albounceex}.  

\vspace{-0.5cm}	\begin{align*}
\scalefont{0.8}
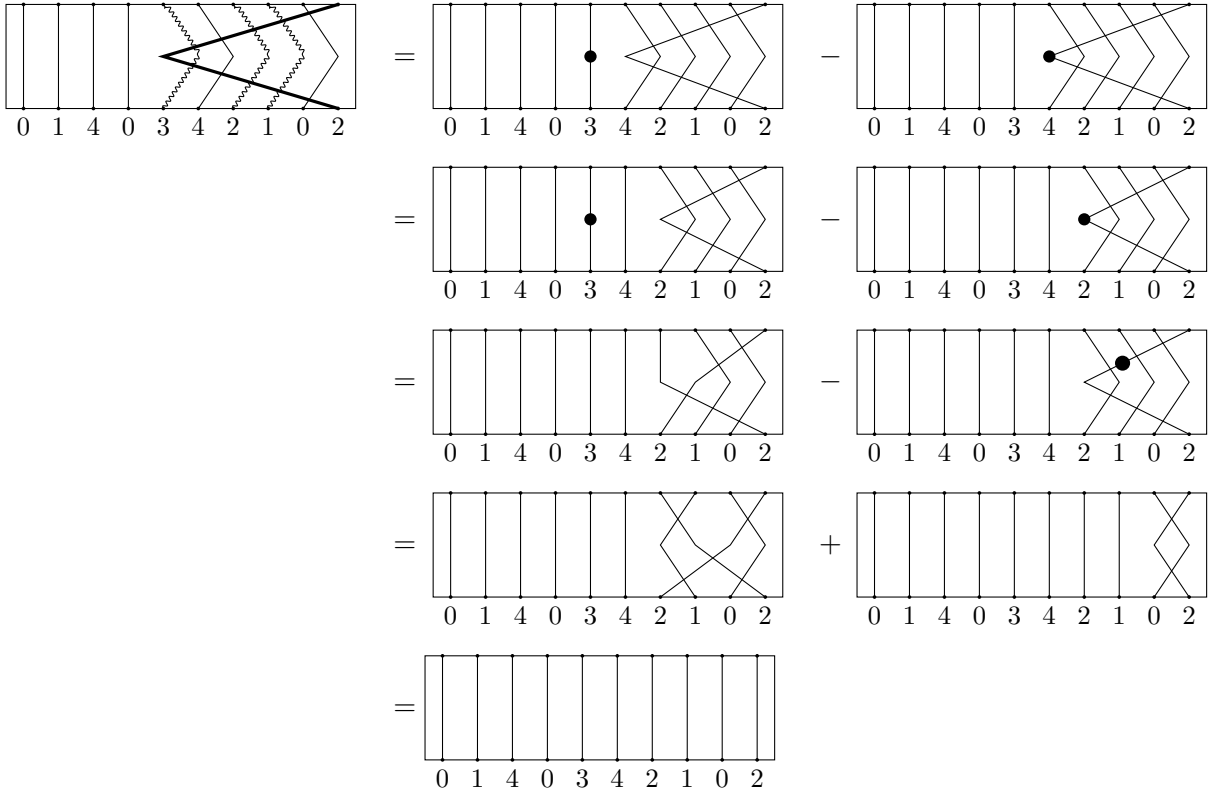
\begin{minipage}{4.3cm} 

    \end{minipage}
\end{align*}
The first and second equalities hold  by case 4 and case 3 of relation \ref{rel1.10}.  The first term in the second line
 and the second term in the third line are both zero by  case 1 of relation \ref{rel1.10}. Thus the third equality  follows by  relation  \ref{rel1.8} and the fourth equality follows from case 1 or relation \ref{rel1.11}.  
    The first term in the fourth line is zero by \cref{resconsider} (the partition $(2^3)$ does not have an addable node of residue 1).    The second term in the fourth line is equal the
 term in the fifth line by case 2 of relation \ref{rel1.10}.  
\end{eg}

 \subsection{Good and bad braids. } 
 Given $w\in \mathfrak{S}_n$, we define a $w$-braid to be any triple 
   $1\leq p<q<r\leq n$ such that  
     $w(p)>w(q)>w(r)$. 
         We recall that an element $w\in \mathfrak{S}_n$ is said to be {\sf  fully-commutative}   if 
      there do   not exist any $w$-braid triples.  
      We define a bad $w$-braid to be a triple 
      $1\leq p<q<r\leq n$ with $i_p = i_r= i_q\pm 1$ 
 such that  
     $w(p)>w(q)>w(r)$.   
     We say that a $w$-braid which is not bad is {\sf good}. 
  We say that $w $ is {\sf residue-commutative} if there do not exist any bad-braid triples.   
 
\begin{lem}\label{ahfdgklhdgajlf}
Suppose   that $w $ is   residue-commutative 
and  let $\w$  be a reduced expression for $w$.  Then $\psi_\w$ is 
independent of the choice of reduced expression and we denote this element simply by $\psi_w$. 
\end{lem}

\begin{proof}
If   $w $ is fully-commutative then 
any two reduced  expressions differ only by the commuting Coxeter relations see \cite[Theorem 2.1]{MR1241505} (in particular, one need not use the braid relation).  Thus the claim follows by the second equality of \ref{rel1.7}.  
 An identical argument shows that  if 
 $w $ is residue-commutative, then 
any two reduced  expressions differ only by the  commuting Coxeter relations   and good braid relations.    The condition for a braid to be good is precisely the commuting case of relation  \ref{rel1.11}.  
  Thus the claim follows by relation \ref{rel1.7} and \ref{rel1.11}.  
 \end{proof}

\subsection{Breadth dilation of permutations }
We will see later on in the paper that the   commutator and hexagonal     generators of \cref{Soergelgens} roughly correspond
to ``dilated" copies of transpositions and braids in the KLR algebra. Similarly, the 
tetrahedron relation roughly corresponds to the equality between two expressions for 
a ``dilated" copy of  $(1,4)(2,3)$.  
In this section, we provide the necessary background results which will allow us to make 
these ideas more precise in \cref{embed,relations}.    
    Given $b>1$, we define the $ b $-dilated transpositions to be the elements 
$$
(i,i+1)_{b} = s_{bi }(s_{bi -1}s_{bi +1})\dots (s_{bi -b+1}s_{bi -b+3}\dots s_{bi -b-3} s_{bi +b-1})	\dots	(s_{bi -1}s_{bi +1})	s_{bi }
$$
for $1\leq i < n$. (The examples in Figure \ref{ghghgh33} should make this definition clear.)   Now,  we note that    $\mathfrak{S}_{n}\cong \langle (i,i+1)_{b} \mid 1\leq i <n\rangle \leq \mathfrak{S}_{bn} .$ 
 We remark that $(i,i+1)_b$ is fully commutative.  
   Given any permutation $w\in \mathfrak{S}_n$ and $\w$ an expression   for $w\in\mathfrak{S}_n$,  we let    $ \w_b$ denote the corresponding expression in the generators  $(i,i+1)_{b} $ 
 of this  $b$-dilated  copy of $ \mathfrak{S}_{n}\leq \mathfrak{S}_{bn} $.    We set 
$ 
 {B}= (  -1 , -2,\dots, -b )^n \in( \ZZ/e\ZZ)^{bn} $ 
and we let $\psi_{\w_b} e_B$ denote the corresponding element in   $ \langle 
e_{B}\psi_{(i,i+1)_{b}}e_{B}\mid 1\leq i <n \rangle
\subseteq \mathcal{H}^\sigma_n$.

  \begin{figure}[ht!]
 $$  
     $$
   
 \!\!\!         \caption{The $2$- $3$- and $4$- dilated  elements  $ e_B \psi_
   {(1,2) _b}e_B$ 
    for $b=2,3,4$.  }
    \label{ghghgh33}
\end{figure}

     \!\!
 \begin{figure}[ht!]
$$   \begin{minipage}{9cm}
 %
 \end{minipage}     
 $$
   
        \!\! 
\caption{The 5-dilated element $e_B\psi_{ (2,3)_5(1,2)_5(2,3)_5}e_B$ for $B=5$.  }  
   \label{nibseg} \end{figure} 
   
   We fix  $\w $ a reduced word for $w\in \mathfrak{S}_n$.  
We say that $D\in \mathcal{H}_{bn}^\sigma$ is a quasi-$b$-dilated expression for
 $\w$ if  for each
$1\leq r <b$,   the 
subexpression consisting solely  of the  $-r$-strands and $-(r+1)$-strands from $D$ forms 
the 2-dilated element $\psi_{\w_2} e_{ (-r,-r-1)^n}$.  
It is easy to see that a 
quasi-$b$-dilated   element  for $\w$   differs from   $\psi_{\w_b}$  simply by undoing some crossings of {\em degree zero}.  
In particular, all quasi-$b$-dilated  expressions for $\w$  (including $\psi_{\w_b}$ itself)  have the same bad braids (in the same order, modulo the commutativity relations).

 \begin{figure}[ht!]
  \begin{minipage}{4cm}
%
%
%
 \begin{tikzpicture}
 [xscale=0.6,yscale=  
  -0.6]
    \draw(2,2.5) rectangle (10,4.5);
        \foreach \i in {2.5,3.5,...,9.5}
  {
   \fill(\i,2.5) circle(1.5pt) coordinate (a\i);
          \fill(\i,4.5)circle(1.5pt)  coordinate (d\i);
    }

\draw(2.5,4.5) --(5.5,2.5);
\draw(3.5,4.5) --(7.5,2.5);
 
\draw(4.5,4.5) --(8.5,2.5);
\draw(5.5,4.5) --(9.5,2.5);
\draw(6.5,4.5) --(2.5,2.5);

\draw(7.5,4.5) --(3.5,2.5);

\draw(8.5,4.5) --(4.5,2.5);
\draw(9.5,4.5) --(6.5,2.5);

  \scalefont{0.8}
\draw (-1+3.5,4.5) node[below] {$-1$};
\draw (-1+4.5,4.5) node[below] {$-2$};
\draw (-1+5.5,4.5) node[below] {$-3$};
\draw (-1+6.5,4.5) node[below] {$-4$};
\draw (-1+7.5,4.5) node[below] {$-1$};
\draw (-1+8.5,4.5) node[below] {$-2$};
\draw (-1+9.5,4.5) node[below] {$-3$};
\draw (-1+10.5,4.5) node[below] {$-4$};
    \end{tikzpicture}    
       \end{minipage}
   \caption{A quasi-$4$-dilated expression for $(1,2)$. 
 This diagram  is obtained from the final diagram  of  \cref{ghghgh33} by undoing a degree zero crossing. }
   \end{figure}

 Finally, we define the nibs of a  permutation $w$ to be 
  the nodes  $1$ and $n$  and $w^{-1}(1)$ and $w^{-1}(n)$ from the top edge and 
  the nodes 
  $1$ and $n$  and $w (1)$ and $w (n)$ from the bottom edge.  
  We define the nib-truncation of $\w$ 
  to be the expression, ${\rm nib}(\w)$, obtained by deleting the 4 pairs of nibs  of $w$ 
   and then deleting the (four) strands connecting these vertices.  
 Similarly, we define  $
{\rm nib}( \psi_{   {\w} }e_{\underline{i}}) =  \psi_{{\rm nib}( {\w})} e_{{\rm nib}(\underline{i})}$ where  the residue sequence ${\rm nib}(\underline{i})\in (\ZZ/e\ZZ)^{bn-4}$ is  inherited  by deleting the $1$st, $n$th, $ w(1)$th and $w(n)$th entries of $\underline{i} \in 
(\ZZ/e\ZZ)^{n}$.  
See \cref{labeler,labeler2} for     examples.  

     \!\!       
  \begin{figure}[ht!]
  $$  \begin{minipage}{9.25cm}
 \end{minipage}     
  $$
   
     \!\! 
\caption{A quasi-$5$-expression element for $ \w=(23)(12)(23)$. Conjugating this diagram by the invertible element  
$(\psi_{10}\psi_{12}\psi_{9}\psi_{11}\psi_{10})e_{
(-1,-2,-3,-4,-5)^3}$   we obtain the  diagram in \cref{nibseg}.  }
\label{labeler}.   \end{figure}

 \begin{figure}[ht!]
 $$  \begin{minipage}{6.5cm}
 %
 \end{minipage}   $$
   \caption{ 
   A diagram obtained by nib-truncation from that in 
   \cref{labeler}.  
    This  diagram is a subdiagram of the hexagonal generator in \cref{braidexample}. 
}
   \label{labeler2}
   \end{figure}

  \!\!\!\!
 \subsection{Freedom of expression} 
We  now prove that the  quasi-dilated elements and their nib-truncations are independent of the choice of reduced expressions.  
For $0\leq q \leq b$, we  define the element $\psi_{[b,q]} 
  $ which breaks the strands into two groups (left and right) according to their residues as follows 
$$
 \psi_{[b,q]}=
\prod_{0\leq p <n} 
 \bigg(
\prod_{1\leq  i \leq q}\psi^{pb+i}_{pq+i}
\bigg) 
\quad \text{where}\quad 
e_B \psi_{[b,q]}\in 
 e_{(-1,\dots,-b)^n}  
\mathcal{H}^\sigma_n
e_{(-1, \dots,-q)^n \boxtimes (-q-1 , \dots, -b)^n} . 
$$ 
   We remark that $\psi_{[b,0]}=\psi_{[b,b]}=1 \in \mathfrak{S}_{bn}$.  
 
\begin{lem}\label{kkkkkkk}
We have that $e_B\psi_{(1,2)_b}\psi_{(1,2)_b}e_B=0$ for $b\geq 1$.  
\end{lem}

\begin{proof}
For $b=1$ the result is immediate by case 1 of relation \ref{rel1.10}.  
Now let $b>1$.  We pull the strand connecting the strand connecting the $1$st top and bottom vertices to the right through the strand connecting the $(b+2)$th top and bottom vertices using case 4 of relation \ref{rel1.10} and hence obtain 
 $$e_B\psi_{[b,b-1]} 
 \Big (\big(
\psi_{(1,2)_{b-1}}
y_{2b-2}
\psi_{(1,2)_{b-1}} \boxtimes \psi_{(1,2)}\psi_{(1,2)}
\big)-\big(
\psi_{(1,2)_{b-1}}
\psi_{(1,2)_{b-1}} \boxtimes \psi_{(1,2)}y_{1} \psi_{(1,2)}
\big)
\Big)
\psi_{[b,b-1]}^\ast e_B$$
and the first (respectively second) terms is zero  by 
 the $(b-1)$th  (respectively  $1$st) inductive step.  \end{proof}

\begin{prop}\label{pfpfpfpfpfpfpfpfpf} 
Let $1\leq b <e$.  The elements $e_{B}\psi_{(i,i+2)_b}e_{B}$
 and 
  ${\rm nib}(e_{B}\psi_{(i,i+2)_b}e_{B})$   are independent of the choice of reduced expression of  $ {(i,i+2)_b}\in \mathfrak{S}_{bn}$.  
\end{prop}
\begin{proof}
For ease of notation we consider the $i=1$ case, the general case is identical up to relabelling of strands. We first consider  $e_{B}\psi_{(1,3)_b}e_{B}$, as the enumeration of strands is easier.   We will refer to two reduced expressions in the KLR algebra as {\em distinct} if they are not trivially equal by the commuting relations (namely,   the latter case of \ref{rel1.7},   case 2 of relation \ref{rel1.10}  and case 3 of relation \ref{rel1.11}). 
There are precisely $b+1$ distinct expressions, $\Omega_q$,  
  of $e_{B}\psi_{(1,3)_b}e_{B}$ as follows    
  \begin{align}\label{sfjkhhlkja}
\Omega_q=
e_B\psi_{[b,q]}\big(
   \psi_{(12)_q} \psi_{(23)_q} \psi_{(12)_q} \boxtimes 
 \psi_{ (23)_{b-q}} \psi_{(12)_{b-q}} \psi_{(23)_{b-q}}
 \big)
 \psi_{[b,q]} ^\ast e_B
\end{align}   for $0\leq q\leq  b$.  See \cref{hh33,4proof} for examples.    We remark that   $
\Omega_0= 
 e_B \psi_{ (23)_{b}} 
  \psi_{ (12)_{b}}  \psi_{ (23)_{b}} e_B$ and $\Omega_b= 
 e_B \psi_{ (12)_{b}}   \psi_{ (23)_{b}}     \psi_{ (12)_{b}}
 e_B$.  We will show that $\Omega_q=\Omega_{q+1}$ for $1\leq q<b$ and hence deduce the result.

\begin{figure}[ht!]
$$   \begin{minipage}{4cm}
 \end{minipage}     $$
   
   \!\! 
\caption{The 3 distinct expressions, $\Omega_0$, $\Omega_1$, and 
$\Omega_2$    for $\psi_{(1,3)_2}$. 
The  $b+1$ distinct  expressions for $\psi_{(1,3)_b}$  are determined by where the central ``fat strand" is broken into ``left" and ``right" parts.  }
\label{hh33}
\end{figure}

\!\!\!\!
\begin{figure}[ht!]
$$    \begin{minipage}{4cm}
 %
 \end{minipage}  $$

\!\! \caption{The 4 equivalent  expressions for  $\Omega_1$ of \cref{hh33}. 
These differ only by applications of case 3 of  relation \ref{rel1.11} (and so the bad braids are all the same).}
\label{4proof}
\end{figure}

\noindent{\bf Step 1.} If $q=0$ proceed to Step 2, otherwise we pull the $(-q)$-strand connecting the $(b+q)$th northern and southern nodes of $\Omega_q$ to the right. We first
use relation \ref{rel1.11} to pull  $(-q)$-strand through 
the crossing of $(1-q)$-strands connecting the 
the $(q-1)$th and $(2b+q-1)$th top and bottom vertices.   We obtain two terms: the first is equal to 
\begin{equation}\label{refff}
e_B\psi_{[b,q]}\big(
 \psi_{[q,q-1]}\big(
 \psi_{ (12)_{q-1}} \psi_{(23)_{q-1}} \psi_{(12)_{q-1}} \boxtimes  \psi_{(12)(23)(12)}\big)
  \psi_{[q,q-1]}^\ast 
   \boxtimes  
 \psi_{ (23)_{b-q}} \psi_{(12)_{b-q}} \psi_{(23)_{b-q}}
 \big)
 \psi_{[b,q]} ^\ast
 e_B
\end{equation} and an error term of strictly smaller length (in which we undo the crossing pair of $(1-q)$-strands). 
If $q=1$, the error term contains a double-crossing of $(r-q)$-strands and so is zero by case 1 of relation \ref{rel1.10}.
If $q>1$,  then we apply   relation \ref{rel1.11} to the error term to obtain two distinct terms; one of which is zero by \cref{kkkkkkk} and the other is zero by case 2 or relation \ref{rel1.10} and the commutativity  relations.

     \medskip
\noindent{\bf Step 2.}  The  output    from Step 1 has a subexpression $  \psi_{(12)(23)(12)}$ which we rewrite as $ \psi_{  (23)(12)(23)}$ using case 3 of relation \ref{rel1.11} (as the three strands are all of the same residue, $-q\in \ZZ/e\ZZ$).  
We also have that    $ \psi_{[b,q]}\psi_{[q,q-1]}=
 \psi_{[b,q-1]}(
1 _{\mathcal{H}^\sigma_{3b-3}}\boxtimes \psi_{[b-q+1,1]})$.
Thus \ref{refff} is equal to  $$
\psi_{[b,q]}\big(
 \psi_{ (12)_{q-1}} \psi_{(23)_{q-1}} \psi_{(12)_{q-1}}   
    \boxtimes  \psi_{[b-q+1,1]} \big( \psi_{(23)(12)(23)}\boxtimes 
 \psi_{ (23)_{b-q}} \psi_{(12)_{b-q}} \psi_{(23)_{b-q}}
 \big)\psi_{[b-q+1,1]}^\ast   \big)
 \psi_{[b,q]} ^\ast
$$
   Now,  by the mirror argument to that used in Step 1, we have that this   equals $$\psi_{[b,{q-1}]}\big(
   \psi_{(12)_{q-1}} \psi_{(23)_{q-1}} \psi_{(12)_{q-1}} \boxtimes 
 \psi_{ (23)_{b-{q+1}}} \psi_{(12)_{b-{q+1}}} \psi_{(23)_{b-{q+1}}}
 \big)
 \psi_{[b,{q-1}]} ^\ast
$$
as required.    
The argument for    ${\rm nib}(e_{B}\psi_{(1,3)_b}e_{B})$    is identical (up to relabelling of strands) except that the $q=0$ and $q=b$ cases  do not appear.  \end{proof}

\!\!

\begin{cor}\label{nibs}
Let  $\x$ be any expression in the Coxeter generators  of $\mathfrak{S}_n$.  
Any  quasi-$b$-dilated expression  of  $\x$  is  independent of the choice of expression $\x$.  
   Similarly, the nib truncations of these elements are independent of the choice of expression $\x$.
\end{cor}
\begin{proof}
By \cref{ahfdgklhdgajlf} it is enough to consider the bad braids in $\psi_\x$.  
If  $\x=\w_b$ for some $w\in \mathfrak{S}_n$, 
then we can resolve each bad braid in $\psi_\x$ and ${\rm nib}(\psi_\x)$ using \cref{pfpfpfpfpfpfpfpfpf}.  
Now, if  $\psi_\x$ is quasi-$b$-dilated then 
 $\psi_\x$ and ${\rm nib}(\psi_\x)$ are  
  obtained from    $\psi_{\w_b}$ and ${\rm nib}(\psi_{\w_b})$
   by undoing some degree zero crossings (thus introducing no new bad braids) and the result follows.  
\end{proof}

 \section{Recasting the diagrammatic Bott--Samelson  generators\\ in the quiver Hecke  algebra } \label{embed} 
 
  We continue with the  notation of \Cref{conventioning}.    The  elements of  the (breadth-enhanced)  diagrammatic Bott--Samelson endomorphism algebras can be thought of as morphisms
  relating pairs of expressions from $\Shl$.    
We have also seen that one can  think of an  element  of  the quiver Hecke algebra as a morphism between pairs of paths in the alcove geometries of \cref{sec3}.
   This will allow us, through the relationship between paths and their colourings described in \cref{sec3},  to define the isomorphism behind Theorem A.  
 In what follows we will define generators 
$$
{\sf adj}_{\al \emp}^{  \emp \al}   \quad 
{\sf spot}_\al^\emp   \quad 
{\sf fork}_{\al \al}^{\emp \al} \quad 
{\sf com}^{\bet\gam}_{\gam\bet}\quad 
{\sf hex} _{\al\bet \al}^{\bet\al\bet}
$$
for $\al,\bet,\gam\in \Pi$ and their duals.  
The hyperplane labelled by 
$\al$ (respectively $\bet$)  is a wall of the dominant chamber if and only if 
 $\SSTP_\al$ (respectively  $\SSTP_\bet$) leaves the dominant chamber.  
By the cyclotomic KLR relation, one  of the above generators is zero if (and only if) 
 one of its  indexing  roots labels a path which leaves the dominant chamber.    
However, one should think of these as generators in the sense of a  
right tensor quotient of a  monoidal category.   
In other words, we still require every  generator for every simple root (even if they are zero) as the {\em left concatenates} of these generators will not be zero, in general.

In order to construct our isomorphism, we must first ``sign-twist" the elements of the KLR algebra.  This twist counts the number of degree $-2$ crossings (heuristically, these are the  crossings which ``intersect an alcove wall").  
For $\w$ an arbitrary reduced expression, we set 
$$
 \Upsilon_\w  	=(-1)^{\sharp\{1\leq p<q\leq n\mid w(p)>w(q), i_p=i_q\}}e_{\underline{i}}   \psi_\w  e_{w({\underline{i}})}   .
$$
While KLR diagrams are usually only defined up to a choice of expression, we   emphasise that each of the generators we define is independent of this choice.  
  Thus there is no ambiguity 
  in defining the  elements $\Upsilon^\SSTP_\SSTQ$ for $w^\SSTP_\SSTQ$   without reference to the underlying expression.  
 In other words: these generators are {\em canonical elements} of $\mathcal{H}_n^\sigma$.  
 Examples of concrete choices of  expression  can be found  in \cite[Section 2.3]{cell4us}.   
 In various proofs it will be convenient to denote by $\SSTT$ and $\southT$ the top and bottom paths of certain diagrams (which we define case-by-case).

\subsection{Idempotents in diagrammatic algebras}
We consider an element  of the quiver Hecke or {diagrammatic Bott--Samelson} endomorphism algebra  to   be a morphism between paths, \color{black}  lifting the ideas of \cref{perms as morphisms}.  \color{black}
The easiest elements to construct are the idempotents 
corresponding to the trivial morphism from a  path to itself.      
 Given $\al$ a simple reflection, 
 we have an 
 associated  path 
     $\SSTP_\al $, a trivial bijection $w^{\SSTP_\al} _{\SSTP_\al} =1\in \mathfrak{S}_{b_\al}$,  and  an idempotent element of the quiver Hecke algebra
 $$ 
e_{\SSTP_\al }   := e_{\res(\SSTP_\al)}
\in \mathcal{H}_{b_\al}^\sigma  
 $$
\color{black} where we reemphasise  that $ e_{\res(\SSTP_\al)}=e_{\res(\SSTP_\al^\flat)}$ (see \cref{idemp-remark-for-the-refsssss}). 
 \color{black}
 Given $\al$ a simple reflection,  we also have  a  Soergel diagram 
${\sf1}_{{\al  }}$ given by a single vertical strand coloured  by  $ \al $.  
We define 
 \begin{equation}\label{idemp-iso}
\Psi (
 {\sf 1}_{ \al  }
)=
 e_{\SSTP_\al }  . \end{equation} 
More generally, given any $\w=s_{ \alpha^{(1)}}s_{ \alpha^{(2)}}\dots  s_{ \alpha^{(k)} }$ any expression of breadth $b(\w)=n$, we have an 
 associated  path 
     $\SSTP_\w $, and  an element of the quiver Hecke algebra
 $$ 
e_{\SSTP_\w }   := e_{\res(\SSTP_\w)}=
e_{\SSTP_{\alpha^{(1)}}}
\otimes
e_{\SSTP_{\alpha^{(2)}}}	
\otimes \dots\otimes e_{\SSTP_{\alpha^{(k)}}}	\in \mathcal{H}_{n{\aatch} }^\sigma 
 $$ and a 
  $(\w,\w)$-Soergel diagram 
$${\sf 1}_{\underline{w}  }= 
{\sf 1}_{\alpha^{({\sf 1})}}
\otimes
{\sf 1}_{\alpha^{(2)}}	
\otimes \dots\otimes {\sf 1}_{\alpha^{(k)}}$$ 
given by $k$ vertical strands,    coloured with   $ \alpha^{(1)}, \alpha^{(2)},\mydots,  \alpha^{(k)} $ from left to right.  
We define 
\begin{equation}\label{idemp-iso2}
\Psi (
 {\sf1}_{ \underline{w} }
)=
 e_{\SSTP_\w }  . \end{equation} 
  
 \begin{eg}\label{idemp-iso3}
Continuing with \cref{diag2,diag22point4,exampleforanton}, we    
let   $$\w= 
s_\emptyset s_\emptyset s_{\color{magenta}\varepsilon_3-\varepsilon_1}
 s_{\color{darkgreen}\varepsilon_2-\varepsilon_3}
 s_{\color{cyan}\varepsilon_1-\varepsilon_2}
 s_{\color{magenta}\varepsilon_3-\varepsilon_1}
 s_{\color{darkgreen}\varepsilon_2-\varepsilon_3}
 s_{\color{cyan}\varepsilon_1-\varepsilon_2}
$$  $$ \w'= 
  s_{\color{magenta}\varepsilon_3-\varepsilon_1}s_\emptyset
 s_{\color{darkgreen}\varepsilon_2-\varepsilon_3}
 s_{\color{cyan}\varepsilon_1-\varepsilon_2} s_\emptyset
 s_{\color{magenta}\varepsilon_3-\varepsilon_1}
 s_{\color{darkgreen}\varepsilon_2-\varepsilon_3}
 s_{\color{cyan}\varepsilon_1-\varepsilon_2}
 .$$   
Recall these path came from ``inserting determinants" into the path in \cref{diag2}.  
We have that 
\begin{align*} \Psi ({\sf1}_{\w})=
e_{
0,1,2,4,0,1, 
\color{magenta} 
3,4,
2,3,
1,2, 
0,
4,
3,
\color{darkgreen}	
0,
2,
1,
\color{cyan}
0,
  1,
4,
\color{magenta}
3,  	4,
2,     3,
1,	2,
0,
4,
3,
\color{darkgreen}
0, 
   2,		
1,\color{cyan}
 0,    1,
 4
 \color{black}}
\\
\Psi ({\sf1}_{\w'})=
e_{
\color{magenta} 
0,1,
4,0,
3,4, 
2,
1,
0,
\color{black} 
2,3,4,
\color{darkgreen}	
4,
1,
0,
\color{cyan}
4,
  0,
3,
\color{black} 
2,3,4,
\color{magenta}
3,  	4,
2,     3,
1,	2,
0,
4,
3,
 \color{darkgreen}
0, 
   2,		
1,\color{cyan}
 0,    1,
 4
 \color{black}}.
\end{align*} 

 \end{eg}

  \begin{rmk}\label{likeidemp}
\color{black} For two paths $\SSTS$ and $\SSTT$, we have that  $\SSTS \sim\SSTT$ if and only if  $\res(\SSTS)=\res(\SSTT)$. 
   Therefore if $\SSTS \sim\SSTT$ then 
   $e_\SSTT=e_{\SSTS}e_\SSTT=e_\SSTS$.  
    \color{black} 
 In particular $e_{\SSTP_\al}=e_{\SSTP_\al}e_{\SSTP_\al^\flat}=e_{\SSTP_\al^\flat}$.

 \end{rmk}

 \begin{rmk}
 \color{black}
We have defined two distinct paths $\SSTP_\al$ 
 and $\SSTP_\al^\flat$ which label the same idempotent, thus  $e_{\SSTP_\emp}\mathcal{H}^\sigma_{b_\al }e_{\SSTP_\al}
 =
 e_{\SSTP_\emp}\mathcal{H}^\sigma_{b_\al }e_{\SSTP_\al^\flat}$.
This apparent redundancy is required  
 because we cannot directly compare $\SSTP_\emp$ and $\SSTP_\al$ as they {\em do not have the same shape} --- however, we can compare 
 $\SSTP_\emp$ and $\SSTP_\al^\flat$ as they {\em do have the same shape}.  Thus $\SSTP_\al^\flat$ is required in order to define the spot-morphism.%
 For the remainder of this section, we will restrict our attention to a subset of morphisms between paths of the same shape which form a set of  monoidal generators of our truncated KLR algebra. 
  \end{rmk}

\renewcommand{\SSTS}{{\sf S}}

\subsection{Local adjustments and isotopy}  We will refer to the passage between alcove paths which differ only 
     by occurrences of $s_\emptyset=1$  (and their associated idempotents) as ``adjustment". 
     We wish to understand the  morphism relating the paths 
     $ 
     \SSTP_{\al}\otimes \SSTP_{\emptyset} 
     $ and $\SSTP_{\emptyset}\otimes \SSTP_{\al}$.
\begin{prop}\label{easy}
The element $\psi^{   \SSTP_{\al\emptyset} }_{\SSTP_{\emptyset\al}} $ is independent of the choice of reduced expression.  
\end{prop}
\begin{proof}
There are precisely three crossings in   
 $ 
  \psi^{   \SSTP_{\al\emptyset} }_{\SSTP_{\emptyset\al}}  
 $ 
 of non-zero degree.   Namely, the $r$-strand (for some $r\in \ZZ/e\ZZ$) connecting the
  $\SSTP_{ \emptyset\al}^{-1}(1,\eps_i)$th top vertex to the  
    $\SSTP_{\al\emptyset}^{-1}(1,\eps_i)$th  bottom vertex 
crosses each of the strands connecting 
  $\SSTP_{ \emptyset\al}^{-1}(q,\eps_{i+1})$th top vertices to the  
    $\SSTP_{\al\emptyset}^{-1}(q,\eps_{i+1})$th  bottom vertices 
for $q=b_\al-1,b_\al,b_\al+1$ (of residues $r+1$, $r$, and $r-1$ respectively) precisely once with degrees $+1$, $-2$, and $+1$ respectively.  
Thus  $ 
  \psi^{   \SSTP_{\al\emptyset} }_{\SSTP_{\emptyset\al}} 
 $ is residue-commutative and the result follows from \cref{ahfdgklhdgajlf}.   
\end{proof}
Thus we are free to define the {\sf KLR}-adjustment to be 
 $$
 {\sf adj}^{\al\emptyset}_{{\emptyset\al}}
 :=
\Upsilon^{   \SSTP_{\al\emptyset} }_{\SSTP_{\emptyset\al}}
 $$
 which is independent 
 of the choice of reduced expression of the permutation.

\!\! 
\begin{figure}[ht!] 
  $$  
 $$  
    
    \!\!\!
    \caption{  We let $h=1$, $\ell=6$, $e=12$, $\sigma=(0,2,4,6,8,10)$ and  $\al=\varepsilon_3-\varepsilon_4$.   
    The adjustment term $  {\sf adj}_{{\al \emptyset}}^{{ \emptyset \al}}
  $ is illustrated.  
    The steps of the path $\SSTP_{\al}$ and $\SSTP_\emptyset$ are coloured pink and black respectively  within   both  $\SSTP_{\al\emptyset}$ (along the top of the diagram) and 
    $\SSTP_{ \emptyset\al}$  (along the bottom of the diagram).   }
\label{adjustmntexample}    \end{figure}

\begin{prop}
We have that 
$$ {\sf adj}^{\emptyset \al}_{\al \emptyset} \circ 
  e_{\SSTP_{\al \emptyset}}
 \circ   {\sf adj}_{\emptyset \al}^{\al \emptyset}   
 = 
   e_{\SSTP_{ \emptyset \al}}
   \qquad\text{and}\qquad
   {\sf adj}_{\emptyset \al}^{\al \emptyset} \circ 
  e_{\SSTP_{ \emptyset\al }}
 \circ   {\sf adj}^{\emptyset \al}_{\al \emptyset}   
 = 
   e_{\SSTP_{   \al\emptyset}}$$
   and so adjustment is an invertible process.
\end{prop}
  \begin{proof}
The paths $\SSTP_{\al\emptyset}$ and $\SSTP_{\emptyset\al}$ satisfy the conditions of \cref{adjust1} and so the result follows.  
  \end{proof}

Finally, we remark that 
the above adjustment can be generalised from the $b_\emptyset=1$ case to the $b_\empg\geq 1$ case as follows.  
 For $\w=s_{\al} s_\empg$ with $\al,\gam\in \Pi $ two (equal, adjacent, or non-adjacent) simple roots, we set 
 $$
 {\sf A}^{\empg\al}_{\al\empg}(q)= \SSTP_{q\emptyset} \otimes 
 \SSTP_{\al}\otimes \SSTP_{(b_\gam-q)\emptyset}
 $$
 for $0\leq q \leq b_\gam$ and we set 
$$
    {\sf adj}_{\al \empg}^{\empg \al}(q )  =
 e_{{\sf A}^{\empg\al}_{\al\empg}(q+1)}
 \left(e_{\SSTP_{q\emptyset}} \otimes  {\sf adj}^{ \emptyset \al }
   _{ \al \emptyset }
 \otimes    e_{\SSTP_{(b_\gam-q-1 )\emptyset}} \right) e_{{\sf A}^{\empg\al}_{\al\empg}(q)}$$
 and we define $$    {\sf adj}^ { \empg \al}_{\al \empg} =   {\sf adj}_{\al \empg}^ { \empg \al}(b_\gam-1 ) \mydots   {\sf adj} _{\al \empg}^ { \empg \al}(1)
    {\sf adj} _{\al \empg}^ { \empg \al}(0).
$$

   \subsection{The KLR-spot diagram}\label{KLRSPOTSECTION}
 We now define the  spot path morphism. 
 Recall that 
\begin{align*}
\SSTP_{\emp}  = 
(\varepsilon_{1 }, \mydots, \varepsilon_{i-1 }
 ,  {\varepsilon_{ i  }},   \varepsilon_{ i+1  },\mydots, \varepsilon_{{\aatch} })^\exx 
\qquad 
\SSTP_{\al }^\flat 
= 
(\varepsilon_{1 }, \mydots, \varepsilon_{i-1 }
 , \widehat{\varepsilon_{ i  }},   \varepsilon_{ i+1  },\mydots, \varepsilon_{{\aatch} })^\exx 
 \boxtimes (\eps_i)^\exx
 \end{align*}
 \color{black} are both paths of the same shape. 
 \color{black}
The permutation $ w^ {\SSTP_\emp}_{\SSTP_\al^\flat}$ is fully-commutative and so we are free to define the {\sf KLR-spot}  to be the elements  $$
{\sf spot}^\emp_\al
 :=
\Upsilon  ^ {\SSTP_\emp}_{\SSTP_\al^\flat }
\qquad
{\sf spot}_\emp^\al
 :=
\Upsilon _ {\SSTP_\emp}^{\SSTP_\al^\flat }
 $$
which  are both independent 
 of   choice of reduced expressions
 \color{black}
  and both belong to
 $e_{\SSTP_\al}\mathcal{H}^\sigma_{b_\al}
 e_{\SSTP_\al}=e_{\SSTP_\al^\flat}\mathcal{H}^\sigma_{b_\al}
 e_{\SSTP_\al^\flat}$.  
 \color{black}

 We wish to inductively pass between the paths 
   $\SSTP_{\al }^\flat $ and $\SSTP_{\emp }$ 
{\color{black}   by means of a visual  timeline (pictured in \cref{pathbendingdot}).   
   This allows us to  factorise} the   KLR-spots    and to simplify our proofs later on.   
   To this end we define 
$$
\Spotq
 = 
 \SSTP_{ q\emptyset }
   \boxtimes
 \REMOVETHESE i {\exx -q}
   \boxtimes  
  \ADDTHIS i  {\exx -q } 
=
  (\varepsilon_1,\varepsilon_2,\mydots, \varepsilon_{\aatch})^{ q}
   \boxtimes
(\varepsilon_{1 }, \mydots, \varepsilon_{i-1 }
 , \widehat{\varepsilon_{ i  }},   \varepsilon_{ i+1  },\mydots, \varepsilon_{{\aatch} })^{\exx -q}
   \boxtimes (\eps_i)^{\exx -q }
 $$
 for $0\leq q \leq b_\al$  and we notice that $\Spotzero=  \SSTP_{{\al}} ^\flat  $ and $ \Spotb=   \SSTP_{{\emp}} $.
  We define 
$ {\sf spot}_{ {{\al}} } ^\emp (q)$ to be the element 
 $ 
  {\sf spot}_{ {{\al}} } ^\emp (q)= \psi ^{\Spotqplus	}_{\Spotq}      $
for  $0\leq q <\exx$   and we  factorise $ {\sf spot}_{\al}^{\emp} $ as follows 
 $$
 {\sf spot}_{\al}^{\emp} :=
 e_{  \SSTP_{\emp }  }	  \circ 
{\sf spot}_{ {{\al}} } ^\emp (\exx-1) \circ 
 \cdots 
 \circ {\sf spot}_{ {{\al}} } ^\emp (1) \circ  {\sf spot}_{ {{\al}} } ^\emp (0)    \circ  e_{\SSTP _{{\al}}^\flat  }.
$$

 \!\!\!
 \begin{figure}[ht!]
$$
  \begin{minipage}{2.6cm}
\end{minipage} 
$$
\caption{
An example  timeline for the KLR spot.  
Fix $\ell=1$ and $h=3$ and $e=5$ and $\al=\varepsilon_3-\varepsilon_1$
 (so that $\exx=3$).  
From left to right we picture  
$ \Spottwo={\sf S}_{{\al}}(3)=\SSTP_{\emp}$, 
$ \Spotone$, 
  $ \Spotzero=\SSTP_{\al }^\flat $.  
We do not picture the $k= 2,1,0$ copies of the path $(+\varepsilon_1,
+\varepsilon_2,+\varepsilon_3)$ at the start of each path, for ease of readability.   
  }
\label{pathbendingdot}
\end{figure}

%

\begin{figure}[ht!]
\begin{equation*}
{\sf spot}^{{\emp}}_{{\al}}\; =\; \begin{minipage}{9cm}
\end{minipage} \end{equation*}
\caption{The element ${\sf spot}_\al^\emp$ of \cref{reasonsbecomeclear}.  We have added the step labels on top and bottom so that one can appreciate that this element is a morphism between paths.   
 However, we remark that while a necessary condition for a product of two KLR diagrams to be non-zero is that their residue sequences must coincide, the same  is not true 
for their step labels (see \cref{likeidemp}).}
 \end{figure}

 \begin{eg}\label{reasonsbecomeclear}
Let $h=3$ and $\ell=1$ and $e=5$ and  $\al=\varepsilon_3-\varepsilon_1$.
We have that $\exx=3$ and 
\begin{align*}
\SSTP_{\al }^\flat   = \Spotzero
&=( \varepsilon_1, \varepsilon_2, 
 \varepsilon_1, \varepsilon_2, 
 \varepsilon_1, \varepsilon_2, 
  \varepsilon_3, \varepsilon_3, \varepsilon_3)
\\
 \Spotone
&=( \varepsilon_1, \varepsilon_2,   \varepsilon_3)   \boxtimes  (
 \varepsilon_1, \varepsilon_2, 
 \varepsilon_1, \varepsilon_2, 
 \varepsilon_3, \varepsilon_3)
\\
 \Spottwo
&= 
( \varepsilon_1, \varepsilon_2,   \varepsilon_3)   \boxtimes 
( \varepsilon_1, \varepsilon_2,   \varepsilon_3)   \boxtimes 
( \varepsilon_1, \varepsilon_2, 
  \varepsilon_3)
 \\
 \SSTP_{\emp}= \Spotthree
&= 
( \varepsilon_1, \varepsilon_2,   \varepsilon_3)   \boxtimes 
( \varepsilon_1, \varepsilon_2,   \varepsilon_3)   \boxtimes 
( \varepsilon_1, \varepsilon_2, 
  \varepsilon_3)
 \end{align*}
which are depicted in \cref{pathbendingdot}.    Of course, $ \Spotthree= \Spottwo $ in this case, but this is only because $\al$ is the affine root $\varepsilon_3-\varepsilon_1$ with $3={\aatch}$.
 
 \end{eg}

 \begin{rmk}
We have that 
 $w^{{\sf S}_{{\al},{q+1}}}_{{\sf S}_{{\al,q}}  }   =  w 
 ^{ q{\aatch}   + \isit }_{   \exx {\aatch}-\exx + q+1	 }      $ 
   for $0\leq q < \exx$,     
    where the sub and superscripts correspond to 
    $$ \Spotq^{-1} (q+1,\eps_i)
    =q{\aatch} +i
    \qquad 
    \Spotqplus ^{-1} (q+1,\eps_{i})=\exx {\aatch}-\exx + q+1
    $$
    and so one can think of the spot morphism as successively   removing each  $+\eps_i$ step from the latter path and adding it to the former.  
 \end{rmk}

 \begin{rmk}The element 
$ e_{ \Spotqplus} {\sf spot}_{ {{\al}} } ^\emp (q)  e_{ \Spotq}  $
 is of degree 1 for $q=0$ and degree 0 for $0<q <b_\al$.  
 The terms with $0<q <b_\al$ are invertible by \cref{adjust1}.   
 Thus one can think of the $q=0$ term as the real substance of ${\sf spot}_{\al}^{\emp} $.  
 One should intuitively think of this degree contribution as coming from the fact that the path ${\SSTS_{0,\al}}$ steps onto and off of a hyperplane but  
${\SSTS_{1,\al}}$  
does not   touch the hyperplane at any point.  
The diagram $ {\sf spot}_{{{\al}} } (0)$ has a   crossing involving the strand  from the
 $  \Spotzero^{-1}(1,\eps_{i})$th node on the   \south       edge   to the 
$   \Spotone^{-1}(1,\eps_i) $th node on the   \north    edge   and the strand from 
the $ 
 	 \Spotone	^{-1}(b_\al,\eps_{i+1})$th node on the   \south       edge   to the
 $  \Spotzero^{-1}(b_\al,\eps_{i+1})$th node  on the   \north    edge.  
  See   \cref{pathbendingdot} for a visualisation. 
\end{rmk}

  \subsection{The KLR-fork   diagram}
We wish to understand the  morphism from 
      $\SSTP_{\emp}\otimes \SSTP_{\al}$ to  $\SSTP_{\al }\otimes\reflectpath$ \color{black}(which are both paths of the same shape, so this makes sense).  
      \color{black}
 The permutation   $w_{ \SSTP_{\al }\otimes\reflectpath }^{\SSTP_{{\emp}}\otimes \SSTP_{{\al}}}$ is {\em not} fully commutative and so we must do a little work prior to our definition.   
 
 \begin{prop}\label{forkkkk}
 The elements $\psi_{ \SSTP_{\al }\otimes\reflectpath }^{\SSTP_{{\emp}}\otimes \SSTP_{{\al}}}$ 
 and 
 $\psi_{\reflectpath \otimes  \SSTP_{\al } }
 ^{    \SSTP_{{\al}}	\otimes \SSTP_{{\emp}}	}$ 
are  independent of the  reduced expressions.   \end{prop}
 \begin{proof}
 We focus on the former case, as the latter is similar.  
The element $w_{ \SSTP_{\al }\otimes\reflectpath }^{\SSTP_{{\emp}}\otimes \SSTP_{{\al}}}$ contains  precisely $b_\al$   crossings of strands with the same residue label:  
 Namely for each $1\leq q \leq b_\al$ the strand connecting the     \north    and   \south       vertices labelled by the integers 
  $$
   \SSTP_{\emp\al}^{-1}(q,\eps_{i})= q{\aatch}+i\qquad 
 (\SSTP_\al \otimes  \reflectpath  )^{-1}(q,\eps_{i})= \exx {\aatch}  +
 (q-1)({\aatch} -1)+\isitone   
  $$
 crosses the strand 
  connecting the     \north    and   \south       vertices labelled by the integers 
   $$
   \SSTP_{\emp\al}^{-1}(b_\al+q,
   \eps_{i+1})
   = 	 	\exx {\aatch}  +(q-1)({\aatch} -1)+ \al(i+1)	\qquad 
 (\SSTP_\al \otimes  \reflectpath  )^{-1}(b_\al+q,\eps_{i+1})= 
 b_\al {\aatch} - b_\al + q .
  $$
The $q$th of these like-labelled crossings forms a braid 
  with a third strand if and only if this third strand connects a      \north    and   \south       node  labelled by the integers 
 $$
     \SSTP_{\emp\al}^{-1}(b_\al+p,
   \eps_{j})
   = 	\exx {\aatch}  +(p-1)({\aatch} -1)+ \al(j) 	\quad 
 (\SSTP_\al \otimes \reflectpath  )^{-1}(b_\al+p,\eps_{j})=
  \exx {\aatch}  +(p-1)({\aatch} -1)+ \al(j) 
  $$
 for $ \al( j)  \neq \al( i+1)$ and  $1\leq p < q$
 or $p=q$ and $\al(j)  <\al(i+1)$.  None of the resulting braids is bad; thus $\psi_{ \SSTP_{\al }\otimes\reflectpath }^{\SSTP_{{\emp}}\otimes \SSTP_{{\al}}}$ is residue-commutative and the result follows.  
\end{proof}

\begin{figure}[ht!]
$$
  \begin{minipage}{2.6cm}
\end{minipage}
$$
\caption{An example of a timeline for the KLR fork.   Fix $\ell=1$ and $h=3$ and $e=5$ and $\al=\varepsilon_3-\varepsilon_1$
 (so that $\exx=3$).  
From left to right we picture  the paths $\Forkzero=
\SSTP_{\al} \otimes    \SSTP_{\al}^\flat $, 
$\Forkone$, 
$\Forktwo$, 
$\Forkthree=\SSTP_{\emp \al}$.      
Notice that we do not picture the $q=0,1,2,3$ copies of the path $(+\varepsilon_1,
+\varepsilon_2,+\varepsilon_3)$ at the start of each path, for ease of readability.  
}
\label{pathbendingfork}
\end{figure}

Thus we are free to define the {\sf KLR-forks} to be the elements  
 $${\sf fork}_{\al\al}^{\emp\al}:=		
 \Upsilon _{ \SSTP_{\al }\otimes\reflectpath }^{\SSTP_{{\emp}}\otimes \SSTP_{{\al}}}
 \qquad
 {\sf fork}_{\al\al}^{\al\emp}:=	
 \Upsilon  _ {\reflectpath \otimes  \SSTP_{\al } }
 ^{    \SSTP_{{\al}}	\otimes \SSTP_{{\emp}}	}$$ 
 which are independent of the choice of reduced expressions. 
 {
 \color{black} 
We reemphasise   that 
$\res(\SSTP_\al)=\res(\SSTP_\al^\flat)$, thus former  element  belongs to 
 $({e_{\SSTP_{\emp} }\color{magenta} \otimes  
e_{\SSTP_{\al} }} )
 \mathcal{H}^\sigma_{b_\al}
 ({e_{\SSTP_{\al} }\color{magenta} \otimes    
e_{\SSTP_{\al} }})
=
({e_{\SSTP_{\emp} }\color{magenta} \otimes   
e_{\SSTP_{\al} }} )
 \mathcal{H}^\sigma_{b_\al}
 ({e_{\SSTP_{\al} }\color{magenta} \otimes  
e_{\SSTP_{\al}^\flat }}) $   (a similar statement holds for the latter element).    }

We wish to inductively pass between the paths   $\SSTP_{\al }\otimes\reflectpath$ and $\SSTP_{\emp\al}$ (respectively  $  \reflectpath\otimes \SSTP_{\al }$ and $\SSTP_{\al\emp}$) 
{\color{black}   by means of a visual  timeline (as in \cref{pathbendingfork}).   
   This allows us to  factorise}
   KLR-forks    and to simplify our proofs later on.    To this end we define 
\begin{align*}
\Forkq&=
 {\SSTP}^{ q\emptyset}
    \boxtimes 
{\sf M}_i^{\exx }
   \boxtimes 
 {\sf P}_{i+1}^{\exx\text{--}q  }   {\color{magenta}\otimes}_\al  
 {\sf M}_i^{b_\al-q }    \boxtimes{\sf P}_i^{\exx }   
\\
{\sf F }_{q,\al\emp} &=
{\sf M}_i^{\exx }
   \boxtimes 
 {\sf P}_{i}^{\exx\text{--}q  }  \boxtimes 
 {\sf M}_i^{b_\al-q }    \boxtimes{\sf P}_{i+1}^{\exx }   
 {\color{magenta}\otimes}_\al   {\SSTP}^{ q\emptyset}
\end{align*}
 and  we remark that 
\begin{align*}
 \Forkzero = 
 \SSTP_{{\al}}\otimes \SSTP_{\al }^\flat  \quad  
 \Forkb   =   \SSTP_{{\emp}}\otimes \SSTP_{{\al}} 
  \quad {\sf F }_{0,\al\emp}= 
\SSTP_{\al }^\flat  \otimes  \SSTP_{{\al}}  \quad  
{\sf F }_{\exx,\al\emp}   =    \SSTP_{{\al}}\otimes  \SSTP_{{\emp}}.
  \end{align*}
We  define 
$  {\sf fork}^{ {{\emp \al}} } _{ {{\al \al}} }(q) =
\Upsilon _{\Forkqplus	}^{ 
\Forkq} 
 $  and 
$  {\sf fork}^{ {{  \al\emp}} } _{ {{\al \al}} }(q) =
\Upsilon _{{\sf F }_{q+1,\al\emp}	}^{ 
{\sf F }_{q,\al\emp}} 
 $ 
 for $0\leq k < \exx$ and we  factorise  the KLR-forks as follows 
  \begin{align*}
 {\sf fork}_{\al \al}^{\emp \al} &=
 e_{  \SSTP_{\emp \al}  }\circ 
{\sf fork}^{{{\emp\al}} }_{{{\al\al}} }  ({\exx-1}) \circ 
 \cdots  \circ 
{\sf fork}^{{{\emp\al}} }_{{{\al\al}} } (1)   \circ 
{\sf fork}^{{{\emp\al}} }_{{{\al\al}} } (0)  \circ  e_{ \SSTP _{\al }	\otimes \SSTP_{\al }^\flat 	}
\\
 {\sf fork}_{\al\al}^{\al\emp} &=
 e_{  \SSTP_{\al\emp}  }\circ 
{\sf fork}^{{{\al\emp}} }_{{{\al\al}} }  ({\exx-1}) \circ 
 \cdots  \circ 
{\sf fork}^{{{\al\emp}} }_{{{\al\al}} } (1)   \circ 
{\sf fork}^{{{\al\emp}} }_{{{\al\al}} } (0)  \circ  e_{\SSTP_{\al }^\flat \otimes \SSTP _{\al }	 	} .  
\end{align*}

\begin{eg} 
Let $h=1$, $\ell=3$, $e=6$, $\sigma=(0,2,4)\in\ZZ^3$  and  $\al=\varepsilon_2-\varepsilon_3$ (thus $\exx=2$).   
We have  
    \begin{align*}
    \SSTP_{\al}  \otimes \SSTP_{\al }^\flat 	&=
    (  \varepsilon_1,  \varepsilon_3,  \varepsilon_1,
  \varepsilon_3,  \varepsilon_3,  \varepsilon_3)  \otimes 
    (  \varepsilon_1,  \varepsilon_3,  \varepsilon_1,
  \varepsilon_3,  \varepsilon_2,  \varepsilon_2)  \\
&=(  \varepsilon_1,  \varepsilon_3,  \varepsilon_1,
  \varepsilon_3,  \varepsilon_3,  \varepsilon_3, 
  \varepsilon_1,  \varepsilon_2,  \varepsilon_1,
    \varepsilon_2,  \varepsilon_3,  \varepsilon_3)
    \\
    \SSTP_{ \emp \al  }&=  ( 
  \varepsilon_1,  \varepsilon_2,  \varepsilon_3,
   \varepsilon_1,  \varepsilon_2,  \varepsilon_3,   \varepsilon_1,  \varepsilon_3,  \varepsilon_1,
  \varepsilon_3,  \varepsilon_3,  \varepsilon_3)
\end{align*} 
are both  dominant paths terminating at   $  (1^4 \mid  1^2 \mid  1^6)   \in 
 \mathscr{P}_{1,3}(12)
$.   
The KLR-fork  diagram is as follows 
    $$
{\sf fork}_{\al \al }^{{\emp  \al } } \; = \;\; \begin{minipage}{10cm} 


  \end{minipage} $$
 
\end{eg}

  The following proposition allows us to see that these two elements  are essentially the same. 
  We will see in the proof that  the ``timelines" for the fork generators  
  allow us to proceed step-by-step   
  (the  steps are  indexed by   $b_\al \geq q \geq 1$).  
 
 \begin{prop}\label{differentforks}
Let $\al\in\Pi $. We have that 
$
{\sf fork}_{\al \al}^{\al \emp}  =
 {\sf adj}^{\al \emp}_{\emp \al} {\sf fork}_{\al \al}^{ \emp \al} 
$.
 \end{prop}
\begin{proof} 
  We note that ${ {\sf A}^{\emp\al}_{\al\emp}(b_\al)}=\SSTP_\emp\otimes \SSTP_\al= {{\sf F}_{b_\al, \emp\al}}$ and 
  ${ {\sf A}^{\emp\al}_{\al\emp}(0)}=  \SSTP_\al\otimes \SSTP_\emp= {{\sf F}_{0,  \al\emp}}$.
 We   claim that 
\begin{equation}\label{label1}{\sf adj}^{\al \emp} _{\emp \al}
 (q-1) \circ
 \Upsilon ^{ {\sf A}^{\emp\al}_{\al\emp}(q)}_{{\sf F}_{q,\al\emp}}
 \circ
{\sf fork}_{\al \al}^{ \emp \al} (q-1)
=
 \Upsilon ^{ {\sf A}^{\emp\al}_{\al\emp}(q-1)}_{{\sf F}_{q-1,\al\emp}}
\end{equation} for $b_\al\geq q\geq 1$.  
The result follows immediately from \cref{forkkkk} once we have proven the claim.   
We label the top and bottom vertices of the lefthand-side of \cref{label1} by the paths 
$\northT_q={ {\sf A}^{\emp\al}_{\al\emp}(q)}$ and $\southT_q={{\sf F}_{q,\emp\al}}$ respectively.
  We remark $\res({{\sf F}_{q,\emp\al}})=\res({{\sf F}_{q,\al\emp }}) $
 (as these paths are obtained from each other by reflection) and so this labelling makes sense.   

We now prove the claim.  There are two strands in the concatenated diagram which do not respect step-labels.  
Namely, the $r_q$-strands (for some $r_q\in \ZZ/e\ZZ$)  connecting the 
  $ 
 \northT_q^{-1}(q,\eps_{i})$
 and  $
  \southT_q^{-1}(b_\al+q,\eps_{i+1})
  $ top and bottom vertices and the strand connecting the
  $
  \northT_q^{-1}(b_\al+q,\eps_{i+1})
$ 
  and $ 
  \southT_q^{-1}(q,\eps_{i})
$   
 top and bottom  vertices. 
   There are four crossings of non-zero degree in the product, all of which involve the former, ``distinguished", $r_q$-strand.  Namely,  the distinguished $r_q$-strand passes from 
   $ 
 \northT_q^{-1}(q,\eps_{i})$ to the left through the 
 latter  $r_q$-strand and then through the vertical  $(r_q+1)$-strand connecting 
 the $\northT^{-1}(b_\al+q,\eps_{i+1})$
 and $\southT^{-1}(b_\al+q,\eps_{i+1})$ vertices before then passing back again  through both these strands and terminating at   $
  \southT_q^{-1}(b_\al+q,\eps_{i+1})
  $.  (The distinguished strand crosses several other strands in the process, but the crossings are  
of degree zero and so can be undone trivially, by case 2 of relation \ref{rel1.10}.)
Using case 4 of relation \ref{rel1.10}, we pull the distinguished  $r_q$-strand rightwards through the $(r_q-1)$-strand 
and hence change the sign and obtain a dot on the $r_q$-strand 
(the term with a dot on the $(r_{q}+1)$-strand  is zero by case 1 of relation \ref{rel1.10} and the commutativity relations).  
Using relation \ref{rel1.8}, we   pull the dot on the distinguished  strand 
rightwards through the crossing of $r_q$-strands and hence undo this crossing,  
kill the dot, and change the sign again (the other term is again zero by case 1 of relation \ref{rel1.10} and the commutativity relations). 
 The resulting diagram has no double-crossings and respects step labels and thus is equal to the righthand-side of \cref{forkkkk}, as required.   
 \end{proof}

   \subsection{The KLR hexagon diagram} \label{adjustgen}
 We now define the  hexagon in the KLR algebra.  
 We let $\al , \bet  \in \Pi  $  label non-commuting reflections.   
We assume, without loss of generality, that $j=i+1$.  We have two cases to consider: if $b_\al \geq b_\bet$ then we must deform the path  
$ 
  \SSTP_{\al \bet \al}$
into the path
$ 
 \SSTP_{\emp-\empb}\otimes \SSTP_{ \bet \al \bet}$ and if 
 $b_\al \leq b_\bet$ then we must deform the path  
$ 
 \SSTP_{\empb-\emp}\otimes  \SSTP_{\al \bet \al}$
into the path
$ 
  \SSTP_{ \bet \al \bet}$, where here 
  $\emp -\empb:=\emptyset^{b_\al-b_\bet}$.

 \begin{prop}\label{analrsusususus}
 The elements $\psi^{ \SSTP_{\al \bet \al}}_{ \SSTP_{\emp-\empb}\otimes\SSTP_{  \bet \al  \bet}}$ and $ \psi^{  \SSTP_{\empb-\emp}\otimes \SSTP_{\al \bet \al}}_{ \SSTP_{  \bet \al  \bet}}$  are independent of the choice of reduced expressions for $b_\al\geq b_\bet$ and $b_\bet\geq b_\al$, respectively.   \end{prop}

\begin{proof}
 We consider the first case as the second is similar.   
 The bad triples of $\psi^{ \SSTP_{\al \bet \al}}_{ \SSTP_{\emp-\empb}\otimes\SSTP_{  \bet \al  \bet}}$ 
are precisely the triples labelled by the integers 
 $$
\SSTP^{-1}_{\al \bet \al} (q,\eps_i)<
\SSTP^{-1}_{\al \bet \al} (b_{\al\bet} +q\pm1,\eps_{i+2})<
\SSTP^{-1}_{\al \bet \al} (b_\al +q,\eps_{i+1}) 
$$
for $1\leq q \leq b_\al  $, where the first and third steps have residue $r_q\in \ZZ/e\ZZ$ and the second   has residue  
$r_{q\pm1}=r_q\mp1\in \ZZ/e\ZZ$.  
 Thus it is enough to consider the subexpression, $\psi_\w$, formed from the  union of the $(r_q,r_{q+1})$-strands for $0\leq q \leq b_\al $ enumerated above.
 We set   $\SSTT= \SSTP _{\al \bet \al} $ and 
 $\southT=  \SSTP_{\empb-\emp}\otimes  \SSTP_{\al \bet \al}$ and we let 
 $$
{\sf t}_i(q)=
\northT^{-1}  (q,\eps_i) \quad
{\sf t}_{i+1}(q)=
\northT^{-1}  (b_\al +q,\eps_{i+1})  \quad
{\sf t}_{i+2}(q)=
\northT^{-1}  (b_{\al\bet} +q ,\eps_{i+2})  
$$
$$
 {\sf b}_i(q)=
\southT^{-1}  (q,\eps_i)  \quad
{\sf b}_{i+1}(q)=
\southT^{-1}  (b_\al +q,\eps_{i+1}) \quad
{\sf b}_{i+2}(q)=
\southT^{-1}  (b_{\al\bet} +q ,\eps_{i+2}) 
$$ 
 for   $0\leq q\leq  b_\al+1$.      
We have  that 
$${\sf t}_i(q) < {\sf t}_i(q+1)< {\sf t}_{i+2}(q)<{\sf t}_{i+2}(q+1)<
 {\sf t}_{i+1}(q)<{\sf t}_{i+1}(q+1)$$
$${\sf b}_i(q) > {\sf b}_i(q+1)> {\sf b}_{i+2}(q)>{\sf b}_{i+2}(q+1)
> {\sf b}_{i+1}(q)>{\sf b}_{i+1}(q+1)$$
for $1\leq q \leq b_\al$ 
and  
$$
 {\sf t}_i(1)< {\sf t}_{i+2}(0)< {\sf t}_{i+1}(1)
 \qquad\quad 
 {\sf t}_i(b_\al)< {\sf t}_{i+2}(b_\al+1)< {\sf t}_{i+1}(b_\al)
   $$
$$
 {\sf b}_i(1)> {\sf b}_{i+2}(0)> {\sf b}_{i+1}(1)
 \qquad\quad 
 {\sf b}_i(b_\al)> {\sf b}_{i+2}(b_\al+1)> {\sf b}_{i+1}(b_\al).  
   $$
     Thus the subexpression $\psi_\w$ is the nib truncation of a 
     quasi-$({b_\al+2})$-expression  for $  w=(13)\in \mathfrak{S}_3$,  which is independent of the choice of expression by \cref{nibs}.  Thus the result follows.  
  \end{proof}

  We are now free to define the {\sf KLR-hexagon} to be the element 
 $$ 
 {\sf hex}^{\al\bet\al}_{\bet\al\bet}:=
 \Upsilon ^{ \SSTP_{\al \bet \al}}_{ \SSTP_{\emp-\empb}\otimes\SSTP_{  \bet \al  \bet}} \qquad\text{or}\qquad
  {\sf hex}^{\al\bet\al}_{\bet\al\bet}:=
  \Upsilon ^{  \SSTP_{\empb-\emp}\otimes \SSTP_{\al \bet \al}}_{ \SSTP_{  \bet \al  \bet}}
  $$ 
 for $b_\al\geq b_\bet$ or $b_\al \leq b_\bet$ respectively, which are independent of the choice of reduced expressions.  See \cref{braidexample} for an example.     
We wish to inductively pass between the paths   ${ \SSTP_{\al \bet \al}}$ and ${ \SSTP_{\emp-\empb}\otimes\SSTP_{  \bet \al  \bet}} $  {\color{black}   by means of a visual  timeline (as in \cref{Steinberg}).   
   This allows us to  factorise}   the   KLR-hexagon    and to simplify our proofs later on.   
First assume that $b_\al \geq b_\bet$.   We define ${\sf H}_{q,\al\bet\al}$ to be the path 
\begin{align*}
 \begin{cases}
\SSTP_{q\emptyset}\boxtimes   \REMOVETHESE{i}{\exx}
  \boxtimes \ADDTHIS{i+1}{\exx } 
  {{\; \color{magenta}\otimes}}_\al\;  \REMOVETHESE {i+1} {b_\bet-q}  
\boxtimes \ADDTHIS{i+2} {b_\bet} 
{{\; \color{cyan}\otimes}}_\bet\;  
 \Pdiptwo^q
\boxtimes \REMOVETHESE { i  } {b_\al-q}
\boxtimes \ADDTHIS {i+1} {\exx }&0\leq q \leq b_\bet \\[6pt]
\SSTP_{q\emptyset}\boxtimes   \REMOVETHESE{i}{\exx}
  \boxtimes \ADDTHIS{i+1}{ b_{\al\bet}-q} 
{{\; \color{magenta}\otimes}}_\al\;
  \ADDTHIS{ i+2} {b_\bet} 
{{\; \color{cyan}\otimes}}_\bet\;  
 \Pdiptwo^{b_\bet}
\boxtimes \REMOVETHESE { i  } {b_\al-q}
\boxtimes \ADDTHIS {i+1} {\exx }&b_\bet \leq q \leq b_\al \\[6pt]
\SSTP_{\emp}\boxtimes
   \REMOVETHESE{i}{\exx}
  \boxtimes \ADDTHIS{i+1}{ b_{\al\bet}-q} 
{{\; \color{magenta}\otimes}}_\al\;
  \ADDTHIS{ i+2} {b_\bet} 
{{\; \color{cyan}\otimes}}_\bet\;  
 \Pdiptwo^{b_\bet} 
\boxtimes \ADDTHIS {i+1} {\exx }
\boxtimes \ADDTHIS {i} {q-\exx }
&b_\al \leq q \leq  b_{\al\bet}\\[6pt]
\end{cases}
\end{align*}
      This is demonstrated in the first 5 paths in \cref{Steinberg}.  
 We now come from the opposite side to meet in the middle.   We define ${\sf H} _{q, \bet\al\bet}$ to be the path 
\begin{align*}
\begin{cases}
\SSTP_{q\emptyset}\boxtimes
\REMOVETHESE{i+1}{b_\bet-q} \boxtimes 
 \REMOVETHESE {i,i+1} {q}  
\boxtimes
\ADDTHIS{i+2}{b_\bet}
{{\; \color{cyan}\otimes}}_\bet\; 
  \REMOVETHESE{i} {b_\al-q} 
\boxtimes \ADDTHIS{i+1} {b_\al}
{{\; \color{magenta}\otimes}}_\al\;  
\REMOVETHESE{i+1}{b_\bet} \boxtimes 
\ADDTHIS{i+2}{b_\bet}&0\leq q \leq b_\bet \\[6pt]
\SSTP_{\empb }  \boxtimes
  \REMOVETHESE{i} {q-b_\bet} \boxtimes \REMOVETHESE {i,i+1} {b_\bet}  \boxtimes  \ADDTHIS{i+2}{b_\bet}
{{\; \color{cyan}\otimes}}_\bet\;
 \REMOVETHESE{i} {b_\al-q} 
\boxtimes \ADDTHIS{i+1} {b_\al}
{{\; \color{magenta}\otimes}}_\al\;  
\REMOVETHESE{i+1}{b_\bet} \boxtimes 
\ADDTHIS{i+2}{b_\bet}&b_\bet  \leq q \leq b_\al  
\\[5pt]  
\SSTP_{\empb }  \boxtimes  
\REMOVETHESE {i } {  q-b_\bet }   \boxtimes  \ADDTHIS{i+2}{b_\bet}
 {{\; \color{cyan}\otimes}}_\bet\; 
\REMOVETHESE{i ,i+2}{q-b_\al} \boxtimes   \ADDTHIS{i+1} {b_\al}
{{\; \color{magenta}\otimes}}_\al\;  
\REMOVETHESE{i+1}{ b_{\al\bet}-q} \boxtimes 
\ADDTHIS{i+2}{b_\bet}&b_\al     \leq q  \leq b_{\al\bet}
\end{cases}
\end{align*}
  This is demonstrated in the final 5 paths in \cref{Steinberg}. 
  While the definitions seems technical, one can intuitively think of this process as ``flattening" the path layer-by-layer 
  {\color{black}   by means of the timeline depicted in}
   \cref{Steinberg}.     We see that ${\sf H} _{{b_{\al\bet}},\al \bet \al}=\SSTP_{\emp-\empb}\boxtimes {\sf H} _{{b_{\al\bet}},\bet \al \bet}$.

 \begin{figure}[ht!]
 $$  \begin{minipage}{4.05cm}
\end{minipage}$$
 \caption{
 An example of a timeline for the KLR hexagon.   
 Mutating from $  \SSTP_{\al\bet\al}$ to  $\SSTP_{\emp-\empb}\otimes \SSTP_{ \bet\al\bet}$ for $b_\al\geq b_\bet$  (again we do not picture the determinant paths).  
 Steps in the procedure should be read from left-to-right along successive rows
 (the paths are $ {\sf H}_{0,\al\bet\al} $, ${\sf H}_{1,\al\bet\al} $, 
 ${\sf H}_{2,\al\bet\al}$, ${\sf H}_{3,\al\bet\al}$, ${\sf H}_{4,\al\bet\al} $, ${\sf H}_{4,\al\bet\al} =\SSTP_\emptyset\boxtimes 
  {\sf H}_{4,\bet\al\bet} $, $  {\sf H}_{3,\bet\al\bet} $, 
  $  {\sf H}_{2,\bet\al\bet} $, $  {\sf H}_{1,\bet\al\bet} $, $  {\sf H}_{0,\bet\al\bet}$).
}
 \label{Steinberg}\end{figure}

  
  We now assume that $b_\al \leq b_\bet$.   We define 
  ${\sf H}_{q,\al\bet\al}$ to be the path 
\begin{align*}
 \begin{cases}
 \SSTP_{q\emptyset}\boxtimes   \REMOVETHESE{i}{\exx}
  \boxtimes \ADDTHIS{i+1}{\exx } 
  {{\; \color{magenta}\otimes}}_\al\;  \REMOVETHESE {i+1} {b_\bet-q}  
\boxtimes \ADDTHIS{i+2} {b_\bet} 
{{\; \color{cyan}\otimes}}_\bet\;  
 \Pdiptwo^q
\boxtimes \REMOVETHESE { i  } {b_\al-q}
\boxtimes \ADDTHIS {i+1} {\exx }
  &0\leq q \leq  b_\al \\[6pt]
  \SSTP_{\emp }\boxtimes 
   \REMOVETHESE{i}{\exx}
  \boxtimes 
  \ADDTHIS{i+1}{\exx } 
  {{\; \color{magenta}\otimes}}_\al\; 
   \REMOVETHESE {i+1} {b_\bet-q}  
\boxtimes   \REMOVETHESE {i,i+1} {q-b_\al}  
\boxtimes
 \ADDTHIS{i+2} {b_\bet} 
{{\; \color{cyan}\otimes}}_\bet\;  
 \Pdiptwo^{b_\al}
 \boxtimes
 \ADDTHIS {i+1} {q }  \boxtimes
 \ADDTHIS {i} {q-b_\al }
&b_\al \leq q \leq b_\bet 
\\[6pt]
\SSTP_{\emp }\boxtimes 
   \REMOVETHESE{i}{\exx}
  \boxtimes 
  \ADDTHIS{i+1}{b_{\al\bet-q}}
  {{\; \color{magenta}\otimes}}_\al\; 
   \REMOVETHESE {i,i+1} {b_\bet-b_\al}  
\boxtimes
 \ADDTHIS{i+2} {b_\bet} 
{{\; \color{cyan}\otimes}}_\bet\;  
 \Pdiptwo^{b_\al}
 \boxtimes
 \ADDTHIS {i+1} {b_\al }
  \boxtimes
 \ADDTHIS {i} {q-b_\al }&b_\bet \leq q \leq b_{\al  \bet}
\end{cases}
 \end{align*}
    We now come from the opposite side to meet in the middle.   We define ${\sf H} _{q, \bet\al\bet}$ to be the path  \begin{align*}
\begin{cases}
 \SSTP_{q\emptyset}\boxtimes
\REMOVETHESE{i+1}{b_\bet-q} \boxtimes 
 \REMOVETHESE {i,i+1} {q}  
\boxtimes
\ADDTHIS{i+2}{b_\bet}
{{\; \color{cyan}\otimes}}_\bet\; 
  \REMOVETHESE{i} {b_\al-q} 
\boxtimes \ADDTHIS{i+1} {b_\al}
{{\; \color{magenta}\otimes}}_\al\;  
\REMOVETHESE{i+1}{b_\bet} \boxtimes 
\ADDTHIS{i+2}{b_\bet}&0\leq q \leq b_\al \\[6pt]
\SSTP_{q\emptyset}\boxtimes
\REMOVETHESE{i+1}{b_\bet-q} \boxtimes 
 \REMOVETHESE {i,i+1} {q}  
\boxtimes
\ADDTHIS{i+2}{b_\bet}
{{\; \color{cyan}\otimes}}_\bet\; 
   \ADDTHIS{i+1} {b_\al}
{{\; \color{magenta}\otimes}}_\al\;  
\REMOVETHESE{i+1}{b_{\al\bet}-q} \boxtimes 
\ADDTHIS{i+2}{b_\bet}&b_\al\leq q \leq b_\bet \\[6pt]
\SSTP_{\empb}\boxtimes
\REMOVETHESE{i}{q-b_{\bet}} \boxtimes 
 \REMOVETHESE {i,i+1} {b_{\bet\bet}-q}  
\boxtimes
\ADDTHIS{i+2}{b_\bet}
{{\; \color{cyan}\otimes}}_\bet\; 
   \ADDTHIS{i+1} {b_\al}
{{\; \color{magenta}\otimes}}_\al\;  
\REMOVETHESE{i+1}{b_{\al\bet}-q} \boxtimes 
\ADDTHIS{i+2}{b_\bet}&b_\bet\leq q \leq b_{\al\bet} \\[5pt]  \end{cases}
\end{align*} 
 With our paths in place, this  allows us to
define 
$$
{\sf hex} ^{\al\bet\al}(q) =\Upsilon^{{\sf H}_{q,\al\bet\al}}_{{\sf H}_{q+1,\al\bet\al}}
\qquad
{\sf hex} _{ \bet\al\bet}(q) =\Upsilon_{{\sf H}_{q,\bet\al\bet}}^{{\sf H}_{q+1,\bet\al\bet}}
$$
and we set 
$$
 {\sf hex}^{ \al\bet\al}
=
\prod_{ b_{\al\bet}> q \geq 0}  
  {\sf hex} ^{ \al\bet\al}(q)
  \qquad
  {\sf hex}_{\bet\al\bet}
=
\prod_{ 0 \leq q \leq b_{\al\bet} }
  {\sf hex} _{\bet\al\bet}(q) 
    $$
which allows us   to  factorise  the hexagon generators as follows 
$${\sf hex}^{\al\bet\al}_{\bet\al\bet}
=
\begin{cases}
{\sf hex}^{\al\bet\al}  (e_{\SSTP^{\emp-\empb}}\otimes {\sf hex}_{\bet\al\bet})
  &\text{for }b_\al \geq b_\bet 
\\
(e_{\SSTP^{\empb-\emp}}\otimes{\sf hex}^{\al\bet\al}  ) {\sf hex}_{\bet\al\bet}
    &\text{for }b_\al \leq b_\bet 
\end{cases}
$$ and, finally, we define $$
   {\sf hex}^{\empb\al\bet\al 	 }_{ \emp\bet\al\bet 	 }= 
  \begin{cases}
 e_{\SSTP_{\emp }}\otimes  {\sf hex}^{\al\bet\al}_{\bet\al\bet}	  	&\text{if }b_\al\le b_\bet \\
   e_{\SSTP_{\empb }}\otimes  {\sf hex}^{\al\bet\al}_{\bet\al\bet}   &\text{if }b_\al\geq  b_\bet  
  \end{cases}  
$$
 the latter notation will be useful when we wish to consider products of such hexagons without assuming $b_\al\geq b_\bet$ or vice versa.   
  Finally, the following shorthand will come in useful when addressing some of the   relations in \cref{relations}.   
 Recall that adjustment is invertible.  With this in mind, we set 
     $$
 {\sf hex}^
 { \underline{v}  \bet\al\bet   \underline{w}  \emp }
 _{  \underline{v} \al\bet\al  \underline{w}  \empb}=   
 {\sf adj}_{{ \underline{v} \emp  \bet\al\bet   \underline{w}   }}
^{{ \underline{v}  \bet\al\bet   \underline{w}  \emp }}
\big( e_{\SSTP_{\underline{v}}}\otimes 
{\sf hex}^
 {    \emp  \bet\al\bet } 
 _{   \empb \al\bet\al }\otimes  e_{\SSTP_\w  }   \big)
 {\sf adj}^{{ \underline{v} \empb \al\bet   \al \underline{w}   }}
_{{ \underline{v}  \al\bet   \al \underline{w}   \empb }}
  =\Upsilon
  _{{ \underline{v}  \al\bet   \al \underline{w}   \empb }}
  ^
 { \underline{v}  \bet\al\bet   \underline{w}  \emp }     $$
 where the second equality follows by removing the resulting double-crossings using   \cref{adjust1} in each case.    
Independence of the   reduced expression    follows from  residue-commutativity of adjustment.  Alternatively,   the reader is invited to make   minor modifications  to the proof  of \cref{analrsusususus}.

  \begin{figure}[ht!] 
 $$  
  $$

\!\!\!\!
    \caption{  Let $h=3$, $\ell=1$, $e=5$ and $\al=\eps_3-\eps_1$, $\bet=\eps_1-\eps_2$.   
    We depict the element ${\sf hex}^{\bet\al\bet}_{\al\bet\al}$ and   highlight   the dilated word ${\rm nib}(1,3)_5$ in bold. 
     The reader should compare the 11 highlighted strands with  the  diagram from $\mathfrak{S}_11$ depicted  in \cref{labeler2}.     (We have drawn all bad-crossing so that they bi-pass on the right.)    }
\label{braidexample}    \end{figure}

     \subsection{The commuting strands diagram}
 \newcommand{\jay}{k}
  \newcommand{\kay}{j}

Let $\gam,\bet\in \Pi$ be roots labelling commuting reflections (in terms of \cref{conventioning}, this is equivalent to $|k -j |>1$).    
We wish to understand the  morphism relating the paths 
    $\SSTP_{\gam }\otimes \SSTP_{\bet}$ to $\SSTP_{\bet}\otimes \SSTP_\gam$.  
 We  suppose without loss of generality that $b_\gam   \geq b_\bet $.  

\begin{prop}\label{resssssss}
 The element $\psi^{\SSTP_{\gam }\otimes \SSTP_{\bet} }
  _{\SSTP_{\bet}\otimes \SSTP_\gam}$ is independent of the choice of reduced expression
\end{prop}
\begin{proof}
 There are precisely $b_{\gam\bet}$ like-labelled crossings. The first $b_\gam$ of these connect  the
 $
 \SSTP^{-1}_{\gam\bet}(q,\eps_j)
$th and $
 \SSTP^{-1}_{\gam\bet}(b_\gam+q,\eps_{j+1})
 $th
 northern vertices to the  
 $  \SSTP^{-1}_{\bet\gam}(q,\eps_j)$th 
 and $
 \SSTP^{-1}_{\bet\gam}(b_\bet+q,\eps_{j+1})
$th 
southern vertices for $1\leq q \leq b_{\gam}$.   
The latter $b_\bet$ of these   connect  the
   $
 \SSTP^{-1}_{\gam\bet}(b_\bet+q,\eps_{k+1}) 
$th and $
 \SSTP^{-1}_{\gam\bet}( q,\eps_{k})
 $th
 northern vertices to the  
 $
 \SSTP^{-1}_{\bet\gam}(b_\bet+q,\eps_{k+1}) 
$th and $
 \SSTP^{-1}_{\bet\gam}( q,\eps_{k})
 $th
southern vertices for $1\leq q \leq b_{\gam}$.

For $k\neq {\aatch}$ (respectively $k={\aatch}$) each of the first $1\leq q \leq b_\gam$ (respectively $1<q \leq b_\gam$) like-labelled crossings forms a braid with precisely one other strand, namely
  that connecting the 
$\SSTP^{-1}_{\gam\bet}(b_\bet+q,\eps_{k+1})$th
top vertex to the 
$\SSTP^{-1}_{\bet\gam}(b_\bet+q,\eps_{k+1})$th bottom vertex for $1\leq q \leq b_\gam$ (respectively $1\leq q < b_\gam$).  
This strand is of non-adjacent residue (by our assumption that $\gam$ and $\bet$ label commuting reflections).  
The latter  $b_\bet$ cases can be treated similarly.  

Thus   each of the  braids involving a like-labelled crossing (either totalling $b_{\bet\gam}$ if $k,j\neq {\aatch}$ or  $b_{\bet\gam}-1$ otherwise) is residue-commutative.  
Thus $\psi^{\SSTP_{\gam\bet} }
  _{\SSTP_{\bet\gam}}$ is residue commutative and the result follows.  
\end{proof}

   Thus we are free to define the {\sf KLR-commutator} to be the element 
 $$ 
 {\sf com}^{\gam\bet}_{\bet\gam}:= 
\Upsilon   ^{\SSTP_{\gam }\otimes \SSTP_{\bet} }
  _{\SSTP_{\bet}\otimes \SSTP_\gam} $$ 
 which is independent of the choice of reduced expression.    
We wish to inductively pass between the paths   ${\SSTP_{\gam }\otimes \SSTP_{\bet} }$
 and ${\SSTP_{\bet}\otimes \SSTP_\gam} $     by means of a visual  timeline (as in \cref{chopin}).   
    
\!\!\!
  \begin{figure}[ht!]
  $$
   \begin{minipage}{3.4cm}
\end{minipage}
 $$  
  \caption{
  An example   timeline for the KLR commutator. 
 We mutate from $  \SSTP^{\gam\bet}$ to  $\SSTP_{ \bet\gam} $
  for $ b_\gam=4,  b_\bet=3 $.    
Reading from left-to-right along successive rows  the paths are $\SSTP^{-1,\gam\bet} $, 
$\SSTP^{0,\gam\bet} $, $\SSTP^{1,\gam\bet} $, $\SSTP^{2,\gam\bet} $, $\SSTP^{3,\gam\bet} =
 \SSTP_{2, \bet\gam} $, $ \SSTP_{1, \bet\gam} $, 
 $ \SSTP_{0, \bet\gam} $, $ \SSTP_{-1, \bet\gam} $. 
We draw   paths in the projection onto  $\RR\{\eps_j+\eps_{j+1},\eps_k+\eps_{k+1}\}$.   
  }
\label{chopin}   \end{figure}

 We define 
\begin{align*}
{\sf C}^{q,\gam\bet}	&=
\begin{cases}
   \REMOVETHESE {k } {b_\gam }  \boxtimes 
\ADDTHIS {k +1} {b_\gam }
{\color{darkgreen}\otimes} _\gam 
\REMOVETHESE {j } {b_\bet }  \boxtimes \ADDTHIS {j +1} {  b_\bet  }
  	&\text{for }q=-1
\\[4pt]
   \REMOVETHESE {k } {b_\gam }  
{\color{darkgreen}\otimes} _\gam  
\REMOVETHESE {j } {b_\bet }  \boxtimes \ADDTHIS {j +1} {  b_\bet  }
   \boxtimes 
\ADDTHIS {k  } {b_\gam }	&\text{for }q=0
\\[4pt]
  \SSTP_{ q\emptyset} \boxtimes 
   \REMOVETHESE {k } {b_\gam-q }  
{\color{darkgreen}\otimes} _\gam 
 \REMOVETHESE {k +1,j } { q}   \boxtimes 
  \REMOVETHESE {j  } {b_\bet -q}   \boxtimes 
  \ADDTHIS {j +1} {  b_\bet  } \boxtimes 
\ADDTHIS {k  } {b_\gam }
  	&\text{for }0< q \leq   b_\bet 
\end{cases}
\\ 
 {\sf C}_{q, \bet\gam}	&=
\begin{cases}
 \SSTP_{\empb}  
{\color{cyan}\otimes} _\bet  
  \REMOVETHESE {k } { q- b_\bet  }  
  \boxtimes  \REMOVETHESE {k ,j +1} {b_\bet   }  
  \boxtimes 
  \REMOVETHESE {k } {b_\gam-q  }  
  \boxtimes \ADDTHIS {j } {  b_\bet  } \boxtimes 
\ADDTHIS {k +1} {b_\gam }
  	&\text{for } 	b_\gam \geq   q>  b_\bet	\\[4pt]		 \SSTP_ {q\emptyset} \boxtimes 	  \REMOVETHESE {j } {b_\bet -q}  
\; {\color{cyan}\otimes} _\bet  \;
 \REMOVETHESE {k ,j +1} {q  }  
  \boxtimes 
  \REMOVETHESE {k } {b_\gam-q  }  
  \boxtimes \ADDTHIS {j } {  b_\bet  } \boxtimes 
\ADDTHIS {k +1} {b_\gam }
  	&\text{for }   b_\bet \geq q >0 
			\\[4pt]
  	  \REMOVETHESE {j } {b_\bet }  
{{\color{cyan}\otimes} _\bet }
 \REMOVETHESE {k } {b_\gam  }   \boxtimes \ADDTHIS {j } {  b_\bet  } \boxtimes 
\ADDTHIS {k +1} {b_\gam }
  	&\text{for } q=0 
		\\[4pt]
 	  \REMOVETHESE {j } {b_\bet }   \boxtimes 
  \ADDTHIS {k +1} {  b_\bet  }
{{\color{cyan}\otimes} _\bet }
 \REMOVETHESE {k } {b_\gam } 
  \boxtimes 
 \ADDTHIS {k +1} {b_\gam } 
  	&\text{for }q=-1	 
 \end{cases}
\end{align*}
and we note that ${\sf C}_{b_\gam,\bet\gam}={\sf C} ^{ b_\bet,\gam\bet}$ 
(to see this, note that the definition of the former contains  a tensor product $\otimes_\gam$ and the latter contains a tensor product $\otimes_\bet$ and this explains the differences in     the subscripts).  
We now define
$$
  {\sf com}^{ q, {{\gam\bet}} }= 
\Upsilon^  {{\sf C}^{ q, {{\gam\bet}} }}_{{\sf C}^{ q+1, {{\gam\bet}} }}
\qquad 
{\sf com}_{ q, {{\bet\gam}} }= 
 \Upsilon _ {{\sf C}_{ q, {{\bet\gam}} }}^{{\sf C}_{ q+1, {{\bet\gam}} }}.  
  $$
This allows us to factorise  
$$ {\sf com}^{ {{\gam\bet}} }_{ {{\bet\gam}} }= 
{\sf com}^{ {{\gam\bet}} }
{\sf com}_{ {{ \bet\gam}} }  
 \qquad {\sf com}^{ {{\gam\bet}} }
=
\prod_{ -1 \leq q < b_{ \bet} }
 {\sf com}^{ q, {{\gam\bet}} }
\qquad
{\sf com}_{{{\bet\gam}} }
=
\prod_{   b_{ \gam}> q \geq -1 }
 {\sf com}_{ q, {{ \bet\gam}} }  . $$
   The following notation will come in useful in \cref{relations}
     $$
 {\sf com}^
 { \underline{v}  \gam\bet  \underline{w}   }
 _{  \underline{v} \bet\gam  \underline{w} }=   
   e_{\SSTP_{\underline{v}}}\otimes 
 {\sf com}^
 { \underline{v}  \gam\bet  \underline{w}   }
 _{  \underline{v} \bet\gam  \underline{w} }   
\otimes  e_{\SSTP_\w  }  .
    $$

  \subsection{The isomorphism}  
  Finally, we now explicitly state the    isomorphism.  
  Our notation has been chosen so as to make this almost tautological at this point.  
 We suppose that $\al$ and $\bet$ (respectively $\bet$ and $\gam$) label non-commuting  (respectively commuting) reflections. 
   We define 
   \begin{align}\label{thisistheisomorphismfor the intro}
   \Psi : \mathscr{S}^{\rm br}_{\underline{h}} ({n  ,\sigma})
  \xrightarrow {\  \ \ \ }     {\sf f}_{ n,\sigma}\left(\mathcal{H}_{n }^\sigma/\mathcal{H}_{n }^\sigma {\sf y}_{\aatchpair }\mathcal{H}_{n }^\sigma 		\right)     {\sf f}_{ n,\sigma}
  \end{align}to be the map defined on generators ({\color{black} and extended using vertical concatenation and contextualised horizontal  concatenation}) as follows 
 $$\Psi ({\sf 1}_\al)=e_{\SSTP_\al}
\quad 
\Psi ({\sf 1}_\emptyset )=e_{\SSTP_\emptyset}
\quad
\Psi ({\sf 1}_{\al\emptyset}^{ \emptyset\al})
= {\sf adj}_{\al\emptyset}^{\emptyset\alpha}
\quad 
\Psi ({\sf SPOT}_\al^\emp)={\sf spot}_\al^\emp
$$
$$
\Psi ({\sf FORK}_{\al\al}^{\emp\al})={\sf fork}_{\al\al}^{\emp\al}
\quad
\Psi ({\sf HEX}_{\al\bet\al}^{\bet\al\bet})={\sf hex}_{\al\bet\al}^{\bet\al\bet}
\quad
\Psi ({\sf COM}_{\bet\gam}^{\gam\bet})={\sf com}_{\bet\gam}^{\gam\bet}   
$$
  and we extend this to the   flips of these diagrams through their horizontal axes.  
  
  \begin{rmk}\color{black} 
  We note that our use of {\em contextualised} horizontal concatenation implies that 
   \cref{idemp-iso2} holds (see also \cref{idemp-iso3}).
  \end{rmk}

      \section{Recasting the diagrammatic Bott--Samelson relations \\
         in the  quiver Hecke  algebra }\label{relations}

 The purpose of this section is to recast  Elias--Williamson's diagrammatic  relations   of \cref{soergel} in the setting of the quiver Hecke algebra, thus verifying that the map $\Psi _n$ is indeed    a     (graded) $\ZZ$-algebra homomorphism. 
       We have already provided timelines which discretise  each
    Soergel generator (which we think of as  a continuous morphism between paths with a unique singularity, where the strands cross).   
We will verify most of the  Soergel relations  via  a  similar discretisation    
 process which factorises the Soergel relation  into simpler steps; we again record this is a visual timeline.  
 We check each relation in turn, but leave it as an exercise for the reader to verify the flips of these relation through their vertical axes (the flips through horizontal axes follow immediately from  the duality, $\ast$).  
 We continue with the  notations  of \Cref{conventioning}.   
 Our    relations fall into three  categories:
\begin{itemize}[leftmargin=*]
\item[$\bullet$] Products involving only hexagons, commutators, and adjustment generators.  Simplifying such products is  an inductive process.  At each step, one simplifies  a non-minimal expression (in the concatenated diagram) to a minimal one  {\em without changing the underlying permutation}.  
This typically involves   a single ``distinguished" strand which    double-crosses some other strands; these double-crossings can be undone using \cref{adjust1}.  (This preserves the parity of like-labelled crossings.)
\item[$\bullet$] \color{black} Products involving a fork or spot generator. 
 Such generators  reflect  one of the indexing paths in an irreversible manner.  
Simplifying such products is  an inductive process.  At each step, one  rewrites   a {\em single pair of crossing strands} (in the  concatenated  permutation)  which {\em do not respect step-labels of the reflected paths}. By undoing this crossing using relation \ref{rel1.8}, we obtain the scalar $-1$ times a new diagram  {\em which does} respect the new step-labelling for the reflected paths.  
 (Thus changing  the   parity of like-labelled crossings and also changing the scalar $\pm1$.)
 \color{black}
 
\item[$\bullet$] Doubly spotted Soergel diagrams (such as the  Demazure relations) for which we argue separately.    
\end{itemize}
In each of the former two cases, we will decorate the top and bottom of the concatenated diagram with  paths $\SSTT$ and $\southT$ (which we define case-by-case) and use the step-labelling from these paths  to keep track of crossings of  strands in the diagram.

 \subsection{The double fork} 
This leftmost  relation in \ref{rel1}   is incredibly  simple to verify, and so    there is no need to record this in  a timeline.  
For $\al \in \Pi$, we must verify that
\begin{equation}\label{doublefork}
\Psi \left(\; \begin{minipage}{2.35cm}\; 
\end{minipage}\right) 
\end{equation}
Thus we need to check that
 \begin{equation}\label{checkme} 
   \left( e_{ \SSTP_{\al} }  \otimes     {\sf fork}_{\al \al}^{ \al \emp}  \right)
 \circ  
 \left(  {\sf fork}^{\al \al}_{\emp \al}
  \otimes   
e_{ \SSTP_{\al} } \right)
=
\left(  {\sf fork}^{ \al \al}_{ \emp \al} 
    \otimes    e_{\SSTP_\emp}\right)
 \circ  
\left( e_{\SSTP_\emp}\otimes 
 {\sf fork}_{\al \al}^{ \al \emp}\right).   
\end{equation}
The permutation  underlying  
 $e_{ \SSTP_{\al} }  \otimes     {\sf fork}_{\al \al}^{ \al \emp}  $ 
 is the element $w^{\SSTT }_{\southT }$  indexed by the pair of paths
 $$
\northT=  \SSTP _\al \otimes \SSTP_\al \otimes \SSTP_\emp
\quad
\text{and}
\quad
\southT=  \SSTP _\al \otimes \SSTP_\al ^\flat \otimes \SSTP_\al
 $$
 which differ only by permuting the
 final  $(\exx {\aatch} +\exx)$   steps.  
The permutation    underlying $  
  {\sf fork}^{\al \al}_{\emp \al}
  \otimes   
e_{ \SSTP_{\al} }$ is the element $w^{\SSTT'}_{\southT'}$ indexed by the pair of    paths 
$$
\northT'=\SSTP_\al \otimes \SSTP_\al^\flat  \otimes \SSTP_\al 
\quad
\text{and}
\quad
\southT'=\SSTP_\emp\otimes \SSTP_\al   \otimes \SSTP_\al .
$$
which differ only by permuting the first   
 $(b_{\al\al} {\aatch} -\exx)$ steps.    
These elements of $\mathfrak{S}_{3b_\al \aatch }$ commute as they permute disjoint subsets of $1,\dots, 3b_\al \aatch $.  
Thus  the elements     $ 
     {\sf fork}^{\al\al}_{\emp\al}
  \otimes   
e_{ \SSTP_{\al} } $
and 
$  e_{ \SSTP_{\al} }  \otimes    {\sf fork}_{\al\emp}^{ \al \al}  
     $ 
 commute   by relation \ref{rel1.7} (and the result follows immediately). 

\begin{rmk}
  The reader might wonder why the element $w^{\SSTT }_{\southT }$ appears to permute  a greater number of strands than 
$  w^{\SSTT' }_{\southT' }$.  
This is  because   our distinguished choice of $\SSTP_\al$ has a total of   $(b_\al {\aatch}-b_\al)$ steps below (or on) the $\al$-hyperplane and $b_\al$ steps above the hyperplane.  
\end{rmk}

\subsection{The one-colour zero relation}  
We now consider the  rightmost  relation in \ref{rel1}. 
For $\al \in \Pi$, we must verify that 
\begin{align}\label{verifydotzero}
 \Psi \left( 
 \begin{minipage}{1.6cm}\; \begin{tikzpicture} [scale=1.4]
 \clip (0,-0.75) rectangle (1,0.75); 
{ \foreach \i in {0,1,2,...,100}
  {
    \path  (-1,-0.75)++(0:0.125*\i cm)  coordinate (a\i);
        \path  (-1,-0.75)++(90:0.125*\i cm)  coordinate (b\i);
     \draw [densely dotted](a\i)--++(90:5);
        \draw [densely dotted](b\i)--++(0:10);
 }
}
   \draw[ black, line width=0.06cm] 
  (-0,-0.75) --++(0:0.5);
  \draw[ magenta, line width=0.06cm] 
  (0.5,-0.75) --++(0:0.5);
    \draw[ black, line width=0.06cm] 
  (-0,0.75) --++(0:0.5);
  \draw[ magenta, line width=0.06cm] 
  (0.5,0.75) --++(0:0.5); 
  \draw[magenta,line width=0.12cm](0.75,-1)--(0.75,0.75);
   \draw[line width=0.12cm,magenta ](0.75,0.9/2) to [out=-150,in=90] (0.25,0pt);
   \draw[line width=0.12cm,magenta ](0.75,-0.9/2) to [out=150,in=-90] (0.25,0pt);
 \end{tikzpicture}\end{minipage} \; \right)  
 =
  {\sf fork}_{\al \al}^{\emp \al}  \circ  {\sf fork}^{\al \al}_{\emp \al}  
 =0
 \end{align} 
For $  b_\al >q\geq 1$ the paths ${{\sf F}_{q,{\emp\al}}}  $ and 
     $ {{\sf F}_{q-1,{\emp\al}}}  $ are   concatenates  of a single   $\al$-crossing  path and   and  a single $\al$-bouncing path.  
      By \cref{adjust1}   we have that 
$$    {\sf fork}_{{\al\al}}^{\emp\al}(q) e_{{\sf F} _{q-1,\emp\al}  }  {\sf fork}_{{\emp\al}}^{\al\al}(q)     = e_{{\sf F} _{q,\emp\al}  }
$$
for $1\leq q <\exx$.  
We apply this from the centre of the product 
${\sf fork}_{\al \al}^{\emp \al}  \circ  {\sf fork}^{\al \al}_{\emp \al}$ which is equal to 
$$
e_{\SSTP ^{ \emp \al}}
     {\sf fork}_{{\al \al}}^{\emp \al}(\exx-1)
 \cdots 
     {\sf fork}_{{\al \al}}^{\emp \al}(0)   
     e_{\SSTP ^{ \al \al}}%
     \circ   e_{\SSTP ^{ \al \al}}
     {\sf fork}^{{\al \al}}_{\emp \al}(0) 
 \cdots 
     {\sf fork}^{{\al \al}}_{\emp \al}(\exx-1)  e_{\SSTP ^{ \emp \al}}
     $$   until we obtain 
\begin{align}\label{righthandds}
{\sf fork}_{\al \al}^{\emp \al}  \circ   {\sf fork}^{\al \al}_{\emp \al}  &=
  e_{{\sf P}_{  \emp \al } }
  {\sf fork}_{{\al \al}}^{\emp \al}(\exx-1)   e_{{\sf F}_{{\exx-1},{\emp\al}} } 
  {\sf fork}^{{\al \al}}_{\emp \al}(\exx-1)  
  e_{{\sf P}_{   \emp \al} }.  
   \end{align}
This  is illustrated in \cref{asfdkhjalksdfhalk}.


 \!\!
 \begin{figure}[ht!] 
$$ 
\scalefont{0.85}\begin{minipage}{6cm} 
\end{minipage}
    $$

 \!\!\caption{Let $h=1$,   $\ell=3$,   $\sigma=(0,2,4)$ and $e=6$.  
The  lefthand-side is 
${\sf fork}_{\al \al}^{\emp \al}  {\sf fork}^{\al \al}_{\emp \al}$; we apply \cref{adjust1} to undo the highlighted strands  (compare the highlighted strands 
  with the highlighted strands of the first diagram of  \cref{transeg}).  
The thick double-crossing of strands in the rightmost diagram    is zero by the first case of relation \ref{rel1.10} (after applying   commutativity relations).}
\label{asfdkhjalksdfhalk}
 \end{figure}
  \!\!

We cannot apply \cref{adjust1}   to the pair of paths   ${\sf F}_{{\exx-1},{\emp\al}}$ and ${\sf F}_{{\exx-2},{\emp\al}}$ because the  former path passes through the $\al$-hyperplane once, whereas the latter passes through/bounces the $\al$-hyperplane twice.    There is a pair of double-crossing  $r$-strand (for some     $r\in \ZZ/e\ZZ$) between the   
$\SSTP_{\emp\al}^{-1}(b_\al,\eps_i)$th
and 
 $\SSTP_{\emp\al}^{-1}(b_{\al\al},\eps_{i+1})$th   
\north    and   \south       vertices   in 
the diagram 
$$   e_{{\sf P}^{   \emp \al} } {\sf fork}_{{\al \al}}^{\emp \al}(\exx-1)   e_{{\sf F}_{{\exx-1},{\emp\al}} } 
  {\sf fork}^{{\al \al}}_{ \emp \al}(\exx-1)   e_{{\sf P}^{   \emp \al} }$$
 This double-crossing of $r$-strands is not intersected by any strand of adjacent residue.  
Therefore the product is zero    by the commutativity relations   and the first case  of relation \ref{rel1.10}, as required. 

\subsection{Fork-spot contraction} 
  We now consider the second  relation  depicted    in \ref{rel1}, namely
\begin{equation}\label{fork-dot-proof} 
(  {\sf spot}^{{\emp}}_{{\al}}     \otimes      e_{\SSTP_{\al}})
   \circ  {\sf fork}_{\emp \al}^{\al \al} 
  =
e_{\SSTP_\emp} \otimes e_{\SSTP _ \al}
 \end{equation}  for $\al \in \Pi$. 
For $0\leq q \leq b_\al$, we   define the
spot-fork path to be 
$${\sf FS}_{{q,\al}} 
  =
 {\SSTP}_{ q\emptyset}
    \boxtimes 
{\sf M}_i^{\exx }
\; {\color{magenta}\otimes }_\al \;
 {\sf P}_{i+1}^{\exx-q  }  \; {\color{magenta}\otimes }_\al \; 
 {\sf M}_i^{b_\al-q }    \boxtimes{\sf P}_{i +1 }^{\exx }
=  {\SSTP}^{ q\emptyset}
    \boxtimes 
{\sf M}_i^{\exx }
\boxtimes 
 {\sf P}_{i}^{\exx-q  }   \boxtimes 
 {\sf M}_i^{b_\al-q }    \boxtimes{\sf P}_{i+1}^{\exx }    
$$
  which  is   obtained from ${\sf F}_{q,\emp\al}$   by reflection by $s_\al$ (see \cref{pathbendingfork2}).  We note that 
$   {\sf FS}_{{b_\al,\al}} =    \SSTP_\emp \otimes \SSTP_\al 
 $
and      ${\sf FS}_{0,\al}=\SSTP_\al ^\flat\otimes \SSTP_\al$.  
Thus these spot-fork paths allow us to  iteratively   prove \cref{fork-dot-proof}, as we will see below.

 \vspace{-0.1cm}
\begin{figure}[ht!]
$$
    \begin{minipage}{2.3cm}
\end{minipage}
  $$

\caption{An example of a timeline for the KLR spot-fork relation,  with $\ell=1$, $h=3$, $e=5$ and $\al=\varepsilon_3-\varepsilon_1$.  
From left to right we picture  the paths ${\sf FS}_{0,\al}=
\SSTP_{\al} ^\flat \otimes   \SSTP_{\al}  $, 
${\sf FS}_{1,\al}$, 
${\sf FS}_{2,\al}$, 
${\sf FS}_{3,\al}=\SSTP_{\emp  }\otimes \SSTP_\al $.      
 }
\label{pathbendingfork2}
\end{figure}
 \vspace{-0.2cm}

   The following example illustrates all of the important  ideas  in  the proof of this relation (in particular, it  illustrates our  iterative approach using the fork-spot paths,   examples of which are depicted in \cref{pathbendingfork2}).  These ideas will be used repeatedly when we consider  (more complicated) relations in the remainder of this section. 
   
 
\begin{eg}   \label{alhigsflhiufgkhjldfgshkjdfgjdghkf}
  \color{black} 
We set  $\sigma=(0,2,4)$ and $e=6$.  We will consider the following product
$$  
$$where we have emphasised the factorisation of spot  and fork by recording 
the steps within these paths at top and bottom and 
the corresponding labelled $\Upsilon^\SSTP_\SSTQ$ elements for each layer of the righthand-side.
We have also recorded the residues of paths (at the very top and   bottom: $0,2,4,\dots$).

Notice that the path at the bottom of the spot-strand KLR-diagram is not the same as the path at the top of the fork KLR-diagram -- however, the residue sequences are identical (simply trace through the residues on strands).  
We start at the middle of the product   --- that is we first compute 
$$\Upsilon^{\SSTS_{1,\al}\otimes \SSTP_\al     }_{\SSTP_\al^\flat \otimes \SSTP_\al  }
\circ \Upsilon^{\SSTP_\al  \otimes \SSTP_\al  }_{{\sf F}_{1,\emp\al}   }$$
as follows: we first place the diagrams on top of each other
 recording the paths
 $\SSTS_{1,\al}\otimes \SSTP_\al     $ 
 and ${{\sf F}_{1,\emp\al} \otimes \SSTP_\al   }$ at the top and bottom of the diagram  (notice that the permutation is not  step-preserving) 
 and we  highlight the strands in the product  which
 have crossings of non-zero degree
 $$ \Upsilon^{\SSTS_{1,\al}\otimes \SSTP_\al     }_{\SSTP_\al^\flat \otimes \SSTP_\al  }
\circ \Upsilon^{\SSTP_\al  \otimes \SSTP_\al  }_{{\sf F}_{1,\emp\al}   }
=\begin{minipage}{10cm} 
\end{minipage} $$
We apply relation \ref{rel1.11}  to obtain two terms: the term in which we undo this highlighted  braid  
 and the other term which is equal to zero by \cref{resconsider}. 
We  relabel  the bottom of the (non-zero) diagram by the folded fork path, ${\sf FS}_{1,\emp \al}$, and hence obtain
$$\Upsilon^{\SSTS_{1,\al}\otimes \SSTP_\al     }_{\SSTP_\al^\flat \otimes \SSTP_\al  }
\circ \Upsilon^{\SSTP_\al  \otimes \SSTP_\al  }_{{\sf F}_{1,\emp\al}  }
=\begin{minipage}{10cm}  %
\end{minipage}$$
which we now observe is a step-preserving KLR diagram.
We trivially undo the double-crossings in the above diagram (using \cref{adjust1}) and hence obtain   $$
\Upsilon^{\SSTS_{1,\al}\otimes \SSTP_\al     }_{\SSTP_\al^\flat \otimes \SSTP_\al  }
\circ \Upsilon^{\SSTP_\al  \otimes \SSTP_\al  }_{{\sf F}_{1,\emp\al}   }=
\Upsilon^{\SSTS_{1,\al}\otimes \SSTP_\al  } _  {{\sf SF}_{1,\emp\al}     }.$$    
We now  insert this back into the larger product (see also  \cref{forkspot}) and hence obtain the following (not-step-preserving) KLR diagram of 
$$({\sf spot}_\al^\emp (1) \otimes {\sf e}_{\SSTP_\al})\circ
\Upsilon^{\SSTS_{1,\al}\otimes \SSTP_\al  } _  {{\sf SF}_{1,\emp\al}     }\circ
{\sf fork}_{\emp\al}^{\al\al}(1)
$$ which is equal to $$\begin{minipage}{10cm}  
\end{minipage}$$
  where we have highlighted the wiggly strands from the previous step (to facilitate comparison) and we have 
  emboldened the unique pair of crossing strands of the same residue.  
  The rightmost wiggly strand and the pair of bold strands 
are the only strands  have crossings of non-zero degree.  We  apply the same argument as above to undo this braid (we do not need to relabel the bottom of the diagram in this case, as the 
final fork-spot path 
 is equal to $\SSTP_\emp \otimes \SSTP_\al$) and we hence obtain 
$$\begin{minipage}{9cm}  
\end{minipage}$$
which we now observe is a step-preserving KLR diagram.
 We trivially undo the double-crossings
    (using \cref{adjust1}) and hence obtain  
  $$({\sf spot}_\al^\emp (1) \otimes {\sf e}_{\SSTP_\al})\circ
\Upsilon^{\SSTS_{1,\al}\otimes \SSTP_\al  } _  {{\sf SF}_{1,\emp\al}     }\circ
{\sf fork}_{\emp\al}^{\al\al}(1)
={\sf e}_{\SSTP_\emp\otimes \SSTP_\al}
  $$
as required.  
  
\end{eg} 
 \color{black} 

What the above example illustrates  is that 
we start at {\em the middle} of the product on the lefthand-side 
 which is labelled by two {\em distinct paths} which have the same residue sequence, that is we start at the middle term in the product
$$
 \left(   {\sf spot}^{{\emp}}_{{\al}}     \otimes      e_{\SSTP_{\al}}  \right)
(e_{\SSTP_\al^\flat \otimes \SSTP_\al}\circ 
  e_{\SSTP_\al  \otimes \SSTP_\al ^\flat})
\left({\sf fork}_{\emp \al}^{\al \al}  \right)
  $$
where we note that  $e_{\SSTP_\al^\flat \otimes \SSTP_\al} =  e_{\SSTP_\al  \otimes \SSTP_\al ^\flat}$.  
Each iterative stage (of which there are two in \cref{alhigsflhiufgkhjldfgshkjdfgjdghkf})
 simply transforms a non-step-preserving 
 KLR-permutation into a step-preserving one (by undoing all non-zero-degree crossings and relabelling). Thus the (seemingly technical) spot-fork paths become incredibly natural, as does their ``timeline" construction (each stage corresponds to one KLR braid which we undo).  
Most beautifully of all: one should emphasise that the spot-fork path is simply the reflection of the fork path through the $\al$-hyperplane (what else?!). 
This brings us to the general case:
 \color{black}

 \begin{prop}\label{keytoforkdot}
For $\al\in \Pi $ and $0\leq q <\exx$ we have that 
\begin{equation}\label{forkspot} ({\sf spot}_  {\al}^\emp (q)  \otimes    e_{\SSTP_{\al}})
 \circ 	
\Upsilon^{  {\sf S}_{{q,\al}}\otimes \SSTP_\al }_{{\sf FS}_{{q,\al}} }
 	\circ
 {\sf fork}^{\al \al}_{\emp \al}(q) = 
\Upsilon^{  {\sf S}_{{q+1,\al}}\otimes \SSTP_\al  }_{{\sf FS}_{{q+1,\al}} } . 
 \end{equation}
\end{prop}

\begin{proof}We first note that the righthand-side is residue commutative (one can reindex the proof of \cref{forkkkk}).  
We decorate the top and bottom edges of the concatenated product on the lefthand-side of \cref{forkspot} with the tableaux $\northT_q={\sf S}_{{q,\al}} \otimes \SSTP_\al$ and 
$\southT_q={\sf FS}_{q,\al} $ respectively  for $0\leq q <b_\al$. 
For each $0\leq q< b_\al$, the product on the lefthand-side of \cref{forkspot} 
has a single pair of 
strands whose crossing if of degree $-2$:  Namely,  
the strand   $Q_1$ from 
connecting the
$ 
\southT^{-1}_q(q+1,\eps_i)\text{th}$   \south       node to the 
$\northT^{-1}_q(b_\al+q+1,\eps_{i+1})\text{th}$
    \north    node
   and the strand $Q_2$ connecting the 
   $ 
\southT^{-1}_q(b_\al+q+1,\eps_{i+1})\text{th}$   \south       node to the 
$\northT^{-1}_q(q+1,\eps_{i})\text{th}$
    \north    node.  The strands $Q_1$ and $Q_2$ are both of the same residue, $r_q\in\ZZ/e\ZZ$ say, and they cross each other exactly once.  
  This crossing of $r_q$-strands is bi-passed on the left by the  $(r_q+1)$-strand connecting the 
$  \southT^{-1}_q(b_\al+q,\eps_i)\text{th}$   \south       node
to the 
$  \northT^{-1}_q(b_\al+q,\eps_i)\text{th}$   \north    node.  
We pull the $(r_q-1)$-strand through this crossing, using relation \ref{rel1.11}.  We hence obtain two terms: the term in which we undo this braid is equal to the righthand-side of \cref{forkspot} 
 and the other term is equal to zero by \cref{resconsider}. 
     \end{proof}
  
 \Cref{fork-dot-proof} holds by   iteratively applying \cref{keytoforkdot} a total of $b_\al$ times, as in  \cref{alhigsflhiufgkhjldfgshkjdfgjdghkf}. 

       \subsection{The spot and commutator }
   \renewcommand{\bet}{{{{\color{ao(english)}\boldsymbol\gamma}}}}
 \renewcommand{\gam}{{{\color{cyan}\boldsymbol\beta}}}
\renewcommand{\empb}{{{{\color{ao(english)}\boldsymbol\clock}}}}
\renewcommand{\empg}{{{{\color{cyan}\boldsymbol\clock}}}}

Let $\gam,\bet\in\Pi$ label  two commuting   reflections, 
  we  now verify  the leftmost relation in \ref{rel6}, namely that 
\begin{equation}\label{theresultfollows1111}
   {\sf com}_{\gam \bet}^{\bet \gam}
({\sf spot}_{\empg}^{\gam} \otimes e_{\SSTP_{\bet}}  )
=
\Upsilon^{\SSTP_\bet \otimes \SSTP_\gam^\flat  }
_{\SSTP_\empg \otimes  \SSTP_\bet}
=
(e_{\SSTP_{\bet}} \otimes   {\sf spot}_{\empg}^{\gam})
{\sf adj}^{\bet \empg}_{\empg \bet}   
\end{equation}   where the righthand equality is immediate.
   We now set about proving the lefthand-equality. We assume     that $b_\gam \leq b_\bet$ (the other case is similar, but has fewer steps).  We define 
\begin{align*} 
{\sf SC} _{ q,  \gam\bet}	&=
\begin{cases}
 \SSTP_\empg  
\boxtimes 
  \REMOVETHESE {\jay } { q- b_\gam  }  
  \boxtimes  \REMOVETHESE {\jay ,\kay +1} {b_\gam   }  
  \boxtimes 
  \REMOVETHESE {\jay } {b_\bet-q  }  
  \boxtimes \ADDTHIS {\kay } {  b_\gam  } \boxtimes 
 \ADDTHIS {\jay +1} {b_\bet }
  	&\text{for }  b_\bet\geq   q >b_\gam     \\[4pt]
%
 \SSTP_{q\emptyset}  \boxtimes 	  \REMOVETHESE {\kay } {b_\gam -q}  
\boxtimes 
 \REMOVETHESE {\jay ,\kay +1} {q  }  
  \boxtimes 
  \REMOVETHESE {\jay } {b_\bet-q  }  
  \boxtimes \ADDTHIS {\kay } {  b_\gam  } \boxtimes 
 \ADDTHIS {\jay +1} {b_\bet }
  	&\text{for } b_\gam \geq q >0
 			\\[4pt]
 	  \REMOVETHESE {\kay } {b_\gam }  
\boxtimes 
 \REMOVETHESE {\jay } {b_\bet  }   \boxtimes \ADDTHIS {\kay } {  b_\gam  } \boxtimes 
 \ADDTHIS {\jay +1} {b_\bet }
  	&\text{for } q=0 
		\\[4pt]
		  \REMOVETHESE {\kay } {b_\gam }   \boxtimes 
  \ADDTHIS {\kay +1} {  b_\gam  }
\boxtimes 
 \REMOVETHESE {\jay } {b_\bet } 
  \boxtimes 
 \ADDTHIS {\jay +1 } {b_\bet } 
  	&\text{for }q=-1	\\[4pt]
 \end{cases}
\end{align*}
which is obtained from ${\sf  C} _{ q,  \gam\bet}$   by reflection through $s_\gam$. We invite the reader to draw an example  of  the   timeline by reflecting the final four paths of \cref{chopin}  through $s_\gam$. 


 \begin{prop} For $0\leq q  < b_\gam$, we have   that  
\begin{equation}\label{SC1}
 {\sf com}_{\gam\bet }(q )  \circ  
\Upsilon _{\SSTS_{{ q },\gam} \otimes\SSTP_\bet}
^{{\sf SC}_{q,\gam\bet }  }
\circ ( {\sf spot }^{\gam  }_{\empg}(q )\otimes e_{\SSTP_\bet} )
=
\Upsilon _{\SSTS_{{ q+1},\gam} \otimes\SSTP_\bet}
^{{\sf SC}_{q+1,\gam\bet }}
\end{equation} 
 (note that $\Upsilon _{\SSTS_{{ 0 },\gam} \otimes\SSTP_\bet}
^{{\sf SC}_{0,\gam\bet }  }
={\sf com}_{\gam\bet}(-1)$)  and  for $b_\gam \leq q < b_\bet$, we have that 
\begin{equation}\label{SC2}
 {\sf com}_{\gam\bet }(q )  \circ 
\Upsilon _{\SSTP_\empg \otimes\SSTP_\bet}
^{{\sf SC}_{q,\gam\bet }  }
=  
\Upsilon_{\SSTP_\empg \otimes\SSTP_\bet}
^{{\sf SC}_{{ q+1},\gam\bet } } .
\end{equation} 
 \end{prop}

\begin{proof}
All these elements are residue commutative (by reindexing the proof of \cref{resssssss}).  
  We prove \cref{SC1,SC2} by    induction on $0\leq  q<b_\bet$ (the $q=-1$ case is trivial).  Label the   \north    and   \south       frames of the concatenated diagrams on the lefthand-side  of   \cref{SC1,SC2} by the paths $\northT_{q+1}={\sf SC}_{{q+1},\gam\bet}$ and $\southT_{q+1}={\sf S} _{{q+1},\empg}\otimes \SSTP_\bet$.  
The concatenated diagram on the lefthand-side of both  \cref{SC1} and \cref{SC2}
  has a single crossing which does not   preserve  step labels.  
Namely the strands connecting the $\northT_q^{-1} (q+1,\eps_\kay)$th  and $\northT_q^{-1} (b_\gam+q+1,\eps_{\kay+1})$th 
\north vertices to the 
$\southT^{-1}_q(q+1,\eps_\kay)$th and 
  $\southT_q^{-1} (b_\gam+q+1,\eps_{\kay+1})$th 
  \south vertices form an $r_q$-crossing, for some $r_q\in \ZZ/e\ZZ$ say, and  these strands permute  the labels 
$+\eps_\kay $ and $+\eps_{\kay +1}$.     
This crossing is bi-passed on the left by a strand connecting the 
  $\northT_q^{-1} (b_\gam+q,\eps_{\kay+1})$th 
\north and 
  $\southT_q^{-1} (b_\gam+q,\eps_{\kay+1})$th 
 \south vertices.   
We undo this triple using case 2 of  relation \ref{rel1.11} and hence obtain the righthand-side of \cref{SC1,SC2}.  
\end{proof}

In order to deduce that  \cref{theresultfollows1111} holds, we   observe that 
$$
 {\sf com} ^{\bet \gam} \circ 
({\sf com}_{ {{\gam\bet}} }	(	{\sf spot}^{\gam}_\empg \otimes e_{\SSTP_\bet})) 
=
   {\sf com} ^{\bet \gam} \circ 
\Upsilon_{\SSTP_\empg \otimes\SSTP_\bet}
^{{\sf SC}_{b_\bet,\gam\bet } }
=\Upsilon^{\SSTP_\bet \otimes  \SSTP_\gam^\flat }
_{\SSTP_\empg \otimes  \SSTP_\bet}
$$
as the lefthand-side of the final equality is minimal and respects step-labels.

   \renewcommand{\gam}{{{{\color{ao(english)}\boldsymbol\gamma}}}}
 \renewcommand{\bet}{{{\color{cyan}\boldsymbol\beta}}}
\renewcommand{\empb}{{{{\color{cyan}\boldsymbol\clock}}}}
\renewcommand{\empg}{{{{\color{ao(english)}\boldsymbol\clock}}}}

      \subsection{The spot-hexagon}
For $\al,\bet\in \Pi$ labelling two non-commuting reflections,  we now check the rightmost relation in  \ref{rel3}, namely that 
\begin{equation}\label{labelmebicthc}
\Psi \left(\; \begin{minipage}{2.3cm}
\end{minipage}\right) 
\end{equation}
(and we leave it the reader to check the  reflection of this relation through its vertical axis).  
In other words, we need to check that 
\begin{align*}
( e_{\SSTP_{\emp }}
\otimes 
{\sf spot}_{\bet}^\empb\otimes e_{\SSTP_{\al \bet}}){\sf hex}^{\emp\bet\al\bet}_{\empb\al\bet\al}
\end{align*}
is equal to 
\begin{align*}
{\sf adj}_{\empb\al \bet \emp}^{\emp\empb \al \bet}
(e_{\SSTP_{\empb\al \bet}}\otimes {\sf spot}_{\al}^\emp )
+  e_{\SSTP_\empb}
\otimes (
 {\sf fork}_{\al\al}^{\emp\al}\otimes 
 {\sf spot}^{\bet}_\empb)
  {\sf adj}_{\al\empb\al}^ {\al\al\empb} 
 (  e_{\SSTP_{\al}}\otimes  {\sf spot}_{\bet}^\empb\otimes  e_{\SSTP_{\al}})).   
\end{align*}
  We set $j=i+1$ so that  $\al=\eps_i-\eps_{i+1}$, $\bet=\eps_{i+1}-\eps_{i+2}$.    We will begin by considering the lefthand-side of the equation.  
In order to do this, we need to use the reflections of the braid ${\sf H}_{q,\bet\al\bet} $-paths for $0\leq q \leq b_{\al\bet}$  through the first $\bet$-hyperplane which they come across (namely the hyperplane whose strand we are putting a spot on top of) and we remark that this path will have the {\em same residue sequence}
as the original ${\sf H}_{q,\bet\al\bet} $-paths,  
 but different step labelling.  
We define ${\sf SH} _{q, \bet\al\bet}$ to be the path 
$$
\begin{cases}
  \SSTP_{q\emptyset}\boxtimes 
 \REMOVETHESE{i+1}{b_\bet-q} \boxtimes 
 \REMOVETHESE {i,i+1} {q}  
\boxtimes
\ADDTHIS{i+1}{b_\bet}
   \boxtimes \REMOVETHESE{i} {b_\al-q} 
 \boxtimes \ADDTHIS{i+1} {b_\al}
{{\; \color{magenta}\otimes}}_\al\;  
\REMOVETHESE{i+1}{b_\bet} \boxtimes 
\ADDTHIS{i+2}{b_\bet}&0\leq q \leq b_\bet \\[6pt]
\SSTP_{\empb }\boxtimes 
   \REMOVETHESE{i} {q-b_\bet} \boxtimes \REMOVETHESE {i,i+1}{b_\bet} 
     \boxtimes  \ADDTHIS{i+1}{b_\bet}
\boxtimes  
 \REMOVETHESE{i} {b_\al-q} 
\boxtimes \ADDTHIS{i+1} {b_\al}
{{\; \color{magenta}\otimes}}_\al\;  
\REMOVETHESE{i+1}{b_\bet} \boxtimes 
\ADDTHIS{i+2}{b_\bet}&   b_\bet \leq q \leq b_\al  \\[5pt]  
\SSTP_{\empb }  \boxtimes  
\REMOVETHESE {i } {  q-b_\bet }   \boxtimes  \ADDTHIS{i+2}{b_\bet}
\boxtimes 
\REMOVETHESE{i ,i+2}{q-b_\al} \boxtimes   \ADDTHIS{i+1} {b_\al}
{{\; \color{magenta}\otimes}}_\al\;  
\REMOVETHESE{i+1}{ b_{\al\bet}-q} \boxtimes 
\ADDTHIS{i+2}{b_\bet}&b_\al     \leq q  \leq b_{\al\bet}\end{cases} 
$$  
for $b_\al \geq b_\bet$  (the $b_\al< b_\bet$ case is similar).  
See \cref{yyyyyyy3} for an example.  
\begin{figure}[ht!]
$$  \begin{minipage}{2.8cm}
\end{minipage}
$$

\caption{
An example of the tableaux ${\sf SH} _{q,\bet\al\bet}$ 
for $ 0\leq   q \leq b_{\al\bet} $. 
The reader should compare these reflected paths with the final five paths of 
 \cref{Steinberg}.   }

 \label{yyyyyyy3}
\end{figure}

  \begin{prop}\label{aprooftorefer}
 We have that  
\begin{align}   \label{notlessthanb2}
  \big(e_{\SSTP_\emp}\otimes  {\sf spot}_{\bet}^\empb 
   \otimes e_{\SSTP_{\al \bet}}\big)\;
     {\sf hex}^{\emp \bet\al\bet} 
 = \Upsilon^{\SSTP_{\emp  \empb \al\bet}}  
  _{\SSTP_\emp\otimes {\sf SH}_{b_{\al\bet}, \bet\al\bet}}  
\end{align}
  \end{prop}

\begin{proof} 
First, we remark that the righthand-side of 
\cref{notlessthanb2} is residue-commuting and so makes sense.  
For  $0\leq q<b_{\al\bet}$, we claim that 
\begin{align}\label{lessthanb}
  \big(e_{\SSTP_\emp}\otimes  {\sf spot}_{\bet}^\empb(q )
   \otimes e_{\SSTP_{\al \bet}}\big)\;
 \Upsilon^{\SSTP_{\emp} \otimes \SSTS_{q ,  \bet} \otimes\SSTP_{ \al\bet}}_{\SSTP_\emp\otimes {\sf SH}_{q,\bet\al\bet}}  
    {\sf hex}^{\emp \bet\al\bet}(q )  
 &= \Upsilon^ {\SSTP_{\emp} \otimes \SSTS_{q+1,  \bet} \otimes\SSTP_{ \al\bet}} _{\SSTP_\emp\otimes {\sf SH}_{q+1,\bet\al\bet}}  
\end{align}
and we will we label the top and bottom of these diagrams 
 according to the paths 
 $ \SSTT_q={\SSTP_{\emp} \otimes \SSTS_{q+1,  \bet} \otimes\SSTP_{ \al\bet}}$  and  $\southT_q=
 \SSTP_{\emp} \otimes{\sf SH}_{q+1,\bet\al\bet}  $ respectively 
(with the convention that ${\sf S}_{q,\bet}  =\SSTP_\empb$ for $q\geq b_\bet$).  Again, this element is residue-commuting and so there is no ambiguity here.     
In   the concatenated diagram  on the lefthand-side  of 
  \cref{lessthanb},  
there   is a single pair of strands,  
 $Q$ and $Q'$
 whose crossing if of degree $-2$   (of residue   $r_q\in \ZZ/e\ZZ$, say); 
  these strands 
   connect  the 
$$
\northT^{-1}_q(b_\al+q+1,\eps_{i+1}) \quad  \northT^{-1}_q(b_{\al\bet} +q+ 1,\eps_{i+2})
$$ 
top vertices and the 
$$  \southT^{-1}_q(b_\al+q+1,\eps_{i+1}) \quad 
   \southT^{-1}_q(b_{\al\bet}+q+1,\eps_{i+2})
$$
 bottom vertices (thus  crossing one another).  
This  crossing of   $r_q$-strands,  $Q$ and $Q'$,  is bi-passed on the left by the $(r_q+1)$-strand from 
   $\northT^{-1}_q(b_{\al\bet}+q,\eps_{i+2})$ 
to 
   $\southT^{-1}_q(b_{\al\bet}+q,\eps_{i+2})$.

  Applying  case 2 of relation \ref{rel1.11} to the 
  concatenated diagram  we obtain two terms: 
  the term with the crossing is bi-passed on   the right is zero by \cref{resconsider};
   the term in which we undo the crossing    is equal to
the righthand-side of \cref{lessthanb} (since  the resulting diagram is minimal).  
An example is given in \cref{alabelforafigre23}.  
    \end{proof}

  \begin{figure}[ht!]
 $$  
 \begin{minipage}{7.2cm}
 \end{minipage}
    $$  
    
 \!\!
\caption{The product  $\big(  {\sf spot}_{\bet}^\empb(0 )
   \otimes e_{\SSTP_{\al \bet}}\big)      {\sf hex}^{ \bet\al\bet}(0 )  
$   in  the proof of \cref{aprooftorefer} 
for  $h=1$, $\ell=5$, $\kappa=(0,2,4,6,8)$,  $e=10$ and $\al=\eps_2-\eps_3$, $\bet=\eps_3-\eps_4$. 
The   \north    path is ${\sf S}_{1,\al} \otimes {\sf M}_{2}$ and the   \south       path is 
${\sf SH}_{1,\bet\al\bet}  \otimes {\sf M}_{2}$ (the prefix $\SSTP_\emp$ and the remainder of  the postfix $\SSTP_{\al}
 ={\sf M}_2^{b_\al}\boxtimes {\sf P}_3^{b_\al}$   would not fit).   }
\label{alabelforafigre23}

\end{figure}

We now wish to show that  
$$
  \Upsilon^ {\SSTP_{\emp  \bet \al\bet}}_{\SSTP_\emp\otimes {\sf SH}_{b_{\al\bet},\bet\al\bet}}  
 {\sf hex}_{\empb\al\bet\al}$$
 is equal to 
 $${\sf adj}_{\empb\al \bet \emp}^{\emp\empb \al \bet}
(e_{\SSTP_{\empb\al \bet}}\otimes {\sf spot}_{\al}^\emp )
+  e_{\SSTP_\empb}
\otimes (
 {\sf fork}_{\al\al}^{\emp\al}\otimes 
 {\sf spot}^{\bet}_\empb)
  {\sf adj}_{\al\empb\al}^ {\al\al\empb} 
 (  e_{\SSTP_{\al}}\otimes  {\sf spot}_{\bet}^\empb\otimes  e_{\SSTP_{\al}})).   
$$ 
In what  follows, we assume that $b_\al\geq b_\bet$.    In order to consider the first term, we   use the reflections of the  ${\sf H}_{q,\al\bet\al} $-paths for $0\leq q \leq b_{\al\bet}$  through the final  $\al$-hyperplane which they come across (namely the hyperplane whose strand we are putting a spot on top of) and we remark that this path will have the {\em same residue sequence} as the original ${\sf H}_{q,\al\bet\al} $-paths   but with a different step labelling.  
We define $ {\sf  S_\al H}_{q,\al\bet\al}$ to be the path 
 \begin{align*}\begin{cases}
\SSTP_{q\emptyset}\boxtimes   \REMOVETHESE{i}{\exx}
  \boxtimes \ADDTHIS{i+1}{\exx } 
  {{\; \color{magenta}\otimes}}_\al\;  \REMOVETHESE {i+1} {b_\bet-q}  
\boxtimes \ADDTHIS{i+2} {b_\bet} 
{{\; \color{cyan}\otimes}}_\bet\;  
 \Pdiptwo^q
\boxtimes \REMOVETHESE { i  } {b_\al-q}
\;  {\color{magenta}\otimes} _\al  \; \ADDTHIS {i+1} {\exx } &0\leq q \leq b_\bet \\[4pt]
\SSTP_{q\emptyset}\boxtimes   \REMOVETHESE{i}{\exx}
  \boxtimes \ADDTHIS{i+1}{ b_{\al\bet}-q} 
{{\; \color{magenta}\otimes}}_\al\;
  \ADDTHIS{ i+2} {b_\bet} 
{{\; \color{cyan}\otimes}}_\bet\;  
 \Pdiptwo^{b_\bet}
\boxtimes \REMOVETHESE { i  } {b_\al-q}
\;  {\color{magenta}\otimes} _\al \;  \ADDTHIS {i+1} {\exx }&b_\bet \leq q \leq b_\al \\[4pt]
\SSTP_{\emp}\boxtimes
   \REMOVETHESE{i}{\exx}
  \boxtimes \ADDTHIS{i+1}{ b_{\al\bet}-q} 
{{\; \color{magenta}\otimes}}_\al\;
  \ADDTHIS{ i+2} {b_\bet} 
{{\; \color{cyan}\otimes}}_\bet\;  
 \Pdiptwo^{b_\bet} 
 \; {\color{magenta}\otimes} _\al \;   \ADDTHIS {i+1} {\exx }
 \boxtimes \ADDTHIS {i} {q-\exx }
&b_\al \leq q \leq  b_{\al\bet}\\[4pt]
\end{cases}
\end{align*}  
   In order to consider the second term, we   need 
 the reflections of the   ${\sf H}_{q,\al\bet\al} $-paths for $0\leq q \leq b_{\al\bet}$  through the first $\bet$-hyperplane   which they come across.   We define   $ {\sf  S_\bet H}_{q,\al\bet\al}$ to be the path
  \begin{align*}
 \begin{cases}
\SSTP_{q\emptyset}\boxtimes   \REMOVETHESE{i}{\exx}
  \boxtimes \ADDTHIS{i+1}{\exx } 
  {{\; \color{magenta}\otimes}}_\al\;  \REMOVETHESE {i+1} {b_\bet-q}  
\boxtimes \ADDTHIS{i+2} {b_\bet} \boxtimes 
 \Pdiptwo^q
\boxtimes \REMOVETHESE { i  } {b_\al-q}
\boxtimes \ADDTHIS {i+1} {\exx } &0\leq q \leq b_\bet \\[4pt]
\SSTP_{q\emptyset}\boxtimes   \REMOVETHESE{i}{\exx}
  \boxtimes \ADDTHIS{i+1}{ b_{\al\bet}-q} 
{{\; \color{magenta}\otimes}}_\al\;
  \ADDTHIS{ i+2} {b_\bet} \boxtimes 
 \Pdiptwo^{b_\bet}
\boxtimes \REMOVETHESE { i  } {b_\al-q}
\boxtimes  \ADDTHIS {i+1} {\exx }&b_\bet \leq q \leq b_\al \\[4pt]
\SSTP_{\emp}\boxtimes
   \REMOVETHESE{i}{\exx}
  \boxtimes \ADDTHIS{i+1}{ b_{\al\bet}-q} 
{{\; \color{magenta}\otimes}}_\al\;
  \ADDTHIS{ i+2 } {b_\bet} 
\boxtimes  \Pdiptwo^{b_\bet} 
\boxtimes    \ADDTHIS {i+1} {\exx }
 \boxtimes \ADDTHIS {i} {q-\exx }
&b_\al \leq q \leq  b_{\al\bet}\\[4pt]
\end{cases}
\end{align*}
See \cref{yyyyyyy4} for an example of the  ${\sf S_{\al}H}_{q,\al\bet\al}$  paths.  We leave it as an exercise for the reader to draw the ${\sf S_{\bet}H}_{q,\al\bet\al}$ paths.
Finally, for the purposes of the proof we will also need the following ``error path"  
$$
{\sf  e  S_\bet H}_{\al\bet\al} 
=\SSTP_{\emp}\boxtimes
   \REMOVETHESE{i}{\exx}
{{\; \color{magenta}\otimes}}_\al\;
  \ADDTHIS{ i+2} {b_\bet-1} 
  \boxtimes  \Pdiptwo 
  \boxtimes   \ADDTHIS{ i+2} { \ } 
\boxtimes  \Pdiptwo^{b_\bet-1} 
\boxtimes    \ADDTHIS {i+1} {\exx }
 \boxtimes \ADDTHIS {i} {b_\bet }
$$
which  one should compare with the final path (the $b_{\al\bet}$th case) above.  
 One should repeat   the above definitions for the $b_\al<b_\bet$ case.

\begin{figure}[ht!]
$$  \begin{minipage}{2.8cm}
\end{minipage}
$$

\caption{
An example  of the paths ${\sf S_{\al}H}_{q,\al\bet\al}$ 
for $b_{\al\bet}\geq q \geq 0 $.  } \label{yyyyyyy4}

\end{figure}

    \begin{prop}\label{weshowthat}
  We have that 
\begin{align}\label{thepointo2}
\Upsilon^ {\SSTP_{\emp  \empb \al\bet}}_{\SSTP_\emp\otimes {\sf SH}_{b_{\al\bet}, \bet\al\bet}}  
{\sf hex} _{\empb\al\bet\al}
=
  \Upsilon^ {\SSTP _{\emp } \otimes  \SSTP_ \empb 
 \otimes  \SSTP_\al\otimes  \SSTP_\bet} 
 _{\SSTP_{\empb}\otimes  \SSTP_{\al}\otimes  \SSTP_\bet \otimes  \SSTP_\al^\flat	}+
  \Upsilon^ {\SSTP_{\emp}\otimes \SSTP_{  \empb}\otimes\SSTP_{ \al}^\flat \otimes \SSTP_{\bet}^\flat}
 _{\SSTP_\empb\otimes \SSTP_\al  \otimes \SSTP_\bet^\flat \otimes \SSTP_\al 	}
.
\end{align}
     \end{prop}
 \begin{proof}
 First, we remark that both terms on the righthand-side of 
\cref{thepointo2} are residue-commuting.  We suppose $b_\al \geq b_\bet$ as the other case is similar.  
We observe that  
$$\Upsilon^ {\SSTP_{\emp  \empb \al\bet}}_{\SSTP_\emp\otimes {\sf SH}_{b_{\al\bet}, \bet\al\bet}}  
=
 \Upsilon^ {\SSTP _{\empb} \otimes  \SSTP_ \empb 
 \otimes  \SSTP_\al\otimes  \SSTP_\bet} 
 _{\SSTP_{\empb}\otimes   {\sf S_\al H}_{b_{\al\bet},\bet\al\bet}  } 
 = \Upsilon^ {\SSTP_{\emp}\otimes \SSTP_{  \empb}\otimes\SSTP_{ \al}^\flat \otimes \SSTP_{\bet}^\flat}
_{\SSTP_{\empb}\otimes     {\sf S_\bet H}_{b_{\al\bet},\bet\al\bet}  }
 $$
as the underlying permutations (and residue sequences) are all identical.    We set 
$$\northT  _{ \al}= {\SSTP _{\emp  \empb }
 \otimes  \SSTP_\al\otimes  \SSTP_\bet} 
 \qquad 
\northT_\bet ={\SSTP _{\emp  \empb } 
 \otimes  \SSTP_\al^\flat \otimes  \SSTP_\bet^\flat} 
\qquad 
$$
$$
\southT  _{ q,\al}={\SSTP_{\empb}\otimes   {\sf S_\al H}_{q+1, \al\bet\al}  } 
\qquad
\southT  _{ q,\bet}=  {\SSTP_{\empb}\otimes     {\sf S_\bet H}_{q , \al\bet\al}  }
$$
for $b_{\al\bet}> q \geq 0$.  
We first consider the $q=b_{\al\bet}-1$ case.  
 The concatenated diagram $$\Upsilon^ {\SSTP_{\emp  \empb \al\bet}}_{\SSTP_\emp\otimes {\sf SH}_{b_{\al\bet}, \bet\al\bet}}  
  (e_{\SSTP_{  \empb}} \otimes {\sf hex} _{ \al\bet\al} (b_{\al\bet}-1 ))
$$
 contains a single   like-labelled crossing of $r_{b_{\al\bet}-1}$-strands  
connecting  the pair 
$$
\SSTT_\al^{-1}(b_{\al\bet\al}+1	,\eps_{i+1})
 =
\SSTT_\bet^{-1}( b_{\al\bet}+1	,\eps_{i})
 \qquad \SSTT_\al^{-1}(  2b_{\al\bet}+1	 ,\eps_{i+2})
 =
\SSTT_\bet^{-1}( b_{\al\bet\al}+1	,\eps_{i+1 })
$$
of top vertices to the pair of 
  $$
\southT_\al^{-1}( 2b_{\al\bet}+1	,\eps_{i+2})
=
 \southT_\bet^{-1}(b_{\al\bet\al}+1	,\eps_{i+1 })
 \qquad  \southT_\al^{-1}( b_{\al\bet\al}+1	,\eps_{i+1})
 =
\southT_\bet^{-1}( b_{\al\bet}+1	,\eps_{i})
$$
 These  $r_{b_{\al\bet}-1}$-crossing strands are bi-passed on the left by the 
  $r_{b_{\al\bet}}$-strand  connecting the 
$$\SSTT_\al^{-1}(2b_{\al\bet}	,\eps_{i+2})=\SSTT_\bet^{-1}(2b_{\al\bet},\eps_{i+2})\qquad \southT_\al^{-1}(2b_{\al\bet} 	,\eps_{i+2})
=
\southT_\bet^{-1}(2b_{\al\bet} 	,\eps_{i+2})
$$
top and bottom vertices.   
  We apply   case 2 of relation \ref{rel1.11} to   the 
this triple of  strands      and hence obtain
\begin{equation}\label{1both}\Upsilon^ {\SSTP_{\emp  \empb \al\bet}}_{\SSTP_\emp\otimes {\sf SH}_{b_{\al\bet}, \bet\al\bet}}  
 {\sf hex} _{  \empb\al\bet\al} (b_{\al\bet}-1 )
= \Upsilon ^ {\SSTP _{\emp } \otimes  \SSTP_ \empb 
 \otimes  \SSTP_\al\otimes  \SSTP_\bet} 
_{\SSTP_{\empb}\otimes     {\sf S_\al H}_{b_{\al\bet}-1,\bet\al\bet}  }
+
 \Upsilon^ {\SSTP_{\emp}\otimes \SSTP_{  \empb}\otimes\SSTP_{ \al}^\flat \otimes \SSTP_{\bet}^\flat} 
 _{\SSTP_{\empb}\otimes   		{\sf eS_\bet H}_{\al\bet\al	}	} 
\Upsilon^  {\SSTP_{\empb}\otimes   		{\sf eS_\bet H}_{\al\bet\al		} }
 _{\SSTP_{\empb}\otimes   {\sf S_\bet H}_{b_{\al\bet}-1,\bet\al\bet}  } 
\end{equation}  where in the first term  we have undone the
 triple-crossing and in the second  ``error"   term   the $r_{b_{\al\bet}}$-strand bi-passes the crossing to the right (and is labelled by the ``error path").  
 We are now ready to consider the $b_{\al\bet}-1>q\geq 0$ cases --- which we do separately for $\al$ and $\bet$, in turn.  
 
 \medskip 
\noindent {\bf Case $\al$. }  We first consider the first  term on the righthand-side of \cref{1both}.  We claim that 
 \begin{align}\label{both}
  \Upsilon^ {\SSTP _{\emp } \otimes  \SSTP_ \empb 
 \otimes  \SSTP_\al\otimes  \SSTP_\bet} 
  _{\SSTP_{\empb}\otimes   {\sf S_\al H}_{q+1, \al\bet\al}  } 
 {\sf hex} _{  \empb\al\bet\al} (q) 
&=
\Upsilon^ {\SSTP _{\emp  \empb }
 \otimes  \SSTP_\al\otimes  \SSTP_\bet} 
 _{\SSTP_{\empb}\otimes   {\sf S_\al H}_{q , \al\bet\al}  } 
  \end{align}
for $b_{\al\bet}-1 > q \geq 0$.   For each  $b_{\al\bet}> q \geq b_\al$ the concatenated diagram  in \cref{both} contains a single   like-labelled crossing of $r_q$-strands (for some  $r_q\in \ZZ/e\ZZ$ say) 
connecting  the pair 
$$
\SSTT_\al^{-1}(2b_{\bet}+3b_{\al}-q	,\eps_{i+1})
 \qquad \SSTT_\al^{-1}(3b_{\bet}+3b_{\al}-q	,\eps_{i+2})
 $$
of top vertices to the pair of 
  $$
\southT_\al^{-1}(3b_{\bet}+3b_{\al}-q	,\eps_{i+2})
 \qquad  \southT_\al^{-1}(3b_{\bet}+3b_{\al}-q	,\eps_{i+1})
 $$
bottom vertices, respectively.  
     For $b_{\al\bet}-1>q \geq b_\al$ the aforementioned (unique) pair of crossing $r_q$-strands in 
$$\Upsilon^ {\SSTP _{\emp  \empb } 
 \otimes  \SSTP_\al\otimes  \SSTP_\bet} 
 _{\SSTP_{\empb}\otimes   {\sf S_\al H}_{q+1, \al\bet\al}  } 
{\sf hex} _{  \empb\al\bet\al} (q) 
=
\Upsilon^ {\SSTP _{\emp  \empb }
 \otimes  \SSTP_\al\otimes  \SSTP_\bet} 
 _{\SSTP_{\empb}\otimes   {\sf S_\al H}_{q , \al\bet\al}  } $$
 is bi-passed on the left by the $r_{q+1}$-strand connecting 
 $ \SSTT_\al^{-1}(3b_{\bet}+3b_{\al}-q-1	,\eps_{i+2})$ and 
  $ \southT_\al^{-1}(3b_{\bet}+3b_{\al}-q-1	,\eps_{i+2})$ top and bottom vertices.  
  Applying case 2 of relation \ref{rel1.11} we undo this triple crossing (the other term is zero by \cref{resconsider}) as required.   Now for $b_\al > q \geq 0$ the concatenated product on the lefthand-side of \cref{both} is both minimal and step-preserving and so the claim follows.  

\medskip
\noindent {\bf Case $\bet$}.  We now consider the second term on the right of \cref{1both}.  We have that 
$$  \Upsilon^  {\SSTP_{\empb}\otimes   		{\sf eS_\bet H}_{\al\bet\al}		} 
 _{\SSTP_{\empb}\otimes   {\sf S_\bet H}_{q+1, \al\bet \al }  } 
{\sf hex} _{  \empb\al\bet\al} (q) = 
  \Upsilon^  {\SSTP_{\empb}\otimes   		{\sf eS_\bet H}_{\al\bet\al}		} 
 _{\SSTP_{\empb}\otimes   {\sf S_\bet H}_{q, \al\bet\al}  } 
$$
for $b_{\al\bet}-1> q \geq b_ \al$ as the lefthand-side is minimal and step-preserving.  
Now, we claim that 
\begin{align}\label{step11}
\Upsilon^ {\SSTP _{\emp } \otimes  \SSTP_ \empb 
 \otimes  \SSTP_\al^\flat \otimes  \SSTP_\bet^\flat } 
 _{\SSTP_{\empb}\otimes   		{\sf eS_\bet H}_{\al\bet\al	}	} 
\Upsilon^  {\SSTP_{\empb}\otimes   		{\sf eS_\bet H}_{\al\bet\al}		} 
 _{\SSTP_{\empb}\otimes   {\sf S_\bet H}_{b_\al, \al\bet\al}  } 
{\sf hex} _{  \empb\al\bet\al} (b_\al-1) 
&=
\Upsilon^ {\SSTP _{\emp } \otimes  \SSTP_ \empb 
 \otimes  \SSTP_\al^\flat \otimes  \SSTP_\bet^\flat } 
  _{\SSTP_{\empb}\otimes   {\sf S_\bet H}_{b_{\al}-1, \al\bet\al}  } \quad
\intertext{and that  } 
\label{step12}
\Upsilon^ {\SSTP _{\emp } \otimes  \SSTP_ \empb 
 \otimes  \SSTP_\al^\flat \otimes  \SSTP_\bet^\flat } 
  _{\SSTP_{\empb}\otimes   {\sf S_\bet H}_{q+1,\bet\al\bet}  } 
{\sf hex} _{  \empb\al\bet\al} (q) 
&=
\Upsilon^ {\SSTP _{\emp } \otimes  \SSTP_ \empb 
 \otimes  \SSTP_\al^\flat \otimes  \SSTP_\bet^\flat } 
  _{\SSTP_{\empb}\otimes   {\sf S_\bet H}_{q, \al\bet\al}  } 
\end{align}
for $b_\al-1>q \geq 0$.  For each $b_\al \geq q \geq 0$ the concatenated diagram on 
the lefthand-side of \cref{step11,step12}
contains  a crossing pair of $r_q$-strands connecting the 
$$\SSTT_\bet^{-1}(b_\bet+q+1,\eps_i) \quad 
 \SSTT_\bet^{-1}(2b_{\bet}+b_{\al}+q+1	,\eps_{i+2}) 
$$
and
$$ \southT_\bet^{-1}(2b_{\bet}+b_{\al}+q+1	,\eps_{i+2}) \quad  
\southT_\bet^{-1}(b_\bet+q+1,\eps_i)  
$$
top and bottom vertices, respectively (note that this crossing does not respect step labels).  
 This $r_q$-crossing is bi-passed on the right by the $(r_{q}-1)$-strand connecting the 
$$
\SSTT_\bet^{-1}(b_\bet +q+2,\eps_{i }) 
\qquad
\southT_\bet^{-1}(b_\bet  +q+2,\eps_{i }) 
$$
top and bottom vertices.  We undo this triple-crossing using  case 1 of relation \ref{rel1.11} (the other term is zero by \cref{resconsider}).  The concatenated product is minimal and step-preserving, as required.  \hfill\qedhere
   \end{proof}

 Finally, in order to deduce  \cref{labelmebicthc},  we observe that 
$$
{\sf adj}_{\empb\al \bet \emp}^{\emp\empb \al \bet}
(e_{\SSTP_{\empb\al \bet}}\otimes {\sf spot}_{\al}^\emp )
=
 \Upsilon^ {\SSTP _{\emp  \empb  \al}   \otimes  \SSTP_\bet} 
 _{\SSTP_{\empb \al \bet } \otimes  \SSTP_\al^\flat	} 
$$
$$ e_{\SSTP_\empb}
\otimes ((
 {\sf fork}_{\al\al}^{\emp\al}\otimes 
 {\sf spot}^{\bet}_\empb)
  {\sf adj}_{\al\empb\al}^ {\al\al\empb} 
 (e_{\SSTP_{\al}}\otimes  {\sf spot}_{\bet}^\empb\otimes  e_{\SSTP_{\al}}))   
=\Upsilon^ {\SSTP_{\emp \empb}\otimes\SSTP_{ \al}^\flat \otimes \SSTP_{\bet}^\flat}
 _{\SSTP_{\empb \al}  \otimes \SSTP_\bet^\flat \otimes \SSTP_\al 	}
 $$
as  the concatenated diagrams are minimal, step-preserving, and residue-commutative.

   \subsection{The fork-hexagon}
For $\al,\bet\in \Pi$ labelling two non-commuting reflections,  we now check the leftmost relation in  \ref{rel3}, namely that 
 \begin{equation}\label{aimofthis}
 (e_{\SSTP_{\empb \emp}} \otimes {\sf hex}^{\emp\bet\al\bet}_{\empb\al\bet\al})
(e_{\SSTP_{\empb\empb}}\otimes {\sf fork}_{\al\al}^{\emp\al}\otimes e_{\SSTP_{\bet\al}})
  {\sf adj}_{\empb\al\empb\al\bet\al}^{\empb\empb\al \al\bet\al}
(e_{\SSTP_{\empb\al}}\otimes {\sf hex}_{\emp\bet\al\bet}^{\empb\al\bet\al})
\end{equation} is equal to
 \begin{equation}\label{aimofthis2} {\sf adj}^{\empb\emp\emp  \bet\al\bet}_{\emp\emp \bet\al\empb\bet} 
(e_{ {\SSTP_{ \emp\emp\bet\al}}} \otimes {\sf fork}^{\empb\bet}_{\bet\bet})
 (e_{\SSTP_{\emp}}\otimes {\sf hex}^{\emp \bet\al\bet }_{\empb \al\bet\al }\otimes e_{\SSTP_{ \bet}})
 {\sf adj}^{\emp\empb \al\bet\al   \bet} _{\empb \al \emp\bet\al  \bet}
\end{equation}
Unlike earlier sections, we find that neither    of \ref{aimofthis} or \ref{aimofthis2} is of minimal length.    We again set $j=i+1$.  
 First assume that $b_\al \geq b_\bet$.     For   \ref{aimofthis}, we must simplify the middle of the diagram.   We define ${\sf FH}_{q,\al\bet\al}$ to be the path 
\begin{align*}
 \begin{cases}
\SSTP_{q\emptyset}\boxtimes   \REMOVETHESE{i}{\exx}
  \boxtimes \ADDTHIS{i}{\exx } 
 \boxtimes         \REMOVETHESE {i+1} {b_\bet-q}  
\boxtimes \ADDTHIS{i+2} {b_\bet} 
{{\; \color{cyan}\otimes}}_\bet\;  
 \Pdiptwo^q
\boxtimes \REMOVETHESE { i  } {b_\al-q}
\boxtimes \ADDTHIS {i+1} {\exx }
&0\leq q \leq b_\bet \\[6pt]
\SSTP_{q\emptyset}\boxtimes   \REMOVETHESE{i}{\exx}
  \boxtimes \ADDTHIS{i }{ b_{\al\bet}-q} 
\boxtimes 
  \ADDTHIS{ i+2} {b_\bet} 
{{\; \color{cyan}\otimes}}_\bet\;  
 \Pdiptwo^{b_\bet}
\boxtimes \REMOVETHESE { i  } {b_\al-q}
\boxtimes \ADDTHIS {i+1} {\exx }&b_\bet \leq q \leq b_\al \\[6pt]
\SSTP_{\emp}\boxtimes
   \REMOVETHESE{i}{\exx}
  \boxtimes \ADDTHIS{i}{ b_{\al\bet}-q} 
\boxtimes   \ADDTHIS{ i+2} {b_\bet} 
{{\; \color{cyan}\otimes}}_\bet\;  
 \Pdiptwo^{b_\bet} 
\boxtimes \ADDTHIS {i+1} {\exx }
\boxtimes \ADDTHIS {i} {q-\exx }
&b_\al \leq q \leq  b_{\al\bet}\\[6pt]
\end{cases}
\end{align*} 
 We have that $ {\sf FH}_{q,\al\bet\al}\sim  {\sf  H}_{q,\al\bet\al}$ because the former is obtained from the latter by reflection through the first $\al$-hyperplane it crosses, this is depicted in \cref{yyyyyyy355}.   
Similarly, we define ${\sf FH} _{q, \bet\al\bet}$ to be the path 
 \begin{align*}
\begin{cases}
 \SSTP_{q\emptyset}\boxtimes
\REMOVETHESE{i+1}{b_\bet-q} \boxtimes 
 \REMOVETHESE {i,i+1} {q}  
\boxtimes
\ADDTHIS{i+2}{b_\bet}
{{\; \color{cyan}\otimes}}_\bet\; 
  \REMOVETHESE{i} {b_\al-q} 
\boxtimes \ADDTHIS{i+1} {b_\al}
{{\; \color{magenta}\otimes}}_\al\;  
\REMOVETHESE{i+1}{b_\bet} \;{\color{cyan}\otimes_\bet}\;
\ADDTHIS{i+2}{b_\bet}&0\leq q \leq b_\bet \\[6pt]
\SSTP_{\empb }  \boxtimes
  \REMOVETHESE{i} {q-b_\bet} \boxtimes \REMOVETHESE {i,i+1} {b_\bet}  \boxtimes  \ADDTHIS{i+2}{b_\bet}
{{\; \color{cyan}\otimes}}_\bet\;  \REMOVETHESE{i} {b_\al-q} 
\boxtimes \ADDTHIS{i+1} {b_\al}
{{\; \color{magenta}\otimes}}_\al\;  
\REMOVETHESE{i+1}{b_\bet} \;{\color{cyan}\otimes_\bet}\;
\ADDTHIS{i+2}{b_\bet}&b_\bet  \leq q \leq b_\al  
\\[5pt]  
\SSTP_{\empb }  \boxtimes  
\REMOVETHESE {i } {  q-b_\bet }   \boxtimes  \ADDTHIS{i+2}{b_\bet}
 {{\; \color{cyan}\otimes}}_\bet\; 
\REMOVETHESE{i ,i+2}{q-b_\al} \boxtimes   \ADDTHIS{i+1} {b_\al}
{{\; \color{magenta}\otimes}}_\al\;  
\REMOVETHESE{i+1}{ b_{\al\bet}-q} \;{\color{cyan}\otimes_\bet}\;
\ADDTHIS{i+2}{b_\bet}&b_\al     \leq q  \leq b_{\al\bet}
\end{cases}
\end{align*}
 We have that $ {\sf FH}_{q,  \bet\al\bet}\sim  {\sf  H}_{q, \bet\al\bet}$ because the former is obtained from the latter by reflection through the final $\bet$-hyperplane it crosses.  We note that ${\sf FH}_{b_{\al\bet}, \al\bet\al}=\SSTP_{\emp-\empb}\boxtimes{\sf FH}_{b_{\al\bet},  \bet\al\bet}$.  
  One can define the   paths  $ {\sf FH}_{q,\al\bet\al}$ and 
  $ {\sf FH}_{q, \bet\al\bet}$  for $b_\al<b_\bet$ in an entirely analogous fashion.

 \begin{figure}[ht!]
$$   \begin{minipage}{2.8cm}
\end{minipage}
$$
\caption{
An example of the tableaux ${\sf FH}_{q, \al\bet\al}$ 
for $b_{\al\bet\geq q \geq 0}$. 
We note that ${\sf FH}_{b_{\al\bet}, \al\bet\al}={\sf FH}_{b_{\al\bet},  \bet\al\bet}$.  
The reader should compare these reflected paths with the first  five paths   of \cref{Steinberg}.   }

 \label{yyyyyyy355}
\end{figure}

 \begin{prop}
The element $\Upsilon ^{\SSTP_{\empb\emp\emp  \bet\al\bet}  }
   _{   \SSTP_{\empb\al \emp \bet\al } \otimes \SSTP_\bet^\flat}$ is independent of the choice of reduced expression.  
 \end{prop}
\begin{proof}We proceed as in the proof of \cref{analrsusususus}.  
 We set $\northT={\SSTP_{\empb\emp\emp  \bet\al\bet}  }$ and 
$\southT= {   \SSTP_{\empb\al \emp \bet\al } \otimes \SSTP_\bet^\flat}$.  
 For $0\leq q \leq b_\al+1$, we set   
 $$
 \begin{array}{lll}
 {\sf t}_i(q)= \northT ^{-1}(b_{\al\bet}+q,\eps_i) 
&
 {\sf t}_{i+1}(q)=  \northT^{-1}(b_{\al\al}+q,\eps_{i+1}) \\
 {\sf t}_{i+2}(q)= \northT^{-1}(b_{\al\bet\al}+q,\eps_{i+2})
\\
  {\sf b}_i(q)=  \southT ^{-1}(b_{\al\bet}+q,\eps_{i+1}) &
   {\sf b}_{i+1}(q)=  \southT^{-1}(b_{\al\al\al}+q,\eps_{i+1}) \\
 {\sf b}_{i+2}(q)= \southT^{-1}(b_{\al\bet\al}+q,\eps_{i+2}).
\end{array}
 $$ 
We have  that 
$${\sf t}_i(q) < {\sf t}_i(q+1)< {\sf t}_{i+2}(q)<{\sf t}_{i+2}(q+1)<
 {\sf t}_{i+1}(q)<{\sf t}_{i+1}(q+1)$$
$${\sf b}_i(q) > {\sf b}_i(q+1)> {\sf b}_{i+2}(q)>{\sf b}_{i+2}(q+1)
> {\sf b}_{i+1}(q)>{\sf b}_{i+1}(q+1)$$
for $1\leq q \leq b_\al$ 
and  
$$
 {\sf t}_i(1)< {\sf t}_{i+2}(0)< {\sf t}_{i+1}(1)
 \qquad\quad 
 {\sf t}_i(b_\al)< {\sf t}_{i+2}(b_\al+1)< {\sf t}_{i+1}(b_\al)
   $$
$$
 {\sf b}_i(1)> {\sf b}_{i+2}(0)> {\sf b}_{i+1}(1)
 \qquad\quad 
 {\sf b}_i(b_\al)> {\sf b}_{i+2}(b_\al+1)> {\sf b}_{i+1}(b_\al).  
   $$
     Thus the subexpression $\psi_\w$ is the nib truncation of a 
     quasi-$({b_\al+2})$-expression  for $  w=(13)$,  which is independent of the choice of expression by \cref{nibs}.  Thus the result follows.  
%
%
\end{proof}


   \begin{prop} \label{alhjfkasfhljkfadhsjkafds}
We have that     $$(e_{\SSTP_{\empb \emp}} \otimes {\sf hex}^{\emp\bet\al\bet}_{\empb\al\bet\al})
(e_{\SSTP_{\empb\empb}}\otimes {\sf fork}_{\al\al}^{\emp\al}\otimes e_{\SSTP_{\bet\al}})
 {\sf adj}_{\empb\al\empb{\al\bet\al}}^{\empb\empb\al {\al\bet\al}}
(e_{\SSTP_{\empb\al}}\otimes {\sf hex}_{\emp\bet\al\bet}^{\empb\al\bet\al})
=\Upsilon ^{\SSTP_{\empb\emp\emp  \bet\al\bet}  }
   _{   \SSTP_{\empb\al \emp \bet\al } \otimes \SSTP_\bet^\flat}
$$
     \end{prop}

\begin{proof}
For  $0\leq q < b_{\bet\al }$, we claim that  
\begin{equation*}
 (e_{\SSTP_{  \emp}} \otimes {\sf hex} _{\empb\al\bet\al}(q))
(e_{\SSTP_{ \empb}}\otimes {\sf fork}_{\al\al}^{\emp\al}\otimes e_{\SSTP_{\bet\al}})
  {\sf adj}_{\al\empb\al\bet\al}^{\empb\al\al\bet\al} 
(e_{\SSTP_{\empb\al}}\otimes {\sf hex} ^{\empb\al\bet\al}(q))
    =\Upsilon ^{\SSTP_{\empb\emp} \otimes {\sf H}_{q,\al\bet\al}  }
   _{   \SSTP_{\al\empb}\otimes{\sf  FH}_{q,\bet\al\bet} }
\end{equation*}
and the statement of the proposition will immediately  follow.  We now prove our claim.     
 We set $\northT_q=  {\SSTP_{\empb\emp}} \otimes {\sf hex}_{\al\bet\al}(q)$ and   $\southT_q=  {\SSTP_{\al\empb}} \otimes {\sf  FH}_{q,\al\bet\al} $.  
We consider the  strand, $Q$,  
 from $\northT^{-1}_q(b_{\al\bet}   +q,\eps_i)$ on the   \north    
 edge to
 $\southT^{-1}_q(b_{\al\bet\al}   +q,\eps_{i+1})$ on the   \south       edge of the diagram
 $$ (e_{\SSTP_\emp} \otimes {\sf hex}_{\al\bet\al}(q))
   \circ 
(   {\sf fork}^{\emp\al}_{\al\al}\otimes e_{\SSTP_{\al\bet}})
   \circ 
   (e_{\SSTP_\al }\otimes {\sf hex}^{\al\bet\al}(q))
  $$
for $0\leq q < b_{\al\bet}$.  We wish to consider  the non-zero degree crossings of the $r_q$-strand $Q$ within the diagram.  
 These are with the strands $ \mathscr{Q}_1, \mathscr{Q}_2, \mathscr{Q}_3,  \mathscr{Q}_4,  \mathscr{Q}_5,$ $  \mathscr{Q}_6,  \mathscr{Q}_7 $ connecting the 
 $$
\northT^{-1}_{q+1}( b_{\al\bet} +q-1,\eps_i)
, \;
\northT^{-1}_{q+1}( b_{\al\bet\al} +q ,\eps_{i+1})
, \;
\northT^{-1}_{q+1}( b_{\al
\bet\al} +q+1 ,\eps_{i+1})
, \;
\northT^{-1}_{q+1}( b_{\al\bet\al} +q+2,\eps_{i+1})
$$
$$
\northT^{-1}_{q+1}( b_{\al\bet\al\bet} +q+1,\eps_{i+2})
, \;
\northT^{-1}_{q+1}( b_{\al\bet\al\bet} +q+2,\eps_{i+2})
, \;
\northT^{-1}_{q+1}( b_{\al\bet\al\bet} +q+3,\eps_{i+2})
 $$
   \north    vertices (which are ordered in increasingly from left to right) to the 
 $$
 \southT^{-1}_{q+1}( b_{\al\bet} +q ,\eps_{i})
, \;
\southT^{-1}_{q+1}(  b_{\al\bet\al}+q,\eps_{i+1}) 
, \;
\southT^{-1}_{q+1}( b_{\al\bet} +q+1 ,\eps_{i})
, \;
\southT^{-1}_{q+1}( b_{\al\bet\al} +q+2,\eps_{i+1})
$$
$$
\southT^{-1}_{q+1}( b_{\al\bet\al\bet} +q+1,\eps_{i+2})
, \;
\southT^{-1}_{q+1}( b_{\al\bet\al\bet} +q+2,\eps_{i+2})
, \;
\southT^{-1}_{q+1}( b_{\al\bet\al\bet} +q+3,\eps_{i+2})
 $$ 
   \south       vertices, respectively.  
 The residues of these strands are $r_q+1,r_q+1, r_q , r_q-1$ for the first row and and $r_q+1,r_q,r_q-1$ or the second row.   
 We have that 
 $$
 \northT^{-1}_{q+1}( b_{\al\bet} +q-1,\eps_i)
<
\northT^{-1}_{q+1}( b_{\al\bet\al} +q ,\eps_{i+1})
$$
$$ \southT^{-1}_{q+1}( b_{\al\bet} +q ,\eps_{i})>\southT^{-1}_{q+1}(  b_{\al\bet\al}+q,\eps_{i+1}) $$
 and so the   pair of strands $\mathscr{Q}_1$ and $\mathscr{Q}_2$ form a crossing of $(r_q+1)$-strands. 
 The strand $Q$ crosses   $\mathscr{Q}_1$ and $\mathscr{Q}_2$  exactly once each.    
  The remaining 5 strands are all vertical lines  (in other words their   \north    and   \south       vertices coincide).  
 The strand $Q$ crosses each of these 
 vertical strands twice.  (Thus the total degree contribution of these crossings is zero.) 
 
 We undo the  crossing of  $Q$ with  the triple of strands $\mathscr{Q}_5, \mathscr{Q}_6, \mathscr{Q}_7$  as in the proof of \cref{adjust1}.    
Pull the $Q$ strand through  $\mathscr{Q}_4$ using case 4 of relation \ref{rel1.10} at the expense of acquiring a dot on $Q$ (the other term is zero by case 1 of relation \ref{rel1.10}
) we then pull the dot on $Q$ upwards through the crossing of $Q$ and 
$\mathscr{Q}_3$ using relation \ref{rel1.8} and obtain two terms:  
the first term, in  which   the dot has passed through the crossing,  is zero   by case  1 of relation \ref{rel1.10};  
 in the second term, in which  we  undo one (of the two)  crossings between $Q$ and $\mathscr{Q}_3$, is equal to 
 $\psi ^{\SSTP_\empb \otimes {\sf H}_{q,\al\bet\al}  }
   _{   \SSTP_\al\otimes{\sf  FH}_{q,\bet\al\bet}}$ as required.     
   
Now suppose $b_\al \leq q <b_{\al\bet}$.   The $r_q$-strand connecting the  
$\southT^{-1}(4b_\al+2b_\bet-q,\eps_{i+1})$ and 
$\northT^{-1}(4b_\al+2b_\bet-q,\eps_{i+1})$
top and bottom nodes double-crosses the $(r_q+1)$-
$r_q$- and $(r_q-1)$- strands connecting the 
$\northT^{-1}(4b_\al+3b_\bet-q-1,\eps_{i+2})$,  
$\northT^{-1}(4b_\al+3b_\bet-q,\eps_{i+2})$, 
$\northT^{-1}(4b_\al+3b_\bet-q+1,\eps_{i+2})$  
top vertices to the 
$\southT^{-1}(4b_\al+3b_\bet-q-1,\eps_{i+2})$,  
$\southT^{-1}(4b_\al+3b_\bet-q,\eps_{i+2})$, 
$\southT^{-1}(4b_\al+3b_\bet-q+1,\eps_{i+2})$  
bottom vertices.  We undo these double-crossings as in the proof of \cref{adjust1}.   \end{proof}

    \begin{prop}  
We have that 
\begin{equation}\label{again again} {\sf adj}^{\empb\emp\emp  \bet\al\bet}_{\emp\emp \bet\al\empb\bet} 
(e_{ {\SSTP_{ \emp\emp\bet\al}}} \otimes {\sf fork}^{\empb\bet}_{\bet\bet})
 {\sf hex}^{\emp\emp \bet\al\bet \bet }_{\emp\empb \al\bet\al  \bet}
 {\sf adj}^{\emp\empb \al\bet\al   \bet} _{\empb \al \emp\bet\al  \bet}
=\Upsilon ^{\SSTP_{\empb\emp\emp  \bet\al\bet}  }
   _{   \SSTP_{\empb\al \emp \bet\al } \otimes \SSTP_\bet^\flat}.
\end{equation}
     \end{prop}
      \begin{proof} For $0 \leq q \leq b_{\al\bet}$, we claim that  
 \begin{align}\label{align2}
   \Upsilon
   _{\SSTP_\emp\otimes  {\sf H}_{q,\empb \al\bet\al}   \otimes \SSTP_\bet}
   ^{	 \SSTP_{\emp\emp\empb\bet\al \bet}	}
  (e_{\SSTP_{\emp}}\otimes 	{\sf hex}_{ \empb\al\bet\al}(q+1)\otimes e_{\SSTP_\bet})
 &=
   \Upsilon_{\SSTP_\emp\otimes {\sf H}_{q ,\empb\al\bet\al}  \otimes \SSTP_\bet}^{	 \SSTP_{\emp\emp\empb\bet\al \bet}	}.
 \end{align}
 We decorate the top and bottom edges of the concatenated diagram in  \cref{align2} by  the paths 
    $\northT= \SSTP_{{\emp\emp\empb \bet\al \bet}}$ and 
    $\southT_{q+1}= {\SSTP_\emp\otimes {\sf H}_{q +1,\empb\al\bet\al}  \otimes \SSTP_\bet}$. 
         For each  $  0\leq q <b_\bet$  the   strand    (of residue $r_q\in \ZZ/e\ZZ$, say)
  connecting the 
    \north    $\northT^{-1}(b_{\al\bet}+q,\eps_{i+1})$)th  and 
  $\southT^{-1}_q(b_{\al\bet}+q,\eps_{i+1}) $th  
   \south       vertices (both of which are equal to $ (b_{\al\bet}+q){\aatch} + \emptyset(i+1) $)
  of the concatenated diagram  
  has  double-crossings of non-zero degree with three    strands 
 of residues $r_q+1$, $r_q$ and $r_q-1$ connecting the
  $\northT^{-1}(b_{\bet\al\bet} -1+q,\eps_{i+2})$th,    $\northT^{-1}(b_{\bet\al\bet} +q,\eps_{i+2})$th,   and  $\northT^{-1}(b_{\bet\al\bet} +q+1,\eps_{i+2})$th     \north    vertices to the 
    $\southT_q^{-1}(b_{\bet\al\bet} -1+q,\eps_{i+2})$th,    $\southT^{-1}_q(b_{\bet\al\bet} +q,\eps_{i+2})$th,   and  $\southT_q^{-1}(b_{\bet\al\bet} +q+1,\eps_{i+2})$th     \south       vertices
   respectively; we undo these crossings using \cref{adjust1}.  
 Now, for  $b_\bet \leq q <b_{\al\bet}$ the claim is immediate as the concatenated diagram is step-preserving and has minimal length.  
Finally, we substitute \cref{align2} into \cref{again again}  and the resulting diagram is again  step-preserving and has minimal length and  the result follows.  
             \end{proof}

\subsection{The tetrahedron relation}\label{tetratetra}
\renewcommand{\bet}{{{\color{cyan}\boldsymbol\beta}}} 
 \renewcommand{\gam}{{{{\color{ao(english)}\boldsymbol\gamma}}}}

We now check that the image of  relation \ref{dragrace} holds in the quiver Hecke algebra.
 Our aim is to  show that 
$$
 {\sf hex}^{ \gam\al\gam\bet\al\gam{\empb\emp\empb}}_{ \al\gam\al\bet\al\gam\empb\empg\empb  }  
 {\sf hex}^{  \al\gam\al\bet\al  \gam\empb\empg\empb}_{  \al\gam \bet\al\bet  \gam\emp\empg\empb}  
 {\sf com}^{ \al\gam\bet\al \bet\gam\emp \empg\empb}_{\al\bet\gam\al \gam\bet\emp \empg\empb}    {\sf hex}^{   \al\bet	\gam\al\gam\bet\emp\empb\empg}_{  \al\bet	 \al\gam\al\bet\empg\empb\empg} 
    {\sf hex}^{  \al\bet	 \al\gam\al \bet\empg\empb\empg} _{  \bet \al\bet\gam\al\bet\empg\emp\empg } 
	  {\sf com}^{\bet\al\bet\gam\al\bet\empg\emp \empg}_{\bet\al \gam\bet\al\bet\empg\emp \empg}  
 $$
is equal to 
$$
 {\sf com}_{ \gam\al\bet\gam \al\gam\empb\emp\empb}^{  \gam\al\gam\bet \al\gam\empb\emp\empb}  
  {\sf hex}^{ \al\bet	\gam \al\gam\empb\emp\empb\gam}
 _{  \gam\al\bet  \al\gam\al\empb\empg\empb}
  {\sf hex}^ 
  {  \gam\al\bet  \al\gam\al\empb\empg\empb}
 _{   \bet  \al\bet\gam\al\emp \empg\empb\gam}
    	{\sf com}^{ \gam\bet\al \bet\gam\al{\emp \empg\empb }}_{\bet\gam\al \gam\bet\al{\emp \empg\empb }} 
   {\sf hex}^ 
 { 	  \bet \gam \al \gam\bet\al \emp \empg\empb	}
 _{ 	 \empg \bet		 \al \gam\al \bet\al\empg \empb	 }
   {\sf hex}^ 
  {   \bet		 \al \gam\al \bet\al\empg \empb	 \empg	}
_  {   \bet		 \al \gam  \bet\al \bet\empg \emp 	 \empg	}.
$$

\begin{prop}\label{length5}
The element   $\psi^{ 
  \SSTP_{\gam\al \gam\bet\al\gam\empb\emp\empb }
}_{
  \SSTP_{\bet\al\gam\bet \al\bet\empg\emp\empg }
 }$ 
 is independent of the choice of reduced expression.  
\end{prop}

\begin{proof}
For notational ease, we let $j=i+1$ and $k=i-1$ and we 
decorate the top and bottom edges with 
$ 
\SSTT=\SSTP_{\gam\al\gam\bet\al\gam\empb\emp\empb}
$ and $
{\sf B}=\SSTP_{\bet\al\gam\bet\al\bet\empg\emp\empg}
$  respectively.  
 For each  
$b_\bet \leq  q \leq b_{\al\bet}+1$,   we consider  the collection of  permutations $w_q$ formed from the $r_q$-strands connecting  each of  the 
\begin{align*} 
{\sf B}_{i-1}(q) = {\sf B}^{-1}( q,\eps_{i-1})\quad\quad&
{\sf B}_{i}(q) ={\sf B}^{-1}(b_{\gam }+ q,\eps_{i}) 
\\
{\sf B}_{i+1}(q) ={\sf B}^{-1}(b_{ \al  \gam}+q,\eps_{i+1})
  \quad\quad & {\sf B}_{i+1}(q)  ={\sf B}^{-1}(b_{ \al \bet\gam}+q,\eps_{i+2})
 \intertext{bottom vertices   
to  }
{\sf T}_{i-1}(q) = {\sf T}^{-1}( q,\eps_{i-1})\quad\quad&
{\sf T}_{i}(q) ={\sf T}^{-1}(b_{\gam }+ q,\eps_{i}) 
\\
{\sf T}_{i+1}(q) ={\sf T}^{-1}(b_{ \al  \gam}+q,\eps_{i+1})\quad\quad
   & {\sf T}_{i+1}(q) ={\sf T}^{-1}(b_{ \al \bet\gam}+q,\eps_{i+2})
\end{align*}
  top vertices  respectively.   By definition $r_q=r_{q+1}+1$ for $b_\bet \leq  q < b_{\al\bet}+1$. 
We let $\w$ denote the subexpression consisting of all strands from    (the union of)  the 
 $w_q$-subexpressions for $b_\bet \leq  q \leq b_{\al\bet}+1$.  
One can verify, simply by looking at the paths ${\sf T}$ and ${\sf B}$ (and their residue sequences) that any bad-crossing  in $w$ belongs to $\psi_{{\rm nib}(\w)}$.  
We have that 
\begin{align*}
&{\sf B}_{i-1}(q)<{\sf B}_{i-1}(q+1) <{\sf B}_{i+2}(q)<{\sf B}_{i+2}(q+1)
<{\sf B}_{i+1}(q)<{\sf B}_{i+1}(q+1)<{\sf B}_{i}(q)<{\sf B}_{i}(q+1) 
\\ &{\sf T}_{i-1}(q)>{\sf T}_{i-1}(q+1) >{\sf T}_{i+2}(q)>{\sf T}_{i+2}(q+1)
>{\sf T}_{i+1}(q)>{\sf T}_{i+1}(q+1)>{\sf T}_{i}(q)>{\sf T}_{i}(q+1).
\end{align*}
for  $b_\bet < q < b_{\al\bet}$.  
In other words,   the  $r_q$-strands for $b_\bet < q \leq b_{\al\bet}$  form a 
$\psi_{(1,4)(2,3)_{b_\al}}$ braid (and thus this subexpression is quasi-dilated  and  of breadth $b_\al$).   
We now restrict to the case $q=b_\bet$, as the  $q={b_{\al\bet}+1}$ is similar.  
We have that 
\begin{align*}
 {\sf B}_{i-1}(b_\bet+1)<
 {\sf B}_{i+2}(b_\bet)<{\sf B}_{i+2}(b_\bet+1)
<{\sf B}_{i+1}(b_\bet)<{\sf B}_{i+1}(b_\bet+1) 
< {\sf B}_{i}(b_\bet+1) 
\\
 {\sf T}_{i-1}(b_\bet+1)>
   {\sf T}_{i+2}(b_\bet)>{\sf T}_{i+2}(b_\bet+1)
>{\sf T}_{i+1}(b_\bet)>{\sf T}_{i+1}(b_\bet+1) 
>{\sf T}_{i}(b_\bet+1) .
\end{align*}
(We have not considered the strands connecting ${\sf B}_{i-1}(b_\bet)$ and $ {\sf T}_{i-1}(b_\bet)$ 
or ${\sf B}_{i}(b_\bet)$ and $ {\sf T}_{i}(b_\bet)$ as these were removed under the nib truncation map.) 
Thus $\psi_{{\rm nib}(\w)}$ is independent of the choice of expression by  \cref{nibs} and  the result follows.   See \cref{alabelforafigre22222} for an example.  
  \end{proof}

 \!\!\!\!\!
\begin{figure}[ht!]
$$
  
$$
\caption{The element  $\psi^
{   \SSTP_{\gam\al \gam\bet\al\gam\empb\emp\empb }
}_{
  \SSTP_{\bet\al\gam\bet \al\bet\empg\emp\empg }
 }$ 
 for $p=5$, $h=3$,   $\ell=1$ and $\al=\eps_2-\eps_3$, $\bet=\eps_3-\eps_4$, $\gam=\eps_1-\eps_2$. The thick black $4$-strands   form a $\w=s_3s_2s_1s_3 s_2 s_3$ braid.  Together with the wiggly strands, these form a  
 subexpression ${\rm nib}\psi_{\w_3}$  containing all bad crossings.  
 }
\label{alabelforafigre22222}\end{figure}

\begin{prop}\label{yyyyyyyyyy}
We have that $ \Upsilon^{ 
  \SSTP_{\gam\al\gam\bet \al\gam\empb\empb\emp\empb}
}_{
  \SSTP_{ \bet\al\gam\bet \al\bet\empb\empg\emp\empg} 
 } $ is equal to both 
$$
 {\sf hex}^{ \gam\al\gam\bet\al\gam{\empb\emp\empb}}_{ \al\gam\al\bet\al\gam\empb\empg\empb  }  
 {\sf hex}^{  \al\gam\al\bet\al  \gam\empb\empg\empb}_{  \al\gam \bet\al\bet  \gam\emp\empg\empb}  
 {\sf com}^{ \al\gam\bet\al \bet\gam\emp \empg\empb}_{\al\bet\gam\al \gam\bet\emp \empg\empb}    {\sf hex}^{   \al\bet	\gam\al\gam\bet\emp\empb\empg}_{  \al\bet	 \al\gam\al\bet\empg\empb\empg} 
    {\sf hex}^{  \al\bet	 \al\gam\al \bet\empg\empb\empg} _{  \bet \al\bet\gam\al\bet\empg\emp\empg } 
	  {\sf com}^{\bet\al\bet\gam\al\bet\empg\emp \empg}_{\bet\al \gam\bet\al\bet\empg\emp \empg}  
 $$
and 
$$
 {\sf com}_{ \gam\al\bet\gam \al\gam\empb\emp\empb}^{  \gam\al\gam\bet \al\gam\empb\emp\empb}  
  {\sf hex}^{ \al\bet	\gam \al\gam\empb\emp\empb\gam}
 _{  \gam\al\bet  \al\gam\al\empb\empg\empb}
  {\sf hex}^ 
  {  \gam\al\bet  \al\gam\al\empb\empg\empb}
 _{   \bet  \al\bet\gam\al\emp \empg\empb\gam}
    	{\sf com}^{ \gam\bet\al \bet\gam\al{\emp \empg\empb }}_{\bet\gam\al \gam\bet\al{\emp \empg\empb }} 
	    {\sf hex}^ 
 { 	  \bet \gam \al \gam\bet\al \emp \empg\empb	}
 _{ 	 \empg \bet		 \al \gam\al \bet\al\empg \empb	 }
   {\sf hex}^ 
  {   \bet		 \al \gam\al \bet\al\empg \empb	 \empg	}
_  {   \bet		 \al \gam  \bet\al \bet\empg \emp 	 \empg	}.
$$

\end{prop}
\begin{proof}
We set $k=i-1$, $j=i+1$.  We will prove the first  equality as the second is very similar (for more details, see \cref{remarkremarkbelow}).    We   proceed from the centre of the diagram, considering the first pair of hexagons (on top and bottom of a pair of commutators), the second pairs of hexagons (on top and bottom of the previous product) and then finally the last commutator (below the previous product).  

\smallskip
\noindent{\bf Step 1}. We  add the first pair of  hexagonal generators symmetrically as follows 
 \begin{equation}\label{length7}{\sf hex}^{  \al\gam\al\bet\al  \gam\empb\empg\empb}_{  \al\gam \bet\al\bet  \gam\emp\empg\empb}  
(	e_{\SSTP_{ \al}}\otimes
{\sf com}^{ \gam\bet\al \bet\gam}_{\bet\gam\al \gam\bet}  
 \otimes e_{\SSTP_{\emp \empg\empb}}  	)
 {\sf hex}^ {   \al\bet	\gam\al\gam\bet\emp\empb\empg} _{  \al\bet	 \al\gam\al\bet\empg\empb\empg}  
=\Upsilon^{\SSTP_{  \al\gam\al\bet\al  \gam\empb\empg\empb}} _{ \SSTP_{ \al\bet	 \al\gam\al\bet\empg\empb\empg}}.    
\end{equation} The only points  worth bearing in mind are $(i)$
double-crossings strands of non-adjacent residue can be undone trivially and $(ii)$ that 
  the implicit adjustments in the definitions of 
 ${\sf hex}^{  \al\gam\al\bet\al  \gam\empb\empg\empb}_{  \al\gam \bet\al\bet  \gam\emp\empg\empb}    $ and 
$ {\sf hex}^ {   \al\bet	\gam\al\gam\bet\emp\empb\empg} _{  \al\bet	 \al\gam\al\bet\empg\empb\empg}   $ will give rise to  (a total of 
$|b_\al-b_\bet|+ |b_\al-b_\gam|+ |b_\bet-b_\gam|$) double-crossings  which can be undone as in the proof of  \cref{adjust1}.

\smallskip
\noindent{\bf Step 2}.    We now add the next pair of hexagonal generators symmetrically to the diagram, $\Upsilon ^ {\SSTP_{\al\gam\al\bet\al  \gam\empb\empg\empb}}_{ \SSTP_{ \al\bet	 \al\gam\al\bet\empg\empb\empg}}$, output by the previous step in the procedure.  
   We first note   
   that 
$$
 {\sf adj}^{	\SSTP_\empg\otimes{\sf H}_{ 0, \al\gam\al	}  \otimes \SSTP_{\bet\al  \gam 	\empb\empb	}	}
 _{  \SSTP_{\al\gam\al\bet\al  \gam\empb\empg\empb }}
\circ 
\Upsilon 
^ {\SSTP_{\al\gam\al\bet\al  \gam\empb\empg\empb}}
 _{  \SSTP_{  \al\bet	 \al\gam\al\bet\empg\empb\empg}}
 \circ {\sf adj}
 ^{   \SSTP_{ \al\bet	 \al\gam\al\bet\empg\empb\empg}}
 _{\SSTP_\empb\otimes	{\sf H}_{ 0, \al\bet\al	} \otimes \SSTP_{ \gam\al\bet\empg\empg}}
=\Upsilon ^{\SSTP_\empg\otimes 	{\sf H}_{0, \al\gam\al	}  \otimes \SSTP_{\bet\al  \gam 	\empb\empb	}	}
_{\SSTP_\empb \otimes 	{\sf H}_{0,  \al\bet\al	}   \otimes \SSTP_{\gam\al\bet\empg\empg	} }
$$ again  by  (a total of 
$|b_\bet - b_\gam|$ applications of)  \cref{adjust1}. 
     We claim that 
\begin{equation}\label{hj1} \big( 	{\sf hex}^{\empg \al\gam\al	}(q)\otimes e_{	\SSTP_{	\bet\al  \gam\empg\empb }}  \big)
\Upsilon ^{  	\SSTP_\empg\otimes 	{\sf H}_{q
 ,  \al\gam\al	}\otimes 	\SSTP_{ \bet\al  \gam\empg\empg  }
 }
 _{ 
\SSTP_\empb \otimes 	{\sf H}_{q, \al\bet \al	} 	\otimes \SSTP_{\gam\al\bet\empg\emp  } }
 \big( 
 {\sf hex}_{\empb \al\bet \al	}(q)	\otimes e_{ \SSTP_{\gam\al\bet\empg\empg } }  \big)
 =
\Upsilon ^{ 
 	\SSTP_\empg \otimes	{\sf H}_{ q+1, \al\gam\al	} \otimes 	\SSTP_{	\bet\al  \gam\empb\empb }}
 _{ 
   \SSTP_\empb \otimes 	{\sf H}_{ q+1, \al\bet \al	} 	\otimes \SSTP_{\gam\al\bet\empg\empg  } }
\end{equation}
for $   0\leq q <\max\{b_{\bet},b_{\gam}\}+b_\al$.  
  For  
  $0 \leq q \leq 		b_\al+ 
  |b_{\bet}- b_{\gam}| $
the concatenated diagram  on the lefthand-side of  \cref{hj1} 
contains a distinguished strand  connecting the 
   $  {\sf T}^{-1}(\min\{b_{\bet},b_{\gam}\}+  q+1,\eps_{i  }) $ top and 
   $  {\sf B}^{-1}(\min\{b_{\bet},b_{\gam}\}+  q+1,\eps_{i  }) $ 
  bottom vertices.  
   For 
  $0  \leq q \leq 		b_\al 	+|b_\bet-b_\gam|   	 $   the distinguished strand passes from left to right and back again,  thus admitting a  double-crossing with each of the  
   $(r_q-1)$-, $r_q$-, $(r_q+1)$-strands  connecting the 
      $$
      \begin{array}{lll}
  {\sf T}^{-1}(\min\{b_{\bet},b_{\gam}\}+    b_{\al}+ q ,\eps_{i+1  }) 
&
    {\sf T}^{-1}(\min\{b_{\bet},b_{\gam}\}+b_{\al}+ q +1,\eps_{i+1  }) 
  \\
      {\sf T}^{-1}(\min\{b_{\bet},b_{\gam}\}+ b_{\al}+ q+2,\eps_{i+1  }) 
      \end{array}
$$
top vertices to the 
 $$
      \begin{array}{lll}
  {\sf B}^{-1}(\min\{b_{\bet},b_{\gam}\}+b_{\al}+ q  ,\eps_{i+1  }) 
  &
    {\sf B}^{-1}(\min\{b_{\bet},b_{\gam}\}+b_{\al}+ q+1 ,\eps_{i+1  }) 
\\
      {\sf B}^{-1}(\min\{b_{\bet},b_{\gam}\}+b_{\al}+ q+2,\eps_{i+1  }) 
      \end{array}
$$
bottom vertices.
For   $|b_\bet-b_\gam|     \leq q \leq 		b_\al 	+|b_\bet-b_\gam|   	 $   the distinguished strand {\em also} admits a  double-crossing with each of the 
 $(r_q-1)$-, $r_q$-, $(r_q+1)$-strands  connecting the
 $$
       \begin{array}{lll}
       {\sf T}^{-1}(\min\{b_{\bet},b_{\gam}\}+b_{\al\bet}+ q ,\eps_{i+2  }) 
&
    {\sf T}^{-1}(\min\{b_{\bet},b_{\gam}\}+b_{\al\bet}+ q+1 ,\eps_{i+2  }) 
  \\
      {\sf T}^{-1}(\min\{b_{\bet},b_{\gam}\}+b_{\al\bet}+ q+2,\eps_{i+2  })
     \end{array} $$
  top vertices to the     
       $$
     \begin{array}{lll}  {\sf B}^{-1}(\min\{b_{\bet},b_{\gam}\}+b_{\al\bet}+ q ,\eps_{i+2  }) 
&
    {\sf B}^{-1}(\min\{b_{\bet},b_{\gam}\}+b_{\al\bet}+ q+1 ,\eps_{i+2  }) 
  \\
      {\sf B}^{-1}(\min\{b_{\bet},b_{\gam}\}+b_{\al\bet}+ q+2,\eps_{i+2  })
      \end{array}$$
      bottom vertices.  
 Note     we have broken these strands into two triples.  
  For  $0  \leq q \leq 		b_\al 	+|b_\bet-b_\gam|   	 $ we undo the 
 double-crossing of the distinguished strand with the former triple using a single application of \cref{adjust1}.  
  For  $|b_\bet-b_\gam|    \leq q \leq 		b_\al 	+|b_\bet-b_\gam|   	 $ we undo the 
 double-crossing of the distinguished strand with the latter triple and then the former triple as in the proof of  \cref{adjust1}.  
Thus \cref{hj1} follows.  
If $b_\bet > b_\gam$ (respectively $b_\gam >b_\bet$) 
we must now multiply on the bottom (respectively top) by the remaining 
terms to obtain a minimal, step-preserving diagram.  
We hence deduce that  
$$\big( 	{\sf hex}_{\empg \al\gam\al	} \otimes e_{	\SSTP_{	\bet\al  \gam\empg\empb }}  \big)
\Upsilon^{\SSTP_{  \al\gam\al\bet\al  \gam\empb\empg\empb}} _{ \SSTP_{ \al\bet	 \al\gam\al\bet\empg\empb\empg}}      \big( 
 {\sf hex}_{\empb \al\bet \al	} 	\otimes e_{ \SSTP_{\gam\al\bet\empg\empg } }  \big)
 =
\Upsilon ^{ 
 		{\sf H}_{ b_{\al\gam}, \empg \al\gam\al	} \otimes 	\SSTP_{	\bet\al  \gam\empb\empb }}
 _{ 
   	{\sf H}_{b_{\al\bet},\empb \al\bet \al	} 	\otimes \SSTP_{\gam\al\bet\empg\empg  } }.
$$
We now multiply   on the top and bottom by the  other ``halves" of the hexagonal generators   to get 
\begin{equation}\label{length6}
 {\sf hex}^{ \gam\al\gam\bet\al\gam\empb\emp\empb}_{ \al\gam\al\bet\al\gam\empb\empg\empb  }  
\Upsilon^{\SSTP_{  \al\gam\al\bet\al  \gam\empb\empg\empb}} _{ \SSTP_{ \al\bet	 \al\gam\al\bet\empg\empb\empg}}       {\sf hex}^{  \al\bet	 \al\gam\al \bet\empg\empb\empg} _{  \bet \al\bet\gam\al\bet\empg\emp\empg } 
 = \Upsilon ^
  {\SSTP_{ \gam\al\gam\bet\al  \gam\empb\emp \empb}}
  _
  { \SSTP_{ \bet	 \al\bet\gam\al\bet\empg\emp \empg}} 
 \end{equation} 
where here the hexagonal terms are minimal and step-preserving, but we must 
again  undo  any double-crossings arising from adjustments as in the proof of \cref{adjust1}. 
We emphasise that the righthand-side of \cref{length6} is independent of the choice of reduced expression, which can be shown in a similar fashion to \cref{length5}.  

\medskip
\noindent{\bf Step 3}.  For $0\leq q < b_{\bet\gam} $, we claim that 
$$ \Upsilon^{ 
  \SSTP_{\gam\al \gam\bet\al\gam\empb\emp\empb}
}_{{\SSTP_{ \bet\al}}\otimes 	{\sf C}^{q,\bet\gam}   \otimes e_{\SSTP_{\al\bet\empg\emp \empg}}}
   (	e_{\SSTP_{ \bet\al}}\otimes 	{\sf com}^{\bet\gam} (q) \otimes e_{\SSTP_{\al\bet\empg\emp \empg}} )
 =
 \Upsilon^{ 
  \SSTP_{ \gam\al \gam\bet\al\gam\empb\emp\empb}
}_{{\SSTP_{ \bet\al}}\otimes 	{\sf C}^{q+1,\bet\gam}   \otimes {\SSTP_{\al\bet\empg\emp \empg}}}
$$  and for $  b_{\bet\gam}\geq q >0 $, we claim that 
$$ \Upsilon^{ 
  \SSTP_{ \gam\al \gam\bet\al\gam\empb\emp\empb}
}_{{\SSTP_{ \bet\al}}\otimes 	{\sf C}_{q,\gam\bet}   \otimes e_{\SSTP_{\al\bet\empg\emp \empg}}}
  (	e_{\SSTP_{ \bet\al}}\otimes 	{\sf com}_{ \gam\bet} (q) \otimes e_{\SSTP_{\al\bet\empg\emp \empg}} )
 =
 \Upsilon^{ 
  \SSTP_{ \gam\al \gam\bet\al\gam\empb\emp\empb}
}_{{\SSTP_{ \bet\al}}\otimes 	{\sf C} _{ \gam\bet} (q-1) \otimes  {\SSTP_{\al\bet\empg\emp \empg}}}.
$$    We consider the former product, as the latter is similar.   If $b_\gam> b_\bet$, then the   concatenated diagram is minimal and step-preserving.  
If $b_\gam\leq  b_\bet$ then    the $r_q$-braid connecting the strands
$$
  {\sf T}^{-1}( q+1,\eps_{i-1})   \quad
 {\sf T}^{-1}(b_{\gam}+ q+1,\eps_{i}) 
\quad {\sf T}^{-1}(b_{  	 \al\gam		}+q+1,\eps_{i+1})    \quad 
 {\sf T}^{-1}(b_{\al\bet\gam}+q+1,\eps_{i+2}) 
 $$
$$
  {\sf B}^{-1}( q+1,\eps_{i-1})   \quad
 {\sf B}^{-1}(b_{\gam}+ q+1,\eps_{i}) 
\quad {\sf B}^{-1}(b_{  	 \al\gam		}+q+1,\eps_{i+1})    \quad 
 {\sf B}^{-1}(b_{\al\bet\gam}+q+1,\eps_{i+2}) 
 $$
 top and bottom vertices form the non-minimal expression    $(s_2 s_1 s_3   s_2  s_3) s_3$ 
(the   bracketed term belongs to the multiplicand 
$ \Upsilon ^
  {\SSTP_{ \gam\al\gam\bet\al  \gam\empb\emp \empb}}
  _
  { \SSTP_{ \bet	 \al\bet\gam\al\bet\empg\emp \empg}} $   and so can be   chosen  
arbitrarily, we have chosen the simplest form for what follows).  
  The $r_q$-strand with label $\eps_i$ double-crosses the $(r_{q}-1) $-strand connecting the 
$  {\sf T}^{-1}(b_{  	 \al\gam		}+q+2,\eps_{i+1})$  and 
$  {\sf B}^{-1}(b_{   \al\gam		}+q+2,\eps_{i+1})$ top and bottom vertices.   We undo this double-crossing at the expense of placing a KLR dot on the  $r_q$-strand (the other term is zero, by case 1 of \cref{rel1.10}). 
We then pull this dot through the $r_q$-crossing labelled by the $\eps_i$ and $\eps_{i+2}$ strands   
 and hence  undoing the bottommost crossing (the other, dotted, term is zero, again 
by case 1 of \cref{rel1.10}).  
Thus our $r_q$-braid now forms 
the   non-minimal expression $ s_2 s_1 s_3   s_2   s_3$.  The 
  $r_q$-crossing of strands connecting the 
  $$
   {\sf T}^{-1}(b_{  	 \al\gam		}+q+1,\eps_{i+1})     ,\;
 {\sf T}^{-1}(b_{\bet}+ q+1,\eps_{i}) 
,\;  {\sf B}^{-1}(b_{\bet}+ q+1,\eps_{i})  ,\;  {\sf B}^{-1}(b_{  	 \al\gam		}+q+1,\eps_{i+1})$$
top and  bottom vertices is bi-passed on the left by the 
 $(r_q+1)$-strand connecting the 
$  {\sf T}^{-1}(b_{   \al\gam		}+q ,\eps_{i+1})$  and 
$  {\sf B}^{-1}(b_{   \al\gam		}+q ,\eps_{i+1})$ vertices.  
We pull this  $(r_q+1)$-strand through this crossing using relation \ref{rel1.11} and hence obtain the diagram in which the crossing is undone (at the expense of another term, which is zero by \cref{resconsider}).  Thus our $r_q$-braid now forms 
the minimal expression $ s_2 s_1 s_3   s_2   $, and the diagram is minimal and step-preserving, as required.  
\end{proof}

\begin{rmk}\label{remarkremarkbelow}
The reader should note that in \cref{dragrace}, the righthand-side is obtained by first flipping the lefthand-side through the horizontal and vertical axes and then swapping the $\bet$ and $\gam$ labels. 
The  ``very similar" proof of the second equality in \cref{yyyyyyyyyy}  amounts to rewriting the above argument but with indices of the crossing-strands determined by the horizontal and vertical flips and recolouring (swap mentions of $b_\bet$ and $b_\gam$) of the indices in the proof above.  
\end{rmk}

\subsection{The tricoloured  commutativity relations}
We now verify the two relations depicted in \ref{S7}.  Namely, we will show that 
\begin{align}\label{popopopopo1}
\begin{split}
\Upsilon^{\empb \al\bet\al{\color{orange}\boldsymbol\delta}}_{{   {\color{orange}\boldsymbol\delta}\emp \bet \al\bet }}&={\sf hex}^{\empb \al\bet\al{\color{orange}\boldsymbol\delta}}_{\emp  \bet\al\bet{\color{orange}\boldsymbol\delta}}
{\sf com}_{\emp  \bet\al{\color{orange}\boldsymbol\delta}\bet}
^{\emp  \bet\al\bet{\color{orange}\boldsymbol\delta} }
{\sf com}^{\emp  \bet\al{\color{orange}\boldsymbol\delta}\bet}
_{\emp  \bet {\color{orange}\boldsymbol\delta}\al\bet}
{\sf com}^{\emp  \bet {\color{orange}\boldsymbol\delta}\al\bet}
_{\emp    {\color{orange}\boldsymbol\delta}\bet\al\bet}
{\sf adj}^{ \emp {\color{orange}\boldsymbol\delta} \bet \al\bet }
_{   {\color{orange}\boldsymbol\delta}\emp \bet \al\bet }
 \\
&= {\sf com}^{\empb\al\bet\al{\color{orange}\boldsymbol\delta}}
_{ \empb  \al\bet {\color{orange}\boldsymbol\delta}\al }
{\sf com}^{ \empb  \al\bet {\color{orange}\boldsymbol\delta}\al }
_{ \empb  \al  {\color{orange}\boldsymbol\delta}\bet\al }
{\sf com}^ { \empb  \al  {\color{orange}\boldsymbol\delta}\bet\al }_
{ \empb  {\color{orange}\boldsymbol\delta}  \al \bet\al }
 {\sf adj}^
 { \empb  {\color{orange}\boldsymbol\delta}  \al \bet\al }
 _{    {\color{orange}\boldsymbol\delta}\empb  \al \bet\al }
 {\sf hex}^{    {\color{orange}\boldsymbol\delta}\empb  \al \bet\al }_{    {\color{orange}\boldsymbol\delta}\emp   \bet\al  \bet}
 \end{split}
 \intertext{and we have that}\label{popopopopo2}
   \Upsilon^{{   \bet\gam {\color{orange}\boldsymbol\delta} }}_{   {\color{orange}\boldsymbol\delta}  \gam \bet} &= {\sf com}_{   \bet  {\color{orange}\boldsymbol\delta}\gam }
^{   \bet\gam {\color{orange}\boldsymbol\delta} }
{\sf com}^{   \bet  {\color{orange}\boldsymbol\delta}\gam }
_{     {\color{orange}\boldsymbol\delta} \bet\gam }
{\sf com} ^{     {\color{orange}\boldsymbol\delta} \bet\gam }
_{     {\color{orange}\boldsymbol\delta}  \gam \bet}
  =
 {\sf com}_{  \gam \bet  {\color{orange}\boldsymbol\delta} }
^{   \bet\gam {\color{orange}\boldsymbol\delta} }
{\sf com}^{  \gam \bet  {\color{orange}\boldsymbol\delta} }
_{   \bet   {\color{orange}\boldsymbol\delta} \gam }
{\sf com} ^{   \bet   {\color{orange}\boldsymbol\delta} \gam }
_{     {\color{orange}\boldsymbol\delta}  \gam \bet}.
\end{align}
We suppress mention of crossing which  can be undone using the commutativity KLR relations in what follows.

Consider the former product in \cref{popopopopo1}.  
For $1\leq q \leq   b_{\color{orange}\boldsymbol\delta} $  the strand connecting the $\SSTP_{\empb \al\bet\al{\color{orange}\boldsymbol\delta}}^{-1}( q  ,\eps_{j})$ 
 and $\SSTP^{-1}_{{   {\color{orange}\boldsymbol\delta}\emp \bet \al\bet }}( q  ,\eps_{j})$ northern and southern vertices double-crosses the strands connecting each of the 
 $\SSTP_{\empb \al\bet\al{\color{orange}\boldsymbol\delta}}^{-1}(b_{\bet }+p,\eps_{j+1})$ 
 and $\SSTP^{-1}_{{   {\color{orange}\boldsymbol\delta}\emp \bet \al\bet }}(b_{\bet } +p,\eps_{j+1})$
northern and southern vertices for $p=q-1, q , q+1$.   
Now consider the latter product of \cref{popopopopo1}. 
For $1\leq q \leq   b_{\color{orange}\boldsymbol\delta} $  the strand connecting the $\SSTP_{\empb \al\bet\al{\color{orange}\boldsymbol\delta}}^{-1}(b_{\al\bet\al}+q ,\eps_{j})$ 
 and $\SSTP^{-1}_{{   {\color{orange}\boldsymbol\delta}\emp \bet \al\bet }}(b_{\al\bet\al}+q ,\eps_{j})$ northern and southern vertices double-crosses the strands connecting each of the 
 $\SSTP_{\empb \al\bet\al{\color{orange}\boldsymbol\delta}}^{-1}(b_{\al\bet\al\bet}+p,\eps_{j+1})$ 
 and $\SSTP^{-1}_{{   {\color{orange}\boldsymbol\delta}\emp \bet \al\bet }}(b_{\al\bet\al\bet}+p,\eps_{j+1})$
northern and southern vertices for $p=q-1, q , q+1$.   
For each $1\leq q \leq   b_{\color{orange}\boldsymbol\delta} $ we can undo these crossings using \cref{adjust1}.

Consider the former product in \cref{popopopopo2}.  
For $1\leq q \leq   \min\{b_\bet,b_{\color{orange}\boldsymbol\delta}\} $  the strand connecting the $\SSTP_{\bet\gam{\color{orange}\boldsymbol\delta}}^{-1}( q ,\eps_{k})$ 
 and $\SSTP^{-1}_{{\color{orange}\boldsymbol\delta} \gam\bet}( q ,\eps_{k})$ northern and southern vertices double-crosses the strands connecting each of the 
 $\SSTP_{\bet\gam{\color{orange}\boldsymbol\delta}}^{-1}( b_\gam+p,\eps_{k+1})$ 
 and $\SSTP^{-1}_{{\color{orange}\boldsymbol\delta} \gam\bet}( b_\gam+p,\eps_{k+1})$
northern and southern vertices for $p=q-1, q , q+1$.   
Now consider
 the latter product in \cref{popopopopo2}.  
For $0\leq q <   \min\{b_\bet,b_{\color{orange}\boldsymbol\delta}\} $  
the strand connecting the 
$\SSTP_{\bet\gam{\color{orange}\boldsymbol\delta}}^{-1}(b_{\bet\gam{\color{orange}\boldsymbol\delta}}-q ,\eps_{k})$ 
 and
  $\SSTP^{-1}_{{\color{orange}\boldsymbol\delta} \gam\bet}(b_{\bet\gam{\color{orange}\boldsymbol\delta}}-q   ,\eps_{k})$
   northern and southern vertices double-crosses the strands connecting each of the 
 $\SSTP_{\bet\gam{\color{orange}\boldsymbol\delta}}^{-1}(b_{\bet\gam\gam{\color{orange}\boldsymbol\delta}} -p,\eps_{k+1})$ 
 and $\SSTP^{-1}_{{\color{orange}\boldsymbol\delta} \gam\bet}(b_{\bet\gam\gam{\color{orange}\boldsymbol\delta}} -p,\eps_{k+1})$
northern and southern vertices for $p=q+1, q , q-1$.   
For each $0\leq q <   \min\{b_\bet,b_{\color{orange}\boldsymbol\delta}\} $ we can undo these crossings using \cref{adjust1}.  
  
Thus we obtain the desired equalities and  the image of relation \ref{S7} holds.

      \subsection{The fork and commutator  }   \label{fokcom}
Let $\gam,\bet\in\Pi$ label  two commuting   reflections, 
  we   now verify the middle  relation depicted  in \ref{rel6}, namely that 
\begin{align*}
    \Upsilon^{\SSTP_{\bet \empg  \gam   }}_{\SSTP_\gam \otimes \SSTP_\gam^\flat \otimes \SSTP_{\bet}}		&=
   (e_{\SSTP_{\bet}} \otimes {\sf  fork}^{\empg\gam }_{\gam\gam} )
(   {\sf com}_{\gam \bet}^{\bet\gam}\otimes e_{\SSTP_\gam})
(  e_{\SSTP_\gam}\otimes  {\sf com}_{\gam\bet}^{\bet\gam} )
\\ &=  ({\sf adj}^{\bet\empg  }_{\empg \bet  }\otimes e_{\SSTP_\gam}) 
 (  e_{\SSTP_\empg}\otimes  {\sf com}_{\gam \bet}^{\bet \gam} )
     ( {\sf fork}^{\empg\gam}_{\gam\gam} \otimes e_{\SSTP_{\bet}} ) 
     \end{align*}
as   both products produce       minimal, step-preserving, and residue commutative elements (after   undoing any
   double-crossings of non-adjacent residue using the commutativity relations).

\subsection{Naturality of adjustment}
For each generator, we must check the corresponding  adjustment naturality relation pictured in  \cref{nature2,nature}.  For the unique one-sided naturality relation,
 $({\sf spot}_\al^\emp \otimes e_{\SSTP_\emp}){\sf adj}^{\al\emp}_{\emp\al}= 
     e_{\SSTP_\emp } \otimes {\sf spot}_\al^\emp $, this follows by a generalisation of the proof of \cref{differentforks}.  The remaining relations all follow from \cref{adjust1}.

   \subsection{Cyclicity}   
   Given $\al,\bet\in \Pi$ labelling a pair of non-commuting reflections, 
   we now verify relation \ref{rel4}, namely that 
 \begin{equation}\label{cycliticityt}
\Psi \left(\; \begin{minipage}{3.5cm}
\end{minipage}\right). \end{equation}
The lefthand-side of \cref{cycliticityt} is equal to  
$$
\big(e_{\SSTP_{\al \emp \bet\al}  }
\otimes
 ({\sf spot}^\empb_\bet 
\otimes 
e_{\SSTP_{\empb}}
)
{\sf fork}^{\bet\empb}_{\bet\bet}\big)
  {\sf hex}^{ \al\emp\bet\al\bet\bet} _{ \al\empb\al\bet\al\bet } 
 ((  {\sf adj}^{\al  \empb \al  }_{ \empb \al \al  }
 (e_{\SSTP_{\empb  }}\otimes ( {\sf fork}^{\al\al}_{\emp\al} (e_{\SSTP_\emp}\otimes {\sf spot}^\al_\emp))))\otimes e_{\SSTP_{\bet\al\bet}} )
$$ which   is minimal and  step-preserving 
and so   is equal to 
  $
  \Upsilon ^
 {\al   \emp\bet\al\bet\empb }_
  {\empb\emp\emp \bet\al\bet} $  
(which is independent of the choice of reduced expression by simply re-indexing the  proof of 
\cref{analrsusususus}).  
The righthand-side  of \cref{cycliticityt} is equal to 
\begin{equation}\label{cylicwork}
{\sf adj}^{\al   \emp   \bet\al \empb\empb
}
_{ \emp \empb  \al\bet\al\empb}
\big(e_{ \SSTP_{ 
 \emp }}\otimes {\sf hex}_{\emp \bet\al\bet  }^{\empb \al\bet\al }
 \otimes e_{\SSTP_ {\empb }}\big)
(e_{\SSTP_{\emp\emp }}\otimes 
 {\sf adj}^{\bet\al\bet \empb}_{  \empb \bet\al\bet }).
  \end{equation}
It will suffice to show that 
\begin{equation}\label{suff}(   {\sf hex}_{ \bet\al\bet} 
 \otimes e_{\SSTP_ {\emptyset }} )
 {\sf adj}^{\bet\al\bet \emptyset}_{ \emptyset\bet\al\bet } =
\Upsilon	^{{\sf H}_{b_{\al\bet},\bet\al\bet}\otimes \SSTP_\emptyset}_ {\emptyset\empb\al\bet\al}		\end{equation}
as $b_\bet$ applications of this will simplify \cref{cylicwork} so that it   is minimal and step-preserving.  
The lefthand-side of \cref{suff}
contains an $r$-strand from ${\sf H}_{q,\al\bet\al}^{-1}(q+1,\eps_{i+1})$ 
to $\SSTP^{-1}_{\emptyset\empb\al\bet\al}(q+1,\eps_{i+1})$ which double-crosses the strands connecting the  top and bottom vertices 
$${\sf H}_{q,\al\bet\al}^{-1}(b_\al + q ,\eps_{i}) \quad 
 {\sf H}_{q,\al\bet\al}^{-1}(b_\al + q+1 ,\eps_{i})\quad 
 {\sf H}_{q,\al\bet\al}^{-1}(b_\al + q +2,\eps_{i}) $$   
$$\SSTP^{-1}_{\emptyset\empb\al\bet\al}(b_\al + q ,\eps_{i})\quad  \SSTP^{-1}_{\emptyset\empb\al\bet\al}(b_\al + q+1 ,\eps_{i})\quad
\quad
 \SSTP^{-1}_{\emptyset\empb\al\bet\al}(b_\al + q+2 ,\eps_{i}),$$
   respectively.    
We undo these double-crossings as in the proof of  \cref{adjust1} to obtain $\Upsilon ^
 {\al   \emp\bet\al\bet\empb }_
  {\emp\emp\empb\bet\al\bet}  $.  

\subsection{Some results concerning doubly-spotted Soergel diagrams} 
  The remainder of this section is dedicated to proving  results involving the  ``doubly-spotted" Soergel  diagrams.  These proofs are of a different flavour to the ``timeline" proofs considered above.  
We shall see that each Soergel spot diagram roughly corresponds to  ``half" of a KLR dotted diagram.  This idea is easiest to see through its manifestation in the grading (Soergel spots have degree 1, whereas KLR dots have degree 2).  We have that 
\begin{align}\notag
\Psi \!\left(\begin{minipage}{0.7cm} \begin{tikzpicture}[scale=1.3]
    \clip (0,-0.25-1/8) rectangle (0.5,0.75); 
{\clip (0,-0.25-1/8) rectangle (0.5,0.75); \foreach \i in {0,1,2,...,40}
  {
    \path  (-1,-0.75)++(0:0.125*\i cm)  coordinate (a\i);
        \path  (-1,-0.75)++(90:0.125*\i cm)  coordinate (b\i);
     \draw [densely dotted](a\i)--++(90:3);
        \draw [densely dotted](b\i)--++(0:10);
 }
} 
  \draw[ black, line width=0.06cm] 
  (0,-0.25-1/8) --++(0:0.5); 
  \draw[ black, line width=0.06cm] 
  (0,0.75) --++(0:0.5);
 \draw[magenta,line width=0.12cm](0.25,0)--(0.25,0.75/2);
\fill[magenta] (0.25,1/2)  circle (4pt); 
\fill[magenta] (0.25,-1/8)  circle (4pt); 
  \end{tikzpicture}
 \end{minipage}  \right) 
  & =     e_{\SSTP_{\emp}} \Bigg(\!\prod_{   \exx> q \geq 0 } 
 \psi
 ^{ q{\aatch}   + \isit }_{   \exx {\aatch}-\exx + q+1	 }   \Bigg) e_{\SSTP_{\al}}
 \Bigg(\!\prod_{ 0\leq q <\exx} 
\psi
 _{ q{\aatch}   + \isit } ^ {   \exx {\aatch}-\exx + q+1	 }  \Bigg)     e_{\SSTP_{\emp}} 
 \\ 
 \label{tod-tod}
 &=     e_{\SSTP_{\emp}}      \big(
  y_{\exx {\aatch} - {\aatch}+\emptyset(i+1)}-y_{\isit   }
       \big) 
   e_{\SSTP_{\emp}} 
\intertext{by  relation \ref{rel1.10}; this is easily seen from the fact that the only crossings of non-zero degree 
 are a double-crossing of strands which begin and end at the   
 $\SSTP_{\emp}^{-1}(b_\al,\eps_{i+1})=(\exx {\aatch} -{\aatch} +\emptyset (i+1))$  and $
 \SSTP_{\emp}^{-1}(1,\eps_{i})={\isit  }$  points on the   \north    and   \south       edges of the diagram (and application of case 3 of relation \ref{rel1.10}).  
Arguing similarly, one has that }\notag
 \Psi \left(\;\begin{minipage}{0.7cm} \begin{tikzpicture}[scale=1.3]
   \clip (0,-0.25-1/8) rectangle (0.5,0.75); 
{\clip (0,-0.25-1/8) rectangle (0.5,0.75); \foreach \i in {0,1,2,...,40}
  {
    \path  (-1,-0.75)++(0:0.125*\i cm)  coordinate (a\i);
        \path  (-1,-0.75)++(90:0.125*\i cm)  coordinate (b\i);
     \draw [densely dotted](a\i)--++(90:3);
        \draw [densely dotted](b\i)--++(0:10);
 }
}
   \draw[ magenta, line width=0.06cm] 
  (0,-0.25-1/8) --++(0:0.5); 
  \draw[ magenta, line width=0.06cm] 
  (0,0.75) --++(0:0.5);
 \draw[magenta,line width=0.12cm](0.25,1)--(0.25,0.75/2);
  \draw[magenta,line width=0.12cm](0.25,-1)--(0.25,-0/2);
\fill[magenta] (0.25,1/2-1/16)  circle (4pt); 
\fill[magenta] (0.25,-1/16)  circle (4pt); 
  \end{tikzpicture}
 \end{minipage}  \right)\!&=         e_{\SSTP_{\al}}  
  \Bigg(\prod_{ 0\leq q <\exx} 
\psi
 _{ q{\aatch}   + \isit } ^ {   \exx {\aatch}-\exx + q+1	 }  \Bigg)  
e_{\SSTP_{\emp}}
 \Bigg(\prod_{   \exx> q \geq 0 } 
 \psi
 ^{ q{\aatch}   + \isit }_{   \exx {\aatch}-\exx + q+1	 }   \Bigg)         e_{\SSTP_{\al}}  
 \\   \label{dot-dot}
&  =
          e_{\SSTP_{\al}}     \big( y_{\exx {\aatch} -\exx-{\aatch} +1 + \isitone}-y_{ \exx {\aatch} - \exx + 1  }     \big)
        e_{\SSTP_{\al}}  .
\end{align}
 \begin{prop}\label{jumpydots} Let  $\al=\varepsilon_i-\varepsilon_{i+1},\gam=\eps_k-\eps_{k+1} \in \Pi$ with $b_\al>1$ and   $0\leq q < b_\gam$.   We have that 
\begin{align}\label{jumpy8}
 y_{   \emptyset(i+1)}
e_{\SSTP_{\emptyset\emptyset} }=
 y_{ {\aatch} +    \emptyset(i+1)} e_{\SSTP_{\emptyset\emptyset} }
\qquad &
 y_{   \emptyset(i)  }	 
e_{\SSTP_{\emptyset\emptyset} }=
 y_{  {\aatch} +   \emptyset(i)  }	 
e_{\SSTP_{\emptyset\emptyset} }
\\
\label{jumpy1}
y_{{\aatch} + \gam(  i+1)   } e_{\SSTP_{\emptyset \gam}} = 
 y_{ \emptyset(  i+1)}  e_{\SSTP_{\emptyset \gam}}
\qquad &
y_{{\aatch} + \gam(  i)   } e_{\SSTP_{\emptyset \gam}} = 
 y_{ \emptyset(  i)}  e_{\SSTP_{\emptyset \gam}}
\\
\label{jumpy2} 
   y_{q({\aatch}-1) +  \gam(i+1)  } e_{\SSTP_{\gam}}= 
y_{(q+1)({\aatch} -1) + \gam(i+1)   } e_{\SSTP_{\gam}}  
\qquad
   &y_{q({\aatch}-1) +  \gam(i)  } e_{\SSTP_{\gam}}= 
y_{(q+1)({\aatch} -1) + \gam(i)   } e_{\SSTP_{\gam}}  
 \end{align}   
whenever the indices are defined (cross reference \cref{whatabouti}).   \end{prop}
\begin{proof}
We prove both cases of \cref{jumpy8}, the other pairs of cases are similar.  Our assumption that $b_\al>1$ implies that the residues of the $i$th and 
$(i+1)$th strands are {\em non-adjacent} and similarly that the  $({\aatch}+ \emptyset(i))$th and 
$({\aatch}+ \emptyset(i+1))$th strands are {\em non-adjacent}    (this is not true if $b_\al=1$).  Therefore we have that 
$$ 
0=
\psi^{   i}_{  {\aatch} +i}
e_{ \SSTP_{\emptyset \emptyset}}
\psi_{   i}^{  {\aatch} +i}
= (y_{   i}-y_{   {\aatch} +i})e_{\SSTP_{\emptyset \emptyset }},
\quad 0=
\psi^{  {\aatch} + \emptyset( i+1)}_{ \emptyset( i+1)}
e_{ \SSTP_{\emptyset \emptyset}}
\psi_{  {\aatch} + \emptyset( i+1)}^{ \emptyset( i+1)}
= (y_{   i}-y_{   {\aatch} +i})e_{\SSTP_{\emptyset \emptyset }}
$$
where  in both cases,  the first and second equalities follow from 
 \cref{resconsider} and  the final case of relation \ref{rel1.10}.  
 \end{proof}

   \begin{prop}\label{jumpydots1}
Let  $\al=\varepsilon_i-\varepsilon_{i+1},\gam=\eps_k-\eps_{k+1} \in \Pi$ with $b_\al=1$ and $0\leq q <b_\gam$.  We have that 
\begin{align*}
(y_{   i   }	-y_{   i+1 })	
e_{\SSTP_{\emptyset\emptyset} }&=
(y_{  {\aatch} +   i   }	-y_{ {\aatch} +    i+1  })
e_{\SSTP_{\emptyset\emptyset} }
\\ 
( y_{{\aatch} + \gam(  i+1)   } -y_{{\aatch} + \gam(  i)   })e_{\SSTP_{\emptyset \gam}}
&= 
( y_{  i+1 } - y_{   i 		}  ) e_{\SSTP_{\emptyset \gam}}
\\
(y_{ q({\aatch}-1)+   \gam(i)  }	-y_{ q({\aatch}-1)+ \gam(i+1)})	
e_{\SSTP_{\gam} }&=
(y_{ (q+1)({\aatch}-1)+   \gam(i)  }	-y_{ (q+1)({\aatch}-1)+    \gam(i+1)})
e_{\SSTP_{\gam} }
 \end{align*} 
whenever the indices are defined (cross reference \cref{whatabouti}).   \end{prop}

   \begin{proof}
   We prove the first equality as the other cases are similar.  
Since  $b_\al= 1$, we have that 
$\emptyset(i)=i$ and     $\emptyset(i+1)=i+1$  (in other words, $i\neq {\aatch}$) and are of adjacent  residue.  
We have that 
   \begin{align*}
 (y_{{\aatch}+i+1}-y_{{\aatch}+i })e_{\SSTP_{\emptyset\emptyset}}
 &= e_{\SSTP_{\emptyset\emptyset}} \psi _{{\aatch}+i}^{{\aatch}+i+1}   \psi ^{{\aatch}+i} _{{\aatch}+i+1}e_{\SSTP_{\emptyset\emptyset}}
   \\ 
   &=
(e_{\SSTP_{\emptyset\emptyset}}   \psi _{{\aatch}+i}^{{\aatch}+i+1}  )
\psi^{{\aatch}+i}_{i+2}
\psi_{{\aatch}+i}^{i+2}
( \psi _{{\aatch}+i+1}  ^{{\aatch}+i}
e_{\SSTP_{\emptyset\emptyset}})
\\
&=
(e_{\SSTP_{\emptyset\emptyset}}   \psi _{{\aatch}+i} ^{{\aatch}+i+1} 
\psi^{{\aatch}+i}_{i+2} )
(   \psi_{i+1}\psi_i	\psi_{i+1} 	+			\psi_i \psi_{i+1}\psi_i	)
(\psi_{{\aatch}+i}^{i+2}
  \psi^{{\aatch}+i} _{{\aatch}+i+1} 
e_{\SSTP_{\emptyset\emptyset}})\\
   &=
(e_{\SSTP_{\emptyset\emptyset}}   \psi ^{{\aatch}+i+1}  _{{\aatch}+i}
\psi^{{\aatch}+i}_{i+2} )
 \psi_{i}\psi_{i+1}	\psi_{i }	 
(\psi_{{\aatch}+i }^{i+2}
  \psi ^{{\aatch}+i} _{{\aatch}+i+1} 
e_{\SSTP_{\emptyset\emptyset}}) \\
&=
 (e_{\SSTP_{\emptyset\emptyset}}    \psi_{i}\psi ^{{\aatch}+i+1}  _{{\aatch}+i}
 )
 \psi^{{\aatch}+i}_{i+2}
 \psi_{i+1} \psi_{{\aatch}+i }^{i+2} 
( 
  \psi ^{{\aatch}+i} _{{\aatch}+i+1} \psi_{i} 
e_{\SSTP_{\emptyset\emptyset}}) \\
&=
 (e_{\SSTP_{\emptyset\emptyset}}    \psi_{i}\psi ^{{\aatch}+i+1}  _{{\aatch}+i}
 )
  \psi^{ i+1}_{{\aatch}+i-1}
 \psi_{{\aatch}+i-1} \psi_{i+1 }^{{\aatch}+i-1})
  \psi ^{{\aatch}+i} _{{\aatch}+i+1} \psi_{i} 
e_{\SSTP_{\emptyset\emptyset}}) \\
&=
 (e_{\SSTP_{\emptyset\emptyset}}    \psi_{i} 
 \psi^{ i+1}_{{\aatch}+i-1} )
  \psi ^{{\aatch}+i+1}  _{{\aatch}+i}
 \psi_{{\aatch}+i-1}   \psi ^{{\aatch}+i} _{{\aatch}+i+1} 
 ( \psi_{i+1 }^{{\aatch}+i-1} 
  \psi_{i} 
e_{\SSTP_{\emptyset\emptyset}}) \\
&=
 (e_{\SSTP_{\emptyset\emptyset}}    \psi_{i} 
 \psi^{ i+1}_{{\aatch}+i-1} )
(   1 +  
 \psi_{{\aatch}+i-1} \psi _{{\aatch}+i }   \psi_{{\aatch}+i-1} 
)( \psi_{i+1 }^{{\aatch}+i-1} 
  \psi_{i} 
e_{\SSTP_{\emptyset\emptyset}}) \\
&=
 (e_{\SSTP_{\emptyset\emptyset}}    \psi_{i} 
 \psi^{ i+1}_{{\aatch}+i-1} )
 ( \psi_{i+1 }^{{\aatch}+i-1} 
  \psi_{i} 
e_{\SSTP_{\emptyset\emptyset}}) \\
&=
 e_{\SSTP_{\emptyset\emptyset}}    \psi_{i} 
  \psi_{i} 
e_{\SSTP_{\emptyset\emptyset}}  \\
&=
 e_{\SSTP_{\emptyset\emptyset}}    (y_{i+1} - y_{i}) e_{\SSTP_{\emptyset\emptyset}}  
\end{align*}
  where the first equality holds by the third case of relation \ref{rel1.10}, 
  the second holds by the second case of 
    relation \ref{rel1.10} (the commuting version), 
  the third holds by case 2 of relation \ref{rel1.11},
  the fourth holds by \cref{resconsider}, and the fifth to  the seventh  by
   the second case of 
    relation \ref{rel1.10} (the commuting version),
    and the  eighth  by the first case of \ref{rel1.11}, 
    and the ninth by  \cref{resconsider}, the tenth by 
       the second case of 
    relation \ref{rel1.10} (the commuting version),
and the eleventh 
by the third case of relation \ref{rel1.10}.  
\end{proof}

        \subsection{The barbell    and commutator }
  \renewcommand{\bet}{{{{\color{ao(english)}\boldsymbol\gamma}}}}
 \renewcommand{\gam}{{{\color{cyan}\boldsymbol\beta}}}
\renewcommand{\empb}{{{{\color{ao(english)}\boldsymbol\clock}}}}
\renewcommand{\empg}{{{{\color{cyan}\boldsymbol\clock}}}}

For $\gam,\bet\in\Pi$ labelling two commuting reflections,  we check that 
   \begin{equation}
 \Psi \left(\;\begin{minipage}{1.2cm} 
\end{minipage}\right)  
   \end{equation}    In other words,   $$ 
(      {\sf spot}^{\empg}_\gam  {\sf spot}_{\empg}^\gam) \otimes e_{\SSTP_\bet}  
={\sf adj}^{\empg\bet}_{\bet\empg}(  e_{\SSTP_\bet}   \otimes         (      {\sf spot}^{\empg}_\gam  {\sf spot}_{\empg}^\gam)  )  
{\sf adj}_{\empg\bet}^ {\bet\empg}  . $$ 
         This relation is very simple to check.  
       We have that  
  \begin{align*}
   {\sf adj}^{\empg\bet}_{\bet\empg}(  e_{\SSTP_\bet}   \otimes         (      {\sf spot}^{\empg}_\gam  {\sf spot}_{\empg}^\gam))  
{\sf adj}_{\empg\bet}^ {\bet\empg}   
&=    
  {\sf adj}^{\empg\bet}_{\bet\empg} 
(          y_{b_{\bet\gam} {\aatch} - {\aatch}+\emptyset(\kay +1)}-y_{b_\bet {\aatch}+\kay    })
   {\sf adj}_{\empg\bet}^{\bet\empg}
e_{\SSTP_{\empg\bet}}
\\
&= (          y_{b_{\bet\gam} {\aatch} - {\aatch}+1-b_\bet +\bet(\kay +1)}-y_{b_\bet {\aatch}+\bet(\kay)    })
 {\sf adj}^{\empg\bet}_{\bet\empg} 
   {\sf adj}_{\empg\bet}^{\bet\empg}
e_{\SSTP_{\empg\bet}}
\\
&= (          y_{b_{\bet\gam} {\aatch} - {\aatch}+1-b_\bet+\bet(\kay +1)}-y_{b_\bet {\aatch}+\bet(\kay)    })
    e_{\SSTP_{\empg\bet}}\\
    &=   
 (          y_{b_\gam {\aatch} - {\aatch}+\emptyset(\kay +1)}-y_{\kay    })e_{\SSTP_{\empg\bet}}   
\intertext{     
      where the first equality follows from \cref{dot-dot}, the second equality follows from the commuting cases of 
       relations \ref{rel1.8} and  \ref{rel1.7}, 
    the third  equality 
       follows from \cref{adjust1},
 the fourth equality follows from applying   \cref{jumpydots,jumpydots1}. 
Again by \cref{dot-dot}, we have that }
           ( {\sf spot}^{\empg}_\gam 
  {\sf spot}_{\empg}^\gam )
      \otimes e_{\SSTP_\bet} 
    &  =
  (          y_{b_\gam {\aatch} - {\aatch}+\emptyset(\kay +1)}-y_{ \kay    })e_{\SSTP_{\empg\bet}}
\end{align*}         
        as required.

\color{black}
   
   \renewcommand{\empb}{{{{\color{cyan}\boldsymbol\clock}}}}
\renewcommand{\empg}{{{{\color{ao(english)}\boldsymbol\clock}}}}

   \renewcommand{\gam}{{{{\color{ao(english)}\boldsymbol\gamma}}}}
 \renewcommand{\bet}{{{\color{cyan}\boldsymbol\beta}}}

\subsection{The one   colour Demazure relation}    
We now verify \ref{rel2}, namely that 
 \begin{equation}\label{anotherlabeldemazure} 
 \Psi \left(\; \begin{minipage}{1.7cm}\; 
\end{minipage}\right)  .
 \end{equation}
for $\al \in \Pi$. In other words, we must check that    
\begin{align*}
  &
(
 {\sf spot}^\emp _\al{\sf spot}_\emp ^\al
 ) \otimes e_{\SSTP_\al }
 +
   {\sf adj}^{\emp \al }_{\al\emp  } 
 (
 e_{\SSTP_\al} \otimes   {\sf spot}^\emp _\al{\sf spot}_\emp ^\al 
 ) 
     {\sf adj}_{\emp \al }^{\al\emp  }  
= 
2
(
 e_{\SSTP_\emp }		\otimes 	{\sf spot}_\emp ^\al{\sf spot}^\emp _\al
 )  
 \end{align*}
     Substituting \cref{tod-tod,dot-dot}  into the  above, we must show that 
\begin{align}\label{subitallbackin}
\begin{split}& \quad \; e_{\SSTP_{\emp\al}} 
  \left(y_{\exx {\aatch} - {\aatch}+\emptyset(i+1)}-y_{\isit   }
+{\sf adj}^{\emp\al}_{ \al\emp}
 ( y_{b_{\al \al} {\aatch} -{\aatch} +  \emptyset(i+1) }-y_{\exx {\aatch} +\isit })  
 {\sf adj}_{\emp\al}^{ \al\emp} \right)
e_{\SSTP_{\emp\al}}  \\
&= 2 e_{\SSTP_{\emp\al}}( 
 y_{b_{\al \al} {\aatch} -\exx-{\aatch} +1 + \isitone}-y_{ b_{\al \al} {\aatch} - \exx + 1  }
)e_{\SSTP_{\emp\al}}.
\end{split}
\end{align}
This leads us to consider the effect of passing dots through the adjustment terms. 
 \begin{prop}\label{conclus0}
Let  $\al \in \Pi $.  We have that 
\begin{align} \label{Referbacj}
{\sf adj}^{\emp\al}_{\al  \emp} 
 y_{b_{\al } {\aatch} +\isit  }   
{\sf adj}_{\emp\al}^{\al  \emp} 
& 
 =
 y_{b_{\al\al} {\aatch} -b_\al+1 }   
 e_{\SSTP_{\emp\al}	} 
  \end{align}
 \begin{align}
 \label{Referbacj3}
 {\sf adj}^{\emp\al}_{\al  \emp}(b_\al-2 )  \mydots {\sf adj}^{\emp\al}_{\al  \emp}(0 )  
  y_{b_{\al\al} {\aatch} -{\aatch} +\emptyset(i+1) } 
   {\sf adj}^{\emp\al}_{\al  \emp}(0 )  
 \mydots {\sf adj}^{\emp\al}_{\al  \emp}(b_\al-2 ) 
 &=
y_{b_{\al\al} {\aatch} -{\aatch} +  \emptyset(i+1) } 
 e_{\SSTP_{\emp\al}}.
 \end{align}
  \end{prop}

 \begin{proof}
  By  the commuting case of relation \ref{rel1.7}, we have that the lefthand-sides of \cref{Referbacj,Referbacj3} are equal to 
$ y_{b_{\al\al} {\aatch} -b_\al+1 }   
{\sf adj}^{\emp\al}_{\al  \emp} 
{\sf adj}_{\emp\al}^{\al  \emp} 
 $  and $y_{b_{\al\al} {\aatch} -{\aatch} +  \emptyset(i+1) } 
 {\sf adj}^{\emp\al}_{\al  \emp}(b_\al-2 )  \mydots   {\sf adj}^{\emp\al}_{\al  \emp}(b_\al-2 ) $ respectively.  
The result then follows by \cref{adjust1}.   
 \end{proof}
 
In \cref{Referbacj3} we   pulled the dot through most of the adjustment term; in   \cref{equationproof} below, we pull the dot through the final adjustment term.  
  \Cref{equationproof2} has an almost identical proof  and so we record it here, for convenience.

  \begin{prop}\label{conclus1}
Let  $\al \in \Pi $.  We have that  
\begin{equation}\label{equationproof}
{\sf adj}^{\emptyset\al}_{\al  \emptyset}
y_{b_{\al } {\aatch}  +  \emptyset(i+1) }
{\sf adj}_{\emptyset\al}^{\al  \emptyset}
=
   \left(  y_i 
 +
  y_{b_{\al } {\aatch} - \exx+1 + \isitone  }
  -   y_{b_{\al } {\aatch}+ {\aatch} - \exx+1 }
   \right)e_{\SSTP_{\emptyset\al}}  
\end{equation}
\begin{equation}\label{equationproof2}
{\sf adj}^{\emptyset\al}_{\al  \emptyset}
y_{b_{\al } {\aatch} -b_\al +1 }
{\sf adj}_{\emptyset\al}^{\al  \emptyset}
=
    y_{b_{\al } {\aatch}+ {\aatch} - \exx+1 }
 e_{\SSTP_{\emptyset\al}} .  
\end{equation}

\end{prop} 

\begin{proof} 
 We first prove \cref{equationproof}.  
The dotted strand in the concatenated diagram on the left of \cref{equationproof} connects  the $i=\SSTP_{\emptyset\al}^{-1}(1,\eps_i)$  \north   and bottom vertices, by way of the $ b_\al {\aatch}+\emptyset (i+1) =
\SSTP_{\al\emptyset}^{-1}(1,\eps_{i+1})$   vertex in the centre of the diagram.  We suppose this dotted strand is of residue $r\in\ZZ/e\ZZ$,  say.  This dotted strand crosses a single strand of the same residue: namely, the strand connecting the $ \SSTP_{\emptyset\al}^{-1}(b_\al+1, \eps_{i+1} )$th vertices  on the   \north    and   \south       edges.  
By relation \ref{rel1.8}, we can pull the dot upwards along its strand and through this crossing at the expense of an error term.  We thus obtain 
\begin{equation} \label{onlylovecanbreak}
{\sf adj}^{\emptyset\al}_{\al  \emptyset}
y_{b_{\al } {\aatch}  +  \emptyset(i+1) }
{\sf adj}_{\emptyset\al}^{\al  \emptyset}
=
e_{\SSTP_{\al\emptyset}}\big(y_i  \psi_{\SSTP_{\al\emptyset}}^{\SSTP_{ \emptyset\al}}
  \psi^{\SSTP_{\al\emptyset}}_{\SSTP_{ \emptyset\al}}
\big)e_{\SSTP_{\al\emptyset}}+e_{\SSTP_{\al\emptyset}}\big(
\psi^{\SSTP_\emptyset\otimes {\sf S}_{0,\al} }
 _{{\sf S}_{1,\al} \otimes \SSTP_\emptyset}  
\psi^{\SSTS_{1,\al}  \otimes \SSTP_\emptyset}
_{\SSTS_{0,\al}  \otimes \SSTP_\emptyset}  
\psi^{\SSTP_{\al\emptyset}}_{\SSTP_{ \emptyset\al}}\big)e_{\SSTP_{\al\emptyset}} 
\end{equation}
(we note that ${\sf S}_{0,\al}=\SSTP_\al^\flat  $).  
The first term in \cref{onlylovecanbreak}  is equal to 
 $y_i e_{\SSTP_{\emptyset\al}} $   by  \cref{adjust1} (and this is equal to the leftmost term on the righthand-side of \cref{equationproof}).
 We now consider the latter term.  We label the top and bottom edges by $\northT=\SSTP_\emptyset\otimes \SSTP_\al^\flat $ and 
 $\southT= \SSTP_\emptyset\otimes \SSTP_\al$.    There is a unique crossing of  strands of the same residue in the diagram 
 $$e_{\SSTP_{\emptyset\al}}\big(\psi^{\SSTP_\emptyset\otimes {\sf S}_{0,\al}}
 _{{\sf S}_{1,\al}\otimes \SSTP_\emptyset}\circ 
\psi^{\SSTS_{1,\al} \otimes \SSTP_\emptyset}
_{\SSTS_{0,\al} \otimes \SSTP_\emptyset}\circ 
\psi^{\SSTP_{\al\emptyset}}_{\SSTP_{ \emptyset\al}}\big) e_{\SSTP_{\emptyset\al}} $$
 namely the $r$-strands connecting the 
 $ i = \SSTT^{-1}(1,\eps_i)
 $  and $ \southT ^{-1}(b_\al+1, \eps_{i+1}) $  vertices on the   \north    and   \south       edges of the diagram.  
 This crossing of strands is bi-passed on the left 
 by the $(r+1)$-strand connecting the $ \northT^{-1}(b_\al, \eps_{i+1}) =
 \southT^{-1}(b_\al, \eps_{i+1})  $    \north    and   \south       vertices.  
 We pull this $(r+1)$-strand to the right through the crossing $r$-strands using case 2 of relation \ref{rel1.11} (and the commuting relations).  We hence undo this crossing and obtain
 $$
e_{\SSTP_{\emptyset\al}}\big( \psi^{\SSTP_\emptyset\otimes {\sf S}_{0,\al}}
 _{{\sf S}_{1,\al}\otimes \SSTP_\emptyset}  
\psi^ {{\sf S}_{1,\al}\otimes \SSTP_\emptyset}_{\SSTP_\emptyset\otimes {\sf S}_{0,\al}}
 \big)e_{\SSTP_{\emptyset\al}}$$ 
(the other term depicted in \cref{rel1.11} is zero by \cref{resconsider}). Now, this diagram contains a double-crossing of the $r$-strand 
connecting the 
$( \SSTP_\emptyset\otimes \SSTP_\al^\flat)^{-1} (b_\al+1, \eps_{i+1}) $ 
   \north    and   \south       vertices and the $(r-1)$-strand connecting the 
  $( \SSTP_\emptyset\otimes \SSTP_\al^\flat)^{-1}( 2, \eps_{i })$   \north    and \south vertices.  
  We undo this double-crossing using case 4 of relation \ref{rel1.10} (and the commutativity relations) to obtain 
 \begin{equation}\label{yy} e_{\SSTP_{\emptyset\al}}    (y_{b_{\al } {\aatch} - \exx+1 + \isitone  }
  -   y_{b_{\al } {\aatch}+ {\aatch} - \exx+1 }
    )e_{\SSTP_{\emptyset\al}}  
\end{equation} and so \cref{equationproof} follows.  Regarding  the  enumeration above, we note that 
  $$  
   ( \SSTP_\emptyset\otimes \SSTP_\al^\flat)^{-1}(b_\al+1, \eps_{i+1}) = 
 {b_{\al } {\aatch} - \exx+1 + \isitone  }
$$
$$ ( \SSTP_\emptyset\otimes \SSTP_\al^\flat)^{-1} ( 2, \eps_{i }) = b_{\al } {\aatch}+ {\aatch} - \exx+1
. $$ 
 Now  we turn to   \cref{equationproof2}.  We push the KLR-dot upwards along its strand using  \ref{rel1.8} 
  to obtain 
\begin{equation} \label{onlylovecanbreak2}
e_{\SSTP_{\emptyset\al}}\big(
y_{b_{\al } {\aatch} - \exx+1 + \isitone  }   \psi_{\SSTP_{\al\emptyset}}^{\SSTP_{ \emptyset\al}}
  \psi^{\SSTP_{\al\emptyset}}_{\SSTP_{ \emptyset\al}}
\big)e_{\SSTP_{\emptyset\al}} -e_{\SSTP_{\emptyset\al}}\big(
\psi^{\SSTP_\emptyset\otimes {\sf S}_{0,\al}}
 _{{\sf S}_{1,\al}\otimes \SSTP_\emptyset}  
\psi^{\SSTS_{1,\al} \otimes \SSTP_\emptyset}
_{\SSTS_{0,\al} \otimes \SSTP_\emptyset}  \circ 
\psi^{\SSTP_{\al\emptyset}}_{\SSTP_{ \emptyset\al}}\big)e_{\SSTP_{\emptyset\al}}.  
\end{equation} 
The first term is equal to  $
y_{b_{\al } {\aatch} - \exx+1 + \isitone  } e_{\SSTP_{\emptyset\al}} $ (again this follows by  \cref{adjust1}).  The second term is  identical to the second term in \cref{onlylovecanbreak} and so is equal to \cref{yy} but with negative coefficient.  Thus we can rewrite \cref{onlylovecanbreak2} in the form 
 $$e_{\SSTP_{\emptyset\al}}  \left(y_{b_{\al } {\aatch} - \exx+1 + \isitone  }
 -  
 (y_{b_{\al } {\aatch} - \exx+1 + \isitone  }
  -   y_{b_{\al } {\aatch}+ {\aatch} - \exx+1 }
   ) \right )e_{\SSTP_{\emptyset\al}}  
$$
and \cref{equationproof2} follows.  
\end{proof}

We now gather together our conclusions from \cref{conclus0,conclus1} 
(shifting the indexing where necessary) in order to prove \cref{anotherlabeldemazure}.   
 We have that  $({\sf spot}^\emp_\al 
{\sf spot}_\emp^\al )\otimes e_{\SSTP_\al}$ is equal to 
$$ 
  e_{\SSTP_{\emp\al}}(y_{\exx {\aatch} - {\aatch}+\emptyset(i+1)}-y_{\isit   }) e_{\SSTP_{\emp\al}}
  $$
 and     ${\sf adj}_{\al\emp}^{\emp \al}(e_{\SSTP_\al}\otimes {\sf spot}^\emp_\al
{\sf spot}_\emp^\al ){\sf adj}^{\al\emp}_{\emp \al}$ is equal to 
 $$  
-  e_{\SSTP_{\emp\al}} y_{b_{\al\al} {\aatch} -\exx +1  }e_{\SSTP_{\emp\al}} +   e_{\SSTP_{\emp\al}}\big( 
      y_{\exx {\aatch} -{\aatch}+\isit  }  
+
  y_{b_{\al\al} {\aatch} - \exx -{\aatch} +1 + \isitone  }
  -   {\color{black}y_{b_{\al\al} {\aatch} - \exx+1}}
  \big)e_{\SSTP_{\emp\al}}
$$
By    \cref{jumpydots,jumpydots1}, we have that 
$$
y_{\exx {\aatch} - {\aatch}+\emptyset(i+1)}e_{\SSTP_{\emp\al}} =  y_{b_{\al\al} {\aatch} - \exx -{\aatch} +1 + \isitone  } e_{\SSTP_{\emp\al}}
$$
for $b_\al\geq 1$  and   by \cref{jumpydots}  we have that 
$$
y_{\isit   }e_{\SSTP_{\emp\al}}
=y_{\exx {\aatch} -{\aatch}+\isit  }  e_{\SSTP_{\emp\al}}
$$  
for   $b_\al>1$ (we note  that this latter  statement is vacuous if $b_\al=1$  as the subscripts are equal).   
The former pair of  terms sum up  and the latter cancel, so we  obtain 
$$ 
({\sf spot}^\emp_\al 
{\sf spot}_\emp^\al )\otimes e_{\SSTP_\al} + {\sf adj}_{\al\emp}^{\emp \al}(e_{\SSTP_\al}\otimes {\sf spot}^\emp_\al
{\sf spot}_\emp^\al ){\sf adj}^{\al\emp}_{\emp \al}
=
2   e_{\SSTP_{\emp\al}}(y_{\exx {\aatch} - {\aatch}+\emptyset(i+1)} 
  -   {\color{black}y_{b_{\al\al} {\aatch} - \exx+1}}
)   e_{\SSTP_{\emp\al}} 
$$  
Hence  \cref{subitallbackin} holds by a further application of \cref{jumpydots,jumpydots1}.


   \subsection{Two colour Demazure}
For $\al,\bet\in \Pi$ labelling two non-commuting reflections,  we now verify relation \ref{rel5}, namely that 
  \begin{equation}\label{anotherlabeldemazure2} 
 \!\!\!\Psi \!\left(\;\!\begin{minipage}{1.9cm}\; 
\end{minipage}
\right) \end{equation} 
We assume that  the rank of $\Phi$ is at least 2.  The reader is invited to check the rank $1$ case separately (here the scalar 2 appears due to certain coincidences in the arithmetic).

  \begin{prop}\label{anotherdotjumper}
  Let $\al\in\Pi$. If   $b_\al>1$,  we have that 
  $$
y_{b_\al {\aatch}+{\aatch} +\emptyset(i+1) }  e_{\SSTP_{\emptyset \al \emptyset } }
=
(y_i {  +} y_  {\emptyset(i+1)} {  -}
  y_{b_\al {\aatch}+{\aatch} - b_\al +1}
)
e_{\SSTP_{\emptyset \al \emptyset } }
\qquad  $$
 $$
  y_{b_\al \aatch  +{\aatch}+ i } e_{\SSTP_{\emptyset\al\emptyset} }
=y_{b_\al \aatch  +{\aatch} } e_{\SSTP_{\emptyset\al\emptyset} }= 
y_{b_\al \aatch +{\aatch}-b_\al+1} e_{\SSTP_{\emptyset \al\emptyset} }
$$
  and if $b_\al=1$ we have that 
 $$
(  y_{b_\al \aatch  +{\aatch}+ i }  -y_{b_\al {\aatch}+{\aatch} +\emptyset(i+1) } ) e_{\SSTP_{\emptyset \al \emptyset } }
=
 ( 2y_{b_\al \aatch +{\aatch}-b_\al+1}
-y_i - y_  {\emptyset(i+1)} )e_{\SSTP_{\emptyset \al \emptyset } }.  
 $$

  \end{prop}

   \begin{proof}
We check the $b_\al>1$ case as the other is similar.  
   The second equality follows  as in the proof of \cref{jumpydots}.    
We now consider the first equality.  We momentarily drop the prefix $\SSTP_\emptyset$ to the path 
   $\SSTP_{\emptyset\al\emptyset}$ for the sake of more manageable indices.  
 Since $b_\al>1$ we can pull the vertical strand connecting the 
$b_\al {\aatch} +\emptyset(i+1)$ top and bottom vertices leftwards until we reach a strand of adjacent residue (namely the $(b_\al {\aatch} -b_\al +2)$th strand)  as follows
$$   e_{\SSTP_{\al \emptyset } }
 =
    e_{\SSTP_{\al \emptyset } }
   \psi^{b_\al {\aatch} +\emptyset(i+1)}_{\color{black} b_\al {\aatch} -b_\al +3} 
      \psi_{b_\al {\aatch} +\emptyset(i+1)}^{\color{black} b_\al {\aatch} -b_\al +3} 
      e_{\SSTP_{\al \emptyset } }
 $$
we can rewrite the   centre of the diagram which using the braid relation as follows,
\begin{align*}
    e_{\SSTP_{\al \emptyset } }
   \psi^{b_\al {\aatch} +\emptyset(i+1)}_{\color{black} b_\al {\aatch} -b_\al +3} 
   & (	\psi_{\color{black} b_\al {\aatch} -b_\al +2}		   \psi_{\color{black} b_\al {\aatch} -b_\al +1}	\psi_{\color{black} b_\al {\aatch} -b_\al +2}-
   \\
  & \qquad	\psi_{\color{black} b_\al {\aatch} -b_\al +1}	\psi_{\color{black} b_\al {\aatch} -b_\al +2}		
   \psi_{\color{black} b_\al {\aatch} -b_\al +1}	
   		)
      \psi_{b_\al {\aatch} +\emptyset(i+1)}^{\color{black} b_\al {\aatch} -b_\al +3} 
      e_{\SSTP_{\al \emptyset } }  
\end{align*}
where the   latter term is zero   by \cref{resconsider} and so this simplifies  to 
$$ 
   e_{\SSTP_{\al \emptyset } }
   \psi^{b_\al {\aatch} +\emptyset(i+1)}_{\color{black} b_\al {\aatch} -b_\al +2} 
 (	 
 	   \psi_{\color{black} b_\al {\aatch} -b_\al +1} 
   		)	
      \psi_{b_\al {\aatch} +\emptyset(i+1)}^{\color{black} b_\al {\aatch} -b_\al +2} 
      e_{\SSTP_{\al \emptyset } } 
$$
now we use the non-commuting version of relation  \ref{rel1.7}  
 together with case 1 of relation \ref{rel1.10} to rewrite the middlemost crossing as a double-crossing with a KLR-dot,
$$-   e_{\SSTP_{\al \emptyset } }
   \psi^{b_\al {\aatch} +\emptyset(i+1)}_{\color{black} b_\al {\aatch} -b_\al +2} 
   (	 
 	   \psi_{\color{black} b_\al {\aatch} -b_\al +1}	 		 
	    	y_{\color{black} b_\al {\aatch} -b_\al +1}   \psi_{\color{black} b_\al {\aatch} -b_\al +1}	 
   		)
      \psi_{b_\al {\aatch} +\emptyset(i+1)}^{\color{black} b_\al {\aatch} -b_\al +2} 
      e_{\SSTP_{\al \emptyset } }      , 
$$
  we   pull the dotted strand leftwards   through the next strand of adjacent residue
 (namely the $((b_\al-1)( {\aatch}-1)  +\al(i+1))$th strand) 
   using the  commutativity relations and case 4 of relation \ref{rel1.10}   to obtain 
\begin{align*}
  e_{\SSTP_{\al \emptyset } }
   \psi^{b_\al {\aatch} +\emptyset(i+1)}_{{\color{black} b_\al {\aatch} -b_\al +2}}&
  		 	(   
				y_{(b_\al-1)( {\aatch}-1)  +\al(i+1)}+
				\\
				&\qquad   
	\psi^{{\color{black} b_\al {\aatch} -b_\al +2}}
 _{{(b_\al-1)( {\aatch}-1)  +\al(i+1)} }
	\psi_{{\color{black} b_\al {\aatch} -b_\al +2}}
 ^{{(b_\al-1)( {\aatch}-1)  +\al(i+1)} }		
   	) 		
      \psi_{b_\al {\aatch} +\emptyset(i+1)}^{{\color{black} b_\al {\aatch} -b_\al +2}}
      e_{\SSTP_{\al \emptyset } }   
\end{align*}
where the first summand is  zero by case 1 of relation \ref{rel1.10} and
the latter term is equal to 
$$
 	    e_{\SSTP_{\al \emptyset } }
   \psi^{b_\al {\aatch} +\emptyset(i+1)} 
 _{{(b_\al-1)( {\aatch}-1)  +\al(i+1)} }
	\psi 
  ^{{(b_\al-1)( {\aatch}-1)  +\al(i+1)} }		
_{b_\al {\aatch} +\emptyset(i+1)}
       e_{\SSTP_{\al \emptyset } }. $$
       Now we concatenate on the left by $\SSTP_\emptyset$ and then multiply by
       $y_{b_\al {\aatch}+{\aatch} +\emptyset(i+1) }  $  to obtain
       \begin{align}\label{ghghghghgdddddd}y_{b_\al {\aatch}+{\aatch} +\emptyset(i+1) }  
  e_{\SSTP_{\emptyset \al \emptyset } }
&=   y_{b_\al {\aatch}+{\aatch} +\emptyset(i+1) }  
     e_{\SSTP_{\emptyset \al \emptyset } }
   \psi^{b_\al {\aatch} +{\aatch} +\emptyset(i+1)} 
 _{{b_\al {\aatch} -b_\al  {  +1} +\al(i+1)} }
	\psi 
  ^{{b_\al {\aatch} -b_\al  {  +1} +\al(i+1)} }		
_{b_\al {\aatch}+{\aatch}  +\emptyset(i+1)}
       e_{\SSTP_{\emptyset \al \emptyset } }  
       \end{align}
which by relation  \ref{rel1.10} 
 is equal to 
\begin{align*}
    e_{\SSTP_{\emptyset\al \emptyset } }\big(
  & \psi^{b_\al {\aatch} +{\aatch} +\emptyset(i+1)} 
 _{{b_\al {\aatch} -b_\al  {  +1} +\al(i+1)} }
 y_{{b_\al {\aatch} -b_\al  {  +1} +\al(i+1)} }
+\\
&\qquad
   \psi^{b_\al {\aatch}+{\aatch}  
   +\emptyset(i+1)} 
   _{b_\al {\aatch}+{\aatch}  -b_\al +2}
 \psi^{b_\al {\aatch}+{\aatch}  -b_\al +1} 
 _{{b_\al {\aatch} -b_\al  {  +1} +\al(i+1)} }
 \big) 	\psi 
  ^{{b_\al {\aatch} -b_\al  {  +1} +\al(i+1)} }		
_{b_\al {\aatch} +{\aatch} +\emptyset(i+1)}
      e_{\SSTP_{\emptyset\al \emptyset } }  .  
\end{align*}
We consider the first term in the sum  first.  
By the commuting relations,  this term is equal to 
\begin{align*}
 e_{\SSTP_{\emptyset \al \emptyset } }\big(
   \psi^{b_\al {\aatch} +{\aatch} +\emptyset(i+1)} 
 _{{{\aatch} +i} }
 y_{{{\aatch} + i} }
\psi 
  ^{{		{\aatch} +i		} }		
_{b_\al {\aatch} +{\aatch} +\emptyset(i+1)}\big )
      e_{\SSTP_{\emptyset \al \emptyset } } 
\intertext{and by \cref{jumpydots} this is equal to }
  e_{\SSTP_{\emptyset \al \emptyset } }\big(
   \psi^{b_\al {\aatch} +{\aatch} +\emptyset(i+1)} 
 _{{{\aatch} +i} }
 y_{{  i } }
\psi 
  ^{{		{\aatch} +\emptyset(i+1)		} }		
_{b_\al {\aatch} +{\aatch} +\emptyset(i+1)}\big )
      e_{\SSTP_{\emptyset \al \emptyset } } 
\end{align*} and now, having  moved this KLR-dot a total of  ${\aatch}$ strands leftward,
 we can  apply the commutativity relations again  to obtain
  \begin{align}\label{pppppq} 
  y_{ i }
   e_{\SSTP_{\emptyset \al \emptyset } }
   \big(
   \psi^{b_\al {\aatch} +{\aatch} +\emptyset(i+1)} 
 _{{b_\al {\aatch} -b_\al  {  +1} +\al(i+1)} }		
\psi 
  ^{{b_\al {\aatch} -b_\al  {  +1} +\al(i+1)} }		
_{b_\al {\aatch} +{\aatch} +\emptyset(i+1)}
      e_{\SSTP_{\emptyset \al \emptyset } } \big )
=  y_{ i }e_{\SSTP_{\emptyset \al \emptyset } }
\end{align} 
where the final equality follows  by \cref{ghghghghgdddddd}.  
We now turn to the second term in the above sum, namely 
$$     e_{\SSTP_{\emptyset\al \emptyset } }   \psi^{b_\al {\aatch}+{\aatch}  +\emptyset(i+1)} 
   _{b_\al {\aatch}+{\aatch}  -b_\al +2}
 \psi^{b_\al {\aatch}+{\aatch}  -b_\al +1} 
 _{{b_\al {\aatch} -b_\al  {  +1} +\al(i+1)} }
 	\psi 
  ^{{b_\al {\aatch} -b_\al  {  +1} +\al(i+1)} }		
_{b_\al {\aatch} +{\aatch} +\emptyset(i+1)}
      e_{\SSTP_{\emptyset\al \emptyset } }    . $$
This has a crossing of like-labelled strands (of residue $r\in \ZZ/e\ZZ$)  connecting the
 $(b_\al {\aatch} + \emptyset(i+1))$th and
  $({ \exx {\aatch}    -b_\al   + 1 })$th   \north    and   \south       vertices.  This crossing is bi-passed on the right by the $(r-1)$-strand connecting the 
  $(b_\al {\aatch}-b_\al+2)$th   \north    and   \south       vertices.  
  We undo this braid using  case 1 of relation  \ref{rel1.11} to obtain 
  $$
 e_{\SSTP_{\emptyset\al \emptyset } } 
   (   \psi^{b_\al {\aatch}+{\aatch}  +\emptyset(i+1)} 
   _{b_\al {\aatch}+{\aatch}  -b_\al +2}
 \psi^{b_\al {\aatch}+{\aatch}  -b_\al +1} 
 _{{b_\al {\aatch} -b_\al  {  +1} +\al(i+1)} })
( \psi_{b_\al {\aatch}+{\aatch}  -b_\al +1} 
 ^{{b_\al {\aatch} -b_\al  {  +1} +\al(i+1)} }
  \psi_{b_\al {\aatch}+{\aatch}  +\emptyset(i+1)} 
   ^{b_\al {\aatch}+{\aatch}  -b_\al +2})
      e_{\SSTP_{ \emptyset  \al \emptyset } } 
  $$  
  where the other term in relation \ref{rel1.11} is zero by \cref{resconsider}.  
This diagram contains a single double-crossing of adjacent residues, which we undo using case 4 of  relation \ref{rel1.10} (and we undo all  the other crossings using the commutativity relation)
 to obtain 
\begin{align}\label{pppppq2} 
 e_{\SSTP_{\emptyset  \al \emptyset } }  (
y_{b_\al \aatch -b_\al +1+\al(i+1)} 
- y_{b_\al {\aatch} - b_\al +1}) e_{\SSTP_{  \al \emptyset } } 
=
e_{\SSTP_{  \al \emptyset } }  (
y_{\emptyset(i+1)}
- y_{b_\al {\aatch} - b_\al +1}) e_{\SSTP_{  \emptyset \al \emptyset } } 
\end{align}
where the final equality follows by  \cref{jumpydots}.   
  The result follows by summing over \cref{pppppq,pppppq2}.  
   \end{proof}

  \begin{prop}\label{conclus3}
Let  $\al=\eps_i-\eps_{i+1}, \bet=\eps_{i+1}-\eps_{i+2} \in \Pi $.  We have that   \begin{align*}
  & \; ({\sf spot}^\empb_\bet
{\sf spot}_\empb^\bet)\otimes e_{\SSTP_{\al\emp}} 
-
  ({\sf adj}^{\empb\al }_{\al\empb }\otimes e_{\SSTP_\emp})
(e_{\SSTP_\al}\otimes 
{\sf spot}^\empb_\bet
{\sf spot}_\empb^\bet\otimes  e_{\SSTP_{\emp }})
({\sf adj}_{\empb\al }^{\al\empb }\otimes e_{\SSTP_\emp})
  \\ 
 =&\; e_{\SSTP_\al}( 
 y_{ i }-  y_{b_{\al\bet} {\aatch}   
 }  )	e_{\SSTP_\al}
  \\
 =&\; e_{\SSTP_{\empb\al} }\otimes 
(
 {\sf spot}^\emp _\al{\sf spot}_\emp ^\al
 ) -e_{\SSTP_\empb }\otimes 
(
{\sf spot}_\emp ^\al{\sf spot}^\emp _\al
 )\otimes e_{\SSTP_\emp  }.
 \end{align*}

 \end{prop} 
\begin{proof}
 Substituting \cref{dot-dot,tod-tod} into the third line, we obtain
 $$   e_{\SSTP_{\empb\al\emp}}( 
y_{b_{\al\bet\al} {\aatch} - {\aatch}+\emptyset(i+1)}-y_{b_{\al\bet}{\aatch}+i  }
 - y_{b_{\al\bet} {\aatch} -\exx-{\aatch} +1 + \isitone}
+y_{ b_{\al\bet} {\aatch} - \exx + 1  } 
)e_{\SSTP_{\empb\al\emp}} .$$
We  apply    \cref{anotherdotjumper}
  to the first term in the sum and then cancelling terms using 
  \cref{jumpydots,jumpydots1}.  
 Substituting \cref{dot-dot,tod-tod} into the first line, we obtain
\begin{equation}\label{follows3}
e_{\SSTP_{\empb\al\emp}} 
  \left(y_{b_\bet {\aatch} - {\aatch}+\emptyset(i+2)}-y_{\emptyset(i+1)   }
-{\sf adj}^{\empb\al\empb}_{ \al\empb\empb}
 ( y_{b_{\al\bet } {\aatch} -{\aatch} +  \emptyset(i+2) }-y_{b_{\al} {\aatch} +\emptyset(i+1) })  
 {\sf adj}_{\empb\al\empb}^{ \al\empb\empb} \right)
e_{\SSTP_{\empb\al\emp}} .
\end{equation}
We have that 
\begin{equation}\label{follows2}
e_{\SSTP_{\empb\al\emp}}{\sf adj}^{\empb\al\empb}_{ \al\empb\empb}
  y_{b_{\al\bet } {\aatch} -{\aatch} +  \emptyset(i+2) }
  {\sf adj}_{\empb\al\empb}^{ \al\empb\empb}  
e_{\SSTP_{\empb\al\emp}}
= y_{b_\al \aatch  - b_\al - {\aatch} +1+\al(i+2)}
= y_{b_\bet {\aatch} - {\aatch}+\emptyset(i+2)} 
\end{equation}
where the first equality follows from 
 the commuting KLR-dot relation \ref{rel1.8} and the latter follows from 
  \cref{jumpydots,jumpydots1}.   
  We also have that 
\begin{align} \nonumber
  {\sf adj}^{\empb\al\empb}_{ \al\empb\empb} y_{b_{\al} {\aatch} +\emptyset(i+1) }  
 {\sf adj}_{\empb\al\empb}^{ \al\empb\empb}  
&=
 e_{\SSTP_{\empb\al\emp}}   \left(  y_{ b_\bet  {\aatch} -{\aatch}+i }
 +
  y_{b_{\al\bet } {\aatch} -{\aatch}- \exx+1 + \isitone  }
  -   y_{b_{\al \bet} {\aatch}  - \exx+1 }
   \right)e_{\SSTP_{\empb\al\emp}}\quad  \\
\label{follows1}   &=
 e_{\SSTP_{\empb\al\emp}}   \left(  y_i 
 +
  y_{\emptyset(i+1)}
  -   y_{b_{\al \bet} {\aatch}   }
   \right)e_{\SSTP_{\empb\al\emp}}
\end{align} where the first equality follows from \cref{conclus1} and the second by \cref{jumpydots,jumpydots1}.
Thus substituting  \cref{follows1,follows2} in to \cref{follows3}, the first equality follows.  
\end{proof}

\subsection{The cyclotomic relation}   \label{cylo2}    
We now verify relation \ref{cyclotomic}. 
We have that 
$\Psi({\sf 1}_{\al})=e_{\SSTP_\al}$ for any $\al\in \Pi$.  
If the $\al$-hyperplane is a wall of the dominant region, then the tableau $\SSTP_\al$ is {\em non-standard} and therefore 
$e_{\SSTP_\al}=0$ by \cref{resconsider}.  
Now, let $\gam\in \Pi$ be arbitrary.  By \cref{tod-tod}, we have that 
$$
\Psi 
\left(\begin{minipage}{0.8cm} \begin{tikzpicture}[scale=1.5]
    \clip (0,-0.25-1/8) rectangle (0.5,0.75); 
{\clip (0,-0.25-1/8) rectangle (0.5,0.75); \foreach \i in {0,1,2,...,40}
  {
    \path  (-1,-0.75)++(0:0.125*\i cm)  coordinate (a\i);
        \path  (-1,-0.75)++(90:0.125*\i cm)  coordinate (b\i);
     \draw [densely dotted](a\i)--++(90:3);
        \draw [densely dotted](b\i)--++(0:10);
 }
} 
  \draw[ black, line width=0.06cm] 
  (0,-0.25-1/8) --++(0:0.5); 
  \draw[ black, line width=0.06cm] 
  (0,0.75) --++(0:0.5);
 \draw[darkgreen,line width=0.12cm](0.25,0)--(0.25,0.75/2);
\fill[darkgreen] (0.25,1/2)  circle (4pt); 
\fill[darkgreen] (0.25,-1/8)  circle (4pt); 
  \end{tikzpicture}
 \end{minipage}  \right)  
 = 
   e_{\SSTP_{\empg}}      \big(
   y_{b_\gam  {\aatch} - {\aatch}+\emptyset(k+1)}-y_{k   }
       \big) 
   e_{\SSTP_{\empg}} 
   =
     e_{\SSTP_{\empg}}      \big(
   y_{\emptyset(k+1)}-y_{k   }
       \big) 
   e_{\SSTP_{\empg}} 
$$
where the latter equality follows from \cref{jumpydots,jumpydots1}.  
If $x\equiv1$ modulo $h$, then 
\begin{equation}\label{asmany2}y_x e_{\SSTP_\empg}
=
e_{\SSTP_\empg}(\psi^x_1 y_1 \psi^1_x)e_{\SSTP_\empg}=0
\end{equation}
 by relation \ref{rel1.12}.  
If not, then by relation \ref{rel1.10} we have that 
\begin{equation}\label{asmany}y_xe_{\SSTP_\empg}  =  y_{x-1}e_{\SSTP_\empg} - 
e_{\SSTP_\empg}\psi_x \psi_x e_{\SSTP_\empg}\end{equation}
 where the latter term is zero by \cref{inwhatfollows} (as $(\eps_1,\dots ,\eps_{x-1},,\eps_{x+1},\eps_{x},\eps_{x+2},\dots  ,\eps_{\aatch}	 ) $  is non-standard for  $b_\gam=1$).  
 Thus the cyclotomic relation holds    (as we can apply \cref{asmany} as many times as necessary and then apply \cref{asmany2}).

\section{Decomposition numbers of cyclotomic Hecke algebras}\label{decompositionnumbersarepKL}

In this section we recall the construction of the graded cellular and ``light leaves" bases for the algebras $\mathscr{S}^{\rm br}_{\underline{h}}(n,\sigma)$, our quotient algebras 
$\mathcal{H}^\sigma_n /\mathcal{H}^\sigma_n {\sf y}_{\aatchpair}\mathcal{H}^\sigma_n $, and their truncations.  
We show that the homomorphism $\Psi$ preserves these $\ZZ$-bases 
(trivially, by definition) and hence deduce that 
$\Psi$ is indeed an isomorphism and hence prove Theorems  A and B of the introduction.

   \subsection{Why is it enough to consider   the truncated algebras?} 
\label{trunddddd}
Thus far in the paper, we have truncated to consider paths which terminate at a point $\la\in \mathscr{P}_{\underline{h}}(n,\sigma)\subseteq \mathscr{P}_{\underline{h}}(n )$.  
This  is, in general, a  proper  co-saturated subset of the principal linkage class of multipartitions for a given $n\in \ZZ_{\geq0}$.  

\begin{thm}[{\cite[Corollary 2.14]{cell4us}}]\label{cellularity breeds contempt}
 For each $\la$, we fix  $\SSTP_\la\in \Std  (\la)$ a choice of  reduced path.   
 The algebra  $\mathcal{H}^\sigma_n /\mathcal{H}^\sigma_n {\sf y}_{\aatchpair}\mathcal{H}^\sigma_n $ is quasi-hereditary with graded cellular basis $$
\{\psi^\northT _{\SSTP_\la}  \psi_\southT ^{\SSTP_\la}\mid 
 \northT,\southT \in \Std  (\la),\la\in {\mathscr P} _{\underline{h}}(n)\}
$$
  with respect to the reverse cylindric  order on ${\mathscr P} _{\underline{h}}(n)$ 
  (see \cite[Definition 1.3]{cell4us}, but for the subset  $\mathscr{P} _{\underline{h}} (n,\sigma) \subseteq \mathscr{P}_{\underline{h}} (n )$ is a refinement of   the opposite of the Bruhat ordering on their alcoves) 
  and the anti-involution, $\ast$, given by flipping a diagram through the horizontal axis.   
\end{thm}
 
 \begin{rmk}
 In {\cite[Corollary 2.14]{cell4us}} it is not explicitly stated that the algebra is quasi-hereditary.  
 However, this is immediate from  the  fact that 
each layer in the cell-filtration has an idempotent  $e_{\SSTP_\la}$ for $\la \in {\mathscr P} _{\underline{h}}(n)$
(and standard facts about cellular algebras).  
 \end{rmk}
 
  \begin{rmk}
In the case of the Hecke algebra of the symmetric group, the  basis of {\cite[Corollary 2.14]{cell4us}} is equivalent (via uni-triangular change of basis with respect to the dominance ordering) to the cellular basis of Hu--Mathas \cite{hm10}.   
   \end{rmk}

\begin{eg}
Let  $\la =(3^n,1^{15})$ with $n\geq 0$.
The first $n=0,1,2,3,4,5$ partitions in this sequence are
 $ (1^{15})$, $ (3,1^{15})$, 
$ (3^2,1^{15})$, $ (3^3,1^{15})$, $ (3^4,1^{15})$ and 
$ (3^5,1^{15})$, all of which label simple modules which 
belong to the principal blocks of their  corresponding group algebras.  In fact, they all label the same point, in the alcove $ s_{\color{magenta}\varepsilon_3-\varepsilon_1}s_{\color{darkgreen}\varepsilon_1-\varepsilon_2}
	s_{\color{cyan}\varepsilon_2-\varepsilon_3}s_{\color{magenta}\varepsilon_3-\varepsilon_1}
	s_{\color{darkgreen}\varepsilon_1-\varepsilon_2}s_{\color{cyan}\varepsilon_2-\varepsilon_3}A_0$, in the projection onto 2-dimensional space in 
\cref{diag2}.  
	However, 
	$\Std_{n,\sigma}(\la)=\emptyset$ for the first five of these partitions.  For $\la =(3^n,1^{15})$ with $n\geq 5$ we have that 	$\Std_{n,\sigma}(\la )\neq \emptyset$.  
Thus, one might be forgiven in thinking that 	our Theorem A only allows us to see   $\la $ for $n\geq 5$.  This is, in fact, not the case as we shall soon see.  
\end{eg}


   \begin{prop}
   \label{dettensor}
Given a partition $\lambda=(\lambda_1 ,\lambda _2,\dots)$, we set 
  ${\rm det}_h(\lambda) =  	(h, \lambda_1 ,\lambda _2,\dots) .$  We have 
 an injective map  of partially ordered sets 
${\rm det}_h:\mathscr{P} _{\underline{h}}(n) 
\hookrightarrow \mathscr{P} _{\aatchpair} (n+{\aatch})  $ given by 
 $${\rm det}_h(\lambda^{(0)},\lambda^{(1)},\dots ,\lambda^{(\ell-1)}) =  
 	({\rm det}_{h_0}(\lambda^{(0)}),{\rm det}_{h_1}(\lambda^{(1)}),\dots,{\rm det}_{h_{\ell-1}}(\lambda^{(\ell-1)}))  $$
and  ${\rm det}_h(\mathscr{P} _{\underline{h}}(n))\subseteq 
 \mathscr{P} _{\aatchpair} (n+{\aatch}) 
 $    is a co-saturated subset.  We have 
  an isomorphism of graded $\ZZ$-algebras
\begin{equation}\label{isomorphimmm2}
\sum_{\begin{subarray}c \southT,\northT \in \Std_n \end{subarray}} e_{\northT}  (\mathcal{H}_n^\sigma /\mathcal{H}_n^\sigma  {\sf y}_{\aatchpair}\mathcal{H}_n^\sigma )
e_{\southT}
  \cong 
\sum_{\begin{subarray}c \southT,\northT \in \Std_n \end{subarray}} e_{\SSTP_\emptyset \otimes \northT}  (\mathcal{H}_{n+{\aatch}} ^\sigma /\mathcal{H}_{n+{\aatch}} ^\sigma  {\sf y}_{\aatchpair}\mathcal{H}_{n+{\aatch}} ^\sigma )e_{\SSTP_\emptyset \otimes  \southT} 
\end{equation}
 where $\Std_n= \cup_{\la\in  \mathscr{P} _{\aatchpair} (n ) } \Std(\la)  $.  
\end{prop}
\begin{proof}
On the level of graded $\ZZ$-modules  the isomorphism, $\phi$ say, is clear.    
The local KLR relations also go through easily.  
We have that 
\begin{equation}\label{cylot}
 \phi(y_1 e_\SSTP)=y_{{\aatch}+1} 
 e_{\SSTP_\emptyset \otimes  \SSTP}=
 y_{ 1} 
 e_{\SSTP_\emptyset \otimes  \SSTP}=0 = y_1 e_\SSTP
 \end{equation}
 where the second equality follows using the same argument as \cref{jumpydots,jumpydots1} and the other equalities   all hold by definition.  
 We further note  that $\SSTP $ is   dominant  path if and only if $\SSTP_\emptyset\otimes \SSTP  $  is   a  dominant  path.   
 Thus the cyclotomic relation follows from  \cref{cylot} and \cref{inwhatfollows}.   
\end{proof}

We wish to only explicitly consider the principal linkage class, but to make deductions for all   regular linkage classes. 
 This is a standard Lie theoretic trick known as the {\sf translation principle}. 
 Given     $\Gamma \subseteq  \mathscr{P}_{\underline{h}}(n)$ any co-saturated subset and  $r \in \ZZ/e\ZZ$ we let 
$$  
e_{\Gamma}= \sum_{
 \begin{subarray}c
 \SSTP\in \Std(\mu)
 \\
 \mu \in \Gamma 
 \end{subarray}} e_\SSTP
 \qquad  \qquad 
 {\sf E}_r = \sum_{i_1,\dots ,i_{n}\in \ZZ/e\ZZ}e(i_1,\dots, i_{n},r) $$
 denote the corresponding   idempotents.  
Given $\la \in \mathscr{P}_{\underline{h}}(n)$ we set $\Lambda=( \Shl \cdot \la) \cap
 \mathscr{P}_{\underline{h}}(n)$.   
 Since every $\la$ belongs to some linkage class, 
 we have that 
 $ \mathscr{P}_{\underline{h}}(n)= \Lambda'\cup \Lambda'' \cup \dots $ 
 and we have a corresponding decomposition 
 $$ 
 \mathcal{H}_n^\sigma /\mathcal{H}_n^\sigma  {\sf y}_{\aatchpair}\mathcal{H}_n^\sigma  =
  \mathcal{H}^{ \Lambda',\sigma}_n
  \oplus 
    \mathcal{H}^{ \Lambda'',\sigma}_n\oplus \dots 
  \quad \text{where}\quad 
 \mathcal{H}^{ \Lambda,\sigma}_n=   e_{ \Lambda}  (\mathcal{H}_n^\sigma /\mathcal{H}_n^\sigma  {\sf y}_{\aatchpair}\mathcal{H}_n^\sigma )e_{ \Lambda}   
 $$
 and similarly for the primed cases.  
Now, we  let $\square$ denote an addable node of the Young diagram multipartition $\la\in \mathscr{P}_{\underline{h}}(n)$, that is 
we suppose that $ \la \cup \square= \la' $ for some $\la' \in \mathscr{P}_{\underline{h}}(n+1)$.

 \begin{prop} 
   \label{decompositionnumbermadeasy}
Suppose that  $\la \in   \mathscr{P}_{\underline{h}}(n)$ and $ \la + \square=\lambda' \in   \mathscr{P}_{\underline{h}}(n+1)$
 are $\sigma$-regular  and  $\square$ is of residue $r\in \ZZ/e\ZZ$ say.  We have an injective map 
 $$
\varphi: \Lambda \hookrightarrow \Lambda' \qquad \varphi(\mu)=\mu+\square 
 $$
 for $\square$ the unique addable node of residue $r\in\ZZ/e\ZZ$.  
The image, $\varphi(\Lambda)$,   is a co-saturated subset of $\Lambda'$.  
We have an isomorphism of graded $\ZZ$-algebras:  
\begin{equation}\label{isomorphimmm}
  \mathcal{H}^{\Lambda,\sigma }_n  
    \cong 
{\sf E}_r(	e_{\varphi(\Lambda)} \mathcal{H}^{\Lambda' }_{n+1} e_{\varphi(\Lambda)} ){\sf E}_r  
\end{equation}
and this preserves the cellular structure.  
 \end{prop}
\begin{proof}
Since both $\la$ and $\la+\square$ are both $e$-regular, there is a bijection between the path bases of the algebras in \cref{isomorphimmm}.  
(Note that if  $\la$ were on a hyperplane and $\la+\square$ in an alcove, then the number of paths would  double.)
Thus we need only check that this $\ZZ$-module homomorphism lifts to an algebra homomorphism.  However this is obvious, as all we have done is add a single strand (of residue $r\in \ZZ/e\ZZ$) to the righthand-side of the diagram and this preserves the multiplication.  
\end{proof}



 Thus  any 
regular block of $\mathcal{H}^\sigma_N / \mathcal{H}^\sigma_N  {\sf y}_{\aatchpair} \mathcal{H}^\sigma_N $ 
is isomorphic to a co-saturated idempotent subalgebra of 
$\mathcal{H}^\sigma_n/ \mathcal{H}^\sigma_n  {\sf y}_{\aatchpair} \mathcal{H}^\sigma_n $  
 for some $n\geq N$. 
 Such truncations preserve decomposition numbers \cite[Appendix]{Donkin}  and much cohomological structure and so it suffices to consider only these truncated algebras (which is precisely what we  have done thus far in the paper!).

\subsection{Bases  
of diagrammatic    algebras} 
For  $\la,\mu\in {\mathscr P}_{\underline{h}}(n,\sigma)$, we choose reduced paths $\SSTP_\w \in \Std_{n,\sigma} (\la)$ and 
$\SSTP_{\underline{v}} \in \Std_{n,\sigma} (\mu)$  which will remain fixed for the remainder of this section. We remind  the reader that this implicitly says that $\la \in  w    A_0$ and $\mu \in v    A_0$.  
We have  shown that the map $$\Psi:
 \mathscr{S}^{\rm br}_{\underline{h}}(n,\sigma) \to
 {\sf f}_{n,\sigma}( \mathcal{H}_n^\sigma
  /
   \mathcal{H}_n^\sigma {\sf y}_{\aatchpair}  \mathcal{H}_n^\sigma){\sf f}_{n,\sigma}$$ is    a   graded $\ZZ$-algebra homomorphism.  It remains to show that this map is an isomorphism.  
Let $\la  \in  \mathscr{P}_{\underline{h}}(n,\sigma)$.   Given any reduced 
path  $\SSTP_\w \in \Std_{n,\sigma}(\la)$ and any (not necessarily reduced)  $\SSTQ \in \Std_{n,\sigma}(\la )$ we will inductively construct   elements   
$$
C^{\SSTP}_{\SSTQ} \in	 {\sf1}_{\SSTP }  \mathscr{S}^{\rm br}_{\underline{h}}(n,\sigma){\sf1}_{\SSTQ}
\qquad
c^{\SSTP}_{\SSTQ} \in 	e_{\SSTP} ( \mathcal{H}_n^\sigma
  /
   \mathcal{H}_n^\sigma {\sf y}_{\aatchpair}  \mathcal{H}_n^\sigma) e_{\SSTQ}   
$$
which provide  (cellular) $\ZZ$-bases of both algebras  which match up under the homomorphism, thus proving that $\Phi$ is  indeed  an isomorphism.  

We can extend a path    $\SSTQ'\in\Std_{n,\sigma}(\la)$   to obtain a new path $\SSTQ$ in one of three possible ways 
$$
\SSTQ=\SSTQ'\otimes \SSTP_\al 	\qquad 	\SSTQ=\SSTQ'\otimes  \reflectpath   		\qquad \SSTQ=\SSTQ'\otimes \SSTP_\emptyset 
$$
for some  $\al \in \Pi$. 
  The first two cases each subdivide into a further two cases based on whether $\al$ is an upper or lower wall of the alcove containing $\la$.    These four cases are pictured in \cref{upanddown} (for $\SSTP_\emptyset$ we refer the reader to \cref{figure1}).  
  Any two reduced paths $ \SSTP_{\w} , \SSTP_{\underline{v}}\in \Std_{n,\sigma} (\la)$  can be obtained from one another 
by some iterated application of hexagon and commutativity permutations.
    We let 
 $${\sf rex}_{\SSTP_{\w} }^{\SSTP_{\underline{v}}  } \qquad {\sf REX}_{\SSTP_{\w} }^{\SSTP_{\underline{v}}  }$$ denote the corresponding path-morphism in the  
   algebras $\mathcal{H}_n^\sigma/\mathcal{H}_n^\sigma{\sf y}_{\aatchpair} \mathcal{H}_n^\sigma$ and $ \mathscr{S}^{\rm br}_{\underline{h}}(n,\sigma)$, respectively (so-named as they permute  reduced expressions).  In the following construction, we will assume that 
   the elements $c^{\SSTP'}_{\SSTQ'}$ and 
   $C^{\SSTP'}_{\SSTQ'}$ exist for any choice of
   reduced path  $\SSTP'$.  
 We then extend $\SSTP'$ using one of the 
$   U_0, U_1,  D_0, $ and $D_1$ paths (which puts a restriction on the form of the reduced expression)
 but then use a ``rex move" to remove 
 obtain  elements  $c^{\SSTP }_{\SSTQ }$ and $C^{\SSTP }_{\SSTQ }$
for $\SSTP$  an arbitrary reduced expression.

\begin{figure}[ht!]
 
$$
   \begin{minipage}{2.6cm}
\end{minipage} 
 $$
\caption{The first (respectively last) two paths   are  $\SSTP_\al$ and $\reflectpath $ originating in an alcove with $\al$ labelling an upper (respectively lower) wall. 
The origin lies below the $\al$-hyperplane.  
We call these paths $  U_0, U_1,  D_0, $ and $D_1$ respectively.  
}
 \label{upanddown}
\end{figure}

\begin{defn}\label{twins} 
Suppose that $\la $ belongs to an alcove which has a hyperplane labelled by $\al$ as an   upper alcove wall.  
Let 
   $\SSTQ'\in\Std_{n,\sigma}(\la)$. 
 If $\SSTQ= \SSTQ' \otimes \SSTP_\al$ then we set $\deg(\SSTQ)=\deg(\SSTQ')$  and  we define 
\begin{align*}
 C_{\SSTQ}^\SSTP= 
 {\sf REX}_{\SSTP'\otimes \SSTP_\al}^{\SSTP}(
C_{\SSTQ'}^{\SSTP'} \;\otimes \; {\sf 1}_{  \al})  
\qquad
 &c_{\SSTQ}^\SSTP= 
 {\sf rex}_{\SSTP'\otimes \SSTP_\al}^{\SSTP}(
c_{\SSTQ'}^{\SSTP'} \;\otimes \; e_{\SSTP_\al}). 
\intertext{
If $\SSTQ= \SSTQ' \otimes  \reflectpath  $ then
we set  $\deg(\SSTQ)=\deg(\SSTQ')+1$ and we define 
}
C_{\SSTQ}^\SSTP= 
 {\sf REX}_{\SSTP'\otimes \SSTP_\emp}^{\SSTP}(
 C_{\SSTQ'}^{\SSTP'} \;\otimes \; {\sf SPOT}_{\al}^\emp)  
\qquad
&c_{\SSTQ}^\SSTP= 
 {\sf rex}_{\SSTP'\otimes \SSTP_\emp}^{\SSTP}(
 c_{\SSTQ'}^{\SSTP'} \;\otimes \; {\sf spot}_{\al}^\emp). 
\end{align*}
Now suppose that $\la $ belongs to an alcove which has a hyperplane labelled by $\al$ as a lower alcove wall. 
Thus we can choose $ \SSTP_{\underline{v}} \otimes \SSTP_\al=\SSTP'\in\Std(\la) $.  
For 
$\SSTQ= \SSTQ' \otimes \SSTP_\al$, 
we set  $\deg(\SSTQ)=\deg(\SSTQ') $   
 and   define 
  \begin{align*}
C_{\SSTQ}^\SSTP= &
{\sf REX}^{		
  \SSTP  }_{\SSTP_{\underline{v}  \emp\emp}}
\big( 
{\sf 1}_{ {\underline{v}}}\otimes 
({\sf SPOT}_\al^\emp \circ{\sf FORK}_{\al\al}^{\al\emp})\big) 
\big(
  C^{\SSTP'}
_{\SSTQ'} 
\otimes {\sf1}_{ \al}\big) 
 \\
\qquad c_{\SSTQ}^\SSTP=  &
{\sf rex}^{		
  \SSTP  }_{\SSTP_{\underline{v} \emp\emp}}
\big( 
{\sf e} _ {\SSTP_{\underline{v}}}\otimes
( {\sf spot}_\al^\emp \circ {\sf fork}_{\al\al}^{\al\emp})\big) 
\big(  c^{\SSTP'}
_{\SSTQ'} 
\otimes e_{\SSTP_\al}\big)  
\end{align*}
  and if $\SSTQ= \SSTQ' \otimes  \reflectpath  $ then 
  then
we set  $\deg(\SSTQ)=\deg(\SSTQ')-1$ and we define   
\begin{align*}
C_{\SSTQ}^\SSTP&= 
{\sf REX}^{		
  \SSTP  }_{  \SSTP_{\underline{v}\al \emp}}
\big( 
{\sf 1}_{ {\underline{v}}}\otimes {\sf FORK}_{\al\al}^{\al\emp}\big) 
\big(
  C^{\SSTP'}
_{\SSTQ'} 
\otimes {\sf 1}_{ \al}\big) 
\\ 
c_{\SSTQ}^\SSTP&=  
{\sf rex}^{		
  \SSTP  }_{ \SSTP_{\underline{v}\al \emp}}
\big( 
{ e}_{\SSTP_{\underline{v}}}\otimes {\sf fork}_{\al\al}^{\al\emp}\big) 
\big(
      c^{\SSTP'}
_{\SSTQ'} 
\otimes e_{\SSTP_\al}\big) .  
 \end{align*}
  \end{defn}

In each of the four cases above,
the path $\SSTP$  is a reduced path by construction (and our assumption that $\SSTP'$ is  reduced).  
We remark that the degree of the path, $\SSTQ$, is   equal to the degree of both the elements $c^\SSTP_\SSTQ$  and $C^\SSTP_\SSTQ$ (recall that $\SSTP$ is a path associated to a  reduced word and so is of degree zero).    

\begin{thm}[Light leaves basis, \cite{MR3555156,antiLW}]\label{PPPPPPSDSDS}
   For each $\la\in \mathscr{P}_{\underline{h}}(n,\sigma)$, we fix an arbitrary 
reduced path   $
 \SSTP_\w \in \Std_{n,\sigma}(\la).
 $ 
 The algebra 
  $ \mathscr{S}^{\rm br}_{\underline{h}}(n,\sigma)$ is   quasi-hereditary  with  graded  integral cellular basis  
 $$
 \{ C^{\SSTP}_{\SSTP_\w} 
  C^{\SSTP_\w}_\SSTQ
 \mid
  \SSTP,\SSTQ \in \Std_{n,\sigma}  (\la), \la \in  \mathscr{P}_{\underline{h}} (n,\sigma)  \}  
 $$
  with respect to the Bruhat ordering $\rhd $ on  $ \mathscr{P}_{\underline{h}} (n,\sigma)$,  the anti-involution $\ast$ given by flipping a diagram through the horizontal axis and the map $\deg :   \Std_{n,\sigma}  (\la) \to \ZZ$.  
\end{thm}


We recalled a  general construction of a cellular basis of $\mathcal{H}^\sigma_n/ \mathcal{H}^\sigma_n {\sf y}_{\aatchpair}\mathcal{H}^\sigma_n$ in \cref{cellularity breeds contempt} subject to  choosing the reduced expressions.
This  provides a cellular basis of  ${\sf f}_{n,\sigma}\mathcal{H}^\sigma_n/ \mathcal{H}^\sigma_n {\sf y}_{\aatchpair}\mathcal{H}^\sigma_n {\sf f}_{n,\sigma}$
by idempotent truncation.  
Choosing our  expressions    
 so as to be compatible with \cref{PPPPPPSDSDS} through the map $\Psi$, we obtain the following.

\begin{thm}[{Light leaves basis, \cite[Theorem 3.12]{cell4us}}] 
   For each $\la\in \mathscr{P}_{\underline{h}}(n,\sigma)$,   choose an arbitrary 
reduced path   $
 \SSTP_\w \in \Std_{n,\sigma}(\la)
 $.   
The algebra 
  ${\sf f}_{n,\sigma} (\mathcal{H}^\sigma_n/ \mathcal{H}^\sigma_n {\sf y}_{\aatchpair}\mathcal{H}^\sigma_n) {\sf f}_{n,\sigma}$ is   quasi-hereditary   with  graded  integral cellular basis  
 $$
 \{ c^{\SSTP}_{\SSTP_\w} 
   c^{\SSTP_\w}_\SSTQ
 \mid
  \SSTP,\SSTQ \in \Std_{n,\sigma}(\la), \la \in  \mathscr{P}_{\underline{h}}(n,\sigma)  \}  
 $$
 with respect to the  Bruhat ordering $\rhd $ on  $ \mathscr{P}_{\underline{h}} (n,\sigma)$,  the anti-involution $\ast$ given by flipping a diagram through the horizontal axis and the map $\deg :   \Std_{n,\sigma}  (\la) \to \ZZ$.  
\end{thm}

\begin{thm}\label{theoremA:2}
Let     $\sigma\in \ZZ^\ell $  and $e\in \ZZ_{>1}$ and suppose that 
$\underline{h}\in \ZZ^\ell_{\geq0}$ is $(\sigma,e)$-admissible.  
We have a canonical isomorphism of  graded $\ZZ$-algebras, 
$$
    {\sf f}_{n,\sigma}^+\left(\mathcal{H}_n^\sigma/\mathcal{H}_n^\sigma {\sf y}_{\aatchpair }\mathcal{H}_n^\sigma 		\right)     {\sf f}_{n,\sigma}^+
\cong {\rm End}_{\mathcal{D}^{\rm asph,\oplus}_{\rm BS}  (  A_{h_0}\times \mydots \times A_{h_{\ell-1}} \backslash \widehat{A}_{h_0+\dots+h_{\ell-1}} )}\left(\oplus_{\underline{w} \in \Lambda({n,\sigma})} B_{\underline{w}} \right).
$$
That is, Theorem A of the introduction holds.

 \end{thm}
\begin{proof}
In \cref{embed} we defined a map from $\mathscr{S}^{\rm br}_{\underline{h}}(n)$ to 
$\mathcal{H}_n^\sigma/\mathcal{H}_n^\sigma{\sf y}_{\aatchpair} \mathcal{H}_n^\sigma$ via the generators of the former algebra.  
In \cref{relations} we showed that this map was a homomorphism by 
verifying that the relations for 
$\mathscr{S}^{\rm br}_{\underline{h}}(n)$ held in the image of the homomorphism.   
Now, the construction of the light leaves bases in   $\mathscr{S}^{\rm br}_{\underline{h}}(n)$ 
(respectively $\mathcal{H}_n^\sigma$) 
is given in terms of the generator (respectively their images).  
Thus the map preserves the $\ZZ$-bases and hence is an isomorphism.  
Thus the result follows from \cref{amit}.
\end{proof}

   An earlier attempt  to solve  the Libedinsky--Plaza conjecture for the classical  blob algebra 
    (the case of $h=1$ and $\ell=2$) 
    has already  led to  a deeper understanding of 
 structure of the diagrammatic Soergel category  \cite{STEEN}.    We remark that their is no obvious intersection between their results and ours (they do not succeed in proving the $h=1$ and $\ell=2$ case, but nor do our results imply theirs).

\subsection{Decomposition numbers of Hecke algebras }\label{whyyyy}
\renewcommand{\SSTU}{{\mathsf{P}_\w}}
\renewcommand{\SSTV}{{\mathsf{P}_\w}}
 
For  $\la,\mu\in {\mathscr P}_{\underline{h}}(n,\sigma)$, we reiterate that we have chosen to fix reduced paths $\SSTP_\w \in \Std_{n,\sigma} (\la)$ and 
$\SSTP_{\underline{v}} \in \Std_{n,\sigma} (\mu)$.  
 We define  one-sided ideals 
\begin{align*}
\mathscr{S}_{n,\sigma} ^{ \trianglerighteq \underline{v} }  =  \mathscr{S} ^{\rm br} _{\underline{h}}(n,\sigma) {\sf1}_{\SSTP_{\underline{v}}}&&
\mathscr{S}_{n,\sigma} ^{ \vartriangleright \w }  =
    \mathscr{S}_{n,\sigma} ^{ \trianglerighteq \w }   \cap \ZZ \{C^\northT_{\SSTP_{\underline{v}}} C^{\SSTP_{\underline{v}}}_\southT \mid 
\northT,\southT \in \Std_{n,\sigma} (\mu),\mu\vartriangleright \la \}   
\\
  \mathcal{H}_{+}   ^{ \trianglerighteq \mu }  = \mathscr{S}^{\rm br}_{\underline{h}}(n) e_{\SSTP_{\underline{v}}}&
&  \mathcal{H}_{+}  ^{ \vartriangleright \la }  =
 \mathcal{H}_{+}   ^{ \trianglerighteq \la }    \cap \ZZ \{c^\northT_{\SSTP_{\underline{v}}} C^{\SSTP_{\underline{v}}}_\southT \mid 
\northT,\southT \in \Std_{n,\sigma} (\mu),\mu\vartriangleright \la \}   
\end{align*}
and we   define the   {\sf standard}  modules of $\mathscr{S}^{\rm br}_{\underline{h}}(n,\sigma)$ and   ${\sf f} _{n,\sigma}(\mathcal{H}^\sigma_n/ \mathcal{H}^\sigma_n {\sf y}_{\aatchpair}\mathcal{H}^\sigma_n){\sf f} _{n,\sigma}$  by considering the resulting subquotients.  The light leaves construction gives us explicit bases of these quotients as follows 
\begin{equation}  \label{identification}
\Delta _\ZZ(\w) = \{ C^\sts_{\SSTP_\w}   +  \mathscr{S}_{n,\sigma} ^{ \vartriangleright \la } 
       \mid \sts \in \Std_{+} (\lambda)\}\qquad 
         {\sf f}_{n,\sigma} {\bf S}_\ZZ(\lambda) = \{ c^\sts_{\SSTP_\w} +  \mathcal{H}  ^{ \vartriangleright \la } 
       \mid \sts \in \Std_{+} (\lambda)\}
\end{equation}
respectively for $\la\in{\mathscr P}_{\underline{h}}(n,\sigma)$. 
  The modules  $  {\sf f}_{n,\sigma} {\bf S}_\ZZ(\lambda) $ are obtained by  truncating   the cell modules (${\bf S}_\ZZ(\lambda) $, say) for the  cellular structure in \cref{cellularity breeds contempt}.  
 For $\Bbbk$ a field,  we define $$ \Delta_\Bbbk(\w)=
 \Delta_\ZZ(\w)\otimes _\ZZ\Bbbk\qquad 
{\sf f}_{n,\sigma}  {\bf S}_\Bbbk(\lambda)=
{\sf f}_{n,\sigma} {\bf S}_\ZZ(\lambda)\otimes _\ZZ\Bbbk.$$ 
We recall that the   cellular structure allows us to define bilinear forms, for each
 $\la \in \mathscr{P}_{\underline{h}}(n)$,  there  are bilinear forms
  $\langle\ ,\ \rangle_{\mathscr{S} }^{   \la} $ 
and $ \langle\ ,\ \rangle_{\mathcal {H} }^{   \la}$ 
  on $\Delta(\lambda) $ and $ {\sf f}_{n,\sigma}{\bf S}_\Bbbk(\la)  $ respectively,   which
are determined by
\begin{align}\label{geoide}
\begin{split}
C ^{\SSTU}_{\SSTP}C ^{\SSTQ}_{\SSTV}&\equiv
  \langle C ^\SSTP_{\SSTU},C ^\SSTQ_{\SSTU}
  \rangle_{\mathscr{S} }^{   \la}\; {\sf 1} _\w \pmod{\mathscr{S}_{n,\sigma}^{\vartriangleright \lambda}}
\\
  c ^{\SSTU}_{\SSTP}c ^{\SSTQ}_{\SSTV}&\equiv
  \langle c ^\SSTP_{\SSTU},c ^\SSTQ_{\SSTU}
  \rangle_{\mathcal {H} }^{   \la} \;  
e_{\SSTP_\w}\pmod{\mathcal{H}_{n,\sigma}^{\vartriangleright \lambda}}
\end{split}  \end{align}
for any $\SSTP,\SSTQ, \SSTU,\SSTV\in \Std(\lambda  )$.  
Factoring out by the radicals of these forms,  we obtain a complete set of non-isomorphic simple modules 
 for $\mathscr{S}^{\rm br}_{\underline{h}}(n,\sigma)$ and $\mathcal{H}^\sigma_n/ \mathcal{H}^\sigma_n {\sf y}_{\aatchpair }\mathcal{H}^\sigma_n $ as follows 
 $$ 
L_\Bbbk(\w) =
\Delta_\Bbbk(\w) /
  \rad(\Delta_\Bbbk(\w) )
 \qquad 
 {\sf f}_{n,\sigma}{\bf D}_\Bbbk(\lambda) =
  {\sf f}_{n,\sigma}{\bf S}_\Bbbk(\lambda) /
  \rad({\sf f}_{n,\sigma} {\bf S}_\Bbbk(\lambda) )
 $$
respectively for $\la \in \mathscr{P}^{+}_{\underline{h}}(n)$. 
Finally, the projective indecomposable modules are as follows,  
\begin{equation}\label{oooooooo}
\mathscr{S}_{n,\sigma} ^{ \trianglerighteq {\underline{v}} } = 
\bigoplus_{ w \leq v}
\dim_t({\sf 1}_{\underline{v}}L_\Bbbk(\w))P_\Bbbk(\w)
\;\;\; \mathcal{H}_{n,\sigma} ^{ \trianglerighteq \mu  } = 
\bigoplus_{ \la \trianglerighteq \mu		}
\dim_t(e_{\SSTP_\mu} {\bf D}_\Bbbk(\la)){\bf P}_\Bbbk(\la).
\end{equation}
The isomorphism, $\Psi$,  preserves standard, simple, and projective modules.

 The categorical (rather than geometric) definition of    $p$-Kazhdan--Lusztig polynomials is given via the {\em diagrammatic character} of \cite[Definition 6.23]{MR3555156}.  This graded character is defined in terms of dimensions of certain weight spaces in the light leaves basis.  
Using the identifications of \cref{identification,oooooooo},  the definition of the anti-spherical $p$-Kazhdan--Lusztig polynomial, ${^pn}_{v,w}(t) $,   is as follows, 
$$
{^pn}_{v,w}(t)
:=
 \dim_t
\Hom_{\mathscr{S}^{\rm br}_{\underline{h}}(n,\sigma)}(  P({  {\underline{v}}}), \Delta(\w)  )  =
 \sum_{k\in\ZZ}\dim [\Delta_\Bbbk(\w):
 L_\Bbbk(\underline{v})\langle k \rangle ] t^k 
   $$ for $v,w\in \Lambda(n,\sigma)$.  
We claim no originality in this observation and refer to  \cite[Theorem 4.8]{MR3591153} for more details.  Through our isomorphism this allows us to see that 
  the    graded decomposition numbers of symmetric groups and more general cyclotomic Hecke algebras are {\em tautologically equal to the associated $p$-Kazhdan--Lusztig polynomials} as follows,  
\begin{align*}
{^pn}_{v,w}(t)
=
 \sum_{k\in\ZZ}\dim [\Delta_\Bbbk(\w):
 L_\Bbbk(\underline{v})\langle k \rangle ] t^k 
 =
  \sum_{k\in\ZZ}\dim_t[{\sf f}_{n,\sigma}{\bf S}_\Bbbk(\la): {\sf f}_{n,\sigma}{\bf D}_\Bbbk(\mu)\langle k \rangle ] t^k 
\end{align*}
for $\la,\mu \in {\mathscr P}_{\underline{h}}(n,\sigma)$ where the equality follows immediately from our isomorphism.  
Finally, we remind the reader that 
 truncation by  ${\sf f}_{n,\sigma}$ is to a co-saturated subset of weights and so preserves the decomposition matrices of these algebras, see for example \cite[Appendix]{Donkin}

\subsection{Counterexamples to Lusztig's conjecture and intersection forms} In \cite{MR3671935}, the counterexamples to Soergel's conjecture are presented in the classical (rather than diagrammatic) language of  intersection forms associated to the fibre  of a Bott--Samelson resolution  of a Schubert varieties.   However, Williamson  emphasises that all his calculations were done using the equivalent diagrammatic setting  of the light leaves basis, which is {\em``explicit and amenable to computation"}.  Moreover, Williamson's counterexamples are dependent on the diagrammatics  because  it  is only {\em``from the diagrammatic approach [that] it is clear that  [the intersection form] $I^\Bbbk_{x,\w,d}$ is defined over $\ZZ$"} in the first place (see Section 3 of \cite{MR3671935} for more details).   
  In terms of the light leaves  cellular basis, Williamson's calculation makes a clever choice of a pair of partitions $\la,\mu$ (equivalently,  words $w,v \in \Shl$ labelling the alcoves containing these partitions) for which there exists a unique element   $\SSTQ\in \Std_{n,\sigma}  (\la)$ such that $\SSTQ\sim\SSTP_{\underline{v}} \in  \Std_{n,\sigma}  (\mu)$.  By highest weight theory, we have that 
$$d_{\la\mu}(t)=\begin{cases}
t^{ \deg(\SSTQ)} &\text{ if }\langle C^\SSTQ_{\SSTU}, C^\SSTQ_{\SSTU}\rangle_{\mathscr{S} }^{   \la}  =0 \in \Bbbk		\\
0									&\text{ otherwise }	
\end{cases}$$  
and Williamson proved  for $\la,\mu \in  \mathscr{P}_{h,1}(n)$ (a pair from ``around the Steinberg weight") that the form is zero for certain primes $p>h$ whereas it is equal to $1$ for $\Bbbk=\CC$ (and hence disproved Lusztig's and James' conjectures).  

Now, clearly the Gram matrices of the bilinear forms  
in \cref{geoide}  are preserved under isomorphism.  
Thus applying our isomorphism (and Brundan--Kleshchev's \cite{MR2551762}) one can view Williamson's counterexamples as being found entirely within the context of the symmetric group.  
More generally, we deduce the following:

\begin{thm}
Theorem B of the introduction  holds.  
\end{thm}

  \appendix

  \section{ Weakly graded monoidal categories }
 \label{Appendix}

In this appendix we describe the framework for constructing the breadth-enhanced diagrammatic Bott--Samelson endomorphism algebras. 
Informally, ``breadth-enhanced'' means that we record and keep track of the ``breadth''  of Soergel diagrams, including the ``blank spaces'' between strands. 
This is contrary to the usual working assumption that Soergel diagrams are defined only up to isotopy. 
We will say a few words for why we have chosen to break this convention in this paper.

Soergel diagrams and KLR diagrams have an important fundamental difference. 
KLR diagrams, which are essentially decorated wiring diagrams, always have the same number of nodes on the top and bottom edges.
By contrast, the top and bottom edges of a Soergel diagram may not have the same number of nodes. 
This basic observation is enough to ensure that a Soergel diagram cannot correspond to only one KLR diagram under the isomorphism in the main theorem.
For example, suppose the isomorphism maps the $\alpha$-coloured spot diagram to a KLR diagram $\mathsf{spot}_{\mathbf{\alpha}}$, with bottom edge $\mathsf{P}$ and top edge $\mathsf{Q}$.
Then the empty Soergel diagram (with no strands at all) should map to the KLR idempotent $e_{\mathsf{Q}}$. 
However it is also clear that the empty Soergel diagram should correspond to the empty KLR diagram. 

The breadth-enhanced diagrammatic Bott--Samelson endomorphism algebra introduces new idempotents, indexed by expressions in the extended alphabet $S \cup \{\emptyset\}$. 
This ensures that the isomorphism is well defined, with each breadth-enhanced Soergel diagram corresponding to a single KLR diagram. 
The breadth of a breadth-enhanced Soergel diagram is simply the number of strands of the corresponding KLR diagram, divided by ${\aatch}$. 
We draw breadth-enhanced Soergel diagrams so that the width is proportional to the breadth. 
In particular, we write $\mathsf{1}_{\emptyset}$ to indicate the empty Soergel diagram of breadth $1$ (i.e.~a ``blank space''), which corresponds to the KLR idempotent $e_{\mathsf{P}_{\emptyset}}$ with ${\aatch}$ strands. 
The breadth-enhanced algebras are Morita equivalent to the usual diagrammatic Bott--Samelson endomorphism algebras, by simply truncating with respect to the idempotents indexed by expressions which do not contain $\emptyset$. 
Thus once we prove the isomorphism for the breadth-enhanced algebras, we immediately obtain an isomorphism for the usual Bott--Samelson algebras.



The machinery for building breadth-enhanced algebras is the notion of a weakly graded monoidal category. 
Weakly graded monoidal categories can be thought of as generalizations of graded monoidal categories, with the grade shifts represented by tensoring with a fixed shifting object. 
The construction of breadth-enhanced algebras is then analogous to defining a graded category from a non-graded category by concentrating the objects in certain fixed degrees. 

We have chosen to write this appendix using the categorical (rather than the algebraic) perspective. 
We hope that this will make the results more applicable and the proofs easier to read. 
All the categories below will be assumed to be small.
We will also use ``monoidal'' to mean ``strict monoidal'' unless stated otherwise.
It is probably possible to generalize everything to arbitrary monoidal categories, but this will not be necessary for our purposes.

\subsection{Definition and examples}

\begin{defn}
  A {\sf weakly graded monoidal category} is a monoidal category $(\mathcal{A},\otimes)$ together with an object in the Drinfeld centre with trivial self-braiding.
This consists of the following data:
  \begin{itemize}
    \item an object $I$ in $\mathcal{A}$ called the {\sf shifting object};
    
    \item for each object $X$ in $\mathcal{A}$, an isomorphism $s_X:X \otimes I \xrightarrow{\sim} I \otimes X$ called a {\sf simple adjustment}
  \end{itemize}
  such that
  \begin{enumerate}[label=(WG\arabic*)]    
    \item \label{item:naturality} the simple adjustments $\{s_X\}$ are the components of a natural isomorphism $s : (-) \otimes I \Rightarrow I \otimes (-)$;
    
    \item \label{item:mult} for any objects $X,Y$ in $\mathcal{A}$ the following diagram  commutes
    \begin{equation*}
      \xymatrix{%
      X \otimes Y \otimes I \ar[rr]^{s_{X \otimes Y}} \ar[rd]_{1_X \otimes s_Y} & & I \otimes X \otimes Y \\
       & X \otimes I \otimes Y \ar[ru]_{s_X \otimes 1_Y}
      }
    \end{equation*}

    \item \label{item:identity} we have $s_{I}=1_{I \otimes I}$.
%
  \end{enumerate}
\end{defn}


\begin{eg} \label{eg:grwkgr}
Suppose $\mathcal{A}^{\bullet}$ is a graded monoidal category, i.e.~a monoidal category whose $\Hom$-spaces are graded modules. 
For the moment, let us drop the assumption of strictness and suppose that $\mathcal{A}^{\bullet}$ is strictly associative, but with non-trivial unitors. 
In the usual way we may construct a new category $\mathcal{A}$ by adding grade shifts and restricting to homogeneous morphisms. 
More precisely, the objects of $\mathcal{A}$ are the formal symbols $X(m)$ for each object $X$ of $\mathcal{A}^{\bullet}$ and each $m \in \mathbb{Z}$, and the $\Hom$-spaces are 
\begin{equation*}
\Hom_{\mathcal{A}}(X(m),Y(n))=\Hom_{\mathcal{A}^{\bullet}}^{n-m}(X,Y) \text{.}
\end{equation*}
It is clear that the grade shift $(1)$ is an autoequivalence of $\mathcal{A}$.
Moreover, the tensor product $X(m) \otimes Y(n)=(X \otimes Y)(m+n)$ gives $\mathcal{A}$ the structure of a monoidal category. 
Now let $\mathbbm{1}$ be the identity object in $\mathcal{A}^{\bullet}$ and set $I=\mathbbm{1}(1)$.
We observe that 
\begin{equation*}
X(m) \otimes \mathbbm{1}=(X \otimes \mathbbm{1})(1) \xrightarrow{\rho_X(1)} X(m+1) \xleftarrow{\lambda_X(1)} (\mathbbm{1} \otimes X)(1)=\mathbbm{1} \otimes X(m) \text{,}
\end{equation*}
and it is straightforward to check that the isomorphisms $s_{X(m)}=\lambda_{X(m)}(1)^{-1} \circ \rho_{X(m)}(1)$ satisfy axioms \ref{item:naturality}--\ref{item:identity}.
Thus $\mathcal{A}$ has the structure of a weakly graded monoidal category. 
\end{eg}

The main result which we will need in the next subsection is a coherence theorem for weakly graded monoidal categories.
Roughly, coherence for weakly graded monoidal categories means that every diagram built up from $s$ and identity morphisms (using composition and tensor products) commutes.
 The precise formulation of coherence requires some combinatorial constructions, which we describe below.
Let $\mathscr{W}$ be the set of non-empty words in the symbols $e$ and $x$. 
We define the following semigroup homomorphisms $\mathrm{length}:\mathscr{W} \rightarrow \mathbb{Z}_{\geq 0}$ and $\mathrm{breadth}:\mathscr{W} \rightarrow \mathbb{Z}_{\geq 0}$ on the generators:
\begin{align*}
\mathrm{length}(e)& =0 & \mathrm{breadth}(e)& =1 \\
\mathrm{length}(x)& =1 & \mathrm{breadth}(x)& =0 \text{.}
\end{align*}
For $w\in \mathscr{W}$ of length $n$, we can associate a functor $w_{\mathcal{A}}:\mathcal{A}^n \rightarrow \mathcal{A}$ by replacing each $e$ with the object $I$, each $x$ with the identity functor $1 _{\mathcal{A}}$, and tensoring the resulting sequence. 
More formally, we fix
\begin{align*}
e_{\mathcal{A}} : \ast & \longrightarrow \mathcal{A} & x_{\mathcal{A}} : \mathcal{A} & \longrightarrow \mathcal{A} \\
\ast & \longmapsto I & A & \longmapsto A
\end{align*}
and inductively define
\begin{align*}
(ew)_{\mathcal{A}} : \mathcal{A}^n & \longrightarrow \mathcal{A} & (xw)_{\mathcal{A}} : \mathcal{A}^{n+1} & \longrightarrow \mathcal{A} \\
(A_1,\dotsc,A_n) & \longmapsto I \otimes w_{\mathcal{A}}(A_1,\dotsc,A_n) & (A_1,\dotsc,A_{n+1}) & \longmapsto A_1 \otimes w_{\mathcal{A}}(A_2,\dotsc,A_{n+1})
\end{align*}
where $n=\mathrm{length}(w)$.

\begin{thm} \label{thm:wkgradingcoherence}
Let $u,v \in \mathscr{W}$ such that $\mathrm{length}(u)=\mathrm{length}(v)$ and $\mathrm{breadth}(u)=\mathrm{breadth}(v)$. 
There is a unique natural isomorphism $u_{\mathcal{A}} \cong v_{\mathcal{A}}$ built up from tensor products and compositions of components of $s$, $s^{-1}$, and the identity.
\end{thm}

We will defer the proof to the end of this appendix.

We call a component of any natural isomorphism arising from Theorem~\ref{thm:wkgradingcoherence} an {\sf adjustment}. 
For two morphisms $f:X \rightarrow Y$ and $g:Z \rightarrow W$ we write $f \sim g$ and say that $f$ and $g$ are {\sf adjustment equivalent} if there exist adjustments
\begin{align*}
q:X   \xrightarrow{\sim} Z \qquad  r:Y   \xrightarrow{\sim} W
\end{align*}
such that $g=r \circ f \circ q^{-1}$. 

\begin{eg}
For any morphism $f:X \rightarrow Y$ in $\mathcal{A}$, we have $f\otimes 1_I \sim 1_I \otimes f$, because
\begin{equation*}
f \otimes 1_I= s_Y^{-1} \circ (1_I \otimes f) \circ s_X
\end{equation*}
by the naturality of simple adjustments.
\end{eg}


%

\subsection{Breadth grading}

Suppose $\mathcal{A}$ is a monoidal category.
Assuming $\mathcal{A}$ is small, the set $\mathrm{Ob}(\mathcal{A})$ has the structure of a monoid.
We call a monoidal homomorphism $b:\mathrm{Ob}(\mathcal{A}) \rightarrow \mathbb{Z}_{\geq 0}$ a {\sf breadth function}.


\begin{defn}
Let $\mathcal{A}$ be a monoidal category with a breadth function $b$.
The {\sf weak grading of $\mathcal{A}$ concentrated in breadth $b$} is the following weakly graded monoidal category $\mathcal{A}[b]$. 
\begin{description}
\item[Objects] The objects of $\mathcal{A}[b]$ are formal free tensor products of objects in $\mathcal{A}$ and a new object $I$. 
In other words, each object $X$ in $\mathcal{A}[b]$ is a formal sequence
\begin{equation*}
I^{\otimes r_0} \otimes X_1 \otimes I^{\otimes r_1} \otimes X_2 \otimes \dotsb \otimes I^{\otimes r_{m-1}} X_{m} \otimes I^{\otimes r_{m}}
\end{equation*}
for some non-negative integers $r_{0},r_{m}$, positive integers $r_1,r_2,\dotsc,r_{m-1}$, and non-identity objects $X_1,X_2,\dotsc,X_m$ in $\mathcal{A}$.
The tensor product on objects in $\mathcal{A}$ extends in the obvious way to objects in $\mathcal{A}[b]$. 
We also extend the breadth function $b$ to a monoidal homomorphism $b:\mathrm{Ob}(\mathcal{A}[b]) \rightarrow \mathbb{Z}_{\geq 0}$ by fixing $b(I)=1$.

\item[Morphisms] For any object $X$ of the above form write $X'$ for the object
\begin{equation*}
X_1 \otimes X_2 \otimes \dotsb \otimes X_m
\end{equation*}
in $\mathcal{A}$. 
We define
\begin{equation*}
\Hom_{\mathcal{A}[b]}(X,Y)=
\begin{cases}
\Hom_{\mathcal{A}}(X',Y') & \text{if $b(X)=b(Y)$,} \\
0 & \text{otherwise.}
\end{cases}
\end{equation*}
Composition and tensor products follow from those in $\mathcal{A}$.

\item[Weak grading] For $X$ an object in $\mathcal{A}[b]$, the natural isomorphism $s_X : X \otimes I \rightarrow I \otimes X$ in $\mathcal{A}[b]$ corresponding to the identity morphism $1_{X'}$ in $\mathcal{A}$ gives $\mathcal{A}[b]$ the structure of a weakly graded monoidal category.
\end{description}
\end{defn}

If $f:X \rightarrow Y$ is a morphism in $\mathcal{A}[b]$, write $f' : X' \rightarrow Y'$ for the corresponding morphism in $\mathcal{A}$. 
It is easy to check that this mapping is functorial. 
We write $b(f)$ for the non-negative integer $b(X)=b(Y)$.

\begin{rmk}
The category $\mathcal{A}[b]$ is the weak graded analogue of the following graded construction. 
For a monoidal category $\mathcal{A}$ with a breadth function $b$, define a grading by setting $\deg f=b(X)-b(Y)$ for each morphism $f:X \rightarrow Y$. 
As in Example~\ref{eg:grwkgr}, we add grade shifts and restrict to homogeneous morphisms to obtain the category $\mathcal{A}\langle b\rangle$.
We may extend the breadth function $b$ to all of $\mathcal{A}\langle b \rangle$ as above. 
For any morphism $g:U \rightarrow V$ in $\mathcal{A}\langle b\rangle$, we have $0=\deg g=b(U)-b(V)$, which allows us to define the breadth of $g$ to be $b(g)=b(U)=b(V)$ as in the weakly graded case.

Our naming convention for $\mathcal{A}[b]$ (``concentrated in breadth $b$'') comes from a special case of the above graded construction. 
If $\mathcal{A}$ is a category of modules over some ring $R$, then we may equivalently construct the grading by considering $R$ to be a graded ring concentrated in degree $0$ and each object $X$ to be concentrated in degree $-b(X)$.
\end{rmk}

As a consequence of our coherence result, there is an alternative presentation of $\mathcal{A}[b]$ in terms of generators and relations. 
First we introduce a way of embedding morphisms from $\mathcal{A}$ into $\mathcal{A}[b]$.

\begin{defn}
Let $f:U \rightarrow V$ be a morphism in $\mathcal{A}$. 
The (left) minimal breadth representative of $f$ is the morphism $g:X \rightarrow Y$ in $\mathcal{A}[b]$ such that $g'=f$ and 
\begin{align*}
X  =I^{\otimes \max(0,b(V)-b(U))} \otimes U \text{,} \qquad 
Y   =I^{\otimes \max(0,b(U)-b(V))} \otimes V \text{.}
\end{align*}
\end{defn}

\begin{thm}
Let $\mathcal{M}$ be the set of all minimal breadth representatives of morphisms in $\mathcal{A}$.
The category $\mathcal{A}[b]$ is generated as a monoidal category by the morphisms
\begin{equation*}
\{1_I\} \cup \{s_X: X \in \mathrm{Ob}(\mathcal{A})\} \cup \mathcal{M}
\end{equation*}
%
%
subject to the following relations:
\begin{itemize}
\item the usual weak grading axioms \ref{item:naturality}--\ref{item:identity};

\item for morphisms
$f: X  
 \longrightarrow Y, 
g: Z 
 \longrightarrow W, 
h: U 
  \longrightarrow V
$ 
in $\mathcal{M}$ such that $f' \circ g'=h'$, we have 
\begin{equation*}
(1_I^{\otimes \max(0,b(g)-b(f))} \otimes f) \circ (1_I^{\otimes \max(0,b(f)-b(g))} \otimes g) \sim 1_I^{\otimes \max(b(f),b(g))-b(h)} \otimes h \text{;}
\end{equation*}
\item for morphisms $f: X  
 \longrightarrow Y, 
g: Z 
 \longrightarrow W, 
h: U 
  \longrightarrow V
$  
in $\mathcal{M}$ such that $f' \otimes g'=h'$, we have 
\begin{equation*}
f \otimes g \sim 1_I^{\otimes b(f)+b(g)-b(h)} \otimes h \text{.}
\end{equation*}
\end{itemize}
\end{thm}

\begin{proof}
Let $\mathcal{B}$ be the monoidal category defined by the above generators and relations. 
It is clear that the same relations hold in $\mathcal{A}[b]$, so there is a functor $\mathcal{B} \rightarrow \mathcal{A}[b]$. 
It is enough to show that this functor is full and faithful. 
Let $X,Y$ be objects in $\mathcal{B}$ such that $b(X)=b(Y)$.
We will show that any morphism $X \rightarrow Y$ can be written in the form
\begin{equation*}
q \circ (1_I^{b(X)-\max(b(X'),b(Y'))} \otimes f) \circ p^{-1} \text{,}
\end{equation*}
where $p,q$ are adjustments and $f$ is a minimal breadth representative.
In other words, we will show that every morphism in $\mathcal{B}$ is adjustment equivalent to the tensor product of a minimal breadth representative and some number of copies of $1_I$.
This automatically gives fullness and faithfulness of the functor above, which proves the result.
Since the generating morphisms of $\mathcal{B}$ are all already of this form, it is enough to show that any composition or tensor product of two morphisms of this form is again of this form.  
Now, consider a composition
\begin{equation*}
q \circ (1_I^{\otimes m} \otimes f) \circ p^{-1} \circ t \circ (1_I^{\otimes n} \otimes g) \circ r^{-1}
\end{equation*}
of two morphisms of the above form.
Both $f$ and $g$ are minimal breadth representatives, so their domains and codomains are ``left-adjusted'', i.e.~of the form $I^{\otimes l} \otimes U$ for some object $U$ in $\mathcal{A}$ and some non-negative integer $l$.
The adjustment $p^{-1} \circ t$ is an isomorphism between $I^{\otimes n} \otimes \cod g$ and $I^{\otimes m} \otimes \dom f$ which are both left-adjusted, so in fact they must be equal.
By Theorem~\ref{thm:wkgradingcoherence} we must have $p=t$, so the composition above equals
\begin{align*}
q \circ (1_I^{\otimes m} \otimes f) \circ (1_I^{\otimes n} \otimes g) \circ r^{-1} & = q \circ (1_I^{\otimes (m-j)} \otimes 1_I^{j} \otimes f) \circ (1_I^{\otimes (n-k)} \otimes 1_I^{\otimes k} \otimes g) \circ r^{-1} \\
& \sim q \circ (1_I^{\otimes (m-j)} \otimes h) \circ r^{-1}
\end{align*}
where $j=\max(0,b(g)-b(f))$, $k=\max(0,b(f)-b(g))$, and $h$ is the minimal breadth representative of $f' \circ g'$.
 Similarly, consider a tensor product of two morphisms of the above form.
We have
\begin{align*}
&\qquad(q \circ (1_I^{\otimes m} \otimes f) \circ p^{-1}) \otimes (t \circ (1_I^{\otimes n} \otimes g) \circ r^{-1})\\
& = (q \otimes t) \circ (1_I^{\otimes m} \otimes f \otimes 1_I^{\otimes n} \otimes g) \circ (p^{-1} \otimes r^{-1}) \\
& \sim (q \otimes t) \circ (1_I^{\otimes (m+n)} \otimes f \otimes g) \circ (p^{-1} \otimes r^{-1}) \\
& \sim (q \otimes t) \circ (1_I^{\otimes (m+n+b(f)+b(g)-b(h))} \otimes h) \circ (p^{-1} \otimes r^{-1}) \text{,}
\end{align*}
where $h$ is the minimal breadth representative of $f' \otimes g'$.
\end{proof}

\subsection{Proof of coherence}

We conclude with the proof of the coherence theorem for weakly graded monoidal categories (Theorem~\ref{thm:wkgradingcoherence}). 
The strategy is broadly similar to Mac Lane's proof of the coherence theorem for monoidal categories \cite[VII.2]{maclane-cwm}. 
This involves first proving the result for a single object $X$ in the category $\mathcal{A}$, and then extending to all of $\mathcal{A}$. 

Now let $\mathscr{S}$ be the set of words in the symbols $\{\sigma_w,\sigma_w^{-1} : w \in \mathscr{W}\} \cup \{\iota_e,\iota_x\}$ defined inductively as follows.
For any $w \in \mathscr{W}$ we have $\sigma_w,\sigma_w^{-1} \in \mathscr{S}$.
Moreover, for any $\alpha \in \mathscr{S}$ and $w \in \mathscr{W}$ we also have $\iota_e\alpha,\iota_x\alpha \in \mathscr{S}$ and $\alpha \iota_e, \alpha\iota_x \in \mathscr{S}$. 
For convenience we write $\iota_w$ for $\iota_{w_1} \iota_{w_2} \dotsm \iota_{w_m}$, where $w=w_1 w_2 \dotsm w_n$ is a word in $\mathscr{W}$.
We inductively define $\dom : \mathscr{S} \rightarrow \mathscr{W}$ and $\cod : \mathscr{S} \rightarrow \mathscr{W}$ as follows:
\begin{align*}
\dom(\sigma_w)& =we & \cod(\sigma_w)& =ew \\
\dom(\sigma_w^{-1})& =ew & \cod(\sigma_w^{-1})& =we \\
\dom(\iota_w \alpha)& =w\dom(\alpha) & \cod(\iota_w \alpha)& =w\cod(\alpha) \\
\dom(\alpha \iota_w)& =\dom(\alpha)w & \cod(\alpha \iota_w)& =\cod(\alpha)w  
\end{align*}
Let $\mathscr{G}$ be the quiver with vertices given by $\mathscr{W}$ and arrows given by $\mathscr{S}$.
It is easy to verify that for any word in $\alpha \in \mathscr{S}$, $\mathrm{length}(\dom(\alpha))=\mathrm{length}(\cod(\alpha))$ and $\mathrm{breadth}(\dom(\alpha))=\mathrm{breadth}(\cod(\alpha))$.
Thus the graph $\mathscr{G}$ has components $\mathscr{G}_{n,k}$ whose vertices $\mathscr{W}_{n,k}$ consist of words of length $n$ and breadth $k$.

Now let $\mathcal{A}$ be a weakly graded monoidal category. 
We fix an object $X$ in $\mathcal{A}$ and set 
\begin{align*}
\mathscr{J}_X(e)& =I & \mathscr{J}_X(x)& =X \\
\mathscr{J}_X(ew)& =I \otimes \mathscr{J}_X(w) & \mathscr{J}_X(xw) &= X \otimes \mathscr{J}_X(w) \\
\mathscr{J}_X(\sigma_w)& =s_{w_X} & \mathscr{J}_X(\sigma_w^{-1})& =s_{w_X}^{-1} \\
\mathscr{J}_X(\iota_w\alpha)& =1_{w_X} \otimes \mathscr{J}_X(\alpha) & \mathscr{J}_X(\alpha\iota_w)& =\mathscr{J}_X(\alpha) \otimes 1_{w_X}
\end{align*}

\begin{prop} \label{prop:jpathsequal}
Let $u,v \in \mathscr{W}$ such that $\mathrm{length}(u)=\mathrm{length}(v)$ and $\mathrm{breadth}(u)=\mathrm{breadth}(v)$.
Suppose $\alpha_1\circ \dotsb \circ \alpha_m$ and $\alpha_1' \circ \dotsb \circ \alpha_{m'}'$ are two paths in $\mathscr{G}$ from $u$ to $v$.
Then
\begin{equation*}
\mathscr{J}_X(\alpha_m) \circ \dotsb \circ \mathscr{J}_X(\alpha_1)=\mathscr{J}_X(\alpha_{m'}') \circ \dotsb \circ \mathscr{J}_X(\alpha_1') \text{.}
\end{equation*}
\end{prop}

\begin{proof}
Let $n=\mathrm{length}(u)=\mathrm{length}(v)$ and $k=\mathrm{breadth}(u)=\mathrm{breadth}(v)$. 
We will pivot on the sink vertex $w^{(n,k)}=e^k x^n$ in the component $\mathscr{G}_{n,k}$.
Every nonempty word in $S$ contains exactly one symbol of the form $\sigma_w$ or $\sigma_w^{-1}$ for $w \in W$. 
Call such words directed or anti-directed respectively.
It is easy to check that for any two directed words $\alpha,\alpha'$ with the same domain and codomain, we must have $\mathscr{J}_X(\alpha)=\mathscr{J}_X(\alpha')$.

We inductively define a function $\rho:W \rightarrow \mathbb{Z}_{\geq 0}$ by
\begin{align*}
\rho(e)  =0 \qquad 
\rho(x)  =0 \qquad 
\rho(ew)  =\rho(w) \qquad 
\rho(xw)  =\rho(w)+\mathrm{breadth}(w) \text{.}
\end{align*}
We also inductively define a function $\mathrm{can}_{n,k}$ mapping words in $W_{n,k}$ to directed paths in $\mathscr{G}_{n,k}$ by
\begin{align*}
\mathrm{can}_{0,1}(e)  &=\emptyset \quad 
\mathrm{can}_{1,0}(x)  =\emptyset \quad 
\mathrm{can}_{n,k}(ew)  =\iota_e \mathrm{can}_{n,k-1}(w) \\ 
\mathrm{can}_{n,k}(xw)& =(\iota_e^{k-1} \sigma_x \iota_x^{n-1}) \circ \dotsb \circ (\iota_e \sigma_x \iota_{e}^{k-2} \iota_{x}^{n-1}) \circ (\sigma_x \iota_{e}^{k-1} \iota_{x}^{n-1}) \circ (\iota_x \mathrm{can}_{n-1,k}(w))
\end{align*}
It can be shown that $\mathrm{can}_{n,k}(w)$ is the longest directed path in $\mathscr{G}_{n,k}$ from $w$ to $w^{(n,k)}$, and that $\rho(w)=\mathrm{length}(\mathrm{can}_{n,k}(w))$. 


\begin{lem} \label{lem:wg-dirpath}
For any $u \in \mathscr{W}_{n,k}$, $\mathscr{J}_X$ maps all directed paths from $u$ to $w^{(n,k)}$ to the same morphism. 
\end{lem}

Before we prove this lemma, we will show that the proposition follows from it almost immediately.
For $\alpha \in \mathscr{S}$ let $\mathrm{inv}(\alpha)$ be the word obtained by switching the symbols $\sigma_w \leftrightarrow \sigma_w^{-1}$.
Clearly $\mathscr{J}_X(\mathrm{inv}(\alpha))=\mathscr{J}_X(\alpha)^{-1}$, and we may write any anti-directed word as the formal inverse of a directed word. 
Let us write the path $\alpha_m \circ \dotsb \circ \alpha_1$ from $u$ to $v$ in this manner, using formal inverses of directed words for any anti-directed word that appears.
For example, if $\alpha_2$ is the only anti-directed word in this path, we write:
\begin{equation*}
\xymatrix{
u \ar[r]^{\alpha_1} & {\bullet}  & {\bullet} \ar[l]_{\mathrm{inv}(\alpha_2)} \ar[r]^{\alpha_3} & {\bullet} \ar@{.}[r]  & {\bullet} \ar[r]^{\alpha_m}  & v
}
\end{equation*}
Now draw canonical paths downwards to $w^{(n,k)}$ underneath each of these objects:
\begin{equation*}
\xymatrix{
u \ar[r]^{\alpha_1} \ar[d] & {\bullet}  \ar[d] & {\bullet} \ar[l]_{\mathrm{inv}(\alpha_2)} \ar[r]^{\alpha_3} \ar[d] & {\bullet} \ar@{.}[r] \ar[d] & {\bullet} \ar[r]^{\alpha_m} \ar[d] & v_X \ar[d] \\
w^{(n,k)} \ar@{=}[r] & w^{(n,k)} \ar@{=}[r] & w^{(n,k)} \ar@{=}[r] & w^{(n,k)} \ar@{.}[r] & w^{(n,k)} \ar@{=}[r] & w^{(n,k)}
} 
\end{equation*}
After applying $\mathscr{J}_X$, each square commutes by the above lemma, so
\begin{equation*}
\mathscr{J}_X(\alpha_m) \circ \dotsb \circ \mathscr{J}_X(\alpha_1)=\mathscr{J}_X(\mathrm{can}_{n,k}(v))^{-1} \circ \mathscr{J}_X(\mathrm{can}_{n,k}(u)) \text{.}
\end{equation*}
Since the right-hand side only depends on $u$ and $v$, we are done.
\end{proof}

\begin{proof}[Proof of Lemma~\ref{lem:wg-dirpath}]
We induct on $\rho(u)$, $n$ and $k$.
Suppose we have two directed paths from $u$ to $w_{n,k}$ which start with $\alpha$ and $\alpha'$ respectively.
 $$
$$As $\rho(w)<\rho(u)$, we are then done by induction.    
Otherwise, suppose $w \neq w'$ and $\alpha \neq \alpha'$. 
It is enough to find some $w'' \in W$ and some paths from $w$ and $w'$ to $w''$ such that the following diamond
  $$%
$$ 
commutes after applying $\mathscr{J}_X$. 
For if so, then $\rho(w'')<\rho(u)$, and by induction the trapezoids in the following diagram

\!\!\!\!\! $$%
$$ 
commute after applying $\mathscr{J}_X$, and therefore the whole diagram commutes.

\medskip\noindent\textbf{Case 1.} 
If $\alpha=\iota_z \beta$ and $\alpha'=\iota_{z'} \beta'$ for some $z,z' \in \mathscr{W}$ and $\beta,\beta' \in \mathscr{S}$, then both $z$ and $z'$ begin with some non-empty word $z''$. 
Thus $u$, $w$, and $w'$ also begin with $z''$, and we can write $\alpha$ and $\alpha'$ as $\iota_{z''}\gamma$ and $\iota_{z''}\gamma'$ respectively.
Let $u'=\dom(\gamma)$, $y=\cod(\gamma)$, and $y'=\cod(\gamma')$, and let $n'$ and $k'$ be the length and breadth of $y$ (or $y'$) respectively.
Since $y$ is a strict subword of $w$, we must have $n'<n$ or $k'<k$. 
Taking $w''=z''w^{(n',k')}$ we obtain the following diamond
 $$
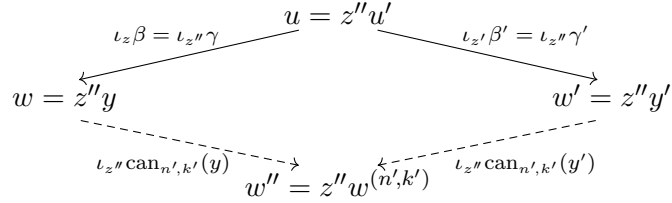
\begin{tikzpicture}[scale=0.9]
  \draw(0,1.25) node {$u=z^{\prime  \prime} u^\prime $};
  
    \draw(-4,0) node {$w=z''y$};    \draw(4,0) node {$w^\prime=z''y' $};
  
   \draw[->] (120:0.9+0.25) --  (176:3.8) node [midway,above ] { \scalefont{0.7}$  \iota_{z }\beta =\iota_{z^{\prime\hspace{-1pt} \prime} }\gamma \quad\quad $}  ; 
      \draw[->] (180-120:0.9+0.25) --  (180-176:3.8) node [midway,above ] { \quad\quad\; \scalefont{0.7}$  \iota_{z^\prime }\beta^\prime =\iota_{z^{\prime\hspace{-1pt} \prime} }\gamma^\prime   $}  ;

   \draw[<-, densely dashed] (-120:1.15) --  (-175:3.8)   node [midway,below ] {   \scalefont{0.7}$  \iota_{z^{\prime\hspace{-1pt} \prime} } {\rm can}_{n^\prime \! ,k^\prime }(y ) \quad\quad\; $}  ; 
      \draw[<-, densely dashed] (-180+120:  1.15) --  (-180+175:3.8)   node [midway,below ] { \quad\quad\; \scalefont{0.7}$  \iota_{z^{\prime\hspace{-1pt} \prime} } {\rm can}_{n^\prime \! ,k^\prime }(y^\prime )  $}  ;

      \draw(0,-1.25) node {$ w^{\prime \prime}  = z^{\prime   \prime}  w^{(n^\prime \!,k^\prime )}$};  
      
   \end{tikzpicture}$$ 
 which commutes after applying $\mathscr{J}_X$ by induction on $n$ and $k$.
A similar proof works if $\alpha=\beta \iota_z$ and $\alpha'=\beta' \iota_{z'}$ for some $z,z'' \in \mathscr{W}$ and $\beta,\beta' \in \mathscr{S}$.

\medskip\noindent\textbf{Cases 2 \& 3.}
The next cases to consider occur when one of $\alpha$ or $\alpha'$ is $\sigma_{y}$ for some $y \in \mathscr{W}$. 
Without loss of generality suppose $\alpha=\sigma_y$. 
If $\alpha'$ is of the form $\iota_{z'} \sigma_{y'}$ for some $y',z' \in \mathscr{W}$ then we must have $y=z'y'$ and thus $u=ye=z'y'e$. 
Taking $w''=ez'y'$ we obtain the following diamond
 $$
$$ 
which commutes after applying $\mathscr{J}_X$ by \ref{item:mult}.
On the other hand, if $\alpha'$ is of the form $\sigma_{y'} \iota_{z'}$ for some $y',z' \in W$, then we must have $z'=z''e$ for some $z'' \in W$, and thus $y=y'ez''$. 
Taking $w''=eey'z''$ we obtain the following diagram

\!\!\!\!\!$$%
$$ 
which commutes after applying $\mathscr{J}_X$, by the naturality of $s$.

\medskip\noindent\textbf{Cases 4 \& 5.}
The last cases are when $\alpha=\sigma_y \iota_z$ and $\alpha'=\iota_{z'} \sigma_{y'}$ for some $y,y',z,z' \in W$, so that $u=yez=z'y'e$. 
Suppose first that $z'$ starts with $ye$. 
Then there is some $z'' \in W$ such that $z'=yez''$. 
Using $yez=z'y'e$ it is also clear that $z=z''y'e$ too.
Taking $w''=eyz''ey'$ we obtain the   diamond
$$\begin{tikzpicture}[scale=0.9]
  \draw(0,1.25) node {$u= y   e z^{\prime \prime  }y' e$};
  
    \draw(-4,0) node {$  ey  z '' y'e= w$};    \draw(4,0) node {$w^\prime= y e  z''e y' $};
  
   \draw[->] (120:0.9+0.25) --  (176:3.8) node [midway,above ] { \scalefont{0.7}$   \sigma_{y} \iota_{{z  ^    {\prime\hspace{-1pt} \prime}} \!y'\!e} $}  ; 
      \draw[->] (180-120:0.9+0.25) --  (180-176:3.8) node [midway,above ] { \quad\quad\; \scalefont{0.7}$  \iota_{ ye{z  ^    {\prime\hspace{-1pt} \prime}}  } \sigma_{y'}     $}  ;

   \draw[<-, densely dashed](-120:1.15) --  (-175:3.8)   node [midway,below ] {    \scalefont{0.7}$\iota   _{ey  } \iota_{z  ^    {\prime\hspace{-1pt} \prime}} \sigma_{y'}$\quad\quad}  ;

         \draw[<-, densely dashed] (-180+120:  1.15) --  (-180+175:3.8)   node [midway,below ] { \quad\quad\; \scalefont{0.7}$ \sigma_{ y   }\iota_{z  ^    {\prime\hspace{-1pt} \prime}} \iota _{e y '  }$}  ;

      \draw(0,-1.25) node {$ w^{\prime \prime}  = e   y  z^{\prime  \prime } e y' $};  
      
   \end{tikzpicture}$$ 
which commutes after applying $\mathscr{J}_X$ by bifunctoriality of the tensor product.
On the other hand, if $ye$ starts with $z'$, then there exists some $y'' \in W$ such that $y=z'y''$. 
This also implies that $y'e$ ends with $z$, so there also exists some $z'' \in W$ such that $z=z''e$. 
This means that $y'=y''ez''$. 
This time we complete the diamond in two steps. 
First, we compose $\iota_{z'}\sigma_{y''ez''}$ with $\sigma_{z'} \iota_{y''ez''}$.
By \ref{item:mult} of a weak grading, this composition equals $\sigma_{z'y''ez''}$. 
Thus we have reduced to a previous case and so we are done.

\!\!\!\!\!\!$$\begin{tikzpicture}[scale=1.2]
  \draw(0,1.6) node {$ u=z'y''ez''e$};
  
    \draw(-3,0.5) node {$ez'y''z''e=w $};    \draw(3,0.5) node {$w'=z'ey''ez'' $};
  
   \draw[->] (111:1.5) --  (165:2.9) node [midway,above ] { \scalefont{0.7}$\sigma_{z'y^{\prime\hspace{-1pt} \prime} } 1_{z^{\prime\hspace{-1pt} \prime} e}$\qquad\quad}  ; 
      \draw[->] (180-111:1.5) --  (180-165:2.9) node [midway,above ] {\qquad \quad \scalefont{0.7}$\iota_{z'} \sigma_{y^{\prime\hspace{-1pt} \prime} ez^{\prime\hspace{-1pt} \prime} } $ }  ; 
      \draw[double] (-3,0.3 )--(-3,-0.22);
      \draw[densely dashed,->] (3,0.3 )--(3,-0.22)			node [midway,right] {\scalefont{0.7}  $	 \sigma_{z'}1_{y^{\prime\hspace{-1pt} \prime}  e z^{\prime\hspace{-1pt} \prime} }$};
      
      \draw(-3,-0.5) node {$ez'y''z''e$};    \draw(3,-0.5) node {$ez'y''ez''  $}			;

 \draw(0,-1.6) node {$w''=eez'y''z'' $};

  \draw[<-] (90-21:-1.5) --  (-165:2.9) node [midway,below ] {  \scalefont{0.7}$\quad {\sigma_{ez'y^{\prime\hspace{-1pt} \prime} z^{\prime\hspace{-1pt} \prime} }}$ \qquad \quad}  ; 
      \draw[->]     (-180+165:2.9)--(111:-1.5)  node [midway,below] { \qquad \quad \scalefont{0.7}$ 1_e\sigma_{z'y^{\prime\hspace{-1pt} \prime} }1_{z^{\prime\hspace{-1pt} \prime} } $}  ;

   \end{tikzpicture} 
   $$  \end{proof}

To extend to the full coherence theorem, we consider objects in a higher category. 

\begin{proof}[{Proof of Theorem~\ref{thm:wkgradingcoherence}}]
Let $\mathrm{Iter}(\mathcal{A})$ be the category of functors of the form $\mathcal{A}^n \rightarrow \mathcal{A}$, where $n$ is a non-negative integer.
It is clear that $\mathrm{Iter}(\mathcal{A})$ is also monoidal, with the tensor product of two functors $F:\mathcal{A}^m \rightarrow \mathcal{A}$ and $G:\mathcal{A}^n \rightarrow \mathcal{A}$ defined to be
\begin{align*}
(F \otimes G): \mathcal{A}^{m+n} \longrightarrow \mathcal{A},   \; \;
(A_1,\dotsc, A_{m+n}) & \longmapsto F(A_1,\dotsc,A_m) \otimes G(A_{m+1},\dotsc,A_{m+n})
\end{align*}
We observe that $w_{\mathcal{A}}$ is precisely $\mathscr{J}_{1_{\mathcal{A}}}(w)$ as defined above, where we consider the identity functor $1_{\mathcal{A}}$ as an object in $\mathrm{Iter}(\mathcal{A})$.
Applying $\mathscr{J}_{1_{\mathcal{A}}}$ to any path between $u$ and $v$ gives a isomorphism in $\mathrm{Iter}(\mathcal{A})$ between $u_{\mathcal{A}}$ and $v_{\mathcal{A}}$, or in other words, a natural isomorphism between the two functors.
Uniqueness of this natural isomorphism follows from Proposition~\ref{prop:jpathsequal}.
\end{proof}

\begin{Acknowledgements*}
The first and third authors    thank the   Institut Henri Poincar\'e for hosting us during the 
  thematic trimester on representation theory.
   The first author thanks Oberwolfach Institute for hosting him during the mini-workshop ``Kronecker, Plethysm, and Sylow branching Coefficients and their Applications to Complexity Theory". 
   The first author was funded by EPSRC grant EP/V00090X/1 and the 
third author was funded by the Royal Commission for the Exhibition of 1851.  
 We would also like to thank Jon Brundan and  Martina Lanini  
     for entertaining and informative discussions and 
  Geordie Williamson and the anonymous referees  for their helpful comments on a previous version of this paper.
 \end{Acknowledgements*}

%
%

\section{List of symbols}

For the convenience of the reader we list the symbols used in the main body of the paper in three categories: those corresponding to the general setup and basic combinatorics; those corresponding to the geometry and choice of paths; and those corresponding to the various algebras of interest. As Appendix A is relatively short and self-contained we omit those symbols here.

 
\begin{table}[h]
	\centering
	\caption{General symbols}
\begin{tabular}{m{6em} m{4em} m{6em} m{4em} m{6em} m{4em}}
		Symbol &  \S\S & Symbol &  \S\S & Symbol & \S\S    \\
		$h$ & 2 & $\ell$& 2 & $e$ & 2 \\
		$\sigma$ & 2 & $\lambda$ & 2.1 & $|\lambda|$ & 2.1\\
		$\lambda^{(i)}$ & 2.1 & ${\rm ct}(r,c,m)$& 2.1 & ${\rm res}(r,c,m)$& 2.1 \\
		$\mathscr{P}_{\ell}(n)$ & 2.1 & $\SSTT_\la$ &2.1 & $\mathscr{P}_{\underline{h}}(n)$ & 2.1\\
		$\Std (\lambda)$ & 2.1 & $\varnothing$ & 2.1& $\mathfrak{S}_{\aatch}$ & 2.2\\
		$\mathfrak{S}_f$ & 2.2 & $\widehat{\mathfrak{S}}_{\aatch}$ & 2.2 & $\mathfrak{S}^f$ & 2.2\\
	    $\underline{w}$ & 2.4 & $\overline{w}$ & 2.4 & $r_{\aatch}(t)$ & 2.4 \\
		$\al(p)$ & 2.4 & $\emptyset(q)$ & 2.4 &	$\underline{i}$ & 3.3 \\
		$s_r$ & 3.3 & $s_r(\underline{i})$ & 3.3 & $w_q^p$ & 3.3 \\
		$(i,i+1)_b$ & 4.3 & $\underline{w}_b$ & 4.3 & $B$ & 4.3 \\ 
		${\rm nib}(\w)$ & 4.3 & ${\rm nib}(\underline{i})$ & 4.3 
	&	${\rm det}_{\underline{h}}$ & 7.1 \\
		  $\Std_n$ & 7.1& $\Gamma$ & 7.1 
	&	$\square$ & 7.1
	\end{tabular}
\end{table}

\begin{table} [h]
 	\caption{Geometry and paths}
	\begin{tabular}{m{6em} m{4em} m{6em} m{4em} m{6em} m{4em}}
		Symbol &  \S\S & Symbol &  \S\S &Symbol & \S\S    \\
		$\varepsilon_i$ & 2.2 & ${\mathbb E}_{\aatchpair}$ & 2.2& $\overline{\mathbb E}_{\underline{h}}$ & 2.1\\
		$\langle \; , \; \rangle$ & 2.2 & $({-})^T$ & 2.2 & $\Phi$ & 2.2\\
		$\Phi_0$ & 2.2	& $\Delta$ & 2.2 & $\Delta_0$ & 2.2 \\ 
		$\alpha$ & 2.2 &$ s_{\alpha,re}$ & 2.2 & $\alpha_0$ & 2.2\\
		$\Pi$ & 2.2 & $S$ & 2.2 & $s_{\alpha}$ & 2.2\\
		$\rho$ & 2.2 & $w\cdot x$ & 2.2 & ${\mathbb E} ({\alpha, re})$ & 2.2\\
		${\mathbb E}^> ({\alpha, re})$ & 2.2 & ${\mathbb E}^< ({\alpha, re})$ & 2.2 & $\overline{\mathbb E}_{\underline{h}}^+$ & 2.2\\
\end{tabular}
\end{table}

\begin{table}[h]
	\caption{Algebras and elements}
	\begin{tabular}{m{6em} m{4em} m{6em} m{4em} m{6em} m{4em}}
		Symbol &  \S\S & Symbol &  \S\S & Symbol & \S\S    \\
		${\sf 1}_{ {\al}}$ & 3.1 & 	${\sf 1}_\w	$ & 3.1 & $\mathscr{S} ({n,\sigma})$ & 3.1\\
		$\mathscr{S}_{\underline{h}} ({n,\sigma})$ & 3.1 & $\ast$ & 3.1 &${\sf 1}_{ {\al}}$ & 3.2 \\
		${\sf 1}_{\emptyset}$ &  3.2& ${\sf 1}_\w$ & 3.2 & ${\sf 1}_{\underline{w'}}^\w$ & 3.2\\
		${\sf 1}_{\emptyset\al}^{\al\emptyset}$ & 3.2 & $ {\sf SPOT}_\al^\emp$ & 3.2 & ${\sf FORK}_{\al\al}^{\emp\al}$ & 3.2\\ 
		${\sf HEX}_{\al\bet\al}^{\bet\al\bet}$ & 3.2 & ${\sf COMM}_{\bet\gam}^{\gam\bet}$ & 3.2& $\mathscr{ S}^{\rm br}({n,\sigma})$ & 3.2\\
		$\mathscr{S}^{\rm br}_{\aatchpair}({n,\sigma})$ & 3.2 & $\Lambda^+(\leq n,\sigma)$ & 3.2 & ${\sf 1}^+_{n,\sigma}$ & 3.2\\
		$\mathcal{H}_n $ & 3.3 & $e_{\underline{i}}$ & 3.3 & $y_i$ & 3.3\\
\end{tabular}
\end{table}

%
%
%
\providecommand{\bysame}{\leavevmode\hbox to3em{\hrulefill}\thinspace}
\providecommand{\MR}{\relax\ifhmode\unskip\space\fi MR }
\providecommand{\MRhref}[2]{%
  \href{http://www.ams.org/mathscinet-getitem?mr=#1}{#2}
}
\providecommand{\href}[2]{#2}

\end{document}